%% file: Two_phase_local_revised.tex
\documentclass[reqno]{amsart}

\usepackage{amsmath,amssymb,amsthm,braket,stmaryrd}

\usepackage{color}

\definecolor{red}{cmyk}{0,0,0,1}
\definecolor{blue}{cmyk}{0,0,0,1}

\newtheorem{theo}{Theorem}[section]
\newtheorem{lemm}[theo]{Lemma}

\numberwithin{equation}{section}

\theoremstyle{definition}
\newtheorem{defi}[theo]{Definition}

\newtheorem{rema}[theo]{Remark}
\newtheorem{prop}[theo]{Proposition}
\newtheorem{assu}[theo]{Assumption}

\newcommand{\ba}{\mathbf{a}}
\newcommand{\bb}{\mathbf{b}}
\newcommand{\bd}{\mathbf{d}}

\newcommand{\bff}{\mathbf{f}}
\newcommand{\bg}{\mathbf{g}}
\newcommand{\bh}{\mathbf{h}}

\newcommand{\bm}{\mathbf{m}}
\newcommand{\bn}{\mathbf{n}}

\newcommand{\bu}{\mathbf{u}}
\newcommand{\bv}{\mathbf{v}}
\newcommand{\bw}{\mathbf{w}}
\newcommand{\bA}{\mathbf{A}}
\newcommand{\bB}{\mathbf{B}}

\newcommand{\bD}{\mathbf{D}}

\newcommand{\bF}{\mathbf{F}}
\newcommand{\bG}{\mathbf{G}}

\newcommand{\bI}{\mathbf{I}}

\newcommand{\bL}{\mathbf{L}}
\newcommand{\bM}{\mathbf{M}}
\newcommand{\bN}{\mathbf{N}}

\newcommand{\bQ}{\mathbf{Q}}

\newcommand{\bT}{\mathbf{T}}
\newcommand{\bV}{\mathbf{V}}

\newcommand{\DV}{\mathrm{Div}\,}
\newcommand{\dv}{\mathrm{div}\,}
\newcommand{\supp}{\mathrm{supp}\,}

\newcommand{\BR}{\mathbb{R}}
\newcommand{\BC}{\mathbb{C}}

\newcommand{\BN}{\mathbb{N}}
\newcommand{\BG}{\mathbb{G}}

\newcommand{\BI}{\mathbb{I}}

\newcommand{\BT}{\mathbb{T}}
\newcommand{\BU}{\mathbb{U}}

\newcommand{\td}{\widetilde{d}}
\newcommand{\tf}{\widetilde{f}}
\newcommand{\tg}{\widetilde{g}}
\newcommand{\tilh}{\widetilde{h}}
\newcommand{\tk}{\widetilde{k}}
\newcommand{\tu}{\widetilde{u}}

\newcommand{\trho}{\widetilde{\rho}}

\newcommand{\tvarphi}{\widetilde{\varphi}}
\newcommand{\tgamma}{\widetilde{\gamma}}

\newcommand{\tzeta}{\widetilde{\zeta}}

\newcommand{\hf}{\widehat{f}}

\newcommand{\hh}{\widehat{h}}

\newcommand{\hW}{\widehat{W}}

\newcommand{\CA}{\mathcal{A}}
\newcommand{\CB}{\mathcal{B}}

\newcommand{\CD}{\mathcal{D}}

\newcommand{\CF}{\mathcal{F}}
\newcommand{\CG}{\mathcal{G}}
\newcommand{\CH}{\mathcal{H}}
\newcommand{\CI}{\mathcal{I}}

\newcommand{\CK}{\mathcal{K}}
\newcommand{\CL}{\mathcal{L}}
\newcommand{\CM}{\mathcal{M}}
\newcommand{\CN}{\mathcal{N}}

\newcommand{\CP}{\mathcal{P}}
\newcommand{\CQ}{\mathcal{Q}}
\newcommand{\CR}{\mathcal{R}}
\newcommand{\CS}{\mathcal{S}}
\newcommand{\CT}{\mathcal{T}}
\newcommand{\CU}{\mathcal{U}}
\newcommand{\CV}{\mathcal{V}}
\newcommand{\CW}{\mathcal{W}}
\newcommand{\CX}{\mathcal{X}}
\newcommand{\CY}{\mathcal{Y}}
\newcommand{\CZ}{\mathcal{Z}}
\newcommand{\fg}{\mathfrak{g}}
\newcommand{\fh}{\mathfrak{h}}
\newcommand{\fp}{\mathfrak{p}}

\newcommand{\pd}{\partial}
\newcommand{\Hol}{\mathrm{Hol}\,}
\newcommand{\RP}{{\BR^N_+}}
\newcommand{\RM}{{\BR^N_-}}
\newcommand{\RZ}{{\BR^N_0}}
\newcommand{\lbrac}{\llbracket}
\newcommand{\rbrac}{\rrbracket}
\newcommand{\ov}{\overline}
\newcommand{\dx}{\,\mathrm{d}x}

\newcommand{\dt}{\,\mathrm{d}t}
\newcommand{\dtheta}{\,\mathrm{d}\theta}
\newcommand{\dtau}{\,\mathrm{d}\tau}

\newcommand{\wh}{\widehat}
\newcommand{\wt}{\widetilde}

\begin{document}
\title[Strong solutions to compressible--incompressible flows]
{Strong solutions to compressible--incompressible \\ two-phase flows with phase transitions}	

\author{Keiichi Watanabe}
\address{Department of Pure and Applied Mathematics, Graduate School of Fundamental Science and Engineering, Waseda University, 3-4-1 Ookubo, Shinjuku-ku, Tokyo, 169-8555, Japan}		
\email{keiichi-watanabe@akane.waseda.jp}
\thanks{This research was partly supported by JSPS Grant-in-Aid for JSPS Fellows \#19J10168 and Top Global University Project of Waseda University.}

\begin{abstract}
We consider a free boundary problem of compressible-incompressible two-phase flows with phase transitions in general domains of $N$-dimensional Euclidean space (e.g. whole space; half-spaces; bounded domains; exterior domains). The compressible fluid and the incompressible fluid are separated by either compact or non-compact sharp moving interface, and the surface tension is taken into account. In our model, the compressible fluid and incompressible fluid are occupied by the Navier-Stokes-Korteweg equations and the Navier-Stokes equations, respectively. This paper shows that for given $T > 0$ the problem admits a unique strong solution on $(0,T)$ in the maximal $L_p - L_q$ regularity class provided the initial data are small in their natural norms.
\end{abstract}

\maketitle
\section{Introduction}
\noindent
The present paper deals with a free boundary problem for compressible-incompressible two-phase flows with phase transitions in the isentropic case. Two immiscible viscous fluids are separated by a sharp interface with taking a surface tension into account. Our problem is formulated as follows: Let $\Omega$ be a domain in $N$-dimensional Euclidean space $\BR^N$ $(N \ge 2)$ surrounded by boundaries $\Gamma_+$ and $\Gamma_-$. In this paper, $\Gamma_+ = \emptyset$ or $\Gamma_- = \emptyset$ are admissible. For $t \ge 0$, the hypersurface $\Gamma_t$ stands a sharp moving interface, which separates $\Omega$ into $\Omega_{t+}$ and $\Omega_{t-}$ such that $\Omega \backslash \Gamma_t = \Omega_{t+} \cup \Omega_{t-}$, $\Omega_{t+} \cap \Omega_{t-} = \emptyset$, $\pd \Omega_{t+} = \Gamma_t \cup \Gamma_+$, and $\pd \Omega_{t-} = \Gamma_t \cup \Gamma_-$. Let $\dot\Omega_t = \Omega_{t+} \cup \Omega_{t-}$, and for any function $f$ defined on $\dot \Omega_t$, we write $f_\pm = f \vert_{\Omega_{t\pm}}$. We consider the following Cauchy problem:
\begin{align}
\label{eq-1.1}
&\left\{\begin{aligned}
\pd_t\varrho_+ + \dv (\varrho_+ \bv_+) & = 0 & \qquad & \text{in $\Omega_{t+}, \enskip t > 0$},\\
\varrho_+ (\pd_t \bv_+ + (\bv_+ \cdot \nabla) \bv_+) - \DV \BT_+ & = 0 & \qquad & \text{in $\Omega_{t+}, \enskip t > 0$},\\
\dv \bv_- & = 0 & \qquad & \text{in $\Omega_{t-}, \enskip t > 0$}, \\
\varrho_- (\pd_t \bv_- + (\bv_- \cdot \nabla) \bv_-) - \DV \BT_- & = 0 & \qquad & \text{in $\Omega_{t-}, \enskip t > 0$}
\end{aligned}\right.
\intertext{with the interfacial boundary conditions on $\Gamma_t$ $(t > 0)$:}
\label{eq-1.2}
&\left\{\begin{aligned}
& V_{\Gamma_t} = \bv_{\Gamma_t} \cdot \bn_t = \cfrac{\lbrac \varrho \bv \rbrac \cdot \bn_t}{\lbrac \varrho \rbrac},\\
& \lbrac \bv \rbrac = \j \bigg \lbrac \frac{1}{\varrho} \bigg \rbrac \bn_t, 
\quad \j \lbrac \bv \rbrac - \lbrac \BT \rbrac \bn_t = - \sigma H_{\Gamma_t} \bn_t,\\
& \lbrac \psi \rbrac + \cfrac{\j^2}{2} \bigg \lbrac \cfrac{1}{\varrho^2} \bigg \rbrac 
- \bigg \lbrac \cfrac{1}{\varrho} (\BT \bn_t \cdot \bn_t) \bigg \rbrac = 0,\\
& (\nabla \varrho_+) \cdot \bn_t \rvert_+ = 0,
\end{aligned}\right.
\end{align}
and the homogeneous Dirichlet boundary conditions on $\Gamma_+$ and $\Gamma_-$:
\begin{equation}
\label{cond-1.3}
\begin{aligned}
\bv_+ = 0 , \quad \nabla \varrho_+ \cdot \bn_+ = 0 \quad \text{on $\Gamma_+$},\qquad \bv_- = 0 \quad \text{on $\Gamma_-$,}
\end{aligned}
\end{equation}
and the initial conditions:
\begin{equation}
\label{cond-1.4}
\begin{aligned}
(\varrho_+, \bv_+)\rvert_{t = 0} = (\rho_{*+} + \rho_{0 +}, \bv_{0 +}) \quad \text{in $\Omega_{0 +}$}, \quad
\bv_-\rvert_{t = 0} = \bv_{0 -} \quad \text{in $\Omega_{0 -}$},
\end{aligned}
\end{equation}
where {\color{red} $\varrho_+$ and $\varrho_- := \rho_{* -}$} are the densities, $\bv_\pm$ the velocity fields, $\psi_\pm$ the Helmholtz free energy functions, and {\color{red} $\rho_{* \pm}$ are positive constants}, $H_{\Gamma_t}$ the $(N - 1)$-times mean curvature of $\Gamma_t$, $\sigma$ a positive constant describing the coefficient of the surface tension, $V_t$ the velocity of evolution of $\Gamma_t$ with respect to $\bn_t$, $\bv_{\Gamma_t}$ the interfacial velocity, $\bn_t$ the outer unit normal to $\Gamma_t$ pointed from $\Omega_{t+}$ to $\Omega_{t-}$, and $\bn_+$ the outer unit normal to $\Gamma_+$. Here, $\j = \varrho_+ (\bv_+ - \bv_\Gamma) \cdot \bn_t = \varrho_- (\bv_- - \bv_\Gamma)\cdot\bn_t$ is the phase flux and $\BT_\pm$ are the Stress tensors defined by
\begin{equation*}
\begin{split}
\BT_+ = & \mu_+ \bD(\bv_+) + (\nu_+ - \mu_+) \dv \bv_+ \bI - \fp_+ \bI + \bigg( \frac{\kappa_+}{2} \lvert \nabla \varrho_+ \rvert^2
+ \kappa_+ \varrho_+ \Delta \varrho_+ \bigg) \bI - \kappa_+ \nabla \varrho_+ \otimes \nabla \varrho_+,\\
\BT_- = & \mu_- \bD(\bv_-) - \fp_- \bI,
\end{split}
\end{equation*}
where $\fp_\pm$ are the pressure fields. Notice that if phase transitions occur on the moving interface $\Gamma_t$, the phase flux $\j$ should be taken arbitrary. Furthermore, the jump of a quantity $\bg(x, t)$ defined on $\dot \Omega_t$ across the interface $\Gamma_t$ is defined by 
\begin{equation*}
\lbrac \bg \rbrac (x_0) := \lim_{\delta \to 0 +} \left(\bg (x_0 + \delta \bn_t(x_0) ) - \bg (x_0 - \delta \bn_t (x_0)) \right)
\end{equation*}
for all $x_0 \in \Gamma_t$, where $\bn_t (x_0)$ is the outer unit normal to $\Gamma_t$ at $x_0$. In addition, we adopt the notations $\bg \rvert_\pm (x_0) = \lim_{\delta \to 0 +} \bg (x_0 \mp \delta \bn_t(x_0))$ for all $x_0 \in \Gamma_t$. The free boundary problem is said to be finding a family of hypersurfaces $\{\Gamma_t\}_{t \ge 0}$ and appropriately smooth functions $\varrho_+$, $\bu_+$, $\bu_-$, and $\fp_-$. We mention that the problem of finding a family of hypersurface $\{\Gamma_t\}_{t \ge 0}$ is equivalent to the problem of finding a family of $\{\Omega_{t +}\}_{t \ge 0}$ and $\{\Omega_{t-}\}_{t \ge 0}$. Notice that our system is thermodynamically consistent model in the sense of second law of thermodynamics, which was derived in the previous paper~\cite{Wat2017}. In particular, the condition $(\nabla \varrho_+)\cdot \bn_t \vert_+ = 0$ not only guarantees the generalized Gibbs-Thomson law and the Stefan law on the interface $\Gamma_t$ but also implies the \textit{interstitial working}: $(\kappa_+ \varrho_+ \dv \bv_+) \nabla \varrho_+$ vanishes in the normal direction of the interface $\Gamma_t$, see Watanabe~\cite{Wat2017}. For further remarks on the interstitial working, the readers may consult the paper by Dunn~\cite{Dun1986} or Dunn and Serrin~\cite{DS1985}. \par
In view of the \textit{Hanzawa transformation}, see Appendix, this paper mainly deals with the following fixed boundary system associated with~\eqref{eq-1.1},~\eqref{eq-1.2}, \eqref{cond-1.3}, and~\eqref{cond-1.4}:
\begin{align}
\label{eq-1.5}
\left\{\begin{aligned}
\pd_t \rho_+ + \rho_{*+} \dv \bu_+ & = f_M (\rho_+, \bu_+, h) &\enskip &\text{ in $\Omega_+ \times(0, T)$}, \\
\rho_{*-} \dv \bu_- = f_d (\bu_-, h) & = \rho_{*-} \dv \bF_d (\bu_-, h) &\enskip &\text{ in $\Omega_- \times (0, T)$}, \\
\rho_{*+}\pd_t\bu_+ -\DV \bT_+ (\gamma_1, \gamma_2, \gamma_3, \rho_+, \bu_+) & = \bff_+ (\rho_+, \bu_+, h) &\enskip &\text{ in $\Omega_+ \times (0, T)$}, \\
\rho_{*-}\pd_t\bu_--\DV \bT_- (\gamma_4, \bu_-, \pi_-) & = \bff_- (\bu_-, h) &\enskip &\text{ in $\Omega_- \times (0, T)$}, \\
\pd_t h - \frac{\langle \rho_{*-} \bu_-, \bn \rangle \rvert_- - \langle \rho_{*+} \bu_+, \bn \rangle \rvert_+}{\rho_{*-} - \rho_{*+}} 
& = d (\rho_+, \bu_+, \bu_-, h) &\enskip &\text{ on $\Gamma \times (0, T)$}, \\
\bB (\gamma_1, \gamma_2, \gamma_3, \gamma_4, \rho_+, \bu_+, \bu_-, \pi_-)
& = \bG (\rho_+, \bu_+, \bu_-, h), &\enskip &\text{ on $\Gamma \times (0, T)$}, \\
\bu_+ = 0, \qquad \langle \nabla \rho_+, \bn_+ \rangle & = 0
&\enskip &\text{ on $\Gamma_+ \times (0, T)$}, \\
\bu_- & = 0 &\enskip &\text{ on $\Gamma_- \times (0, T)$}, \\
(\rho_+, \bu_+, \bu_-, h) \rvert_{t = 0} & = (\rho_{0+}, \bu_{0+}, \bu_{0-}, h_0)
&\enskip &\text{ on $\Omega_+ \times \Omega_+ \times \Omega_- \times \Gamma$}.
\end{aligned}\right.
\end{align}
By abuse of notation, we let $\bT_+$ and $\bT_-$ be ``linearized'' stress tensors defined by
\begin{align*}
\bT_+ (\gamma_1, \gamma_2, \gamma_3, \bu_+, \rho_+)
& :=\gamma_1 \bD (\bu_+) + (\gamma_2 - \gamma_1) (\dv \bu_+) \bI + (- \gamma_{*+} + \rho_{*+} \gamma_3 \Delta) \rho_+ \bI, \\
\bT_- (\gamma_4, \bu_-, \pi_-) & := \gamma_4 \bD(\bu_-) - \pi_- \bI.
\end{align*}
and $\bB (\gamma_1, \gamma_2, \gamma_3, \gamma_4, \rho_+, \bu_+, \bu_-, \pi_-) = \bG (\rho_+, \bu_+, \bu_-, h)$ stands for the following interface conditions on $\Gamma \times (0, T)$:
\begin{align}
\label{def-B}
\left\{\begin{aligned}
\Pi_\bn (\gamma_4 \bD(\bu_-) \bn) \rvert_- - \Pi_\bn (\gamma_1 \bD(\bu_+) \bn) \rvert_+ & = g (\rho_+, \bu_+, \bu_-, h), \\
\langle \bT_- (\gamma_4, \bu_-, \pi_-) \bn, \bn \rangle \rvert_-
- \langle \bT_+ (\gamma_1, \gamma_2, \gamma_3, \bu_+, \rho_+) \bn, \bn \rangle \rvert_+ - \sigma (\langle \Delta_\Gamma \bn, \bn \rangle + \Delta_\Gamma) h & = f^+_B (\rho_+, \bu_+, \bu_-, h), \\
\frac{1}{\rho_{*-}} \langle \bT_- (\gamma_4, \bu_-, \pi_-) \bn, \bn \rangle \rvert_- - \frac{1}{\rho_{*+}} \langle \bT_+ (\gamma_1, \gamma_2, \gamma_3, \bu_+, \rho_+) \bn, \bn \rangle \rvert_+ - \frac{\gamma_{**}^+}{\rho_{*+}} \rho_+ \rvert_+ & = f^-_B (\rho_+, \bu_+, \bu_-, h),\\
\Pi_\bn \bu_- \rvert_- - \Pi_\bn \bu_+ \rvert_+ & = \bh(\bu_+, \bu_-, h), \\
\langle \nabla \rho_+, \bn \rangle \rvert_+ & = k_-(\rho_+, h),
\end{aligned}\right.
\end{align}
that is, $\bB$ and $\bG$ denote the left-hand and right-hand side of \eqref{def-B}, respectively. We denote the pull back of $(\varrho_+, \bv_+, \bv_-, \fp_-)$ by $(\rho_{* +} + \rho_+, \bu_+, \bu_-, \pi_{* -} + \pi_-)$, where $\pi_{* -}$ is a positive constant defined below, and $h$ is a height function --- the unknown moving interface $\Gamma_t$ is parametrized over the fixed hypersurface $\Gamma$ by means of a height function $h(x, t)$. We will give explicit formulas of right-hand members, the nonlinear terms, $f_M, \bff_+, \bff_-, \dots, k_-$ in Appendix. We emphasize that the operator $\bB$ is depending on $\pi_-$ but the function $\bG$ is \textit{independent} of $\pi_-$. Namely, the all nonlinear terms do not include the pressure term $\pi_-$, which plays an important role in formulating the system in a semigroup setting --- the pressure term and the divergence equation can be eliminated due to this fact. To simplify the notation, we have set $\Pi_\bn \bV = \bV - \langle \bV, \bn \rangle \bn$ for any $\bV \in \BR^N$, $\gamma_1 = \mu_+$, $\gamma_2 = \nu_+$, $\gamma_3 = \kappa_+$, $\gamma_4 = \mu_-$, $\gamma_{* +} = \fp'_+ (\rho_{* +})$, and $\gamma_{**}^+ = \pd_{\varrho_+} \psi_+ (\rho_{* +}, 0)$, where $\pd_{\varrho_+} = \pd / \pd \varrho_+$. Here, $\rho_{* -}$ is a positive constant denoting the reference mass density of $\Omega_-$, $\bn$ the outer unit normal to $\Gamma$ pointed from $\Omega_+$ into $\Omega_-$, and $\Delta_\Gamma$ the Laplace-Beltrami operator on $\Gamma$. Since we consider the system \eqref{eq-1.1} in isentropic case, we assume that the fluid occupied in $\Omega_{t +}$ is the barotropic viscous fluid, that is, the pressure $\fp_+$ is a function depending only on the density $\varrho_+$. In this paper, $\fp_+ $ and $\psi_+ $ are given smooth (at least $C^2$) functions. Note that the Helmholtz free energy $\psi_+$ depends not only on the density $\varrho_+$ but also on the square of gradient of density $\lvert \nabla \varrho_+ \rvert^2$ in the Navier-Stokes-Korteweg flow (cf. Dunn~\cite{Dun1986} or Dunn and Serrin~\cite{DS1985}). We will mainly devote to prove a unique solvability of the problem~\eqref{eq-1.5} in the \textit{maximal} $L_p - L_q$ \textit{regularity class}. \par
The thermodynamically consistent model, which takes the Gibbs-Thomson correction (cf. Pr{\"u}ss and Simonett~\cite{PS2016}) into account, including phase transitions has been studied in \cite{PS2012,PSW2014,PSSS2012,PS2016,SY2017} in the case of incompressible-incompressible two-phase flows. On the other hand, as far as the author knows, a free boundary problem of compressible-incompressible flows including phase transitions are few. In this direction, Shibata~\cite{Shi2016b} considered the linearized problem of the compressible-incompressible Navier-Stokes-Fourier system and proved the existence of solutions to the linearized system. However, it seems to be difficult to seek a unique solution to the original free boundary problem of Navier-Stokes-Fourier system with phase transitions based on his result. In fact, the regularity of density of a compressible fluid is not enough to solve the kinetic equation: $\bv_{\Gamma_t}\cdot\bn_t = \lbrac\varrho\bv\rbrac\cdot\bn_t / \lbrac\varrho\rbrac$, which causes the \textit{regularity loss} on the free boundary $\Gamma_t$. To overcome this difficulty, the author proposed the extended system of the Navier-Stokes-Fourier system by taking the Korteweg tensor into account~\cite{Wat2017}. Especially, in the previous paper~\cite{Wat2017}, the author showed the thermodynamically consistency of the extended model adopting the Navier-Stokes-Korteweg equations for compressible fluids and the Navier-Stokes equations for incompressible fluids, and formulate \eqref{eq-1.1},~\eqref{eq-1.2}, and \eqref{cond-1.3}. In addition, the existence of solution operator families for the corresponding generalized resolvent problem was shown in the case of whole space with flat interface. For further historical review or physical backgrounds of our model, the readers may consult the introduction in~\cite{Dun1986,DS1985,PS2012,PSSS2012,PS2016,Wat2017} and references therein. \par
To describe the compressible-incompressible two-phase flows, the Navier-Stokes-Allen-Cahn equations and Navier-Stokes-Cahn-Hilliard equations have been mainly considered, see Anderson, McFadden, and Wheeler~\cite{AMW1998} and Dreyer and Krauss~\cite{DK2010} for some examples. These models are said to be the \textit{diffuse-interface model}, in which we consider the interface as a non-zero thickness. As for these models, recently, Freist\"{u}hler and Kotschote~\cite{FK2017} proved that the Navier-Stokes-Allen-Cahn equations and Navier-Stokes-Cahn-Hilliard equations can be deduced to the Navier-Stokes-Korteweg equations. On the other hand, the \textit{sharp-interface} model regards the interface as a zero thickness. From a mathematical point of view, the sharp interface models are suitable models for a \textit{free boundary problem} because the position of the interface is a priori unknown. Hence, it seems to be reasonable to adopt the coupling system of the Navier-Stokes-Korteweg equations and the Navier-Stokes equations instead of the compressible-incompressible Navier-Stokes-Fourier system --- employ the Navier-Stokes-Korteweg equations instead of the ``usual'' compressible Navier-Stokes equations to describe the motion of viscous compressible fluid --- to investigate the compressible-incompressible two-phase flows. We emphasize that the Navier-Stokes-Korteweg equations were originally introduced to describe the structure of phase transition (cf. Dunn and Serrin~\cite{DS1985}). \par
The aim of this paper is to prove a unique existence of strong solution to \eqref{eq-1.5} with appropriate initial data in \textit{general} domains, which admits bounded or unbounded domains with either compact or non-compact free interface, see Section~\ref{def-uni} below for the definition. The problem \eqref{eq-1.5} reduced locally to so-called \textit{model problems} in a neighborhood of either an interior point or a boundary point by using the partition of unity associated with the domain $\Omega$ and standard localized methods. In the neighborhood of $\Gamma$, $\Gamma_+$, and $\Gamma_-$, \textcolor{red}{the system} \eqref{eq-1.5} is transformed to the problem in the whole space $x_N \in \BR$, the Stokes-Korteweg equations in the half space $x_N > 0$, and the Stokes equations in the half space $x_N < 0$, respectively. Corresponding problems have been studied by the author~\cite{Wat2017}, Saito~\cite{Sai2019}, and Shibata~\cite{Shi2014}, respectively. The essential part of the present paper is eliminating the pressure and divergence equation from the linearized system with the help of a unique solvability of the \textit{weak Dirichlet-Neumann problem}, which is slightly different from the two-phase Stokes equations case~\cite{MS2017}. \par
The rest of this paper is organized as follows: In the next section, we introduce some symbols and definitions needed in this paper. Sect.~\ref{sect-3} is concerned with the equivalence of \eqref{eq-1.5} and its \textit{reduced} equations. Sect.~\ref{sect-4} is devoted to prove the existence of $\CR$-bounded solution operators for the generalized resolvent problem. In Sect.~\ref{sect-7}, we prove a solvability of~\eqref{eq-1.5} in the maximal $L_p - L_q$ regularity class with the help of the Laplace transform and the operator-valued Fourier multiplier theorem. Finally, in Sect.~\ref{sect-8}, we prove the existence of a unique strong solution to~\eqref{eq-1.5} by the Banach fixed point theorem. In \textcolor{red}{the} appendix, we give a definition of the Hanzawa transformation and the explicit formulas of nonlinear terms $f_M, \bff_+, \bff_-, \dots, k_-$.

\section{Preliminaries}
\subsection{Notation}
We set $\pd_i = \pd / \pd x_i$ for $i = 1, \dots, N$ and $\pd_t = \pd / \pd t$. For any vector fields $\ba = {}^\top\!(a_1, \dots, a_N)$, the deformation tensor $\bD (\ba)$ is defined by $\bD (\ba) = \nabla \ba + {}^\top\!(\nabla \ba)$ whose $(j,k)$th components are $\pd_j a_k + \pd_k a_j$. In addition, for any vector fields $\ba = {}^\top\!(a_1, \dots, a_N)$ and $\bb = (b_1, \dots, b_N)$, the notation $\langle \ba, \bb \rangle = \ba \cdot \bb = \sum_{j = 1}^{N} a_j b_j$ denotes the inner product of $\ba$ and $\bb$. For any $N \times N$ matrix $\bM = (M_{jk})_{1 \le j, k \le N}$, the quantity $\DV \bM$ denotes the $N$-vectors with $j$-th component of $\sum_{k}^{N} \pd_k M_{jk}$. For a domain $U\subset \BC$ and Banach spaces $X$ and $Y$, the symbol $\Hol(U, \CL(X,Y))$ denotes the set of all $\CL (X, Y)$-valued holomorphic functions defined on $U$, where $\BC$ is the set of complex numbers and $\CL (X, Y)$ is the set of all bounded linear operators from $X$ to $Y$. For $\varepsilon \in (0, \pi/2)$ and $\lambda_0 > 0$, let $\Sigma_\varepsilon = \{\lambda \in \BC \backslash \{0\} \mid \lvert \arg \lambda \rvert \leq \pi - \varepsilon\}$ and $\Sigma_{\varepsilon, \lambda_0} = \{\lambda \in\Sigma_\varepsilon \mid \lvert \lambda \rvert \geq \lambda_0 \}$. The letter $C$ denotes a constant and $C_{a,b,c,\dots}$ denotes the constant depending on $a$, $b$, $c$, and so forth. In addition, the value of $C$ and $C_{a,b,c,\dots}$ may change from line to line.

\subsection{Function spaces}
For any domain $D \subset \BR^N$, let $B^s_{q, p} (D)$, $L_q (D)$, and $W^m_q (D)$ be inhomogeneous Besov spaces, Lebesgue spaces, and Sobolev spaces on $D$, respectively, and their norms are denoted by $\lVert \, \cdot \, \rVert_{B^s_{q, p} (D)}$, $\lVert \, \cdot \, \rVert_{L_q (D)}$, and $\lVert \, \cdot \, \rVert_{W^m_q (D)}$, respectively. For simplicity, we may write $B^s_{q, q} (D)$ as $W^s_q (D)$ and $L_q (D)$ as $W^0_q (D)$. We denote Lebesgue spaces and Sobolev spaces of $X$-valued functions defined on $\BR$ as $L_p (\BR, X)$ and $W^m_p (\BR, X)$ for $1 < p < \infty$, respectively, and their norm are denoted by $\lVert \, \cdot \, \rVert_{L_p (\BR, X)}$ and $\lVert \, \cdot \, \rVert_{W^m_p (\BR, X)}$, respectively. For a Banach space $X$ the $X$-valued Bessel potential spaces of order $1/2$ are defined by
\begin{align*}
H^{1/2}_p (\BR, X) = & \{f \in L_p (\BR, X) \mid \lVert f \rVert_{H^{1/2}_p (\BR, X)} < \infty \}, \quad \lVert f \rVert_{H^{1/2}_p (\BR, X)} = \lVert \CF^{-1} [(1+ \lvert \xi \rvert^2)^{-1/4} \CF[f]] \rVert_{L_p (\BR, X)},
\end{align*}
where $\CF$ and $\CF^{-1}$ stand the Fourier transform and \textcolor{red}{its inverse, respectively.}
For $\gamma \in \BR$ and $I \subset \BR$, we define function spaces with exponential weights as
\begin{align*}
L_{p, \gamma} (I, X) = & \{f \colon I \to X \mid e^{- \gamma t} f(t) \in L_p (I, X)\}, \\
W^1_{p, \gamma} (I, X) = & \{f \in L_{p, \gamma} (I, X) \mid {\color{red} e^{- \gamma t} \pd^k_t f(t)} \in L_p (I, X) \enskip \text{for $k = 0, 1$}\}.
\end{align*}
For simplicity of notation, in this paper we use the following symbols:
\begin{align*}
H^{1/2}_{p, \gamma} (\BR, X) & = \{f \in L_{p, \gamma} (\BR, X) \mid \lVert e^{- \gamma t} f \rVert_{H^{1/2}_p (\BR, X)} < \infty\}, \\
H^{1, 1/2}_{q, p} (D \times I) & = L_p (I, W^1_q (D)) \cap H^{1/2}_p (I, L_q (D)), \\
H^{1, 1/2}_{q, p, \gamma} (D \times I) & = L_{p, \gamma} (I, W^1_q (D)) \cap H^{1/2}_{p, \gamma} (I, L_q (D)), \\
W^{m + 2, 1}_{q, p} (D \times I) & = L_p (I, W^{m + 2}_q (D)) \cap W^1_p (I, W^m_q (D)), \\ 
W^{m + 2, 1}_{q, p, \gamma} (D \times I)  & = L_{p, \gamma} (I, W^{m + 2}_q (D)) \cap W^1_{p, \gamma} (I, W^m_q (D))
\end{align*}
for $1 < p,q < \infty$, $m = 0, 1$, $\gamma > 0$, a domain $D\subset\Omega$, and a nontrivial time interval $I \subset \BR$.

\subsection{Uniform domains}
\label{def-uni}
The present paper consider the system in the uniform $W^{4,3}_r$ and $W^{4,2}_r$ domains, which are analogous to usual uniform domains, defined by as follows: Let $\Omega_+$ be a connected domain surrounded by boundaries $\Gamma$ and $\Gamma_+$, while $\Omega_-$ be a connected domain surrounded by boundaries $\Gamma$ and $\Gamma_-$, where $\Gamma \neq \emptyset$ but $\Gamma_+ = \emptyset$ or $\Gamma_- = \emptyset$ is acceptable. For $1 < r < \infty$, the domains $\Omega_+$ and $\Omega_-$ are said to be uniform $W^{4, 3}_r$ and uniform $W^{4,2}_r$ domain, respectively, if there exist positive constants $\alpha$, $\beta$, and $K$ such that the following there assertions hold true: (1) For any $x_0 = (x_{01}, \dots, x_{0N}) \in \Gamma$ there exist a number $j$ and a $W^{4 - 1/r}_r$ function $h_\Gamma (x')$ defined on $B'_\alpha(x'_0)$ such that $\lVert h_\Gamma \rVert_{W^{4 - 1/r}_r (B'_\alpha (x'_0))} \leq K$ and
\begin{align*}
\Omega_\pm \cap B_\beta (x_0) = & \{x \in \BR^N \mid \pm x_j > h_\Gamma (x'_j) \enskip (x'_j \in B'_\alpha (x'_{0j}))\} \cap B_\beta (x_0), \\
\Gamma\cap B_\beta (x_0) = & \{x \in \BR^N \mid x_j = h_\Gamma (x'_j) \enskip (x'_j \in B'_\alpha (x'_{0j}))\} \cap B_\beta (x_0).
\end{align*}
(2) If $\Gamma_\pm \neq \emptyset$, for any $x_0 = (x_{01}, \dots, x_{0N}) \in \Gamma_\pm$, there exist a number $j$ and a $W^{(5 \pm 1) / 2 - 1 / r}_r$ function $h_{\Gamma_\pm} (x')$ defined on $B'_\alpha(x'_0)$ such that $\lVert h_{\Gamma_\pm} \rVert_{W^{(5 \pm 1) / 2 - 1 / r}_r (B'_\alpha (x'_0))} \leq K$ and
\begin{align*}
\Omega_\pm \cap B_\beta (x_0) = & \{x \in \BR^N \mid x_j > h_{\Gamma_\pm} (x'_j) \enskip (x'_j \in B'_\alpha (x'_{0j})) \} \cap B_\beta (x_0), \\
\Gamma_\pm \cap B_\beta(x_0) = & \{x \in \BR^N \mid x_j = h_{\Gamma_\pm} (x'_j) \enskip (x'_j \in B'_\alpha (x'_{0j})) \} \cap B_\beta (x_0).
\end{align*}
Here, we have set
\begin{gather*}
x'_j = (x_1, \dots, x_{j - 1}, x_{j + 1}, \dots, x_N), \quad x'_{0j} = (x_{01}, \dots, x_{0(j - 1)}, x_{0(j + 1)}, \dots, x_{0N}), \\
B'_\alpha (x'_{0j}) = \{x'_j \in \BR^{N - 1} \mid \lvert x'_j - x'_{0j} \rvert < \alpha \}, \quad B_\beta (x_0) = \{x \in \BR^N \mid \lvert x - x_0 \rvert < \beta\}.
\end{gather*}
Notice that if the boundary is compact, the uniformness is satisfied without any assumption. Let us give a typical situation of the domains:
\begin{enumerate}\renewcommand{\labelenumi}{(\roman{enumi})}
\item $\Omega_+$ and $\Omega_-$ are bounded domains with $\Gamma_+ \neq \emptyset$ and $\Gamma_- \neq \emptyset$;
\item $\Omega_+$ and $\Omega_-$ are a bounded and exterior domain, respectively, with $\Gamma_- = \emptyset$;
\item $\Omega_-$ is a bounded domain and $\Omega_+$ is its complement in $\BR^N$ assuming $\Gamma_+ = \emptyset$;
\item $\Omega_+$ and $\Omega_-$ are infinite layers with fixed boundaries $\Gamma_+$ and $\Gamma_-$ assuming that $\Gamma$ is non-compact.
\end{enumerate}
Here, the case (ii) and (iii) admit the case $\Gamma_+ = \emptyset$ and $\Gamma_- = \emptyset$, respectively. Our results mentioned below cover the \textit{all} cases of these domains.

\subsection{$\CR$-boundedness}
To establish a maximal regularity property, we have to introduce the concept of $\CR$-\textit{boundedness}, which is the stronger concept than the uniformly boundedness, due to the requirement of boundedness properties in vector-valued function spaces.
\begin{defi}
\label{R-bound}
For Banach spaces $X$ and $Y$, a family of operators $\mathcal T \subset \CL (X, Y)$ is said to be $\CR$-bounded if there exist constants $p \in [1, \infty)$ and $C \in (0, \infty)$ such that for any $m \in \BN$, $T_1, \dots, T_m \in \CT$, and $x_1, \dots, x_m \in X$ the inequality
\begin{align*}
\bigg(\int_0^1 \Big\lVert \sum_{n = 1}^m r_n (t) T_n x_n \Big\rVert_Y^p \dt \bigg)^{1/p} \leq C \bigg(\int_0^1 \Big\lVert \sum_{n = 1}^m r_n (t) x_n \Big\rVert_X^p \dt \bigg)^{1/p}
\end{align*}
holds, where $r_n (t) = \mathrm{sign}\, \sin (2^n \pi t)$ are the Rademacher functions on $[0,1]$. The infimum of $C$ such that the bound holds is said to be \textcolor{red}{the} $\CR$-bound of $\CT$ on $\CL (X, Y)$ denoted by $\CR_{\CL (X, Y)} (\CT)$.
\end{defi}

\begin{rema}
The constant $C$ in Definition \ref{R-bound} depends on $p$ in general. From the \textcolor{red}{Kahane--Khintchine} inequality {\color{red} (cf. Hyt{\"o}nen \textit{et al.}~\cite[Theorem~6.2.4]{HNVW})}, however, the definition of $\CR$-boundedness is independent of $p$. Namely, a family of operators $\CT$ is $\CR$-bounded for any $p \in [1, \infty)$ supposing that $\CT$ is $\CR$-bounded for some $p \in [1, \infty)$.
\end{rema}

The following properties \textcolor{red}{show} that $\CR$-bounds behave like norms (cf. Pr{\"u}ss and Simonett~\cite[Proposition 4.1.6]{PS2016}).
\begin{lemm}
\label{lem-4.2}
Let $1 < p, q < \infty$ and let $X$, $Y$, and $Z$ be Banach spaces. The following properties are valid:
\begin{enumerate}\renewcommand{\labelenumi}{(\arabic{enumi})}
\item Let $\CS$ and $\CT$ be $\CR$-bounded families in $\CL(X,Y)$. 
Then $\CS+\CT$ is also an $\CR$-bounded family in $\CL(X,Y)$ and hold the estimate:
\begin{align*}
\CR_{\CL(X,Y)}(\CS+\CT) \leq \CR_{\CL(X,Y)}(\CS)+\CR_{\CL(X,Y)}(\CT).
\end{align*}
\item Let $\CS$ and $\CT$ be $\CR$-bounded families in $\CL(X,Y)$ and $\CL(Y,Z)$, respectively. Then $\CT \CS = \{T S \mid S \in \CS, \enskip T \in \CT\}$ is an $\CR$-bounded family in $\CL(X, Z)$ satisfying the estimate
\begin{align*}
\CR_{\CL (X, Z)} (\CT \CS) \leq \CR_{\CL (X, Y)} (\CS) \CR_{\CL (Y, Z)} (\CT).
\end{align*}
\end{enumerate}
\end{lemm}

\section{Main result}
\subsection{Technical setup}
Before stating our main result, we first define the solenoidal space and the function space for the pressure. To this end, we introduce the \textit{weak Dirichlet-Neumann problem}.

\begin{defi}
\label{def-1.2}
For $1 < q < \infty$ and $q' = q / (q - 1)$ define the spaces $\hW^1_{q, \Gamma} (\Omega_-)$ and $W^1_{q, \Gamma} (\Omega_-)$ such that
\begin{align*}
\hW^1_{q,\Gamma} (\Omega_-) = \{\theta \in L_{q, \mathrm{loc}} (\Omega_-) \mid \nabla \theta \in L_q (\Omega_-)^N, \enskip \theta \rvert_\Gamma = 0 \}, \quad W^1_{q, \Gamma} (\Omega_-) = \{\theta \in W^1_q (\Omega_-) \mid \theta \rvert_\Gamma = 0 \}.
\end{align*}
The weak Dirichlet-Neumann problem is said to be uniquely solvable on $\wh W^1_{q, \Gamma} (\Omega_-)$ if the following assertion holds: For arbitrary $\bff \in L_q (\Omega_-)^N$, there exists a unique solution $\theta \in \wh W^1_{q, \Gamma} (\Omega_-)$ to the variational equation:
\begin{align}
\label{eq-1.9}
(\nabla \theta, \nabla \varphi)_{\Omega_-} = (\bff, \nabla \varphi)_{\Omega_-} \quad \text{for all $\varphi \in \wh W^1_{q', \Gamma} (\Omega_-)$}
\end{align}
satisfying $\lVert \nabla \theta \rVert_{L_q (\Omega_-)} \leq C_q \lVert \bff \rVert_{L_q (\Omega_-)}$ for some positive constant $C_q$ independent of $f$, $\theta$, and $\varphi$. If the weak Dirichlet-Neumann problem is uniquely solvable, we define a bounded linear operator $K_1 \in \CL (L_q (\Omega_-)^N, \wh W^1_{q, \Gamma} (\Omega_-))$ by $K_1 (\bff) = \theta$ with $\bff \in L_q (\Omega_-)$ and $\theta \in \wh W^1_{q, \Gamma} (\Omega_-)$ given in \eqref{eq-1.9}.
\end{defi}

\begin{rema}
\label{rem-1.3}
Let $\bff \in L_q (\Omega_-)^N$ and $\fp \in W^{1 - 1 / q}_q (\Gamma)$ for $1 < q < \infty$. Then there exists $p \in W^1_q (\Omega_-) + \wh W^1_{q, \Gamma} (\Omega_-)$ which is a unique solution to the variational equation:
\begin{align}
\label{eq-1.10}
(\nabla p, \nabla \theta)_{\Omega_-} = (\bff, \nabla \varphi)_{\Omega_-} \quad \text{ for any $\varphi \in \wh W^1_{q', \Gamma} (\Omega_-)$}
\end{align}
subject to $p = \fp$ on $\Gamma$ with $W^1_q (\Omega_-) + \wh W^1_{q, \Gamma} (\Omega_-) = \{p_1 + p_2 \mid p_1 \in W^1_q (\Omega_-), p_2 \in \wh W^1_{q, \Gamma} (\Omega_-) \}$. In fact, let $\bT_\Gamma \colon W^{1 - 1 / q}_q (\Gamma) \to W^1_q (\Omega_-)$ be a map such that for any $\fp \in W^{1 - 1 / q}_q (\Gamma)$, the function $\bT_\Gamma (\fp) \in W^1_q (\Omega_-)$ satisfies the conditions: $\bT_\Gamma (\fp) = \fp$ on $\Gamma$ and $\lVert \bT_\Gamma (\fp) \rVert_{W^1_q (\Omega_-)} \leq C \lVert \fp \rVert_{W^{1 - 1 / q}_q (\Gamma)}$ with some positive constant $C$ independent of $\fp$. If we define $p$ such that $p = \bT_\Gamma (\fp) + K_1 (\bff - \bT_\Gamma (\fp))$, where $p$ belongs to $W^1_q (\Omega_-) + \wh W^1_{q, \Gamma} (\Omega_-)$ with $\fp \in W^{1 - 1 / q}_q (\Gamma)$, then $p$ satisfies \eqref{eq-1.10} with the estimate
\begin{align*}
\lVert \nabla p \rVert_{L_q (\Omega_-)} \leq C_q \Big( \lVert \fp \rVert_{W^{1 - 1 / q}_q (\Gamma)} + \lVert \bff \rVert_{L_q (\Omega_-)} \Big).
\end{align*}
Especially, $W^1_q (\Omega_-) + \wh W^1_{q, \Gamma} (\Omega_-)$ is the space for the pressure field.		
\end{rema}

\begin{rema}
When $\Omega_-$ is a bounded domain, an exterior domain, a half space, and a bent half space, the weak Dirichlet-Neumann problem is uniquely solvable in $\wh W^1_{q, \Gamma} (\Omega_-)$ for \textit{arbitrary} $q \in (1, \infty)$. As for further examples for the domains such that the the weak Dirichlet-Neumann problem is uniquely solvable, see Shibata~\cite[Example 1.6]{Shi2013} and \cite[Example 1.8]{Shi2014}. Here, let us remark about the case when $\Omega_-$ is an exterior domain: Although the weak Dirichlet-Neumann problem
\begin{align*}
(\nabla \theta, \nabla \varphi)_{\Omega_-} = (\bff_-, \nabla \varphi)_{\Omega_-} \quad \text{for any $\varphi \in \CW^1_{q'} (\Omega_-)$}
\end{align*}
admits a unique solution $\theta \in \CW^1_q (\Omega_-)$ for any $\bff_- \in L_q (\Omega_-)^N$, where $\CW^1_q (\Omega_-)$ is the closure of $C^\infty_0 (\Omega_-)$ with respect to the norm $\lVert \nabla \cdot \rVert_{L_q (\Omega_-)}$, if and only if $N/(N - 1) < q < N$ if $N \geq 3$, and $q = 2$, if $N = 2$ (cf. Galdi~\cite[Theorem 8.4]{Gal2011}), we emphasize that the weak Dirichlet-Neumann problem is uniquely solvable in $\wh W^1_{q, \Gamma} (\Omega_-)$ for \textit{any} $q \in (1, \infty)$, see Pr{\"u}ss and Simonett~\cite[Theorem 7.4.3]{PS2016} and Shibata~\cite[Lemma 3.4]{Shi2018}. These differences arises form the fact that $C^\infty_0 (\Omega_-)$ is \textit{not} dense in $\wh W^1_{q, \Gamma} (\Omega_-)$ with the norm $\lVert \nabla \cdot \rVert_{L_q (\Omega_-)}$ for all $1 < q < \infty$, see Shibata~\cite[Appendix A]{Shi2018}.
\end{rema}

Under the assumption that the weak Dirichlet-Neumann problem is uniquely solvable on $\wh W^1_{q', \Gamma} (\Omega_-)$, we define solenoidal spaces $J_q(\Omega_-)$ by
\begin{align}
\label{def-sol}
J_q (\Omega_-) = \{\bff_- \in L_q (\Omega_-)^N \mid (\bff_-, \nabla \varphi)_{\Omega_-} = 0 \text{ for any $\varphi \in \wh W^1_{q', \Gamma} (\Omega_-)$} \},
\end{align}
where $1 < q < \infty$ and $\Omega_-$ is a uniform $W^{4,2}_r$ domain. We define $\CD\CI_q (\Omega_-)$ as data spaces for the divergence equation: $\rho_{*-} \dv \bu_- = f_d$ in $\Omega_-$ with $\bu_- \cdot \bn_- = 0$ on $\Gamma_-$ such that
\begin{align*}
\CD \CI_q (\Omega_-) = \Set{f_d \in W^1_q (\Omega_-) |
\begin{array}{l}
\text{ there exists a $\fg_d \in L_q (\Omega_-)^N$ such that} \\
\text{ $(f_d, \varphi)_{\Omega_-} = - \rho_{*-} (\fg_d, \nabla \varphi)_{\Omega_-}$ for all $\varphi \in W^1_{q', \Gamma} (\Omega_-)$}
\end{array}}.		
\end{align*}
Let $\CG (f_d) = \{\fh_d \in L_q (\Omega_-)^N \mid \dv \fg_d = \dv \fh_d \}$ and $\bF_d \in \CG (f_d)$. We then see that $\rho_{*-} \dv \bF_d = f_d$ in $\Omega_-$ and $\bF_d \cdot \bn_- = 0$ on $\Gamma_-$. Indeed, for any $\varphi \in C^\infty_0 (\Omega_-)$ we obtain
\begin{align*}
\rho_{*-} (\dv \bF_d, \varphi)_{\Omega_-} = - \rho_{*-} (\bF_d, \nabla \varphi)_{\Omega_-} = (f_d, \varphi)_{\Omega_-},
\end{align*}
which yields that $\rho_{*-} \dv \bF_d = f_d$ in $\Omega_-$. Furthermore, for arbitrary $\tvarphi \in C^1_0 (\Gamma_-)$ if we choose $\varphi \in W^1_{q', \Gamma} (\Omega_-)$ such that $\varphi \rvert_{\Gamma_-} = \tvarphi$, we observe that
\begin{align*}
\rho_{*-} (\bF_d \cdot \bn_-, \tvarphi)_{\Gamma_-} = \rho_{*-} (\dv \bF_d, \varphi)_{\Omega_-} + \rho_{*-} (\bF_d, \nabla \varphi)_{\Omega_-} = (f_d, \varphi)_{\Omega_-} - (f_d, \varphi)_{\Omega_-} = 0,
\end{align*}
which furnishes that $\bF_d \cdot \bn_-$ vanishes on $\Gamma_-$. Setting
\begin{align*}
\lVert f_d \rVert_{\CD \CI_q (\Omega_-)} = \lVert f_d \rVert_{W^1_q (\Omega_-)} + \inf_{\fh_d \in \CG (f_d)} \lVert \fh_d \rVert_{L_q (\Omega_-)}
\end{align*}
for $f_d \in \CD \CI_q (\Omega_-)$, we see that $\CD\CI_q(\Omega_-)$ is a Banach space with norm $\lVert \, \cdot \, \rVert_{\CD \CI_q (\Omega_-)}$. In the present paper, we define that $\bu_- \in W^1_q (\Omega_-)^N$ satisfies
\begin{align}
\label{cond-1.11}
\rho_{*-} \dv \bu_- = f_d \quad \text{ in $\Omega_-$,} \quad \bu_- \cdot \bn_- = 0 \quad \text{ on $\Gamma_-$}
\end{align}
if the identity
\begin{align}
\label{cond-1.12}
(\bu_-, \nabla \varphi)_{\Omega_-} = (\bF_d, \nabla \varphi)_{\Omega_-}
\end{align}
holds for any $\varphi \in \wh W^1_{q', \Gamma} (\Omega_-)$. We remark that \eqref{cond-1.11} and \eqref{cond-1.12} are \textit{not} equivalent, that is, the condition~\eqref{cond-1.11} does not imply the identity \eqref{cond-1.12} because $W^1_{q', \Gamma} (\Omega_-)$ is not dense in $\wh W^1_{q', \Gamma} (\Omega_-)$. From this fact, we shall introduce the following definition.
\begin{defi}
Let $1 < q < \infty$. For $\bu_-, \bF_d \in L_q (\Omega_-)^N$, we define that the identity $\dv \bu_- = \dv \bF_d$ holds in $\Omega_-$ if $\bu_- - \bF_d$ belongs to solenoidal spaces $J_q (\Omega_-)$.
\end{defi}
We finally introduce some technical assumptions for the coefficients of \textcolor{red}{the problem} \eqref{eq-1.5}.
\begin{assu}
\label{as-star}
The coefficients $\gamma_1 = \mu_+$, $\gamma_2 = \nu_+$, $\gamma_3 = \kappa_+$, $\gamma_4 = \mu_-$,
are real valued uniformly continuous functions defined in $\BR^N$.
We assume the following properties:
\begin{enumerate}
\item The coefficients $\rho_{*+}$ and $\rho_{*-}$ satisfy $\rho_{*+} \neq \rho_{*-}$.
\item There exists positive constants $\gamma^-_{k*}$ and $\gamma_{k*}^+$ $(k = 1, \dots, 4)$ such that
$\gamma^-_{k*} \leq \gamma_k (x) \leq \gamma^+_{k*}$ for any $x \in \BR^N$.
\item The coefficients $\gamma_k$ $(k = 1, \dots, 4)$ belong to $W^1_{r, \mathrm{loc}} (\BR^N)$
and $\lVert \nabla \gamma_k \rVert_{L_r (B_R)} \leq C_{r, R}$ with some positive constant
$C_{r, R}$ for any ball $B_R \subset \BR^N$ where $R > 0$ denotes the radius of ball $B_R$.
\item The coefficients $\rho_{*+}$ and $\gamma_k$ $(k = 1, \dots, 3)$ satisfy the condition
\begin{align*}
\bigg(\frac{\gamma_{1}(x) + \gamma_{2}(x)}{2 \rho_{*+}^2 \gamma_{3}(x)} \bigg)^2
\neq \frac{1}{\rho_{*+} \gamma_{3}(x)}, \qquad
\rho_{*+}^3 \gamma_{3}(x) \neq \gamma_{1}(x) \gamma_{2}(x).
\end{align*}
\end{enumerate}
\end{assu}

\begin{rema}
The conditions (2) and (3) are required to employ a localization argument in order to construct a solution to the generalized resolvent problem in general domains, see Maryani and Saito~\cite{MS2017} for an example on the two-phase Stokes equations case. On the other hand, Assumption \ref{as-star} (4) guarantees that we have the three roots with positive real parts different from each other, see \cite{Sai2017a,Sai2019} and \cite[Lemma 5.1]{Wat2017}. Notice that the assumption on $\rho_{* +}$, $\mu_{* +}$, $\nu_{* +}$, and $\kappa_{* +}$ in \cite[Theorem 1.2]{Wat2017} should be corrected as the one given by (4) in Assumption ~\ref{as-star} with $\gamma_1 = \mu_{* +}$, $\gamma_2 = \nu_{* +}$, and $\gamma_3 = \kappa_{* +}$. We remark that this assumption expect to be removed by using the similar argument due to Saito~\cite{Sai2019}.
\end{rema}

\subsection{Main results}
Setting $\dot \Omega = \Omega_+ \cup \Omega_-$, we shall state our main results.
\begin{theo}
\label{th-1.9}
Let $2 < p < \infty$ and $N < q < \infty$ with $2/p + N/q < 1$. In addition, let $T > 0$ and $N < r < \infty$ and $\max (q, q') \leq r$.
Suppose the following assertions:
\begin{enumerate}\renewcommand{\labelenumi}{(\alph{enumi})}
\item The domains $\Omega_+$ and $\Omega_-$ are uniform $W^{4,3}_r$ and $W^{4,2}_r$ domain, respectively.
\item The weak Dirichlet-Neumann problem is uniquely solvable on $\wh W^1_{q, \Gamma} (\Omega_-)$ and $\wh W^1_{q', \Gamma} (\Omega_-)$.
\item Assumption \ref{as-star} holds true..
\item The pressure field $\fp_+ (\rho_+)$ is $C^2$-function defined on $\rho_{* +} /3 \leq \rho_+ \leq 3 \rho_{* +}$ such that $0 \leq \fp'_+ (\rho_+) \leq \pi^*$ with some positive constant $\pi^*$ for any $\rho_{* +} / 3 \leq \rho_+ \leq 3 \rho_{* +}$.
\item The Helmholtz free energy $\psi_+ (\rho_+, \lvert \nabla \rho_+ \rvert^2)$ is $C^2$-function defined on $(\rho_{* +} / 3, 3 \rho_{* +}) \times [0, \infty)$ such that $0 \leq \pd_{\varrho_+} \psi_+ (\rho_+, \lvert \nabla \rho_+ \rvert^2) \leq \psi^*$ with some positive constant $\psi^*$ for any $\rho_{* +} / 3 \leq \rho_+ \leq 3 \rho_{* +}$.
\item There exist positive constants $\pi_{* \pm}$ such that
\begin{align}
\label{cond-GT}
\psi_- (\rho_{* -}) - \psi_+ (\rho_{* +}, 0) = \frac{\pi_{* +}}{\rho_{* +}} - \frac{\pi_{* -}}{\rho_{* -}}, \quad \pi_{* -} - \pi_{* +} = \sigma H_\Gamma,
\end{align}
which represents the Gibbs-Thomson condition and the Young-Laplace law, respectively.
\item The initial data
\begin{align*}
(\rho_{0 +}, \bu_{0 +}, \bu_{0 -}, h_0) \in B^{3 - 2/p}_{q, p} (\Omega_+) \times B^{2(1 - 1/p)}_{q, p} (\Omega_+) \times B^{2(1 - 1/p)}_{q, p} (\Omega_-) \times B^{3 - 1/p - 1/q}_{q, p} (\Gamma)
\end{align*}
satisfies the compatibility conditions:
\begin{align*}
\left\{\begin{aligned}
\dv \bu_{0 -} = f_{d 0} & = \dv \bF_{d 0} &\quad & \text{ in $\Omega_-$}, \\
\Pi_\bn (\gamma_4 \bD(\bu_{0 -}) \bn) \vert_- - \Pi_\bn (\gamma_1 \bD(\bu_{0 +}) \bn) \vert_+ & = g_0 &\quad &\text{ on $\Gamma$}, \\
\Pi_\bn \bu_{0 -} \vert_- - \Pi_\bn \bu_{0 +} \vert_+ & = \bh_0 &\quad &\text{ on $\Gamma$}, \\
\langle \nabla \rho_{0 +}, \bn \rangle \vert_+ & = k_{0 -} &\quad &\text{ on $\Gamma$}, \\
\langle \nabla \rho_{0 +}, \bn_+ \rangle = 0, \quad \bu_{0 +} & = 0 &\quad &\text{ on $\Gamma_+$}, \\
\bu_{0 -} & = 0 &\quad &\text{ on $\Gamma_-$}
\end{aligned}\right. 
\end{align*}
with $\bu_{0 -} - \bF_{d 0} \in J_q (\Omega_-)$, where we have set $f_{d 0} = f_d (\bu_{0-}, h_0)$, $\bF_{d 0} = \bF_d (\bu_{0-}, h_0)$, $g_0 = g (\rho_{0 +}, \bu_{0+}, \bu_{0-}, h_0)$, $\bh_0 = \bh (\rho_{0 +}, \bu_{0+}, \bu_{0-}, h_0)$, and $k_{0 -} = k_- (\rho_{0 +}, h_0)$.
\end{enumerate}
Then there exists a positive constant $\varepsilon_T$ depending on $T$ such that if the initial data satisfy the smallness condition:
\begin{align*}
\lVert \rho_{0 +} \rVert_{B^{3 - 2/p}_{q, p} (\Omega_+)} + \lVert \bu_{0 +} \rVert_{B^{2(1 - 1/p)}_{q, p} (\Omega_+)} + \lVert \bu_{0 -} \rVert_{B^{2(1 - 1/p)}_{q, p} (\Omega_-)} + \lVert h_0 \rVert_{B^{3 - 1/p - 1/q}_{q, p} (\Gamma)} \leq \varepsilon_T,
\end{align*}
there exists a unique solution $(\rho_+, \bu_+, \bu_-, \pi_-, h)$ to the system \eqref{eq-1.5} with
\begin{alignat*}4
\rho_+ & \in L_p ((0,T), W^3_q (\Omega_+)) \cap W^1_p ((0,T), W^1_q (\Omega_+)), &\enskip \bu_\pm & \in L_p ((0,T), W^2_q (\Omega_\pm)^N) \cap W^1_p ((0,T), L_q (\Omega_\pm)^N),\\
\pi_- & \in L_p ((0,T), W^1_q (\Omega_-) + \wh W^1_{q, \Gamma} (\Omega_-)), &\enskip h & \in L_p ((0,T), W^{3 - 1/q}_q (\Gamma)) \cap W^1_p ((0,T), W^{2 - 1/q}_q (\Gamma))
\end{alignat*}
satisfying the estimate $\BI_{p,q} (\rho_+, \bu_+, \bu_-, \pi_-, h, 0; (0,T)) \leq \varepsilon_T$. Here and in the following, for $\delta \in [0, \infty)$ and $0 \leq a < b \leq \infty$ we set
\begin{align*}
& \BI_{p,q} (\rho_+, \bu_+, \bu_-, \pi_-, h, \delta; (a,b)) \\
& := \lVert e^{-\delta t} \pd_t \rho_+ \rVert_{L_p ((a, b), W^1_q (\Omega_+))} + \lVert e^{-\delta t} \rho_+ \rVert_{L_p ((a, b), W^3_q (\Omega_+))} + \lVert e^{-\delta t} \rho_+ \rVert_{L_\infty ((a, b), B^{3 - 2/p}_{q, p} (\Omega_+))} \\
& \quad + \sum_{\ell = \pm} \left(\lVert e^{-\delta t} \pd_t \bu_\ell \rVert_{L_p ((a, b), L_q (\Omega_\ell))} + \lVert e^{-\delta t} \bu_\ell \rVert_{L_p ((a, b), W^2_q (\Omega_\ell))} + \lVert e^{-\delta t} \bu_\ell \rVert_{L_\infty ((a,b), B^{2 (1 - 1/p)}_{q,p} (\Omega_\ell))} \right) \\
& \quad + \lVert e^{-\delta t} \nabla \pi_- \rVert_{L_p ((a, b), L_q (\Omega_-))} + \lVert e^{-\delta t} \pd_t h \rVert_{L_p ((a, b), W^{2 - 1/p}_q (\Gamma))} + \lVert e^{-\delta t} h \rVert_{L_p ((a, b), W^{3 - 1/p}_q (\Gamma))} \\
& \quad + \lVert e^{-\delta t} h \rVert_{L_\infty ((a, b), B^{3 - 1/p - 1/q}_{q, p} (\Gamma))}.
\end{align*}
\end{theo}

\begin{rema}
(1) {\color{red} From the trace method of real interpolation, for $1 < p < \infty$ we have}
\begin{align}
\label{emb-BUC}
W^1_p ((0, T), X_1) \cap L_p ((0, T), X_2) \hookrightarrow \mathrm{BUC}([0, T), (X_1, X_2)_{1 - 1/p, p}),
\end{align}
where $X_1$ and $X_2$ are Banach spaces such that $X_2$ is dense subset of $X_1$ and $\mathrm{BUC} ([0, T), (X_1, X_2)_{1 - 1/p, p})$ denotes the set of all $(X_1, X_2)_{1 - 1/p, p}$-valued uniformly continuous and bounded functions on $[0, T)$ (cf. Amann~\cite[Chapter I\hspace{-.1em}I\hspace{-.1em}I, Theorem 4.10.2]{Ama1995}). Hence, a quadruple $(\rho_+, \bu_+, \bu_-, h)$ is continuous with respect to initial data $(\rho_{0 +}, \bu_{0 +}, \bu_{0 -}, h_0)$, so that the system is a locally well-posed. \\
(2) For given $T >0$, we can find a family of hypersurfaces $\{\Gamma_t\}_{t \ge 0}$ and see that $(\varrho_+, \bv_+, \bv_-, \fp_-)$ is a unique solution to the free boundary problem \eqref{eq-1.1}-\eqref{cond-1.4} for any $t \in (0, T)$ because the Hanzawa transformation is injective, see Appendix below.
\end{rema}

\section{Reduced problem}
\label{sect-3}
\subsection{Eliminating the pressure term and the divergence equation}
\label{sect-3.1}
We first consider the left-hand side of \eqref{eq-1.5} without the lower order terms $\gamma_{*+} \rho_+$, $\langle \Delta_\Gamma \bn, \bn \rangle h$, and $\gamma_{**}^+ \rho_+$. To this end, we define
\begin{align*}
\bT_{0 +} (\gamma_1, \gamma_2, \gamma_3, \bu_+, \rho_+) = \gamma_1 \bD (\bu_+) + (\gamma_2 - \gamma_1) (\dv \bu_+) \bI + \rho_{*+} \gamma_3 \Delta \rho_+ \bI.
\end{align*}
We now decompose the interface condition:
\begin{align*}
\bT_- (\gamma_4, \bu_-, \pi_-) \vert_- - \{\bT_{0 +} (\gamma_1, \gamma_2, \gamma_3, \bu_+, \rho_+)\}\vert_+ - \sigma \Delta_\Gamma h & = f^+_B, \\
\frac{1}{\rho_{*-}} \langle \bT_- (\gamma_4, \bu_-, \pi_-) \bn, \bn \rangle \Big\vert_- - \frac{1}{\rho_{*+}} \langle \bT_{0 +} (\gamma_1, \gamma_2, \gamma_3, \bu_+, \rho_+) \bn, \bn \rangle \Big\vert_+ & = f^-_B
\end{align*}
into
\begin{align*}
\langle \bT_- (\gamma_4, \bu_-, \pi_-) \bn, \bn \rangle \Big\vert_- - \frac{\rho_{*-} \sigma}{\rho_{*-} - \rho_{*+}} \Delta_\Gamma h & = \frac{\rho_{*-} (f^+_B - \rho_{*+} f^-_B)}{\rho_{*-} - \rho_{*+}} =: g_-, \\
\langle \bT_{0 +} (\gamma_1, \gamma_2, \gamma_3, \bu_+, \rho_+) \bn, \bn \rangle \Big\vert_+ - \frac{\rho_{*+} \sigma}{\rho_{* -} - \rho_{* +}} \Delta_\Gamma h & = \frac{\rho_{*+}(f^+_B - \rho_{*-} f^-_B)}{\rho_{* -} - \rho_{* +}} =: g_+.
\end{align*}
We then consider the following linear problem:
\begin{align}
\label{eq-1.6}
\left\{\begin{aligned}
\pd_t \rho_+ + \rho_{*+} \dv \bu_+ & = f_M &\enskip &\text{ in $\Omega_+ \times (0,\infty)$,} \\
\rho_{*-} \dv \bu_- = f_d = \rho_{*-} &\dv \bF_d  &\enskip &\text{ in $\Omega_- \times (0,\infty)$}, \\
\rho_{*+}\pd_t \bu_+ -\DV \bT_{0 +} (\gamma_1, \gamma_2, \gamma_3, \bu_+, \rho_+) & = \bff_+ &\enskip &\text{ in $\Omega_+ \times (0,\infty)$}, \\
\rho_{*-} \pd_t \bu_- - \DV \bT_- (\gamma_4, \bu_-, \pi_-) & = \bff_- &\enskip &\text{ in $\Omega_- \times (0,\infty)$}, \\
\pd_t h - \frac{\langle \rho_{*-} \bu_-, \bn \rangle \vert_- - \langle\rho_{*+} \bu_+, \bn \rangle\vert_+} {\rho_{*-}-\rho_{*+}} & = d &\enskip &\text{ on $\Gamma \times (0,\infty)$}, \\
\bB_0 (\gamma_1, \gamma_2, \gamma_3, \gamma_4, \rho_+, \bu_+, \bu_-) & = \bG_0 &\enskip &\text{ on $\Gamma \times (0,\infty)$}, \\
\langle \bT_- (\gamma_4, \bu_-, \pi_-) \bn, \bn \rangle \Big\vert_- - \frac{\rho_{*-} \sigma}{\rho_{*-} - \rho_{*+}} \Delta_\Gamma h & = g_- &\enskip &\text{ on $\Gamma \times (0,\infty)$}, \\
\color{red}
\bu_+ = 0, \qquad \langle \nabla \rho_+, \bn_+ \rangle &
\color{red}
= k_+ &\enskip &
\color{red}
\text{ on $\Gamma_+ \times (0,\infty)$}, \\
\color{red}
\bu_- & 
\color{red}
= 0 &\enskip &
\color{red}
\text{ on $\Gamma_- \times (0,\infty)$}
\end{aligned}\right.
\end{align}
with $(\rho_+, \bu_+, \bu_-, h) \vert_{t = 0} = (\rho_{0 +}, \bu_+, \bu_{0-}, h_0)$, where $\bB_0 (\gamma_1, \gamma_2, \gamma_3, \gamma_4, \rho_+, \bu_+, \bu_-) = \bG_0$ stands the conditions:
\begin{align*}
\left\{\begin{aligned}
\Pi_\bn (\gamma_4 \bD(\bu_-) \bn) \vert_- - \Pi_\bn (\gamma_1 \bD(\bu_+) \bn) \vert_+ & = g &\enskip &\text{ on $\Gamma \times (0, \infty)$}, \\
\langle \bT_{0 +} (\gamma_1, \gamma_2, \gamma_3, \bu_+, \rho_+) \bn, \bn \rangle \Big\vert_+ - \frac{\rho_{*+}\sigma}{\rho_{*-} - \rho_{*+}} \Delta_\Gamma h
& = g_+ &\enskip &\text{ on $\Gamma \times (0,\infty)$}, \\
\Pi_\bn \bu_- \vert_- - \Pi_\bn \bu_+ \vert_+ & = \bh &\enskip &\text{ on $\Gamma \times (0,\infty)$}, \\
\langle \nabla \rho_+, \bn \rangle\vert_+ & = k_- &\enskip &\text{ on $\Gamma \times (0,\infty)$}, \\
\end{aligned}\right.
\end{align*}
Here, the right-hand members are given functions at this stage, where $k_+$ is an additional function --- we will take $k_+$ as zero if we solve the nonlinear problem. Notice that $\bB_0$ and $\bG_0$ are independent of $\pi_-$. To prove a solvability of \eqref{eq-1.6}, we consider the following resolvent problem:
\begin{align}
\label{eq-4.1*}
\left\{\begin{aligned}
\lambda \wh \rho_+ + \rho_{*+} \dv \wh \bu_+ & = \wh f_M &\enskip &\text{ in $\Omega_+$,} \\
\rho_{*-} \dv \wh \bu_- = \wh f_d = \rho_{*-} & \dv \wh \bF_d  &\enskip &\text{ in $\Omega_-$}, \\
\rho_{*+} \lambda \wh \bu_+ - \DV \bT_{0 +} (\gamma_1, \gamma_2, \gamma_3, \wh \bu_+, \wh \rho_+) & = \wh \bff_+ &\enskip &\text{ in $\Omega_+$}, \\
\rho_{*-} \lambda \wh \bu_- - \DV \bT_- (\gamma_4, \wh \bu_-, \wh \pi_-) & = \wh \bff_- &\enskip &\text{ in $\Omega_-$}, \\
\lambda \wh h - \frac{\langle \rho_{*-} \wh \bu_-, \bn \rangle \vert_- - \langle \rho_{* +} \wh \bu_+, \bn \rangle\vert_+} {\rho_{* -} - \rho_{* +}} & = \wh d &\enskip &\text{ on $\Gamma$}, \\
\bB_0 (\gamma_1, \gamma_2, \gamma_3, \gamma_4, \wh \rho_+, \wh \bu_+, \wh \bu_-) & = \wh \bG_0 &\enskip &\text{ on $\Gamma$}, \\
\langle \bT_- (\gamma_4, \wh \bu_-, \wh \pi_-) \bn, \bn \rangle \Big\vert_- - \frac{\rho_{*-} \sigma}{\rho_{*-} - \rho_{*+}} \Delta_\Gamma \wh h & = \wh g_- &\enskip &\text{ on $\Gamma$}, \\
\color{red}
\wh \bu_+ = 0, \qquad \langle \nabla \rho_+, \bn_+ \rangle &
\color{red}
= \wh k_+ &\enskip &
\color{red}
\text{ on $\Gamma_+$}, \\
\color{red}
\wh \bu_- & 
\color{red}
= 0 &\enskip &
\color{red}
\text{ on $\Gamma_-$}, 
\end{aligned}\right.
\end{align}
where $\bB_0 (\gamma_1, \gamma_2, \gamma_3, \gamma_4, \wh \rho_+, \wh \bu_+, \wh \bu_-) = \wh \bG_0$ denotes the following conditions:
\begin{align*}
\left\{\begin{aligned}
\Pi_\bn (\gamma_4 \bD(\wh \bu_-) \bn) \vert_- - \Pi_\bn (\gamma_1 \bD(\wh \bu_+) \bn) \vert_+ & = \wh g &\enskip &\text{ on $\Gamma$}, \\
\langle \bT_{0 +} (\gamma_1, \gamma_2, \gamma_3, \wh \bu_+, \wh \rho_+) \bn, \bn \rangle \Big\vert_+ - \frac{\rho_{*+}\sigma}{\rho_{*-} - \rho_{*+}} \Delta_\Gamma \wh h & = \wh g_+ &\enskip &\text{ on $\Gamma$}, \\
\Pi_\bn \wh \bu_- \vert_- - \Pi_\bn \wh \bu_+ \vert_+ & = \wh \bh &\enskip &\text{ on $\Gamma$}, \\
\langle \nabla \wh \rho_+, \bn \rangle\vert_+ & = \wh k_- &\enskip &\text{ on $\Gamma$}.
\end{aligned}\right.
\end{align*}
To formulate the problem \eqref{eq-1.6} in the semigroup setting, we have to eliminate the pressure term $\pi_-$ and the divergence equation: $\rho_{*-} \dv \wh \bu_- = \wh f_d$ in \eqref{eq-4.1*}, that is, we deduce the reduced equations equivalent to \eqref{eq-4.1*}. To this end, we follow the idea due to Shibata~\cite{Shi2013,Shi2016a}. \par
Let $\CK_1 (\wh \bu_-)$ be a unique solution to the variational problem
\begin{align}
\label{eq-3.1}
(\nabla \CK_1(\wh \bu_-), \nabla \varphi)_{\Omega_-} = (\DV (\gamma_4 \bD(\wh \bu_-)) - \rho_{*-} \nabla \dv \wh \bu_-, \nabla \varphi)_{\Omega_-} \quad \text{ for any $\varphi \in \wh W^1_{q', \Gamma} (\Omega_-)$,}
\end{align}
subject to $\CK_1 (\wh \bu_-) = \gamma_4 \langle \bD (\wh \bu_-) \bn, \bn \rangle - \rho_{*-} \dv \wh \bu_-$ on $\Gamma$, while $\CK_2 (\wh h)$ is a unique solution to the following variational problem:
\begin{align}
\label{eq-3.2}
(\nabla \CK_2 (\wh h), \nabla \varphi)_{\Omega_-} = 0 \quad \text{ for any $\varphi \in \wh W^1_{q', \Gamma} (\Omega_-)$,}
\end{align}
subject to $\CK_2 (\wh h) = - (\rho_{*-} - \rho_{*+})^{-1} \rho_{*-} \sigma \Delta_\Gamma \wh h$ on $\Gamma$. As we mentioned in Remark \ref{rem-1.3}, the functions $\CK_1 (\wh \bu_-)$ and $\CK_2 (\wh h)$ can be defined by
\begin{align*}
\CK_1 (\wh \bu_-) = K^1_\Gamma (\wh \bu_-) + K_1(\DV(\gamma_4 \bD (\wh \bu_-)) - \rho_{*-} \nabla \dv \wh \bu_- - \nabla K^1_\Gamma (\wh \bu_-)), \quad \CK_2 (\wh h) = - K^2_\Gamma (\wh h) + K_1 (\nabla K^2_\Gamma (\wh h))
\end{align*}
with $K^1_\Gamma (\wh \bu_-) = \bT_\Gamma (\gamma_4 \langle \bD(\wh \bu_-) \bn, \bn \rangle - \rho_{*-} \dv \wh \bu_-)$ and  $K^2_\Gamma (\wh h) = \bT_\Gamma(- \CK_2 (\wh h))$, respectively. We easily see that $\CK_1 (\wh \bu_-)$ and $\CK_2 (\wh h)$ belong to $W^1_q (\Omega_-) + \wh W^1_{q, \Gamma} (\Omega_-)$ satisfying the estimates
\begin{align*}
\lVert \nabla \CK_1 (\wh \bu_-) \rVert_{L_q (\Omega_-)} \leq C \lVert \nabla \wh \bu_- \rVert_{W^1_q (\Omega_-)}, \quad \lVert \nabla \CK_2 (\wh h) \rVert_{L_q (\Omega_-)} \leq C \lVert \wh h \rVert_{W^{3 - 1/q}_q (\Gamma)},
\end{align*}
respectively. We then have the ``reduced'' system:
\begin{align}
\label{eq-3.8}
\left\{\begin{aligned}
\lambda \wh \rho_+ + \rho_{*+} \dv \wh \bu_+ & = \wh f_M &\enskip &\text{ in $\Omega_+$,} \\
\rho_{*+} \lambda \wh \bu_+ - \DV \bT_{0 +} (\gamma_1, \gamma_2, \gamma_3, \wh \bu_+, \wh \rho_+) & = \wh \bff_+ &\enskip &\text{ in $\Omega_+$}, \\
\rho_{*-} \lambda \wh \bu_- - \DV \bT_- (\gamma_4, \wh \bu_-, \CK_1 (\wh \bu_-) + \CK_2 (\wh h)) & = \wh \bff_- &\enskip &\text{ in $\Omega_-$}, \\
\lambda \wh h - \frac{\langle \rho_{* -} \wh \bu_-, \bn \rangle \vert_- - \langle \rho_{* +} \wh \bu_+, \bn \rangle \vert_+} {\rho_{* -} - \rho_{* +}} & = \wh d &\enskip &\text{ on $\Gamma$}, \\
\bB_0 (\gamma_1, \gamma_2, \gamma_3, \gamma_4, \wh \rho_+, \wh \bu_+, \wh \bu_-) & = \wh \bG_0 &\enskip &\text{ on $\Gamma$}, \\
\rho_{*-} (\dv \wh \bu_-) \vert_- & = \wh g_- &\enskip &\text{ on $\Gamma$}, \\
\color{red}
\wh \bu_+ = 0, \qquad \langle \nabla \rho_+, \bn_+ \rangle &
\color{red}
= \wh k_+ &\enskip &
\color{red}
\text{ on $\Gamma_+$}, \\
\color{red}
\wh \bu_- & 
\color{red}
= 0 &\enskip &
\color{red}
\text{ on $\Gamma_-$}.
\end{aligned}\right.
\end{align} \par 
In the following, we shall prove the equivalence between \eqref{eq-4.1*} and \eqref{eq-3.8}. Given $\wh f_M \in W^1_q (\Omega_+)$, $\wh \bff_\pm \in L_q (\Omega_\pm)^N$, $\wh d \in W^2_q (\dot\Omega)$, $\wh g, \wh g_\pm \in W^1_q (\dot \Omega)$, $\wh \bh \in W^2_q (\dot \Omega)^N$, and $\wh k_\pm \in W^2_q(\Omega_+)$. Let $\wh f_d \in W^1_q(\Omega_-)$ be a unique solution to the auxiliary problem
\begin{align*}
\lambda (\wh f_d, \varphi)_{\Omega_-} + (\nabla \wh f_d, \nabla \varphi)_{\Omega_-} = - (\rho_{*-}^{-1} \wh \bff_-, \nabla \varphi)_{\Omega_-}
\end{align*}
for any $\varphi\in W^1_{q',\Gamma} (\Omega_-)$ subject to $\wh f_d = \wh g_-$ on $\Gamma$. Recall that a unique existence of $\wh f_d$ is guaranteed for suitably large $\lambda > 0$ (cf. Shibata~\cite[Sect.~9.6.2]{Shi2016a}). In this case, we see that
\begin{align}
\label{eq-3.9}
\wh \bF_d = \lambda^{-1} (\nabla \wh f_d + \rho_{*-}^{-1} \wh \bff_-).
\end{align}
Let $\wh \rho_+ \in W^3_q (\Omega_+)$, $\wh \bu_\pm \in W^2_q (\Omega_\pm)$, $\wh \pi_- \in W^1_q (\Omega_-) + \wh W^1_{q, \Gamma} (\Omega_-)$, and $\wh h \in W^{3 - 1/q}_q (\Gamma)$ be unique solutions of~\eqref{eq-4.1*} with \eqref{eq-3.9}. From \eqref{cond-1.11}, \eqref{cond-1.12}, and \eqref{eq-3.9}, we have
\begin{align}
\label{eq-3.10}
\begin{aligned}
\rho_{*-} \dv \wh \bu_- & = \wh f_d \in W^1_q (\Omega_-), \\
(\wh \bu_-, \nabla \varphi)_{\Omega_-} & = \lambda^{-1} (\nabla \wh f_d + \rho_{*-}^{-1} \wh \bff_-, \nabla \wh \varphi)_{\Omega_-} &\quad &\text{ for any $\varphi \in \wh W^1_{q', \Gamma} (\Omega_-)$.}
\end{aligned}
\end{align}
On the other hand, by \eqref{eq-3.1} and \eqref{eq-3.2}, for any $\varphi \in \wh W^1_{q', \Gamma} (\Omega_-)$ we observe that
\begin{align*}
& (\rho_{*-}^{-1} \wh \bff_-, \nabla \varphi)_{\Omega_-} \\
& = (\lambda \wh \bu_- - \rho_{*-}^{-1} \DV( \mu_- \bD(\wh \bu) - \wh \pi_- \bI), \nabla \varphi)_{\Omega_-} \\
& = \lambda(\wh \bu_-, \nabla \varphi)_{\Omega_-} - (\nabla \dv \wh \bu_-, \nabla \varphi)_{\Omega_-} - (\rho_{*-}^{-1} \DV (\gamma_4 \bD(\wh \bu_-)) - \rho_{*-} \nabla \dv \wh \bu_-, \nabla \varphi)_{\Omega_-} + (\rho_{*-}^{-1} \nabla \wh \pi_-, \nabla \varphi)_{\Omega_-} \\
& = \lambda(\wh \bu_-, \nabla \varphi)_{\Omega_-} - (\nabla \dv \wh \bu_-, \nabla \wh \varphi)_{\Omega_-} + (\rho_{*-}^{-1} \nabla (\wh \pi_- - (\CK_1 (\wh \bu_-) + \CK_2 (\wh h))), \nabla \varphi)_{\Omega_-},
\end{align*}
which, combined with \eqref{eq-3.10}, furnishes that $(\nabla (\wh \pi_- - (\CK_1 (\wh \bu_-) + \CK_2 (\wh h))), \nabla \varphi)_{\Omega_-} = 0$ for any $\varphi \in \wh W^1_{q', \Gamma} (\Omega_-)$. Furthermore, since $\rho_{*-} \dv \wh \bu_- = f_d$ in $\Omega_-$ with $\wh f_d = \wh g_-$ on $\Gamma$, by \eqref{eq-3.1} and \eqref{eq-3.2}, we see that
\begin{align*}
\pi_- - (\CK_1 (\wh \bu_-) + \CK_2 (\wh h)) & = \mu_- \langle \bD (\wh \bu_-) \bn, \bn \rangle - \frac{\rho_{*-} \sigma}{\rho_{*-} - \rho_{*+}} \Delta_\Gamma \wh h - \wh g_- - \CK_1 (\wh \bu_-) - \CK_2 (\wh h) \\
& = - \wh g_- + \rho_{*-} \dv \wh \bu_- \vert_- = - \wh f_d + \wh f_d = 0
\end{align*}
holds on $\Gamma$. Hence, a uniqueness of solutions to \eqref{eq-4.1*} implies that $\wh \pi_- = \CK_1 (\wh \bu_-) + \CK_2 (\wh h)$, which yields that $\wh \bu_-$ and $\wh h$ satisfy \eqref{eq-3.8}. \par
Conversely, we assume a unique solvability of \eqref{eq-3.8}. Given $\wh f_d \in \CD \CI_q (\Omega_-)$, we assume that
\begin{align}
\label{eq-3.11}
\wh g_- = 0 \quad \text{ on $\Gamma$}, \quad (\rho_{*-}^{-1} \wh \bff_-, \nabla \varphi)_{\Omega_-} = 0 \quad \text{ for any $\varphi \in \wh W^1_{q', \Gamma} (\Omega_-)$.}
\end{align}
Let $\CK (\lambda, \wh f_d) \in W^1_q (\Omega_-) + \wh W^1_{q, \Gamma} (\Omega_-)$ be a solution to the variational problem:
\begin{align}
\label{eq-3.12}
(\nabla \CK(\lambda, \wh f_d), \nabla \varphi)_{\Omega_-} = (\rho_{*-} \lambda \wh \bF_d - \nabla \wh f_d, \nabla \varphi)_{\Omega_-}
\end{align}
for any $\varphi \in \wh W^1_{q', \Gamma} (\Omega_-)$ with $\CK (\lambda,\wh f_d) = - \wh f_d$ on $\Gamma$. Let $\wh \rho_+ \in W^3_q (\Omega_+)$, $\wh \bu_\pm \in W_q^2 (\Omega_\pm)$, and $\wh h \in W^{3 - 1/q}_q (\Gamma)$ be solutions to 
\begin{align}
\label{eq-3.13}
\left\{\begin{aligned}
\lambda \wh \rho_+ + \rho_{*+} \dv \wh \bu_+ & = \wh f_M &\enskip &\text{ in $\Omega_+$,} \\
\rho_{*-} \dv \wh \bu_- = \wh f_d = \rho_{*-} & \dv \wh \bF_d  &\enskip &\text{ in $\Omega_-$}, \\
\rho_{*+} \lambda \wh \bu_+ - \DV \bT_{0 +} (\gamma_1, \gamma_2, \gamma_3, \wh \bu_+, \wh \rho_+) & = \wh \bff_+ &\enskip &\text{ in $\Omega_+$}, \\
\rho_{*-} \lambda \wh \bu_- - \DV \bT_- (\gamma_4, \wh \bu_-, \CK_1 (\wh \bu_-) + \CK_2 (\wh h)) & = \wh \bff_- - \nabla \CK (\lambda, \wh f_d) &\enskip &\text{ in $\Omega_-$}, \\
\lambda \wh h - \frac{\langle \rho_{*-} \wh \bu_-, \bn \rangle \vert_- - \langle\rho_{*+} \wh \bu_+, \bn \rangle\vert_+} {\rho_{*-}-\rho_{*+}} & = \wh d &\enskip &\text{ on $\Gamma$}, \\
\bB_0 (\gamma_1, \gamma_2, \gamma_3, \gamma_4, \wh \rho_+, \wh \bu_+, \wh \bu_-) & = \wh \bG_0 &\enskip &\text{ on $\Gamma$}, \\
\rho_{*-} \dv \wh \bu_- & = \wh g_- + \wh f_d &\enskip &\text{ on $\Gamma$}, \\
\color{red}
\wh \bu_+ = 0, \qquad \langle \nabla \rho_+, \bn_+ \rangle &
\color{red}
= \wh k_+ &\enskip &
\color{red}
\text{ on $\Gamma_+$}, \\
\color{red}
\wh \bu_- & 
\color{red}
= 0 &\enskip &
\color{red}
\text{ on $\Gamma_-$}.
\end{aligned}\right.
\end{align}
From the last condition of \eqref{eq-3.13}, we have
\begin{align}
\label{eq-3.14}
\rho_{*-} \dv \wh \bu_- \vert_- = \wh f_d \quad \text{ on $\Gamma$}.
\end{align}
By \eqref{eq-3.11}, \eqref{eq-3.12}, and \eqref{eq-3.13}, for any $\varphi \in \wh W^1_{q', \Gamma} (\Omega_-)$ we observe that
\begin{equation}
\label{eq-3.15}
\begin{split}
(\rho_{*-} \lambda \wh \bF_d - \nabla \wh f_d, \nabla \varphi)_{\Omega_-} & = (\nabla \CK (\lambda, \wh f_d), \nabla \varphi)_{\Omega_-} \\
& = (\rho_{*-} \lambda \wh \bu_- - \DV(\mu_- \bD(\wh \bu_-) - (\CK_1 (\wh \bu_-) + \CK_2 (\wh h))), \nabla \varphi)_{\Omega_-} \\
& = (\rho_{*-} \lambda \wh \bu_-, \nabla \varphi)_{\Omega_-} - (\rho_{*-} \nabla \dv \wh \bu_-, \nabla \varphi)_{\Omega_-}.
\end{split}
\end{equation}
Since $W^1_{q', \Gamma} (\Omega_-) \subset \wh W^1_{q', \Gamma} (\Omega_-)$, by the divergence theorem of Gauss, the definition of $\CD \CI_q (\Omega_-)$, and~\eqref{eq-3.15}, we have
\begin{align*}
& \lambda (\rho_{*-} \dv \wh \bu_-, \varphi)_{\Omega_-} + (\rho_{*-} \nabla \dv \wh \bu_-, \nabla \varphi)_{\Omega_-} \\
& = - \lambda (\rho_{*-} \wh \bu_-, \nabla \varphi)_{\Omega_-} + \lambda (\rho_{*-} \wh \bu_-, \nabla \varphi)_{\Omega_-} - \lambda (\rho_{*-} \wh \bF_d, \nabla \varphi)_{\Omega_-} + (\nabla \wh f_d, \nabla \varphi)_{\Omega_-} \\
& = \lambda (\wh f_d, \varphi)_{\Omega_-} + (\nabla \wh f_d, \nabla \varphi)_{\Omega_-}
\end{align*}
for any $\varphi \in W^1_{q', \Gamma} (\Omega_-)$. Namely, $\lambda (\wh f_d - \rho_{*-} \dv \wh \bu_-, \varphi)_{\Omega_-} + (\nabla (\wh f_d - \rho_{*-} \dv \wh \bu_-), \nabla \varphi)_{\Omega_-} = 0$ for any $\varphi \in W^1_{q', \Gamma} (\Omega_-)$. From \eqref{eq-3.14}, a uniqueness of solutions implies that $\rho_{*-} \dv \wh \bu_- = \wh f_d$ in $\Omega_-$, which substitute into \eqref{eq-3.15} implies that $(\wh \bF_d, \nabla \varphi)_{\Omega_-} = (\wh \bu_-, \nabla \varphi)_{\Omega_-}$ for any $\varphi \in \wh W^1_{q', \Gamma} (\Omega_-)$ because we may assume $\lambda \neq 0$. Hence, in light of \eqref{cond-1.11} and \eqref{cond-1.12}, $\wh \rho_+$, $\wh \bu_\pm$, $\wh h$, and $\wh \pi_- = \CK_1 (\wh \bu_-) + \CK_2 (\wh h) - \CK(\lambda, \wh f_d)$ satisfy~\eqref{eq-4.1*} assuming \eqref{eq-3.11}.

\subsection{On the $\mathcal{R}$ bounded solution operators for the reduced problem}
In the following, we consider \eqref{eq-3.8} instead of \eqref{eq-4.1*}. We first define function spaces $Y_q$ and $\CY_q$ as follows:
\begin{equation*}
\begin{split} 
Y_q(\Omega_+,\Omega_-,\Gamma) & = \{(\wh f_M, \wh \bff_+,  \wh \bff_-, \wh d, \wh g, \wh g_+, \wh g_-, \wh \bh, \wh k_-, \wh k_+) \mid \wh f_M \in W^1_q (\Omega_+), \enskip \wh \bff_+ \in L_q (\Omega_+)^N, \enskip \wh \bff_- \in L_q (\Omega_-)^N, \\
& \quad \quad \wh d \in W^{2 - 1 / q}_q (\Gamma), \enskip \wh g, \enskip \wh g_+, \wh g_- \in W^1_q (\dot \Omega), \enskip \wh \bh \in W^2_q (\dot \Omega)^N, \enskip \wh k_-, \wh k_+ \in W^2_q (\Omega_+) \}, \\
\CY_q(\Omega_+,\Omega_-,\Gamma) & = \{(F_1, \dots, F_{21}) \mid F_1, F_2 \in L_q (\Omega_+), \enskip F_3 \in L_q (\Omega_+)^N, \enskip F_4 \in L_q (\Omega_-)^N, \enskip F_5 \in W^{2 - 1 / q}_q (\Gamma), \\
& \quad \quad F_6, F_8, F_{10} \in L_q (\dot \Omega), \enskip F_7, F_9, F_{11}, F_{13} \in L_q (\dot \Omega)^N, \enskip F_{12} \in W^1_q (\dot \Omega), \enskip F_{14} \in L_q (\dot \Omega)^{N^2},\\
& \quad \quad F_{15} \in L_q (\dot \Omega)^{N^3}, \enskip F_{16}, F_{19} \in L_q (\Omega_+), \enskip F_{17}, F_{20} \in L_q (\Omega_+)^N, \enskip F_{18}, F_{21} \in L_q (\Omega_+)^{N^2}\}.
\end{split}
\end{equation*}
Furthermore, we set
\begin{align*}
& \lVert (\wh f_M, \wh \bff_+,  \wh \bff_-, \wh d, \wh g, \wh g_+, \wh g_-, \wh \bh, \wh k_-, \wh k_+) \rVert_{Y_q(\Omega_+, \Omega_-, \Gamma)} \\
& = \lVert \wh f_M \rVert_{W^1_q (\Omega_+)} + \sum_{\ell = \pm} (\lVert \wh \bff_\ell \rVert_{L_q(\Omega_\ell)^N} + \lVert \wh g_\ell \rVert_{W^1_q(\dot{\Omega})^N} + \lVert \wh k_\ell \rVert_{W^2_q (\Omega_+)}) + \lVert \wh d \rVert_{W^{2-1/q}_q(\Gamma)} + \lVert \wh g \rVert_{W^1_q(\dot{\Omega})^N} + \lVert \wh \bh \rVert_{W^2_q(\dot{\Omega})}, \\
& \lVert (F_1,\dots,F_{21}) \rVert_{\CY_q(\Omega_+,\Omega_-,\Gamma)} \\
&= \sum_{m = 1}^3 \lVert F_m \rVert_{L_q(\Omega_+)} + \lVert F_4 \rVert_{L_q(\Omega_-)} + \lVert F_5 \rVert_{W^{2-1/q}_q(\Gamma)} + \sum_{n = 6}^{15} \lVert F_n \rVert_{L_q(\dot{\Omega})} + \lVert \nabla F_{12} \rVert_{L_q(\dot{\Omega})} + \sum_{l = 16}^{21} \lVert F_l \rVert_{L_q(\Omega_+)}
\end{align*}
for any $(\wh f_M, \wh \bff_+,  \wh \bff_-, \wh d, \wh g, \wh g_+, \wh g_-, \wh \bh, \wh k_-, \wh k_+) \in Y_q (\Omega_+, \Omega_-, \Gamma)$ and $(F_1, \dots, F_{21}) \in \CY_q (\Omega_+, \Omega_-, \Gamma)$. The following theorem plays an important role in the present paper. We will give the proof in the next section.

\begin{theo}
\label{TH-3.1}
Let $1 < q < \infty$, $N < r < \infty$, and $\max(q, q') \leq r$. Suppose that the assumptions (a)--(c) in Theorem~\ref{th-1.9} holds. Then there exists constant $\varepsilon_* \in (0, \pi/2)$ such that for any $\varepsilon\in (\varepsilon_*, \pi/2)$ there exists a constant $\lambda_* \geq 1$ with the following assertions hold true:
\begin{enumerate}\renewcommand{\labelenumi}{(\arabic{enumi})}
\item For any $\lambda \in \Sigma_{\varepsilon, \lambda_*}$, there exist operators $\CA^+ (\lambda)$, $\CB^\pm (\lambda)$, and $\CH (\lambda)$ with
\begin{align*}
\CA^{\color{red} +} (\lambda) & \in \Hol (\Sigma_{\varepsilon, \lambda_*}, \CL(\CY_q(\Omega_+, \Omega_-, \Gamma), W^3_q (\Omega_+))), \\
\CB^{\color{red} \pm} (\lambda) & \in \Hol (\Sigma_{\varepsilon, \lambda_*}, \CL(\CY_q(\Omega_+, \Omega_-, \Gamma), W^2_q (\Omega_\pm)^N)), \\
\CH (\lambda) & \in \Hol (\Sigma_{\varepsilon, \lambda_*}, \CL(\CY_q(\Omega_+, \Omega_-, \Gamma), W^{3 - 1/q}_q (\Gamma))),
\end{align*}
such that for any $\bF_Y = (\wh f_M, \wh \bff_+,  \wh \bff_-, \wh d, \wh g, \wh g_+, \wh g_-, \wh \bh, \wh k_-, \wh k_+) \in Y_q(\Omega_+, \Omega_-,\Gamma)$, the quadruple 
\begin{equation*}
(\rho_+, \bu_+, \bu_-, h) = (\CA^+ (\lambda) F_\lambda (\bF_Y), \CB^+ (\lambda) F_\lambda (\bF_Y), \CB^- (\lambda) F_\lambda (\bF_Y), \CH (\lambda) F_\lambda (\bF_Y))
\end{equation*}
is a unique solution to~\eqref{eq-3.8}. Here, we have set
\begin{align*}
F_\lambda(\bF_Y) & = (\lambda^{1/2} \wh f_M, \nabla \wh f_M, \wh \bff_+, \wh \bff_-, \wh d, \lambda^{1/2} \wh g, \nabla \wh g, \lambda^{1/2} \wh g_+, \nabla \wh g_+, \lambda^{1/2} \wh g_-, \nabla \wh g_-, \\
& \quad \wh g_-, \lambda \wh \bh, \lambda^{1/2} \nabla \wh \bh, \nabla^2 \wh \bh, \lambda \wh k_-, \lambda^{1/2} \nabla \wh k_-, \nabla^2 \wh k_-, \lambda \wh k_+, \lambda^{1/2} \nabla \wh k_+, \nabla^2 \wh k_+).
\end{align*}
\item There exists a positive constant $c_*$, independent of $\lambda$, such that
\begin{equation*}
\begin{split}
\CR_{\CL(\CY_q(\Omega_+, \Omega_-, \Gamma), W^{3 - i}_q (\Omega_+))} (\{(\tau \pd_\tau)^s (\lambda^{i/2} \CA^+ (\lambda)) \mid \lambda \in \Sigma_{\varepsilon, \lambda_*}\}) & \leq c_*, \\
\CR_{\CL(\CY_q(\Omega_+, \Omega_-, \Gamma), W^{2 - j}_q (\Omega_\pm)^N)} (\{(\tau \pd_\tau)^s (\lambda^{j/2} \CB^\pm (\lambda)) \mid \lambda \in \Sigma_{\varepsilon, \lambda_*}\}) & \leq c_*, \\
\CR_{\CL(\CY_q(\Omega_+, \Omega_-, \Gamma), W^{3 - 1/q - k}_q (\Gamma))} (\{(\tau \pd_\tau)^s (\lambda^k \CH (\lambda)) \mid \lambda \in \Sigma_{\varepsilon, \lambda_*}\}) & \leq c_*,
\end{split}
\end{equation*}	
for $s = 0, 1$, $i = 0, 1, 2, 3$, $j = 0, 1, 2$, and $k = 0, 1$.
\end{enumerate}
\end{theo}

\section{Generalized resolvent problem}
\label{sect-4}
\subsection{Reduced problem with a flat interface}
To prove Theorem~\ref{TH-3.1}, we first consider a flat interface case, that is, we consider the problem \eqref{eq-3.8} with $\Omega_\pm = \BR^N_\pm$, $\Gamma = \BR^N_0$, and $\Gamma_+ = \Gamma_- = \emptyset$, where we have set $\BR^N_\pm := \{x \in \BR^N \mid \pm x_N > 0 \}$ and $\BR^N_0 := \{x \in \BR^N \mid x_N = 0 \}$. Furthermore, we set
\begin{align*}
\hW^1_{q, 0} (\RM) = \{\theta\in L_{p, \mathrm{loc}} (\BR^N_-) \mid \nabla \theta \in L_q (\BR^N_-)^N, \enskip \theta \vert_{x_N = 0} = 0\}, \quad W^1_{q, 0} (\RM) = \{\theta \in W^1_q (\BR^N_-) \mid \theta \vert_{x_N = 0}\}.
\end{align*}
Let
\begin{align*}
\CD\CI_{F, q} (\RM) = \Set{\wh f_d \in W^1_q (\RM) \vert
\begin{array}{l}
\text{ there exists $\wh \fg_d \in L_q (\BR^N_-)^N$ such that} \\
\text{ $(\wh f_d,\varphi)_{\RM} = - \rho_{*-} (\fg_d, \nabla \varphi)_{\RM}$ for all $\varphi\in W^1_{q', 0} (\RM)$}
\end{array}}.		
\end{align*}
In addition, let $\CG_F (f_d) = \{\fh_d \in L_q (\RM)^N \mid \dv \fg_d = \dv \fh_d\}$ and $\bF_d \in \CG_F(f_d)$. We observe that $\rho_{* -} \dv \bF_d = f_d$ in $\RM$. Set $\lVert f_d \rVert_{\CD\CI_{F, q}(\BR^N)} = \lVert f_d \rVert_{W^1_q (\RM)} + \inf_{\fh_d \in \CG_F (f_d)} \lVert \fh_d \rVert_{L_q (\RM)}$ for $f_d \in \CD\CI_{F, q} (\RM)$. We see that $\CD\CI_{F, q} (\Omega_-)$ is a Banach space with norm $\lVert \, \cdot \, \rVert_{\CD\CI_{F, q} (\Omega_-)}$. \par
For any $\bu_- \in W^2_q (\RM)$, let $\CK_{F1} (\bu_-) \in W^1_q (\RM) + \wh W^1_{q, 0} (\RM)$ be a unique solution to the following variational problem:
\begin{align}
\label{eq-4.2}
(\nabla \CK_{F1} (\bu_-), \nabla \varphi)_\RM = (\DV (\gamma_{40} \bD (\bu_-)) - \rho_{* -} \nabla \dv \bu_-, \nabla \varphi)_{\Omega_-} \quad \text{ for any $\varphi \in \hW^1_{q', 0} (\RM)$}
\end{align}
subject to $\CK_{F1} (\bu_-) = \langle \gamma_{40} \bD (\bu_-) \bn_0, \bn_0 \rangle - \rho_{* -} \dv \bu_-$ on $\RZ$, while for $h \in W^{3 - 1/q}_q (\BR^{N - 1})$, $\CK_{F2} (h) \in W^1_q (\RM) +\hW^1_{q, 0} (\RM)$ be a unique solution to the following variational problem:
\begin{align}
\label{eq-4.3}
(\nabla \CK_{F2} (h), \nabla \varphi)_\RM = 0 \quad \text{ for any $\varphi\in \hW^1_{q',0}(\RM)$}
\end{align}
subject to $\CK_{F2} (h) = - (\rho_{* -} - \rho_{* +})^{- 1} \rho_{* -} \sigma \Delta_\RZ h$ on $\RZ$ with $\rho_{* +} \neq \rho_{* -}$. Here, we have the following estimates:
\begin{equation}
\label{est-4.4}
\begin{split}
\lVert \nabla \CK_{F1} (\bu_-) \rVert_{L_q (\RM)} \leq C \lVert \nabla \bu_- \rVert_{W^1_q (\RM)}, \quad \lVert \nabla \CK_{F2} (h) \rVert_{L_q (\RM)} \leq C \lVert h \rVert_{W^{3 - 1/q}_q (\RZ)}.
\end{split}
\end{equation}
According to the previous section, we obtain the reduced equations:
\begin{align}
\label{eq-4.5}
\left\{\begin{aligned}
\lambda \wh \rho_+ + \rho_{*+} \dv \wh \bu_+ & = \wh f_M &\enskip &\text{ in $\BR^N_+$,} \\
\rho_{* +} \lambda \wh \bu_+ - \DV \bT_{0 +} (\gamma_{10}, \gamma_{20}, \gamma_{30}, \wh \bu_+, \wh \rho_+) & = \wh \bff_+ &\enskip &\text{ in $\BR^N_+$}, \\
\rho_{* -} \lambda \wh \bu_- - \DV \bT_- (\gamma_{40}, \wh \bu_-, \CK_{{\color{red} F1}} (\wh \bu_-) + \CK_{{\color{red} F2}} (\wh h)) & = \wh \bff_- &\enskip &\text{ in $\BR^N_-$}, \\
\lambda \wh h - \frac{\langle \rho_{* -} \wh \bu_-, \bn \rangle \vert_- - \langle \rho_{* +} \wh \bu_+, \bn \rangle \vert_+}{\rho_{* -} - \rho_{* +}} & = \wh d &\enskip &\text{ on $\BR^N_0$}, \\
\bB_0 (\gamma_{10}, \gamma_{20}, \gamma_{30}, \gamma_{40}, \wh \rho_+, \wh \bu_+, \wh \bu_-) & = \wh \bG_0 &\enskip &\text{ on $\BR^N_0$}, \\
\rho_{* -} \dv \wh \bu_- \vert_- & = \wh g_- &\enskip &\text{ on $\BR^N_0$}
\end{aligned}\right.
\end{align}
where $\gamma_{10}$, $\gamma_{20}$, $\gamma_{30}$, and $\gamma_{40}$ are positive constants. We now define function spaces $Z_q$ and $\CZ_q$ as follows:
\begin{equation*}
\begin{split} 
& Z_q(D_+, D_-, D_0) \\
& \quad = \{(\wh f_M, \wh \bff_+,  \wh \bff_-, \wh d, \wh g, \wh g_+, \wh g_-, \wh \bh, \wh k_-) \mid \wh f_M \in W^1_q (D_+), \enskip \wh \bff_+ \in L_q (D_+)^N, \enskip \wh \bff_- \in L_q (D_-)^N,  \\
& \quad \qquad \wh d \in W^{2 - 1 / q}_q (D_0), \enskip \wh g, \wh g_+, \wh g_- \in W^1_q (\dot D), \enskip \wh \bh \in W^2_q (\dot D)^N, \enskip \wh k_- \in W^2_q (D_+) \}, \\
& \CZ_q(D_+,D_-,D_0) \\
& \quad = \{(G_1, \dots, G_{18}) \mid G_1, G_2 \in L_q (D_+), \enskip G_3 \in L_q (D_+)^N, \enskip G_4 \in L_q (D_-)^N, \enskip Z_5 \in W^{2 - 1 / q}_q (D_0), \\
& \quad \qquad G_6, G_8, G_{10} \in L_q (\dot D),\enskip G_7, G_9, G_{11}, G_{13} \in L_q (\dot D)^N, \enskip G_{12} \in W^1_q (\dot D), \enskip G_{14} \in L_q (\dot D)^{N^2}, \\
& \quad \qquad G_{15} \in L_q (\dot D)^{N^3}, \enskip G_{16} \in L_q (D_+), \enskip G_{17} \in L_q (D_+)^N, \enskip G_{18} \in L_q (D_+)^{N^2}\}.
\end{split}
\end{equation*}
with $\dot D = D_+ \cup D_-$, $D_\pm \subset \BR^N$, $D_0 \subset \BR^{N - 1}$. In addition, we set
\begin{align*}
& \lVert (\wh f_M, \wh \bff_+,  \wh \bff_-, \wh d, \wh g, \wh g_+, \wh g_-, \wh \bh, \wh k_-) \rVert_{Z_q(D_+, D_-, D_0)} \\
& = \lVert \wh f_M \rVert_{W^1_q (D_+)} + \sum_{\ell = \pm} (\lVert \wh \bff_\ell \rVert_{L_q(D_\ell)^N} + \lVert \wh g_\ell \rVert_{W^1_q(\dot D)^N}) + \lVert \wh d \rVert_{W^{2-1/q}_q(D_0)} + \lVert \wh g \rVert_{W^1_q(\dot D)^N} + \lVert \wh \bh \rVert_{W^2_q(\dot D)} + \lVert \wh k_- \rVert_{W^2_q (D_+)}, \\
& \lVert (G_1,\dots,G_{18}) \rVert_{\CZ_q(D_+,D_-,D_0)} \\
& = \sum_{m = 1}^3 \lVert G_m \rVert_{L_q(D_+)} + \lVert G_4 \rVert_{L_q(D_-)} + \lVert G_5 \rVert_{W^{2-1/q}_q(D_0)} + \sum_{n = 6}^{15} \lVert G_n \rVert_{L_q(\dot D)} + \lVert \nabla G_{12} \rVert_{L_q(\dot D)} + \sum_{l = 16}^{18} \lVert G_l \rVert_{L_q(D_+)}
\end{align*}
for any $(\wh f_M, \wh \bff_+,  \wh \bff_-, \wh d, \wh g, \wh g_+, \wh g_-, \wh \bh, \wh k_-) \in Z_q(D_+,D_-,D_0)$ and $(G_1,\dots,G_{18})\in\CZ_q(D_+,D_-,D_0)$. The purpose of this section is to prove the following theorem.
\begin{theo}
\label{th-4.4}
Let $1 < q < \infty$ and $\rho_{* +} \neq \rho_{* -}$. Assume that $\rho_{* +}$, $\gamma_{10}$, $\gamma_{20}$, and $\gamma_{30}$ satisfy Assumption~\ref{as-star} (c) with $\gamma_k = \gamma_{k0}$ $(k = 1, 2, 3)$. Then there exists constant ${\color{red} \varepsilon_*} \in (0, \pi/2)$ such that for any $\varepsilon\in ({\color{red} \varepsilon_*}, \pi/2)$ there exists a constant $\lambda_1 > 0$ with the following assertions valid:
\begin{enumerate}\renewcommand{\labelenumi}{(\arabic{enumi})}
\item For any $\lambda\in \Sigma_{\varepsilon, \lambda_1}$, there exists operators
\begin{align*}
\CA^+_{F0} (\lambda) & \in \Hol (\Sigma_{\varepsilon, \lambda_1}, \CL (\CZ_q (\RP, \RM,\RZ), W^3_q (\RP))), \\
\CB^\pm_{F0} (\lambda) & \in \Hol (\Sigma_{\varepsilon, \lambda_1}, \CL (\CZ_q (\RP, \RM, \RZ), W^2_q (\BR^N_\pm)^N)), \\
\CH_{F0} (\lambda) & \in \Hol (\Sigma_{\varepsilon, \lambda_1}, \CL (\CZ_q (\RP, \RM, \RZ), W^{3 - 1/q}_q (\RZ))),
\end{align*}
such that for any $\bF_Z := (\wh f_M, \wh \bff_+,  \wh \bff_-, \wh d, \wh g, \wh g_+, \wh g_-, \wh \bh, \wh k_-) \in Z_q (\RP, \RM, \RZ)$, the problem \eqref{eq-4.5} admit a unique solution $(\wh \rho_+, \wh \bu_+, \wh \bu_-, \wh h)$ defined by $\wh \rho_+ = \CA^+_{F0} (\lambda) G_\lambda (\bF_Z)$, $\wh \bu_\pm = \CB^\pm_{F0} (\lambda) G_\lambda (\bF_Z)$, and $\wh h = \CH_{F0} (\lambda) G_\lambda (\bF_Z)$ with
\begin{equation*}
\begin{split}
G_\lambda (\bF_Z) & = (\lambda^{1/2} {\color{red} \wh f_M}, \nabla {\color{red} \wh f_M}, {\color{red} \wh \bff_+}, {\color{red} \wh \bff_-}, {\color{red} \wh d}, \lambda^{1/2} {\color{red} \wh g}, \nabla {\color{red} \wh g}, \lambda^{1/2} {\color{red} \wh g_+}, \nabla {\color{red} \wh g_+}, \lambda^{1/2} {\color{red} \wh g_-}, \\
& \qquad \nabla {\color{red} \wh g_-}, {\color{red} \wh g_-}, \lambda {\color{red} \wh \bh}, \lambda^{1/2} \nabla {\color{red} \wh \bh}, \nabla^2 {\color{red} \wh \bh}, \lambda {\color{red} \wh k_-}, \lambda^{1/2} \nabla {\color{red} \wh k_-}, \nabla^2 {\color{red} \wh k_-}).
\end{split}
\end{equation*}
\item For $s = 0, 1$, $i = 0, 1, 2, 3$, $j = 0, 1, 2$, and $k = 0, 1$, there exists a positive constant $c_1$ such that
\begin{equation*}
\begin{split}
\CR_{\CL (\CZ_q (\RP, \RM, \RZ), W^{3 - i}_q (\RP))} (\{(\tau \pd_\tau)^s (\lambda^{i/2} \CA^+_{F0} (\lambda)) \mid \lambda \in \Sigma_{\varepsilon, \lambda_1} \}) & \leq c_1, \\
\CR_{\CL (\CZ_q (\RP, \RM, \RZ), W^{2 - j}_q (\BR^N_\pm)^N)} (\{(\tau \pd_\tau)^s (\lambda^{j/2} \CB^\pm_{F0} (\lambda)) \mid \lambda \in \Sigma_{\varepsilon, \lambda_1}\}) & \leq c_1, \\
\CR_{\CL (\CZ_q (\RP, \RM, \RZ), W^{3 - 1/q - k}_q (\RZ))} (\{(\tau \pd_\tau)^s (\lambda^k \CH_{F0} (\lambda)) \mid \lambda \in \Sigma_{\varepsilon, \lambda_1}\}) & \leq c_1.
\end{split}
\end{equation*}
Here, the constant $c_1$ is independent of $\lambda$.
\end{enumerate}
\end{theo}

To prove Theorem \ref{th-4.4}, for given $\wh \bff_- \in L_q (\RM)^N$ and $\wh g_- \in W^1_q (\dot \BR^N)$, we consider a function $\wh f_d$ satisfying
\begin{align}
\label{eq-4.6}
(\lambda \wh f_d, \varphi)_\RM + (\nabla \wh f_d, \varphi)_\RM = - (\rho_{* -}^{- 1} \wh \bff_-, \nabla \varphi) \quad \text{ for any $\varphi \in \wh W^1_{q', \Gamma} (\RM)$}
\end{align}
subject to $\wh f_d = \wh g_-$ on $\RZ$. According to the discussion in Sect.~\ref{sect-3.1}, functions $\wh \rho_+$, $\wh \bu_\pm$, and $\wh h$ satisfying
\begin{align}
\label{eq-4.7}
\left\{\begin{aligned}
\lambda \wh \rho_+ + \rho_{*+} \dv \wh \bu_+ & = \wh f_M &\enskip &\text{ in $\BR^N_+$,} \\
\rho_{*-} \dv \wh \bu_- = \wh f_d = \rho_{*-} & \dv \wh \bF_d  &\enskip &\text{ in $\BR^N_-$}, \\
\rho_{*+} \lambda \wh \bu_+ - \DV \bT_{0 +} (\gamma_{10}, \gamma_{20}, \gamma_{30}, \wh \bu_+, \wh \rho_+) & = \wh \bff_+ &\enskip &\text{ in $\BR^N_+$}, \\
\rho_{*-} \lambda \wh \bu_- - \DV \bT_- (\gamma_{40}, \wh \bu_-, \wh \pi_-) & = \wh \bff_- &\enskip &\text{ in $\BR^N_-$}, \\
\lambda \wh h - \frac{\langle \rho_{*-} \wh \bu_-, \bn \rangle \vert_- - \langle\rho_{*+} \wh \bu_+, \bn \rangle\vert_+} {\rho_{*-}-\rho_{*+}} & = \wh d &\enskip &\text{ on $\BR^N_0$}, \\
\bB_0 (\gamma_{10}, \gamma_{20}, \gamma_{30}, \gamma_{40}, \wh \rho_+, \wh \bu_+, \wh \bu_-) & = \wh \bG_0 &\enskip &\text{ on $\BR^N_0$}, \\
\langle \bT_- (\gamma_{40}, \wh \bu_-, \wh \pi_-) \bn, \bn \rangle \Big\vert_- - \frac{\rho_{*-} \sigma}{\rho_{*-} - \rho_{*+}} \Delta_{\BR^N_0} \wh h & = \wh g_- &\enskip &\text{ on $\BR^N_0$}
\end{aligned}\right.
\end{align}
are solutions to \eqref{eq-4.5}. In case of $\wh f_d \equiv 0$ in \eqref{eq-4.6}, the following lemma is applied, see \cite[Theorem 1.2]{Wat2017}.
\begin{lemm}
\label{th-4.5}
Let $1 < q < \infty$ and $\rho_{* +} \neq \rho_{* -}$. Suppose that $\rho_{* +}$, $\gamma_{10}$, $\gamma_{20}$, $\gamma_{30}$ satisfy Assumption \ref{as-star} with $\gamma_k = \gamma_{k0}$ $(k = 1, 2, 3)$. Then there exists constant ${\color{red} \varepsilon_*} \in (0, \pi/2)$ such that for any $\varepsilon\in ({\color{red} \varepsilon_*}, \pi/2)$ there exists a constant $\lambda_2 > 0$ with the following statements satisfy:
\begin{enumerate}\renewcommand{\labelenumi}{(\arabic{enumi})}
\item For any $\lambda \in \Sigma_{\varepsilon, \lambda_2}$, operators
\begin{align*}
\CA^+_{F1} (\lambda) & \in \Hol (\Sigma_{\varepsilon, \lambda_2}, \CL (\CZ_q (\RP, \RM, \RZ), W^3_q (\RP))), \\
\CB^\pm_{F1} (\lambda) & \in \Hol (\Sigma_{\varepsilon, \lambda_2}, \CL (\CZ_q (\RP, \RM, \RZ), W^2_q (\BR^N_\pm)^N)), \\
\CP^-_{F1} (\lambda) & \in \Hol (\Sigma_{\varepsilon, \lambda_2}, \CL(\CZ_q (\RP, \RM, \RZ), \wh W^1_q (\RM))), \\
\CH_{F1} (\lambda) & \in \Hol (\Sigma_{\varepsilon, \lambda_2}, \CL (\CZ_q (\RP, \RM, \RZ), W^{3 - 1/q}_q (\RZ)))
\end{align*}
exists such that for any $\bF_Z\in Z_q(\RP,\RM,\RZ)$, the quintuple
\begin{equation*}
\begin{split}
& (\wh \rho_+, \wh \bu_+, \wh \bu_-, \wh \pi_-, \wh h) \\
& = (\CA_{F1}^+ (\lambda) G_\lambda (\bF_Z), \CB_{F1}^+ (\lambda) G_\lambda (\bF_Z), \CB_{F1}^- (\lambda) G_\lambda (\bF_Z), \CP_{F1}^- (\lambda) G_\lambda (\bF_Z), \CH_{F1} (\lambda) G_\lambda (\bF_Z))
\end{split}
\end{equation*}
is a unique solution to \eqref{eq-4.7} with $\wh f_d \equiv 0$.
\item For $s = 0, 1$, $i = 0, 1, 2, 3$, $j = 0, 1, 2$, and $k = 0, 1$, the estimates
\begin{align*}
\CR_{\CL (\CZ_q (\RP, \RM, \RZ), W^{3 - i}_q (\BR^N_\pm)^N)} (\{(\tau \pd_\tau)^s (\lambda^{i/2} \CA^+_{F1} (\lambda)) \mid \lambda \in \Sigma_{\varepsilon, \lambda_2}\}) & \leq c_2, \\
\CR_{\CL (\CZ_q (\RP, \RM, \RZ), W^{2 - j}_q (\RP))} (\{(\tau \pd_\tau)^s (\lambda^{j/2} \CB^\pm_{F1} (\lambda)) \mid \lambda \in \Sigma_{\varepsilon, \lambda_2}\}) & \leq c_2, \\
\CR_{\CL (\CZ_q (\RP, \RM, \RZ), L_q (\RM)^N)} (\{(\tau \pd_\tau)^s (\nabla \CP^-_{F1}(\lambda)) \mid \lambda \in \Sigma_{\varepsilon, \lambda_2}\}) & \leq c_2, \\
\CR_{\CL (\CZ_q (\RP, \RM, \RZ), W^{3 - 1/q - k}_q (\RZ))} (\{(\tau \pd_\tau)^s (\lambda^k \CH_{F1} (\lambda)) \mid \lambda \in \Sigma_{\varepsilon, \lambda_2}\}) & \leq c_2
\end{align*}			
hold true with some positive constant $c_2$ independent of $\lambda$.
\end{enumerate}
\end{lemm}

To treat the case of $\wh f_d \neq 0$, we consider the divergence equation:
\begin{align}
\label{eq-4.8}
\rho_{* -} \dv \wh \bu_\mathrm{div} = \wh f_d \quad \text{in $\RM$,}
\end{align}
where $\wh f_d$ is a solution to \eqref{eq-4.6}. The solution $\wh f_d$ is given by the following lemma shown by Shibata~\cite[Theorem 9.3.10]{Shi2016a}. 
\begin{lemm}
\label{lem-4.6}
Let $1 < q < \infty$ and $0 < \varepsilon < \pi/2$. Let
\begin{equation*}
\begin{split}
Z^1_q (\RM) & = \{(\bff_-, g_-) \mid \bff_- \in L_q (\RM)^N, \enskip g_- \in W^1_q (\dot \BR^N) \}, \\
\CZ^1_q (\RM) & = \{(G_4, G_{12}, G_{13}) \mid G_4 \in L_q (\RM)^N, \enskip G_{12} \in L_q (\dot \BR^N), \enskip G_{13} \in W^1_q (\dot \BR^N)\}.
\end{split}
\end{equation*}
Then we have the following assertions:
\begin{enumerate}\renewcommand{\labelenumi}{(\arabic{enumi})}
\item There exists an operator $\CD_{F1} \in \Hol (\Sigma_\varepsilon, \CL (\CZ^1_q (\RM), \CD\CI_{F, q} (\RM)))$ such that for any $\lambda \in \Sigma_\varepsilon$ and $(\wh \bff_-, \wh g_-) \in Z^1_q (\RM)$, the problem \eqref{eq-4.6} admits a unique solution $\wh f_d = \CD_{F1} (\lambda)(\wh \bff_-, \wh g_-, \lambda^{1/2} \wh g_-)$. Furthermore, for $s = 0, 1$, $j = 0, 1, 2$, and any $\lambda_0 > 0$, the estimate
\begin{align*}
\CR_{\CL (\CZ^1_q (\RM), W^{1 - j}_q (\RM)^N)} (\{(\tau \pd_\tau)^s (\lambda^{j/2} \CD_{F1}(\lambda)) \mid \lambda \in \Sigma_{\varepsilon, \lambda_0}\}) \leq c_{\lambda_0}
\end{align*}
is valid with some positive constant $c_{\lambda_0}$ independent of $\lambda$.
\item For the function $\wh f_d$ given in (1), there exist an operator $\CD_{F2} \in \Hol (\Sigma_\varepsilon, \CL (\CZ^1_q (\RM), W^2_q (\RM)^N))$ such that for any
$\lambda \in \Sigma_\varepsilon$ and $(\wh \bff_-, \wh g_-) \in Z^1_q (\RM)$, the problem admits a unique solution $\bu_\mathrm{div} = \CD_{F2} (\lambda) (\wh \bff_-, \wh g_-, \lambda^{1/2} \wh g_-)$. In addition, for $s = 0, 1$, $j = 0, 1, 2$ and any $\lambda_0 > 0$. the estimate:
\begin{equation*}
\begin{split}
\CR_{\CL(\CZ^1_q(\RM),W^{2-j}_q(\RM)^N)}(\{(\tau\pd_\tau)^s
(\lambda^{j/2}\CD_{F2}(\lambda))\mid 
\lambda\in\Sigma_{\varepsilon,\lambda_0} \})\leq c_{\lambda_0}
\end{split}
\end{equation*}
holds true, where $c_{\lambda_0}$ is some positive constant independent of $\lambda$.
\end{enumerate}			
\end{lemm}
Let $\bu_\pm$ be a solution to \eqref{eq-4.7} and let $\wh \bw_- = \wh \bu_- - \wh \bu_\mathrm{div}$, where $\wh \bu_\mathrm{div}$ is a solution of \eqref{eq-4.8}. Then $\wh \rho_+$, $\wh \bu_+$, $\wh \bw_-$, $\wh \pi_-$, and $\wh h$ satisfy the equations
\begin{align*}
\left\{\begin{aligned}
\lambda \wh \rho_+ + \rho_{* +} \dv \wh \bu_+ & = \wh f_M &\quad &\text{ in $\RP$,} \\
\rho_{* +} \lambda \wh \bu_+ - \DV \bT_{0 +} (\gamma_{10}, \gamma_{20}, \gamma_{30}, \wh \bu_+, \wh \rho_+) & = \wh \bff_+ &\quad &\text{ in $\RP$}, \\
\dv \wh \bw_- & = 0 &\quad &\text{ in $\RP$}, \\
\rho_{* -} \lambda \wh \bw_- - \DV \bT_- (\gamma_{40}, \wh \bw_-, \wh \pi_-) & = \wh \bff_- + \wh \bff'_- &\quad &\text{ in $\RM$}, \\
\lambda \wh h - \frac{\langle \rho_{* -} \wh \bw_-, \bn_0 \rangle \vert_- - \langle \rho_{* +} \wh \bu_+, \bn_0 \rangle \vert_+}{\rho_{* -} - \rho_{* +}} & = \wh d + \wh d' &\quad &\text{ on $\RZ$}, \\
\Pi_{\bn_0} (\gamma_{40} \bD (\wh \bw_-) \bn_0) \vert_- - \Pi_{\bn_0} (\gamma_{10} \bD (\wh \bu_+) \bn_0) \vert_+ & = \wh g + \wh g' &\quad &\text{ on $\RZ$}, \\
\langle \bT_- (\gamma_{40}, \wh \bw_-, \wh \pi_-) \bn_0, \bn_0 \rangle \Big\vert_- - \frac{\rho_{* -} \sigma}{\rho_{* -} - \rho_{* +}} \Delta_\RZ \wh h & = \wh g_- + \wh g'_- &\quad &\text{ on $\RZ$}, \\
\langle \bT_{0 +} (\gamma_{10}, \gamma_{20}, \gamma_{30}, \wh \bu_+, \wh \rho_+) \bn_0, \bn_0 \rangle \Big\vert_+ - \frac{\rho_{* +} \sigma}{\rho_{* -} - \rho_{* +}} \Delta_\RZ \wh h & = \wh g_+ &\quad &\text{ on $\RZ$}, \\
\Pi_{\bn_0} \wh \bw_- \vert_- - \Pi_{\bn_0} \wh \bu_+ \vert_+ & = \wh \bh + \wh \bh' &\quad &\text{ on $\RZ$}, \\
\langle \nabla \wh \rho_+, \bn_0 \rangle \vert_+ & = \wh k_- &\quad &\text{ on $\RZ$},
\end{aligned}\right.
\end{align*}
where we have set
\begin{align*}
\wh \bff'_- = & - \rho_{* -} \lambda \wh \bu_\mathrm{div} + \DV(\gamma_{40} \bD (\wh \bu_\mathrm{div})), \quad \wh d' = - \frac{\langle \rho_{* -} \wh \bu_\mathrm{div}, \bn_0 \rangle \vert_-}{\rho_{* -} - \rho_{* +}}, \\
\wh g' = & - \Pi_{\bn_0}(\gamma_{40} \bD (\wh \bu_\mathrm{div}) \bn_0) \vert_-, \quad \wh g'_- = - \langle \gamma_{40} \bD (\wh \bu_\mathrm{div}) \bn_0, \bn_0 \rangle \vert_-, \quad \wh \bh' = - \Pi_{\bn_0} \wh \bu_\mathrm{div} \vert_-.
\end{align*}
From Lemmas \ref{th-4.5} and \ref{lem-4.6}, we obtain {\color{red} $\wh \rho_+ = \CA^+_{F1} (\lambda) \bF'$, $\wh \bu_+ = \CB^+_{F1} (\lambda) \bF'$, $\wh \bu_- = \wh \bu_\mathrm{div} + \CB^-_{F1} (\lambda) \bF'$,} $\wh \pi_- = \CP^-_{F1} (\lambda) \bF'$, $\wh h = \CH_{F1} (\lambda) \bF'$ with $\wh \bu_\mathrm{div} = \CD_{F2} (\lambda) (\wh \bff_-, \lambda \wh g_-, \wh g_-)$ and
\begin{align*}
\bF' = (\wh f_M, \wh \bff_+, \wh \bff_- + \wh \bff'_-, \wh d + \wh d', \wh g + \wh g', \wh g_+, \wh g_- + \wh g'_-, \wh \bh + \wh \bh', \wh k_-).
\end{align*}
From the argument in Sect. \ref{sect-3.1} we see that $\wh \pi_- = \CK_{F1} (\wh \bu_-) + \CK_{F2} (\wh h)$, and then $(\wh \rho_+, \wh \bu_+, \wh \bu_-, \wh h)$ is a solution to \eqref{eq-4.5}. Hence, we can define operators $\CA^+_{F0} (\lambda)$, $\CB^\pm_{F0} (\lambda)$, and $\CH_{F0} (\lambda)$ by
\begin{align*}
\CA^+_{F0} (\lambda) (G) & = \CA^+_{F1} (\lambda) (G) + \CA^+_{F1} (\lambda) (\CF'), \\
\CB^\pm_{F0} (\lambda) (G) & = \CB^\pm_{F1} (\lambda) (G) + \CB^\pm_{F1} (\lambda) (\CF'), \\
\CH_{F0} (\lambda) (G) & = \CH_{F1} (\lambda) (G) + \CH_{F1} (\lambda) (\CF'),
\end{align*}
where we have set
\begin{align*}
G & = (G_1, G_2, G_3, G_4, G_5, G_6, G_7, G_8, G_9, G_{10}, G_{11}, G_{12}, G_{13}, G_{14}, G_{15}, G_{16}, G_{17}, G_{18}), \\
\CF' & = (0, 0, 0, \CF^1, \CF^2, \lambda^{1/2} \CF^3, \nabla \CF^3, 0, 0, \lambda^{1/2} \CF^4, \nabla \CF^4, \CF^4, \lambda \CF^5, \lambda^{1/2} \nabla \CF^5, \nabla^2 \CF^5, 0, 0, 0), \\
\CF^1 & = - \rho_{* -} \lambda \CD_{F2} (\lambda) (G_4, G_{10}, G_{12}) + \DV (\gamma_{40} \bD (\CD_{F2} (\lambda) (G_4, G_{10}, G_{12}))), \\
\CF^2 & = - \frac{\langle \rho_{* -} \CD_{F2} (\lambda) (G_4, G_{10}, G_{12}), \bn_0 \rangle \vert_-}{\rho_{* -} - \rho_{* +}}, \quad \CF^3 = - \Pi_{\bn_0} \{\gamma_{40} \bD (\CD_{F2} (\lambda) (G_4, G_{10}, G_{12})) \bn_0\} \vert_-, \\
\CF^4 & = - \langle \gamma_{40} \bD (\CD_{F2} (\lambda) (G_4, G_{10}, G_{12})) \bn_0, \bn_0 \rangle \Big\vert_-, \quad \CF^5 = - \Pi_{\bn_0} \CD_{F2} (\lambda) (G_4, G_{10}, G_{12}) \vert_-.
\end{align*}
By Lemmas \ref{lem-4.2}, \ref{th-4.5}, and \ref{lem-4.6}, operators $\CA^+_{F0}(\lambda)$, $\CB^\pm_{F0} (\lambda)$, $\CH_{F0} (\lambda)$ satisfy the required properties in Theorem \ref{th-4.4}. Hence, we complete the proof of Theorem \ref{th-4.4}. \par
Using a localization argument, we may reduce the problem \eqref{eq-1.5} to the following model problems:
\begin{equation*}
\text{(i) Whole space problem;} \quad \text{(ii) Half space problem;} \quad \text{(iii) Two-phase problem in a whole space}.
\end{equation*}
For a detailed explanation of localization argument, the reader may refer Maryani and Saito~\cite{MS2017} for the two-phase flows case. The first two types of model problems have been studied by Saito~\cite{Sai2019} and Shibata~\cite{Shi2014}, respectively, and the case (iii) is treated in Theorem~\ref{th-4.4}. \textcolor{red}{We will show that these results make us to obtain Theorem \ref{TH-3.1} in the rest of this section.} 

{\color{blue}
\subsection{Reduced probelm with a bent interface}
\input{localized}
}

\section{Maximal $L_p - L_q$ regularity theorem}
\label{sect-7}
\noindent
The purpose of this section is to prove the following theorem.
\begin{theo}
\label{th-1.5}
Let $1 < p, q < \infty$ with $2/p + N/q \neq 1$ and $2/p + N/q \neq 2$. In addition, let $N < r < \infty$ and $\max(q, q') \leq r$. Suppose that the assumptions (a)--(c) in Theorem~\ref{th-1.9} holds. Then there exists a constant $\gamma_0 \geq 1$ such that the following
assertions hold:
\begin{enumerate}\renewcommand{\labelenumi}{(\arabic{enumi})}
\item  Let $\rho_{0 +} \in B^{3 - 2/p}_{q, p}$, $\bu_{0 \pm} \in B^{2(1 - 1/p)}_{q, p} (\Omega_\pm)$, and $h_0 \in B^{3 - 1/p - 1/q}_{q, p} (\Gamma)$. Furthermore, let $f_M$, $\bff_\pm$, $\bF_d$, $f_d$, $d$, $g$, $f^+_B$, $f^-_B$, $\bh$, $k_-$, and $k_+$ be functions in the right-hand members of \eqref{eq-1.5} such that
\begin{equation}
\label{1.11*}
\begin{gathered}
f_M \in L_{p, \gamma} (\BR, W^1_q (\Omega_+)), \quad \bff_\pm \in L_{p, \gamma} (\BR, L_q(\Omega_\pm)^N), \quad \bF_d \in W^1_{p, \gamma} (\BR, L_q (\Omega_-)^N), \\
f_d \in L_{p, \gamma} (\BR, \CD\CI_q (\Omega_-)^N) \cap H^{1/2}_{p, \gamma} (\BR, L_q (\Omega_-)), \enskip d \in L_{p, \gamma} (\BR, W^{2 - 1/q}_q (\Gamma)), \\
g, f^\pm_B, \in H^{1, 1/2}_{q, p, \gamma} (\dot \Omega \times \BR), \enskip \bh\in W^{2, 1}_{q, p,\gamma} (\dot \Omega \times \BR), \enskip k_\pm \in L_{p, \gamma} (\BR, W^2_q (\Omega_+))
\end{gathered}
\end{equation}
for any $\gamma \geq \gamma_0$. Assume that the compatibility conditions:
\begin{align}
\label{1.12*}
\bu_{0 -} - \bF_{d 0} \in J_q (\Omega_-), \qquad \dv \bu_{0 -} = \dv f_d \vert_{t = 0} \quad \text{ in $\Omega_-$}.
\end{align}
In addition, we assume the compatibility conditions:
\begin{align}
\label{1.13*}
\left\{\begin{aligned}
\Pi_\bn (\gamma_4 \bD(\bu_{0 -}) \bn) \vert_- - \Pi_\bn (\gamma_1 \bD(\bu_{0 +}) \bn) \vert_+ & = g &\quad &\text{ on $\Gamma$}, \\
\Pi_\bn \bu_{0 -} \vert_- - \Pi_\bn \bu_{0 +} \vert_+ & = \bh &\quad &\text{ on $\Gamma$}, \\
\langle \nabla \rho_{0 +}, \bn \rangle \vert_+ & = k_- &\quad &\text{ on $\Gamma$}, \\
\langle \nabla \rho_{0 +}, \bn_+ \rangle = k_+, \quad \bu_{0 +} & = 0 &\quad &\text{ on $\Gamma_+$}, \\
\bu_{0 -} & = 0 &\quad &\text{ on $\Gamma_-$}
\end{aligned}\right.
\end{align}
hold when $2/p + N/q < 1$, while we assume the compatibility conditions:
\begin{align}
\label{1.14*}
\left\{\begin{aligned}
\langle \nabla \rho_{0 +}, \bn_+ \rangle = k_+, \quad \bu_{0 +} & = 0 &\quad &\text{ on $\Gamma_+$}, \\
\bu_{0 -} & = 0 &\quad &\text{ on $\Gamma_-$}
\end{aligned}\right.
\end{align}
hold when $1 < 2/p + N/q < 2$. Then the equations \eqref{eq-1.5} admits unique solutions $(\rho_+, \bu_+, \bu_-, \pi_-, h)$ with
\begin{align*}
\rho_+ & \in W^{3, 1}_{q, p, \gamma} (\Omega_+ \times (0, \infty)) \quad \bu_\pm \in W^{2, 1}_{q, p, \gamma} (\Omega_\pm \times (0, \infty)), \\
\pi_- & \in L_{p, \gamma} ((0, \infty), W^1_q (\Omega_-) + \wh W^1_{q, \Gamma}(\Omega_-)), \\
h &\in L_{p, \gamma} ((0, \infty), W^{3 - 1/q}_q (\Gamma)) \cap W^1_{p, \gamma} ((0, \infty), W^{2 - 1/q}_q (\Gamma)).
\end{align*}
\item  The solution $(\rho_+, \bu_+, \bu_-, \pi_-, h)$ satisfies the following
estimate:
\begin{align*}
&\BI_{p, q} (\rho_+, \bu_+, \bu_-, \pi_-, h, \gamma; (0, \infty)) \\
& \leq C \Big\{\lVert \rho_{0 +} \rVert_{B^{3 - 2/p}_{q, p} (\Omega_+)} + \sum_{\ell = \pm} \lVert \bu_{0 -} \rVert_{B^{2(2 - 1/p)}_{q, p} (\Omega_\ell)} + \lVert h_0 \rVert_{B^{3 - 1/p - 1/q}_{q, p} (\Gamma)} + \lVert e^{-\gamma t} f_M \rVert_{L_p (\BR, W^1_q (\Omega_+))} \\
& \quad + \sum_{\ell = \pm} \lVert e^{-\gamma t} \bff_\ell \rVert_{L_p (\BR, L_q (\Omega_\ell))} + \lVert e^{-\gamma t} f_d \rVert_{L_p (\BR, W^1_q (\Omega_-))} + \lVert e^{-\gamma t} f_d \rVert_{H^{1/2}_p (\BR, L_q (\Omega_-))} \\
& \quad + \lVert e^{-\gamma t} \pd_t \bF_d \rVert_{L_p (\BR, L_q (\Omega_-))} + \lVert e^{-\gamma t} d \rVert_{L_p (\BR, W^{2 - 1/q}_q (\Gamma))} + \sum_{\ell = \pm} \lVert e^{-\gamma t} f^\ell_B \rVert_{L_p(\BR, L_q (\dot \Omega))} \\
& \quad + \lVert e^{-\gamma t} (g, \nabla \bh) \rVert_{H^{1/2}_p (\BR, L_q (\dot \Omega))} + \lVert e^{-\gamma t} (\nabla g, \pd_t \bh, \nabla^2 \bh) \rVert_{L_p (\BR, L_q (\dot \Omega))}\\
& \quad + \sum_{\ell = \pm} \Big(\lVert e^{-\gamma t} (f^\ell_B, \nabla k_\ell) \rVert_{H^{1/2}_p (\BR, L_q (\dot \Omega))} + \lVert e^{-\gamma t} (\nabla f^\ell_B, \pd_t k_\ell, \nabla^2 k_\ell) \rVert_{L_p (\BR, L_q (\dot \Omega))} \Big) \Big\}
\end{align*}
with some positive constant $C$ independent of $t$ and $\gamma$.
\end{enumerate}
\end{theo}

\subsection{Generation of $C^0$-analytic semigroup}
To prove Theorem \ref{th-1.5}, we first show a generation of $C^0$-analytic semigroup. For this purpose, we first introduce function spaces $X_q$ and $\CX_q$ defined by
\begin{align*}
X_q(\Omega_+,\Omega_-,\Gamma) & = \{(\wh f_M, \wh \bff_+,  \wh \bff_-, \wh d, \wh g, \wh f^+_B, \wh f^-_B, \wh \bh, \wh k_-, \wh k_+) \mid \wh f_M \in W^1_q (\Omega_+), \enskip \wh \bff_+ \in L_q (\Omega_+)^N, \enskip \wh \bff_- \in L_q (\Omega_-)^N, \\
& \quad \quad \wh d \in W^{2 - 1 / q}_q (\Gamma), \enskip \wh g, \enskip \wh f^+_B, \wh f^-_B \in W^1_q (\dot \Omega), \enskip \wh \bh \in W^2_q (\dot \Omega)^N, \enskip \wh k_-, \wh k_+ \in W^2_q (\Omega_+) \}, \\
\CX_q(\Omega_+,\Omega_-,\Gamma) & = \{(E_1, \dots, E_{21}) \mid E_1, E_2 \in L_q (\Omega_+), \enskip E_3 \in L_q (\Omega_+)^N, \enskip E_4 \in L_q (\Omega_-)^N, \enskip E_5 \in W^{2 - 1 / q}_q (\Gamma), \\
& \quad \quad E_6, E_8, E_{10} \in L_q (\dot \Omega), \enskip E_7, E_9, E_{11}, E_{13} \in L_q (\dot \Omega)^N, \enskip E_{12} \in W^1_q (\dot \Omega), \enskip E_{14} \in L_q (\dot \Omega)^{N^2},\\
& \quad \quad E_{15} \in L_q (\dot \Omega)^{N^3}, \enskip E_{16}, E_{19} \in L_q (\Omega_+), \enskip E_{17}, E_{20} \in L_q (\Omega_+)^N, \enskip E_{18}, E_{21} \in L_q (\Omega_+)^{N^2}\}.
\end{align*}
and we set
\begin{align*}
& \lVert (E_1, \dots, E_{21}) \rVert_{\CX_q (\Omega_+, \Omega_-, \Gamma)} \\
& = \sum_{m = 1}^3 \lVert E_m \rVert_{L_q (\Omega_+)} + \lVert E_4 \rVert_{L_q (\Omega_-)} + \lVert E_5 \rVert_{W^{2 - 1/q}_q (\Gamma)} + \sum_{n = 6}^{15} \lVert E_n \rVert_{L_q (\dot \Omega)} + \lVert \nabla E_{12} \rVert_{L_q (\dot \Omega)} + \sum_{l = 16}^{21} \lVert E_l \rVert_{L_q (\Omega_+)}
\end{align*}
for any $(E_1,\dots,E_{21})\in\CX_q$.
We then consider the generalized resolvent problem
\begin{align}
\label{eq-8.1*}
\left\{\begin{aligned}
\lambda \wh \rho_+ + \rho_{* +} \dv \wh \bu_+ & = \wh f_M &\enskip &\text{ in $\Omega_+ $}, \\
\rho_{*-} \dv \wh \bu_- = \wh f_d & = \rho_{*-} \dv \wh \bF_d &\enskip &\text{ in $\Omega_- $}, \\
\rho_{*+} \lambda \wh \bu_+ - \DV \bT_+ (\gamma_1, \gamma_2, \gamma_3, \wh \rho_+, \wh \bu_+) & = \wh \bff_+ &\enskip &\text{ in $\Omega_+$}, \\
\rho_{*-} \lambda \wh \bu_- -\DV \bT_- (\gamma_4, \wh \bu_-, \wh \pi_-) & = \wh \bff_- &\enskip &\text{ in $\Omega_-$}, \\
\lambda \wh h - \frac{\langle \rho_{*-} \wh \bu_-, \bn \rangle \rvert_- - \langle \rho_{*+} \wh \bu_+, \bn \rangle \rvert_+}{\rho_{* -} - \rho_{* +}} & = \wh d &\enskip &\text{ on $\Gamma$}, \\
\bB (\gamma_1, \gamma_2, \gamma_3, \gamma_4, \wh \rho_+, \wh \bu_+, \wh \bu_-, \wh \pi_-) & = \wh \bG, &\enskip &\text{ on $\Gamma$}, \\
\wh \bu_+ = 0, \qquad \langle \nabla \wh \rho_+, \bn_+ \rangle & = \wh k_+ &\enskip &\text{ on $\Gamma_+$}, \\
\wh \bu_- & = 0 &\enskip &\text{ on $\Gamma_-$}.
\end{aligned}\right.
\end{align}
{\color{red} Using the similar argument as in Sect. \ref{sect-3}, we can replace $\wh \pi_-$ by $\CK_1(\wh \bu_-) + \CK_2(\wh h) + \CK_3(\wh \rho_+, \wh h)$}. Here,  $\CK_3(\wh \rho_+, \wh h)$ is a unique solution to the following variational problem:
\begin{gather*}	
(\nabla \CK_3(\wh \rho_+, \wh h), \nabla \varphi)_{\Omega_-} = 0 \quad \text{ for any $\varphi \in \wh W^1_{q', \Gamma} (\Omega_-)$}, \\
\CK_3 (\wh \rho_+, \wh h) = - \frac{\gamma_{**}^+ \rho_{* +} \rho_{* -}}{\rho_{* -} - \rho_{* +}} \wh \rho_+ -\frac{\rho_{* -} \sigma}{\rho_{* -} - \rho_{* +}} \langle \Delta_\Gamma \bn, \bn \rangle \wh h \quad \text{ on $\Gamma$},
\end{gather*}
where $\langle \Delta_\Gamma \bn, \bn \rangle$ is a given function depending on $\Gamma$. Recalling Remark \ref{rem-1.3}, $\CK_3(\wh \rho_+, \wh h)$ can be defined by
\begin{align*}
\CK_3 (\wh \rho_+, \wh h) = & - K^3_\Gamma (\wh \rho_+, \wh h) + K_1 (\nabla K^3_\Gamma (\wh \rho_+, \wh h))
\end{align*}
with $K^3_\Gamma (\wh \rho_+, \wh h) = \bT_\Gamma (- \CK_3 (\wh \rho_+, \wh h))$ and belongs to $W^1_q (\Omega_-) + \wh W^1_{q, \Gamma} (\Omega_-)$ satisfying
\begin{align*}
\lVert \nabla \CK_3 (\wh \rho_+, \wh h) \rVert_{L_q(\Omega_-)} &\leq C \Big(\lVert \wh \rho_+ \rVert_{W^3_q (\Omega_+)} + \lVert \wh h \rVert_{W^{3 - 1/q}_q (\Gamma)} \Big).
\end{align*}
{\color{blue} Since the system \eqref{eq-8.1*} can be seen as a perturbation from \eqref{eq-4.1*}, we obtain the following theorem, where the proof is analog to the proof of Theorem \ref{TH-3.1}.}
\begin{theo}
\label{TH-8.1*}
Let $1 < q < \infty$, $N < r < \infty$ and $\max (q, q') \leq r$. Suppose that the assumptions (a)--(c) in Theorem~\ref{th-1.9} holds. Then there exists constant ${\color{red} \varepsilon_*} \in (0, \pi/2)$ such that for any $\varepsilon \in ({\color{red} \varepsilon_*}, \pi/2)$ there is a constant $\lambda_{**} \geq 1$ with the following properties hold true:
\begin{enumerate}\renewcommand{\labelenumi}{(\arabic{enumi})}
\item For any $\lambda \in \Sigma_{\varepsilon, \lambda_{**}}$, there exist operators $\widetilde \CA^+ (\lambda)$, $\widetilde \CB^\pm (\lambda)$, and $\widetilde \CH (\lambda)$ with
\begin{align*}
\widetilde \CA^\pm (\lambda) & \in \Hol (\Sigma_{\varepsilon, \lambda_{**}}, \CL(\CX_q(\Omega_+, \Omega_-, \Gamma), W^3_q (\Omega_+))), \\
\widetilde \CB^+ (\lambda) & \in \Hol (\Sigma_{\varepsilon, \lambda_{**}}, \CL(\CX_q(\Omega_+, \Omega_-, \Gamma), W^2_q (\Omega_\pm)^N)), \\
\widetilde \CH (\lambda) & \in \Hol (\Sigma_{\varepsilon, \lambda_{**}}, \CL(\CX_q(\Omega_+, \Omega_-, \Gamma), W^{3 - 1/q}_q (\Gamma))),
\end{align*}
such that for any $\bF_X = (\wh f_M, \wh \bff_+,  \wh \bff_-, \wh d, \wh g, \wh f^+_B, \wh f^-_B, \wh \bh, \wh k_-, \wh k_+) \in Y_q(\Omega_+, \Omega_-,\Gamma)$, the quadruple 
\begin{equation*}
(\rho_+, \bu_+, \bu_-, h) = (\widetilde \CA^+ (\lambda) \widetilde F_\lambda (\bF_X), \widetilde \CB^+ (\lambda) \widetilde F_\lambda (\bF_X), \widetilde \CB^- (\lambda) \widetilde F_\lambda (\bF_X), \widetilde \CH (\lambda) \widetilde F_\lambda (\bF_X))
\end{equation*}
is a unique solution to~\eqref{eq-3.8}. Here, we have set
\begin{align*}
\widetilde F_\lambda (\bF_X) & = (\lambda^{1/2} \wh f_M, \nabla \wh f_M, \wh \bff_+, \wh \bff_-, \wh d, \lambda^{1/2} \wh g, \nabla \wh g, \lambda^{1/2} \wh f^+_B, \nabla \wh f^+_B, \wh f^+_B, \lambda^{1/2} \wh f^-_B, \nabla \wh f^-_B, \\
& \qquad \wh f^-_B, \lambda \wh \bh, \lambda^{1/2} \nabla \wh \bh, \nabla^2 \wh \bh, \lambda \wh k_-, \lambda^{1/2} \nabla \wh k_-, \nabla^2 \wh k_-, \lambda \wh k_+, \lambda^{1/2} \nabla \wh k_+, \nabla^2 \wh k_+).
\end{align*}
\item There exists a positive constant $c_{**}$, independent of $\lambda$, such that
\begin{equation*}
\begin{split}
\CR_{\CL(\CX_q(\Omega_+, \Omega_-, \Gamma), W^{3 - i}_q (\Omega_+))} (\{(\tau \pd_\tau)^s (\lambda^{i/2} \widetilde \CA^+ (\lambda)) \mid \lambda \in \Sigma_{\varepsilon, \lambda_{**}}\}) & \leq c_{**}, \\
\CR_{\CL(\CX_q(\Omega_+, \Omega_-, \Gamma), W^{2 - j}_q (\Omega_\pm)^N)} (\{(\tau \pd_\tau)^s (\lambda^{j/2} \widetilde \CB^\pm (\lambda)) \mid \lambda \in \Sigma_{\varepsilon, \lambda_{**}}\}) & \leq c_{**}, \\
\CR_{\CL(\CX_q(\Omega_+, \Omega_-, \Gamma), W^{3 - 1/q - k}_q (\Gamma))} (\{(\tau \pd_\tau)^s (\lambda^k \widetilde \CH (\lambda)) \mid \lambda \in \Sigma_{\varepsilon, \lambda_{**}}\}) & \leq c_{**},
\end{split}
\end{equation*}	
for $s = 0, 1$, $i = 0, 1, 2, 3$, $j = 0, 1, 2$, and $k = 0, 1$.	
\end{enumerate}
\end{theo}

We now consider the following homogeneous problem:
\begin{align}
\label{eq-7.1}
\left\{\begin{aligned}
\pd_t \rho_+ + \rho_{*+} \dv \bu_+ & = 0 &\enskip &\text{ in $\Omega_+ \times(0, \infty)$}, \\
\rho_{*+}\pd_t\bu_+ -\DV \bT_+ (\gamma_1, \gamma_2, \gamma_3, \rho_+, \bu_+) & = 0 &\enskip &\text{ in $\Omega_+ \times (0, \infty)$}, \\
\rho_{*-}\pd_t\bu_--\DV \bT_- (\gamma_4, \bu_-, \CK_1(\bu_-) + \CK_2(h) + \CK_3(\rho_+, h)) & = 0 &\enskip &\text{ in $\Omega_- \times (0, \infty)$}, \\
\pd_t h - \frac{\langle \rho_{*-} \bu_-, \bn \rangle \rvert_- - \langle \rho_{*+} \bu_+, \bn \rangle \rvert_+}{\rho_{*-} - \rho_{*+}} 
& = 0 &\enskip &\text{ on $\Gamma \times (0, \infty)$}, \\
\bB (\gamma_1, \gamma_2, \gamma_3, \gamma_4, \rho_+, \bu_+, \bu_-, \CK_1(\bu_-) + \CK_2(h) + \CK_3(\rho_+, h))
& = 0, &\enskip &\text{ on $\Gamma \times (0, \infty)$}, \\
\bu_+ = 0, \qquad \langle \nabla \rho_+, \bn_+ \rangle & = 0
&\enskip &\text{ on $\Gamma_+ \times (0, \infty)$}, \\
\bu_- & = 0 &\enskip &\text{ on $\Gamma_- \times (0, \infty)$}, \\
(\rho_+, \bu_+, \bu_-, h) \rvert_{t = 0}= (\rho_{0+}, \bu_{0+}, & \bu_{0-},h_0)
&\enskip &\text{ on $\Omega_+ \times \Omega_+ \times \Omega_- \times \Gamma$}.
\end{aligned}\right.
\end{align}
We introduce an operator $A_q$ and its domain $D_q$ to formulate \eqref{eq-7.1} in the semigroup setting. Let
\begin{align*}
& D_q:=D_q(\Omega_+,\Omega_-,\Gamma)\\
& =\Set{(\rho_+, \bu_+, \bu_-, h) |
\begin{array}{l}
\rho_+ \in W^3_q (\Omega_+), \, \bu_+ \in W^2_q (\Omega_+)^N, \, \bu_- \in J_q (\Omega_-) \cap W^2_q (\Omega_+)^N, \, h \in W^{3 - 1/q}_q (\Gamma), \\
\text{$(\rho_+, \bu_+, \bu_-, h)$ satisfies the fifth, sixth, and seventh condition of \eqref{eq-7.1}}
\end{array}}, \\
& A_q:=A_q(\Omega_+,\Omega_-,\Gamma)=
\begin{pmatrix}
- \rho_{* +} \dv \bu_+ \\
\rho_{* +}^{- 1} \DV \{\gamma_1 \bD (\bu_+) + (\gamma_2 - \gamma_1) \dv \bu_- \bI + (- \gamma_{* +} \nabla + \rho_{* +} \gamma_3 \Delta)\rho_+\bI\} \\
\rho_{*-}^{- 1} \DV \{\gamma_4 \bD (\bu_-) - (\CK_1 (\bu_-) + \CK_2 (h) + \CK_3 (\rho_+, h)) \bI\} \\
\cfrac{\langle \rho_{* -} \bu_-, \bn \rangle \vert_- - \langle \rho_{* +} \bu_+, \bn \rangle \vert_+}{\rho_{* -} - \rho_{* +}}
\end{pmatrix},\\
& B_q := B_q (\Omega_+, \Omega_-, \Gamma) =
\Set{(\rho_+, \bu_+, \bu_-, h) \vert
\begin{array}{l}
\rho_+ \in W^1_q (\Omega_+), \enskip \bu_+ \in L_q (\Omega_+)^N, \\
\bu_- \in J_q (\Omega_-), \enskip h \in W^{2 - 1/q}_q (\Gamma)
\end{array}}.
\end{align*}
Recall that the space $J_q (\Omega_-)$ is the solenoidal space defined \eqref{def-sol} and $\bu_- \in J_q (\Omega_-) \cap W^2_q (\Omega_-)^N$ implies that $\rho_{* -} \dv \bu_- = 0$ in $\Omega_-$ and $\rho_{* -} \dv \bu_- = 0$ on $\Gamma$. Using the symbols given above, the problem \eqref{eq-7.1} can be rewritten as
\begin{align}
\label{eq-7.5}
\pd_t \CU (t) - A_q \CU (t) = 0 \quad (t > 0), \quad \CU (t) \vert_{t = 0} = \CU_0,
\end{align}
where $\CU (t) = (\rho_+, \bu_+, \bu_-, h) \in D_q$ for $t > 0$ and $\CU_0 = (\rho_{0 +}, \bu_{0 +}, \bu_{0 -}, h_0) \in B_q$. The corresponding resolvent problem to \eqref{eq-7.5} is that for any $F \in B_q$ we find $\CU \in D_q$ solving the equation:
\begin{align}
\label{eq-7.6}
\lambda \CU - A_q \CU = F \quad \text{ in $\Omega_+ \times \Omega_- \times \Omega_+ \times \Gamma$}
\end{align}
and possessing the estimate:
\begin{align}
\label{eq-7.7}
\lvert \lambda \rvert \lVert \CU \rVert_{B_q} + \lVert \CU \rVert_{D_q} \leq C \lVert F \rVert_{B_q}
\end{align}
for any $\lambda \in \Sigma_{\varepsilon, \lambda_{**}}$. Here, we have set
\begin{align*}
\lVert \CU \rVert_{B_q} = & \lVert \rho_+ \rVert_{W^1_q (\Omega_+)} + \lVert \bu \rVert_{L_q (\Omega_+)} + \lVert \bu_- \rVert_{L_q (\Omega_-)} + \lVert h \rVert_{W^{2 - 1/q}_q (\Gamma)}, \\
\lVert \CU \rVert_{D_q} = & \lVert \rho_+ \rVert_{W^3_q (\Omega_+)} + \lVert \bu \rVert_{W^2_q (\Omega_+)} + \lVert \bu_- \rVert_{W^2_q (\Omega_-)} + \lVert h \rVert_{W^{3 - 1/q}_q (\Gamma)}
\end{align*}
for $\CU = (\rho_+, \bu_+, \bu_-, h)$. Recalling that \eqref{eq-1.5} and \eqref{eq-8.1*} are equivalent, by Theorem \ref{TH-8.1*} the functions
\begin{align*}
\rho_+ & = \widetilde \CA^+ (\lambda) U_0, \quad \bu_\pm = \widetilde{\CB}^\pm(\lambda) U_0, \quad h = \widetilde \CH (\lambda) U_0 \\
U_0 & = (\rho_{0 +}, \bu_{0 +}, \bu_{0 -}, h_0, 0, 0, 0, 0, 0, 0, 0, 0, 0, 0, 0, 0, 0, 0, 0, 0)
\end{align*}
are solutions to \eqref{eq-7.6}. Since the $\CR$-boundedness implies the uniform boundedness, $\CU$ satisfies the resolvent estimate \eqref{eq-7.7}. Hence, by the standard semigroup theory, we obtain the following theorem.
\begin{theo}
\label{th-7.1}
Let $1 < q < \infty$. Let the assumptions (a)--(c) in Theorem~\ref{th-1.9} hold true. Then the operator $A_q$ defined in \eqref{eq-7.5} generates a $C^0$-analytic semigroup $\{T(t) \}_{t\geq 0}$ on $B_q$ satisfying the following estimates:
\begin{align*}
\lVert T (t) \CU_0 \rVert_{B_q} + t (\lVert \pd_t T (t) \CU_0 \rVert_{B_q} + \lVert T (t) \CU_0 \rVert_{D_q}) & \leq C e^{\gamma t} \lVert \CU_0 \rVert_{B_q}, \\
\lVert \pd_t T (t) \CU_0 \rVert_{B_q} + \lVert T (t) \CU_0 \rVert_{D_q} & \leq C e^{\gamma t} \lVert \CU_0 \rVert_{D_q}.
\end{align*}
\end{theo}

Applying the similar argument as in Shibata and Shimizu~\cite[Theorem 3.9]{SS2008}, by Theorem \ref{th-7.1} we have the following theorem.
\begin{theo}
\label{th-7.2}
Let $1 < p, q < \infty$ with $2/p + N/q \neq 1$ and $2/p + N/q \neq 2$. In addition, let $N < r < \infty$ and $\max (q, q') \leq r$. Let the assumptions (a)--(c) in Theorem~\ref{th-1.9} be valid. Define $\CD_{q, p} (\Omega_+, \Omega_-, \Gamma)$ as a subspace of $B^{3 - 2/p}_{q, p} (\Omega_+) \times B^{2 (1 - 1/p)}_{q, p} (\Omega_+)^N \times (J_q (\Omega_-) \times B^{2 (1 - 1/p)}_{q, p} (\Omega_-)^N) \times B^{3 - 1/p - 1/q}_{q, p} (\Gamma)$ defined by $\CD_{q, p} (\Omega_+, \Omega_-, \Gamma) = (B_q, D_q)_{1 - 1/p, p}$, where $(\cdot, \cdot)_{1 - 1/p, p}$ denotes a real interpolation functor. Then there exists a constant $\gamma>0$ such that the following assertions hold:
\begin{enumerate}\renewcommand{\labelenumi}{(\arabic{enumi})}
\item For any initial data $\CU_0 = (\rho_{0 +}, \bu_{0 +}, \bu_{0 -}, h_0) \in \CD_{q, p} (\Omega_+, \Omega_-, \Gamma)$, the problem \eqref{eq-7.1} admits a unique solution $(\rho_+, \bu_+, \bu_-,
h)$ with
\begin{align*}
\rho_+ & \in W^{3, 1}_{q, p, \gamma} (\Omega_+ \times (0, \infty)) \quad \bu_\pm \in W^{2, 1}_{q, p, \gamma} (\Omega_\pm \times (0, \infty)), \\
h & \in L_{p, \gamma} ((0, \infty), W^{3 - 1/q}_q (\Gamma)) \cap W^1_{p, \gamma} ((0, \infty), W^{2 - 1/q}_q (\Gamma)).
\end{align*}
\item The solution $\CU = (\rho_+, \bu_+, \bu_-, h)$ satisfies the following estimate:
\begin{align*}
\begin{split}
\lVert e^{-\gamma t} \pd_t \CU \rVert_{L_p ((0, \infty), B_q)} + \lVert e^{-\gamma t} \CU \rVert_{L_p ((0, \infty), D_q)} \leq C \lVert \CU_0 \rVert_{\CD_{q, p} (\Omega_+, \Omega_-, \Gamma)}.
\end{split}
\end{align*}
Here, the norm of $\CD_{q, p} (\Omega_+, \Omega_-, \Gamma)$ has been defined by
\begin{align*}
\lVert \CU_0 \rVert_{\CD_{q, p} (\Omega_+, \Omega_-, \Gamma)} = \lVert \rho_{0 +} \rVert_{B^{3 - 2/p}_{q, p} (\Omega_+)} + \lVert \bu_{0 +} \rVert_{B^{2(1 - 1/p)}_{q, p} (\Omega_+)} + \lVert \bu_{0 -} \rVert_{B^{2(1 - 1/p)}_{q, p} (\Omega_-)} + \lVert h_0 \rVert_{B^{3 - 1/p - 1/q}_{q, p} (\Gamma)}.
\end{align*}
\end{enumerate}
\end{theo}

\subsection{Maximal $L_p - L_q$ regularity theorem}
We now prove the maximal $L_p - L_q$ regularity theorem with the help of the operator-valued Fourier multiplier theorem~\cite{Wei2001}. Let {\color{red} $X$ be a Banach space}. We define spaces $\CD(\BR,X)$, $\CD'(\BR,X)$, $\CS(\BR,X)$, and $\CS'(\BR,X)$ by
\begin{align*}
\begin{aligned}
\CD (\BR, X) & \colon \text{the set of all $X$-valued $C^\infty$ functions having compact supports}, \\
\CD' (\BR, X) &\colon \text{the set of all linear bounded operator from $\CD(\BR,\BC)$ to $X$},\\
\CS (\BR, {\color{red} X}) & \colon \text{{\color{red} the Schwartz class of rapidly decreasing smooth functions from $\BR$ into $X$}},\\
\CS' (\BR, X) & \colon \text{{\color{red} $X$-valued tempered distribution}},
\end{aligned}
\end{align*}
respectively. For given $M \in L_{1, \mathrm{loc}} (\BR, \CL(X, Y))$, we define an Fourier multiplier $T_M \colon \CF^{- 1} \CD (\BR, X) \to \CS' (\BR, Y)$ by $T_M \phi = \CF^{- 1} [M \CF [\phi]]$, where $\CF [\phi] \in \CD (\BR, X)$ {\color{blue} and $X$, $Y$ are Banach spaces}.

\begin{theo}
\label{th-7.4}
Suppose that $X$ and $Y$ are $\mathrm{UMD}$ spaces and let $1<q<\infty$. Let $M$ be a function in $C^1 (\BR \backslash \{0\}, \CL (X, Y))$ such that
\begin{align*}
\CR (\{(\rho \pd_\rho)^k M (\rho) \mid \rho \in \BR \backslash \{0\} \}) & = \kappa_k < \infty
\end{align*}
for $k = 0, 1$. Then the operator $T_M$ defined above is extended to a bounded operator from $L_p (\BR, X)$ into $L_p (\BR, Y)$ with norm
\begin{align*}
\lVert T_M \rVert_{\CL (L_p (\BR, X), L_p (\BR, Y))} \leq C (\kappa_0 + \kappa_1),
\end{align*}
where $C > 0$ depends only on $p$, $X$, and $Y$. Here, a Banach space is said to be a $\mathrm{UMD}$ space if the Hilbert transform extends to bounded operator on $L_p(\BR, X)$ for some $1 < p < \infty$.
\end{theo}

We now consider the following problem:
\begin{align}
\label{eq-7.10}
\left\{\begin{aligned}
\pd_t \rho_{1 +} + \rho_{*+} \dv \bu_{1 +} & = f_M &\enskip &\text{ in $\Omega_+ \times \BR$}, \\
\rho_{*-} \dv \bu_{1 -} = f_d & = \rho_{*-} \dv \bF_d &\enskip &\text{ in $\Omega_- \times \BR$}, \\
\rho_{*+}\pd_t \bu_{1 +} -\DV \bT_+ (\gamma_1, \gamma_2, \gamma_3, \rho_{1 +}, \bu_{1 +}) & = \bff_+ &\enskip &\text{ in $\Omega_+ \times \BR$}, \\
\rho_{*-}\pd_t\bu_{1 -} -\DV \bT_- (\gamma_4, \bu_{1 -}, \pi_{1 -}) & = \bff_- &\enskip &\text{ in $\Omega_- \times \BR$}, \\
\pd_t h - \frac{\langle \rho_{*-} \bu_{1 -}, \bn \rangle \rvert_- - \langle \rho_{*+} \bu_{1 +}, \bn \rangle \rvert_+}{\rho_{*-} - \rho_{*+}} & = d &\enskip &\text{ on $\Gamma \times \BR$}, \\
\bB (\gamma_1, \gamma_2, \gamma_3, \gamma_4, \rho_{1 +}, \bu_{1 +}, \bu_{1 -}, \pi_{1 -}) & = \bG, &\enskip &\text{ on $\Gamma \times \BR$}, \\
\bu_{1 +} = 0, \qquad \langle \nabla \rho_{1 +}, \bn_+ \rangle & = k_+ &\enskip &\text{ on $\Gamma_+ \times \BR$}, \\
\bu_{1 -} & = 0 &\enskip &\text{ on $\Gamma_- \times \BR$}.
\end{aligned}\right.
\end{align}
The right-hand members of \eqref{eq-7.10}: $f_M, \bff_+, \bff_-, d, g, g_+, g_-, \bh, k_-, k_+$ are defined on $t \in \BR$. Let $\CL_L$ be the Laplace transform with respect to time variable $t$ defined by $\hf (\lambda) = \CL_L[f] (\lambda) = \int_\BR e^{- \lambda t} f(t) \dt$ for $\lambda = \gamma + i \tau \in \BC$. Applying the Laplace transform to \eqref{eq-7.10} gives
\begin{align}
\label{eq-7.11}
\left\{\begin{aligned}
\lambda \wh \rho_{1 +} + \rho_{*+} \dv \wh \bu_{1 +} & = \wh f_M &\enskip &\text{ in $\Omega_+ $}, \\
\rho_{*-} \dv \wh \bu_{1 -} = \wh f_d & = \rho_{*-} \dv \wh \bF_d &\enskip &\text{ in $\Omega_- $}, \\
\rho_{*+} \lambda \wh \bu_{1 +} - \DV \bT_+ (\gamma_1, \gamma_2, \gamma_3, \wh \rho_{1 +}, \wh \bu_{1 +}) & = \wh \bff_+ &\enskip &\text{ in $\Omega_+$}, \\
\rho_{*-} \lambda \wh \bu_{1 -} -\DV \bT_- (\gamma_4, \wh \bu_{1 -}, \wh \pi_{1 -}) & = \wh \bff_- &\enskip &\text{ in $\Omega_-$}, \\
\lambda \wh h - \frac{\langle \rho_{*-} \wh \bu_{1 -}, \bn \rangle \rvert_- - \langle \rho_{*+} \wh \bu_{1 +}, \bn \rangle \rvert_+}{\rho_{*-} - \rho_{*+}} & = \wh d &\enskip &\text{ on $\Gamma$}, \\
\bB (\gamma_1, \gamma_2, \gamma_3, \gamma_4, \wh \rho_{1 +}, \wh \bu_{1 +}, \wh \bu_{1 -}, \wh \pi_{1 -}) & = \wh \bG, &\enskip &\text{ on $\Gamma$}, \\
\wh \bu_{1 +} = 0, \qquad \langle \nabla \wh \rho_{1 +}, \bn_+ \rangle & = \wh k_+ &\enskip &\text{ on $\Gamma_+$}, \\
\wh \bu_{1 -} & = 0 &\enskip &\text{ on $\Gamma_-$}.
\end{aligned}\right.
\end{align}
From Theorem \ref{TH-8.1*} we have the representation of $\wh \bu_{1 \pm}$, $\wh \rho_{1 +}$, and $\wh \pi_{1 -}$ such that
\begin{align*}
\wh \rho_{1 +} = \widetilde \CA^+ (\lambda) \widetilde F_\lambda (\wh \bF_X), \quad \wh \bu_{1 \pm} = \widetilde \CB^\pm (\lambda) \widetilde F_\lambda (\wh \bF_X), \quad \wh \pi_{1 -} = \widetilde \CP^- (\lambda) \widetilde F_\lambda (\wh \bF_X),\quad \hh_1 = \widetilde \CH (\lambda) \widetilde F_\lambda (\wh \bF_X)
\end{align*}
for $\lambda \in \Sigma_{\varepsilon, \lambda_{**}}$, where $\wh \bF_X = (\wh f_M, \wh \bff_+,  \wh \bff_-, \wh d, \wh g, \wh f^+_B, \wh f^-_B, \wh \bh, \wh k_-, \wh k_+)$. Here, $\widetilde \CP^-$ is the operator belongs to
$\Hol(\Sigma_{\varepsilon, \lambda_{**}}, \CL(\CX_q (\Omega_+, \Omega_-, \Gamma), W^1_q (\Omega_-) + \wh W^1_{q, \Gamma} (\Omega_-)))$ such that
\begin{align*}
\CR_{\CL(\CX_q (\Omega_+, \Omega_-, \Gamma), L_q (\Omega_\pm))} (\{(\tau \pd_\tau)^s (\nabla \widetilde \CP^-(\lambda)) \mid \lambda \in \Sigma_{\varepsilon, \lambda_{**}}\}) \leq c_{**}.
\end{align*}
Let $\CL_L^{- 1}$ be the inverse Laplace transform defined by $\CL_L^{- 1} [f] (t) = (2 \pi)^{- 1} \int_{\BR} e^{\lambda t} f (\tau) \dtau$ for $\lambda = \gamma + i \tau \in \BC$.
Setting
\begin{align*}
\Lambda^{1/2}_\gamma f(t)=\CL_L^{-1}[\lambda^{1/2}\CL_L[f]](t)
=e^{\gamma t}\CF^{-1}[\lambda^{1/2}\CF[e^{-\gamma t}f]](t),
\end{align*}
and using the fact that
$\lambda \widehat{f}_1(\lambda)=\CL_L[\pd_t f](\lambda)$ and $\lambda^{1/2}  \widehat{f}_2(\lambda) = \CF[e^{-\gamma t}\Lambda^{1/2}_\gamma f](\tau)$, we define $\rho_{1+}$, $\bu_{1\pm}$, $\pi_{1-}$, and $h_1$ by
\begin{align*}
\rho_{1 +} (\cdot, t) = & \CL_L^{- 1} [\widetilde \CA^+ (\lambda) \widetilde F_\lambda(\widehat \bF_X)] = e^{\gamma t} \CF^{- 1} [\widetilde \CA^+ (\lambda) \CF[e^{-\gamma t} F (t)] (\tau)], \\
\bu_{1 \pm} (\cdot, t) = & \CL_L^{- 1} [\widetilde \CB^\pm (\lambda) \widetilde F_\lambda(\widehat \bF_X)] = e^{\gamma t} \CF^{- 1} [\widetilde \CB^\pm (\lambda) \CF[e^{-\gamma t} F(t)] (\tau)], \\
\pi_{1 -} (\cdot, t) = & \CL_L^{- 1} [\widetilde \CP^- (\lambda) \widetilde F_\lambda(\widehat \bF_X)] = e^{\gamma t} \CF^{- 1} [\widetilde \CP^- (\lambda) \CF[e^{-\gamma t} F(t)] (\tau)], \\
h_1 (\cdot, t) = & \CL_L^{- 1} [\widetilde \CH (\lambda) \widetilde F_\lambda(\widehat \bF_X)] = e^{\gamma t} \CF^{- 1} [\widetilde \CH (\lambda) \CF[e^{-\gamma t} F(t)] (\tau)]
\end{align*} 
with 
\begin{align*}
F(t) & = (f_M, \bff_+, \bff_-, f_d, \Lambda^{1/2}_\gamma f_d, \pd_t \bF_d, d, \Lambda^{1/2}_\gamma g, \nabla g, \Lambda^{1/2}_\gamma f^+_B, \nabla f^+_B, f^+_B, \Lambda^{1/2}_\gamma f^-_B, \nabla f^-_B,\\
&\quad f^-_B, \pd_t \bh, \Lambda^{1/2}_\gamma \nabla \bh, \nabla^2 \bh, \pd_t k_-, \Lambda^{1/2}_\gamma \nabla k_-, \nabla^2 k_-, \pd_t k_+, \Lambda^{1/2}_\gamma \nabla k_+, \nabla^2 k_+),
\end{align*}
where $\gamma$ is chosen such that $\gamma > \lambda_{**}$ holds, which implies $\lambda = \gamma + i \tau \in \Sigma_{\varepsilon, \lambda_{**}}$ for any $\tau \in \BR$. From the Cauchy theorem, $\rho_{1 +}$, $\bu_{1 \pm}$, $\pi_{1 -}$, and $h_1$ are independent of choice of $\gamma$ whenever $\gamma > \lambda_{**}$ and \eqref{1.11*} satisfied for $\gamma > \lambda_{**}$. Noting that
\begin{align*}
\pd_t\rho_{1+}(\cdot,t)=&\CL_L^{-1}[\lambda \widetilde\CA^+(\lambda) \widetilde F_\lambda(\widehat \bF_X)]
=e^{\gamma t}\CF^{-1}[\lambda \widetilde\CA^+(\lambda)\CF[e^{-\gamma t}F(t)](\tau)],\\
\pd_t\bu_{1\pm}(\cdot,t)=&\CL_L^{-1}[\lambda \widetilde\CB^\pm(\lambda) \widetilde F_\lambda(\widehat \bF_X)]
=e^{\gamma t}\CF^{-1}[\lambda \widetilde\CB^\pm(\lambda)\CF[e^{-\gamma t}F(t)](\tau)],\\
\pd_t\pi_{1-}(\cdot,t)=&\CL_L^{-1}[\lambda \widetilde\CP^-(\lambda) \widetilde F_\lambda(\widehat \bF_X)]
=e^{\gamma t}\CF^{-1}[\lambda \widetilde\CP^-(\lambda)\CF[e^{-\gamma t}F(t)](\tau)],\\
\pd_t h_1(\cdot,t)=&\CL_L^{-1}[\lambda \widetilde\CH(\lambda) \widetilde F_\lambda(\widehat \bF_X)]
=e^{\gamma t}\CF^{-1}[\lambda \widetilde\CH(\lambda)\CF[e^{-\gamma t}F(t)](\tau)]
\end{align*}
and applying Theorem \ref{th-7.4}, we have
\begin{equation}\label{7.12}
\begin{split}
& \lVert e^{-\gamma t} \pd_t \rho_{1 +} \rVert_{L_p ((0, \infty), W^1_q (\Omega_+))} + \lVert e^{-\gamma t} \rho_{1 +} \rVert_{L_p ((0, \infty), W^3_q (\Omega_+))} + \lVert e^{-\gamma t}\pd_t \bu_{1 +} \rVert_{L_p ((0, \infty), L_q (\Omega_+))} \\
& \quad + \lVert e^{-\gamma t} \bu_{1 +} \rVert_{L_p ((0, \infty), W^2_q (\Omega_+))} + \lVert e^{-\gamma t} \pd_t \bu_{1 -} \rVert_{L_p ((0, \infty), L_q (\Omega_-))} + \lVert e^{-\gamma t} \bu_{1 -} \rVert_{L_p ((0, \infty), W^2_q (\Omega_-))} \\
& \quad + \lVert e^{-\gamma t} \nabla \pi_{1 -} \rVert_{L_p ((0, \infty), L_q (\Omega_-))} + \lVert e^{-\gamma t} \pd_t h_1 \rVert_{L_p ((0, \infty), W^{2 - 1/q}_q (\Gamma))} + \lVert e^{-\gamma t} h_1 \rVert_{L_p ((0, \infty), W^{3 - 1/q}_q (\Gamma))} \\
& \leq C \lVert e^{\gamma t} F \rVert_{L_p (\BR, B_q)} \\
&\leq C \Big\{\lVert \rho_{0 +} \rVert_{B^{3 - 2/p}_{q, p} (\Omega_+)} + \sum_{\ell = \pm} \lVert \bu_{0 -} \rVert_{B^{2(2 - 1/p)}_{q, p} (\Omega_\ell)} + \lVert h_0 \rVert_{B^{3 - 1/p - 1/q}_{q, p} (\Gamma)} + \lVert e^{-\gamma t} f_M \rVert_{L_p (\BR, W^1_q (\Omega_+))} \\
& \quad + \sum_{\ell = \pm} \lVert e^{-\gamma t} \bff_\ell \rVert_{L_p (\BR, L_q (\Omega_\ell))} + \lVert e^{-\gamma t} f_d \rVert_{L_p (\BR, W^1_q (\Omega_-))} + \lVert e^{-\gamma t} f_d \rVert_{H^{1/2}_p (\BR, L_q (\Omega_-))} \\
& \quad + \lVert e^{-\gamma t} \pd_t \bF_d \rVert_{L_p (\BR, L_q (\Omega_-))} + \lVert e^{-\gamma t} d \rVert_{L_p (\BR, W^{2 - 1/q}_q (\Gamma))} + \sum_{\ell = \pm} \lVert e^{-\gamma t} f^\ell_B \rVert_{L_p(\BR, L_q (\dot \Omega))} \\
& \quad + \lVert e^{-\gamma t} (g, \nabla \bh) \rVert_{H^{1/2}_p (\BR, L_q (\dot \Omega))} + \lVert e^{-\gamma t} (\nabla g, \pd_t \bh, \nabla^2 \bh) \rVert_{L_p (\BR, L_q (\dot \Omega))}\\
& \quad + \sum_{\ell = \pm} \Big(\lVert e^{-\gamma t} (f^\ell_B, \nabla k_\ell) \rVert_{H^{1/2}_p (\BR, L_q (\dot \Omega))} + \lVert e^{-\gamma t} (\nabla f^\ell_B, \pd_t k_\ell, \nabla^2 k_\ell) \rVert_{L_p (\BR, L_q (\dot \Omega))} \Big) \Big\}
\end{split}
\end{equation}
with some positive constant $C$ independent of $\gamma$ and $t$. \par
We write a solution $(\rho_+, \bu_+, \bu_-, \pi_-, h)$ to \eqref{eq-1.6} in the form of $\rho_+ = \rho_{1 +} + \rho_{2 +}$, $\bu_\pm = \bu_{1 \pm} + \bu_{2 \pm}$, $\pi_- = \pi_{1 -} + \pi_{2 -}$, and $h = h_1 + h_2$. Then $\rho_{2 +}$, $\bu_{2 \pm}$, $\pi_{2 -}$, and $h_2$ enjoy the homogeneous equations \eqref{eq-7.1} with $\rho_{2 +} = \rho_+$, $\bu_{2 \pm} = \bu_\pm$, $\pi_{2 -} = \CK_1 (\bu_-) + \CK_2 (h) + \CK_3 (\rho_+, h)$. Here, $\rho_{* -} \dv \bu_{2 -} = 0$ in $\Omega_- \times (0, \infty)$ means that $\bu_{2 -}$ belongs to $J_q (\Omega_-)$ for any $t > 0$. We know that
\begin{equation}
\label{7.14}
\begin{split}
\sup_{t\in (0,\infty)} e^{-\gamma t} \lVert \rho_{1 +} \rVert_{B^{3 - 2/p}_{q, p} (\Omega_+)} & \leq C \Big(\lVert e^{-\gamma t} \rho_{1 +} \rVert_{L_p ((0, \infty), W^3_q (\Omega_+))} + \lVert e^{-\gamma t} \pd_t \rho_{1 +} \rVert_{L_p ((0, \infty), W^1_q (\Omega_\pm))} \Big), \\
\sup_{t \in (0, \infty)} e^{-\gamma t} \lVert \bu_{1 \pm} \rVert_{B^{2(1 - 1/p)}_{q, p} (\Omega_\pm)} & \leq C \Big(\lVert e^{-\gamma t} \bu_{1 \pm} \rVert_{L_p ((0,\infty), W^2_q (\Omega_\pm))} + \lVert e^{-\gamma t} \pd_t \bu_{1 \pm} \rVert_{L_p ((0, \infty), L_q (\Omega_\pm))} \Big), \\
\sup_{t \in (0, \infty)} e^{-\gamma t} \lVert h_1 \rVert_{B^{3 - 1/p - 1/q}_{q, p} (\Gamma)} & \leq C \Big(\lVert e^{-\gamma t} h_1 \rVert_{L_p ((0, \infty), W^{3 - 1/q}_q (\Gamma))} + \lVert e^{-\gamma t} \pd_t h_1 \rVert_{L_p ((0, \infty), W^{2 - 1/q}_q (\Gamma))} \Big)
\end{split}
\end{equation}
and
\begin{align*}
& (\rho_{0 +} - \rho_{1 +} \vert_{t = 0}, \enskip \bu_{0 +} - \bu_{1 +} \vert_{t = 0}, \enskip \bu_{0 -} - \bu_{1 -} \vert_{t = 0}, \enskip h_0 - h_1 \vert_{t = 0}) \\
& \quad \in B^{3 - 2/p}_{q, p} (\Omega_+) \times B^{2(1 - 1/p)}_{q, p} (\Omega_+) \times B^{2(1 - 1/p)}_{q, p} (\Omega_-) \times B^{3 - 1/p - 1/q}_{q, p} (\Gamma),
\end{align*}
which follows from the embedding \eqref{emb-BUC}. From the compatibility condition \eqref{1.12*} we have
\begin{align*}
(\bu_{0 -} - \bu_{1 -} \vert_{t = 0}, \nabla \varphi)_{\Omega_-} = (\bu_{0 -} - \bF_d \vert_{t = 0}, \nabla \varphi)_{\Omega_-} = 0 \quad \text{ for any $\varphi \in \hW_{q',0} (\Omega_-)$}.
\end{align*}
In addition, if $1 < 2/p + 1/q < 2$, by compatibility condition \eqref{1.13*} we have
\begin{align*}
\langle \nabla (\rho_{0 +} - \rho_{1 +} \vert_{t = 0}), \bn_+ \rangle \vert_+ = \langle\nabla \rho_{0 +}, \bn_+ \rangle \vert_+ - k_+ \vert_{t = 0} & = 0 & \quad & \text{ on $\Gamma_+$}, \\
\bu_{0 +} - \bu_{1 +} \vert_{t = 0} & = 0 & \quad & \text{ on $\Gamma_+$},\\
\bu_{0 -} - \bu_{1 -} \vert_{t = 0} & = 0 & \quad & \text{ on $\Gamma_-$}
\end{align*}
while if $2/p + 1/q < 1$, by compatibility condition \eqref{1.14*} we have
\begin{alignat*}3
\Pi_\bn (\gamma_4 \bD (\bu_{0 -} - \bu_{1 -} \vert_{t = 0})\bn) \vert_- - \Pi_\bn (\gamma_1 \bD (\bu_{0 +} - \bu_{1 +} \vert_{t = 0})\bn) \vert_+ & \\
= \Pi_\bn (\gamma_4 \bD (\bu_{0 -}) \bn) \vert_- - \Pi_\bn (\gamma_1 \bD (\bu_{0 +}) \bn) \vert_+ - g \vert_{t = 0} & = 0 & \quad & \text{ on $\Gamma$}, \\
\Pi_\bn (\bu_{0 -} - \bu_{1 -} \vert_{t = 0}) \vert_- - \Pi_\bn (\bu_{0 +} - \bu_{1 +} \vert_{t = 0}) \vert_+ = \Pi_\bn \bu_{0 -} \vert_- - \Pi_\bn \bu_{0 +} \vert_+ - \bh \vert_{t = 0} & = 0 & \quad & \text{ on $\Gamma$}, \\
\langle \nabla (\rho_{0 +} - \rho_{1 +} \vert_{t = 0}), \bn \rangle \vert_+ = \langle \nabla\rho_{0 +}, \bn \rangle \vert_+ - k_+ \vert_{t = 0} & = 0 & \quad & \text{ on $\Gamma$}, \\
\langle \nabla (\rho_{0 +} - \rho_{1 +} \vert_{t = 0}), \bn_+ \rangle \vert_+ = \langle \nabla\rho_{0 +}, \bn_+ \rangle \vert_+ - k_+ \vert_{t = 0} & = 0 & \quad & \text{ on $\Gamma_+$}, \\
\bu_{0 +} - \bu_{1 +} \vert_{t = 0} & = 0 & \quad & \text{ on $\Gamma_+$}, \\
\bu_{0 -} - \bu_{1 -} \vert_{t = 0} & = 0 & \quad & \text{ on $\Gamma_-$}.
\end{alignat*}
Hence, if $2/p + 1/q \neq 1$ and $2/p + 1/q \neq 2$, we see that
\begin{align*}
(\rho_{0 +} - \rho_{1 +} \vert_{t = 0}, \bu_{0 +} - \bu_{1 +} \vert_{t = 0}, \bu_{0 -} -\bu_{1 -} \vert_{t = 0}, h_0 - h_1 \vert_{t = 0}) \in \CD_{q, p} (\Omega_+, \Omega_-, \Gamma).
\end{align*}
Then, by Theorem \ref{th-7.2}, there exists a positive constant $\gamma'$ such that \textcolor{red}{the system} \eqref{eq-7.11} admits unique solutions $\rho_{2 +}$, $\bu_{2 \pm}$, and $h_2$ with $\pi_{2 -} = \CK_1 (\bu_{2 -}) + \CK_2 (h_2) + \CK_3 (\rho_{2 +}, h_2)$ and
\begin{equation}
\label{7.15}
\begin{split} 
\rho_{2 +} & \in W^{3, 1}_{q, p, \gamma'} (\Omega_+ \times (0, \infty)) \quad \bu_{2 \pm} \in W^{2, 1}_{q, p, \gamma'} (\Omega_\pm \times (0, \infty)), \\
h_2 & \in W^1_{q, \gamma'} ((0, \infty), W^{2 - 1/q}_q (\Gamma)) \cap L_{p, \gamma'} ((0,\infty), W^{3 - 1/q}_q (\Gamma))
\end{split}
\end{equation}
possessing the estimate:
\begin{equation}
\label{7.16}
\begin{split}
& \lVert e^{-\gamma' t} \pd_t \rho_{2 +} \rVert_{L_p((0, \infty), W^1_q (\Omega_+))} + \lVert e^{-\gamma' t} \rho_{2 +} \rVert_{L_p((0, \infty), W^3_q (\Omega_+))} \\
& + \lVert e^{-\gamma' t} \pd_t \bu_{2 +} \rVert_{L_p((0, \infty), L_q(\Omega_+))} + \lVert e^{-\gamma' t} \bu_{2 +} \rVert_{L_p((0, \infty), W^2_q (\Omega_+))} \\
& + \lVert e^{-\gamma' t} \pd_t \bu_{2 -} \rVert_{L_p((0, \infty), L_q(\Omega_-))} + \lVert e^{-\gamma' t} \bu_{2 -} \rVert_{L_p((0, \infty), W^2_q (\Omega_-))} \\
& + \lVert e^{-\gamma' t} \pd_t h_2 \rVert_{L_p((0, \infty), W^{2 - 1/q}_q (\Gamma))} + \lVert e^{-\gamma' t} h_2 \rVert_{L_p((0, \infty), W^{3 - 1/q}_q (\Gamma))} \\
& \leq C \Big(\lVert \rho_{0 +} - \rho_{1 +} \vert_{t = 0} \rVert_{B^{3 - 2/p}_{q, p} (\Omega_+)} + \lVert \bu_{0 +} - \bu_{1 +} \vert_{t = 0} \rVert_{B^{2 (1 - 1/p)}_{q, p} (\Omega_+)} \\
& \quad + \lVert \bu_{0 -} - \bu_{1 -} \vert_{t = 0} \rVert_{B^{2 (1 - 1/p)}_{q, p} (\Omega_-)} + \lVert h_0 - h_1 \vert_{t = 0} \rVert_{B^{3 - 1/p - 1/q}_{q, p} (\Gamma)} \Big).
\end{split}
\end{equation}
Setting $\rho_+ = \rho_{1 +} + \rho_{2 +}$, $\bu_\pm = \bu_{1 \pm} + \bu_{2 \pm}$, $\pi_- =\pi_{1 -} + \CK_1 (\bu_{2 -}) + \CK_2 (h_2) + \CK_3 (\rho_{2 +}, h_2)$, and $h = h_1 + h_2$ and choosing $\gamma_0$ such that $\gamma_0 > \max(\lambda_{**}, \gamma')$, from \eqref{eq-7.10}, \eqref{7.12}, \eqref{7.14}, \eqref{7.15}, and \eqref{7.16}, we see that $\rho_+$, $\bu_\pm$, $\pi_-$, and $h$ are solutions to the problem \eqref{eq-1.6} and satisfy the required estimate. Furthermore, a uniqueness of $\rho_+$, $\bu_\pm$, $\pi_-$, and $h$
follow from the uniqueness of a solution to the generalized resolvent problem~\eqref{eq-8.1*}. We, therefore, complete the proof of Theorem \ref{th-1.5}.

\section{The nonlinear problem}\label{sect-8}
\subsection{Tools for estimating nonlinear terms}
In this subsection, we collect some useful tools and definitions needed later on in the proof of Theorem \ref{th-1.9}. We first introduce the Sobolev embedding theorem:
\begin{equation}
\label{8.1}
\begin{split}
\lVert f \rVert_{L_\infty (\Omega_\pm)} & \leq C \lVert f \rVert_{W^1_q (\Omega_\pm)}, \\
\lVert f g \rVert_{W^1_q (\Omega_\pm)} & \leq C \lVert f \rVert_{W^1_q (\Omega_\pm)} \lVert g \rVert_{W^1_q (\Omega_\pm)}, \\
\lVert f g \rVert_{W^{1 - 1/q}_q (\Gamma)} & \leq C \lVert f \rVert_{W^{1 - 1/q}_q (\Gamma)} \lVert g \rVert_{W^{1 - 1/q}_q (\Gamma)},
\end{split}
\end{equation}
where $N < q < \infty$. On the other hand, if $2/p + N/q < 1$, we have
\begin{equation}
\label{8.2}
\begin{split}
\lVert f (\cdot, t) \rVert_{W^1_\infty (\Omega_\pm)} & \leq C_{p,q} \sup_{t \in (0, T)} \lVert f (\cdot, t) \rVert_{B^{2(1 - 1/p)}_{q,p} (\Omega_\pm)}, \\
\lVert f (\cdot, t) \rVert_{W^2_\infty (\Omega_+)} & \leq C_{p,q} \sup_{t \in (0, T)} \lVert f (\cdot, t) \rVert_{B^{3 - 2/p}_{q,p} (\Omega_+)}, \\
\lVert f (\cdot, t) \rVert_{W^2_\infty (\Omega_\pm)} & \leq C_{p,q} \sup_{t \in (0, T)} \lVert f (\cdot, t) \rVert_{B^{3 - 1/p}_{q,p} (\Omega_\pm)}
\end{split}
\end{equation}
for every $t \in (0, T)$. In fact, by Muramatu~\cite[Theorem 9]{Mur1973}, the embedding $B^{1 + N/p + \varepsilon}_{q, p} (\Omega_\pm) \hookrightarrow W^1_\infty (\Omega_\pm)$ holds for $0 < \varepsilon < 1 - (2/p + N/q)$ provided that $2/p + N/q < 1$. Assuming $2/p + N/q < 1$, we also have $B^{2 (1 - 1/p)}_{q, p} (\Omega_\pm) \hookrightarrow B^{1 + N/p + \varepsilon}_{q, p} (\Omega_\pm)$, which implies $B^{2 (1 - 1/p)}_{q, p} (\Omega_\pm) \hookrightarrow W^1_\infty (\Omega_\pm)$. In addition, we see that $B^{3 - 1/p}_{q, p} (\Omega_\pm) \hookrightarrow B^{3 - 2/p}_{q, p} (\Omega_+) \hookrightarrow W^2_\infty (\Omega_+)$ under the assumption: $2/p + N/q < 1$.

Given a function $\phi$ defined on $\Gamma$, let $H_h$ be a solution to the strong Dirichlet problem:
\begin{align*}
(1 - \Delta) H_h = 0 \quad \text{in $\dot \Omega$}, \qquad H_h = h \quad \text{on $\Gamma$}
\end{align*}
with an initial data $H_h \vert_{t = 0} = H_{h_0}$ such that
\begin{align*}
(1 - \Delta) H_{h_0} = 0 \quad \text{in $\dot \Omega$}, \qquad H_{h_0} = h_0 \quad \text{on $\Gamma$}.
\end{align*}
In the following, we assume
\begin{align*}
\sup_{t \in (0, T)} \lVert H_h (\cdot, t) \rVert_{W^1_\infty (\dot \Omega)} \leq \widetilde \varepsilon,
\end{align*}
whose constant $\widetilde \varepsilon$ is same as in \eqref{2.2}. Here, we have the estimate
\begin{equation}
\label{8.5}
\begin{split}
\lVert H_h (\cdot,t) \rVert_{W^3_q (\dot \Omega)} \leq C \lVert h (\cdot, t) \rVert_{W^{3 - 1/q}_q (\Gamma)}, \quad \lVert \pd_t H_h (\cdot,t) \rVert_{W^2_q (\dot \Omega)} \leq C \lVert \pd_t h (\cdot,t) \rVert_{W^{2 - 1/q}_q (\Gamma)}
\end{split}
\end{equation}
for $t \in (0, T)$. \par
Define the following space:
\begin{equation*}
\BU_{\varepsilon_T} = \Set{(\rho_+, \bu_+, \bu_-, h) |
\begin{aligned}
\rho_+ & \in W^{3, 1}_{q, p} (\Omega_+ \times (0, T)), \quad \bu_\pm \in W^{2, 1}_{q, p} (\Omega_\pm \times (0, T)), \\
h & \in L_{p}((0,T),W^{3-1/q}_q(\Gamma))\cap
W^1_p((0,T),W^{2-1/q}_q(\Gamma)),\\
\rho_+ & \vert_{t = 0} = \rho_{0 +} \quad \text{in $\Omega_+$}, \quad \bu_\pm \vert_{t = 0} = \bu_{0 \pm} \quad \text{in $\Omega_\pm$}, \\
h & \vert_{t = 0} = h_0 \quad \text{on $\Gamma$}, \quad \sup_{t \in (0, T)} \lVert H_h (\cdot, t) \rVert_{W^1_\infty (\dot \Omega)} \leq \widetilde \varepsilon, \\
& \BI_{p,q} (\rho_+, \bu_+, \bu_-, h, 0; (0, T)) \leq \varepsilon_T
\end{aligned}}.
\end{equation*}
Let $(\varrho_+,\bv_+,\bv_-,H_{\phi})\in\BU_{\varepsilon_T}$,
and then we consider the following system:
\begin{align}
\label{8.7}
\left\{\begin{aligned}
\pd_t \rho_+ + \rho_{*+} \dv \bu_+ & = f_M (\varrho_+, \bv_+, H_\phi) &\enskip &\text{ in $\Omega_+ \times(0, T)$}, \\
\rho_{*-} \dv \bu_- = f_d (\bv_-, H_\phi) & = \rho_{*-} \dv \bF_d (\bv_-, H_\phi) &\enskip &\text{ in $\Omega_- \times (0, T)$}, \\
\rho_{*+}\pd_t\bu_+ -\DV \bT_+ (\gamma_1, \gamma_2, \gamma_3, \rho_+, \bu_+) & = \bff_+ (\varrho_+, \bv_+, H_\phi) &\enskip &\text{ in $\Omega_+ \times (0, T)$}, \\
\rho_{*-}\pd_t\bu_--\DV \bT_- (\gamma_4, \bu_-, \pi_-) & = \bff_- (\bv_-, H_\phi) &\enskip &\text{ in $\Omega_- \times (0, T)$}, \\
\pd_t h - \frac{\langle \rho_{*-} \bu_-, \bn \rangle \rvert_- - \langle \rho_{*+} \bu_+, \bn \rangle \rvert_+}{\rho_{*-} - \rho_{*+}} 
& = d (\varrho_+, \bv_+, \bv_-, H_\phi) &\enskip &\text{ on $\Gamma \times (0, T)$}, \\
\bB (\gamma_1, \gamma_2, \gamma_3, \gamma_4, \rho_+, \bu_+, \bu_-, \pi_-)
& = \bG (\varrho_+, \bv_+, \bv_-, H_\phi), &\enskip &\text{ on $\Gamma \times (0, T)$}, \\
\bu_+ = 0, \qquad \langle \nabla \rho_+, \bn_+ \rangle & = 0
&\enskip &\text{ on $\Gamma_+ \times (0, T)$}, \\
\bu_- & = 0 &\enskip &\text{ on $\Gamma_- \times (0, T)$}, \\
(\rho_+, \bu_+, \bu_-, h) \rvert_{t = 0} & = (\rho_{0+}, \bu_{0+}, \bu_{0-}, h_0)
&\enskip &\text{ on $\Omega_+ \times \Omega_+ \times \Omega_- \times \Gamma$}.
\end{aligned}\right.
\end{align}
\par
We extend initial data $\rho_{0 +}$, $\bu_{0 +}$, $\bu_{0 -}$, and $H_{h_0}$ to $x \in \BR^N$. Let $\widetilde \rho_{0+}$, $\widetilde \bu_{0\pm}$, and $\widetilde H_{h_0}$ be extensions of $\rho_{0 +}$, $\bu_{0 \pm}$, and $h_0$ to $x \in \BR^N$, respectively, such that
\begin{alignat*}4
\rho_{0 +} & = \widetilde \rho_{0+} & \quad & \text{in $\Omega_+$}, & \qquad \lVert \widetilde \rho_{0 +} \rVert_{B^{3 - 2/p}_{q,p} (\BR^N)} & \leq C \lVert \rho_{0 +} \rVert_{B^{3 - 2/p}_{q,p} (\Omega_+)}, \\
\bu_{0 \pm} & = \widetilde \bu_{0 \pm} & \quad & \text{in $\Omega_\pm$}, & \qquad \lVert \widetilde \bu_{0 \pm} \rVert_{B^{2(1 - 1/p)}_{q,p} (\BR^N)} & \leq C \lVert \bu_{0 \pm} \rVert_{B^{2(1 -1/p)}_{q,p} (\Omega_\pm)}, \\
H_{h_0} & = \widetilde H_{h_0} & \quad & \text{in $\dot \Omega$}, & \qquad \lVert \widetilde H_{h_0} \rVert_{B^{3 - 1/p}_{q,p} (\BR^N)} & \leq C \lVert H_{h_0} \rVert_{B^{3 - 1/p}_{q,p}(\dot \Omega)}.
\end{alignat*}
Then we define functions $T_{\varrho_+} (t)$, $T_{v \pm} (t)$, and $T_\phi (t)$ as
\begin{align*}
T_{\varrho_+} (t) (\rho_{0 +}) & = e^{- (1 - \Delta)^{3/2}t} (\widetilde \rho_{0+}) = \CF^{- 1} [e^{- (1 + \lvert \xi \rvert^2)^{3/2} t} \CF[\widetilde \rho_{0 +}] (\xi)], \\
T_{v \pm} (t) \bu_{0 \pm} & = e^{- (1 - \Delta) t} \widetilde \bu_{0\pm} = \CF^{- 1} [e^{- (1 + \lvert \xi \rvert^2)} \CF[\widetilde \bu_{0 \pm}] (\xi)], \\
T_\phi (t) H_{h_0} & = e^{- (1 - \Delta)^{3/2}t} \widetilde H_{h_0} = \CF^{- 1} [e^{- (1 + \lvert \xi \rvert^2)^{3/2} t} \CF[\widetilde H_{h_0}] (\xi)].
\end{align*}
We see that $T_{\varrho_+} (0) \rho_{0 +} = \rho_{0 +}$ in $\Omega_+$, $T_{v \pm} (0) \bu_{0 \pm} = \bu_{0 \pm}$ in $\Omega_\pm$, and $T_\phi (0) H_{h_0} = H_{h_0}$ in $\dot \Omega$ satisfying the following estimates:
\begin{equation}
\label{9.8}
\begin{split}
\lVert T_{\varrho_+} (\cdot) \rho_{0 +} \rVert_{W^1_p ((0, \infty), W^1_q (\BR^N))} + \lVert T_{\varrho_+} (\cdot) \rho_{0 +} \rVert_{L_p ((0, \infty), W^3_q (\BR^N))} \leq C \lVert \rho_{0 +} \rVert_{B^{3 - 2/p}_{q,p} (\Omega_+)} & \leq C \varepsilon_T, \\
\lVert T_{v \pm} (\cdot) \bu_{0 \pm} \rVert_{W^1_p ((0, \infty), L_q (\BR^N))} + \lVert T_{v \pm} (\cdot) \bu_{0 \pm} \rVert_{L_p ((0, \infty), W^2_q (\BR^N))} \leq C \lVert \bu_{0 \pm} \rVert_{B^{2(1 - 1/p)}_{q,p} (\Omega_\pm)} & \leq C \varepsilon_T, \\
\lVert T_\phi (\cdot) H_{h_0} \rVert_{W^1_p ((0, \infty), W^2_q (\BR^N))} + \lVert T_\phi (\cdot) H_{h_0} \rVert_{L_p ((0, \infty), W^3_q (\BR^N))} \leq C \lVert h_0 \rVert_{B^{3 - 1/p}_{q,p} (\Gamma)} & \leq C \varepsilon_T
\end{split}
\end{equation}
for some positive constant $C$. \par
Let $\chi(t)\in C^\infty(\BR)$ be cut-off functions such that $\chi (t) = 1$ on $t > - 1$ and {\color{red} $\chi (t) = 0$} on $t < -2$. In addition, let $f(t)$ be a function defined on $t \in (0, T)$, and then we define an extension of $f$ to $t \in \BR$ by
\begin{align*}
e_T [f] (t) =
\begin{cases}
0 & (t \leq 0), \\
f (t) & (0 < t < T), \\
f (2T - t) & (T \leq t < 2T), \\
0 & (t \geq 2T).
\end{cases}
\end{align*}
If $f \vert_{t = 0} = 0$, we see that
\begin{align}
\label{7.10}
\pd_t e_T [f] (t)=
\begin{cases}
0 & (t \leq 0), \\
(\pd_t f) (t) & (0 < t < T),\\
- (\pd_t f) (2T - t) & (T \leq t < 2T), \\
0 & (t \geq 2T).
\end{cases}
\end{align}
We then define the extension of $\varrho_+$, $\bv_\pm$, and $H_\phi$ to $t \in \BR$ by
\begin{align*}
E_{\varrho_+} [\varrho_+] & = e_T [\varrho_+ - T_{\varrho_+} (t) \rho_{0 +}] + \chi(t) T_{\varrho_+} (\lvert t \rvert) \rho_{0 +}, \\
E_{v \pm} [\bv_\pm] & = e_T [\bv_\pm - T_{v \pm} (t) \bu_{0 \pm}] + \chi (t) T_{v \pm} (\lvert t \vert) \bu_{0 \pm}, \\
E_\phi [H_\phi] & = e_T [H_\phi - T_\phi (t) H_{h_0}] + \chi (t) T_\phi (\lvert t \vert) H_{h_0},
\end{align*}
respectively. Here, of course, $E_{\varrho_+} [\varrho_+] = \varrho_+$, $E_{v \pm} [\bv_\pm] = \bv_\pm$, and $E_\phi[\phi] = \phi$ holds for every $t \in (0, T)$. To estimate $H^{1/2}_p$ norm of $\bF_d (\varrho_+, \bv_+, \bv_-, H_\phi)$ and $f_d (\varrho_+, \bv_+,\bv_-,H_\phi)$ we use the following lemmas. The farmer lemma is immediately follows from the complex interpolation methods and the latter has been proven by Shibata~\cite[Proposition 1]{Shi2018}, so that we may omit those proofs.
\begin{lemm}
\label{lem-7.1}
Let $1 < p,q < \infty$ and let $D \subset \BR^N$ be a uniformly $C^2$ domain. In addition, let $f \in W^1_\infty (\BR, L_\infty (D))$ and $g \in H^{1/2}_p (\BR, L_q(D))$. Then, we have the estimate
\begin{align*}
\lVert f g \rVert_{H^{1/2}_p (\BR, L_q (D))} \leq C \lVert f \rVert_{W^1_\infty (\BR, L_\infty (D))} \lVert g \rVert_{H^{1/2}_p (\BR, L_q (D))}.
\end{align*}		
\end{lemm}
\begin{lemm}
\label{lem-7.2}
Let $1 < p,q < \infty$ and let $D$ be uniformly $C^2$ domain. We then have the following properties:
\begin{gather*}
H^1_p (\BR, L_q(D)) \cap L_p (\BR, W^2_q (D)) \subset H^{1/2}_p (\BR, W^1_q (D)), \\
\lVert f \rVert_{H^{1/2} (\BR, W^1_q (D))} \leq  C \{\lVert f \rVert_{L_p (\BR, W^2_q (D))} + \lVert \pd_t f \rVert_{L_p (\BR, L_q (D))}\}.
\end{gather*}		
\end{lemm}

\subsection{Estimating the nonlinear terms}
We finally prove Theorem \ref{th-1.9} with the help of the Banach fixed point argument. To use Theorem \ref{th-1.5}, we now extend the right-hand side of \eqref{8.7} to $t \in \BR$. We first consider $\bff_+ (\varrho_+, \bv_+, H_\phi)$ and $\bff_- (\bv_-, H_\phi)$. Let $\ov \bff_+ (\varrho_+, \bv_+, H_\phi)$ and $\ov \bff_- (\bv_-, H_\phi)$ be the zero extension of $\bff_+ (\varrho_+, \bv_+, H_\phi)$ and $\bff_- (\bv_-, H_\phi)$ to all of $\BR$, respectively. Recalling the representation formulas $\bff_+$ and $\bff_-$, which are given in Appendix, using \eqref{8.1} and \eqref{8.5} and choosing $\varepsilon_T > 0$ so small such that $\varepsilon_T \leq \rho_{* +} / 3$, we have the estimate
\begin{align*}
& \lVert \bff_+ (\varrho_+, \bv_+, H_\phi) \rVert_{L_q (\Omega_+)} \\
& \leq C \Big\{\lVert \varrho_+ (\cdot, t) \rVert_{W^1_q (\Omega_+)} \lVert \pd_t \bv_+ (\cdot, t) \rVert_{L_q (\Omega_+)} + \lVert \bv_+ (\cdot, t) \rVert^2_{W^1_q (\Omega_+)} + \lVert \pd_t H_\phi (\cdot, t) \rVert_{W^1_q (\Omega_+)} \lVert \bv_+ (\cdot, t) \rVert_{W^1_q (\Omega_+)} \\
& \quad + \lVert H_\phi (\cdot, t) \rVert_{W^2_q (\Omega_+)} \Big(\lVert \pd_t \bv_+ (\cdot, t) \rVert_{L_q (\cdot, t)} + \lVert \bv_+ (\cdot, t) \rVert_{W^2_q (\Omega_+)} \Big) + \lVert \varrho_+(\cdot, t) \rVert_{W^1_q (\Omega_+)} \lVert \varrho_+ (\cdot, t) \rVert_{W^3_q (\Omega_+)} \\
& \quad + \lVert \varrho_+ (\cdot, t) \rVert^2_{W^2_q (\Omega_+)} \Big\},
\end{align*}
which yields
\begin{equation}
\label{est-1}
\begin{split}
& \lVert \ov \bff_+ (\varrho_+, \bv_+, H_\phi) \rVert_{L_p (\BR, L_q (\Omega_+))} \\
& \leq C \Big\{\lVert \varrho_+ \rVert_{L_\infty ((0, T), W^1_q (\Omega_+))} \lVert \pd_t \bv_+ \rVert_{L_p ((0, T), L_q (\Omega_+))} + \lVert \bv_+ \rVert_{L_\infty ((0, T), W^1_q (\Omega_+))} \lVert \bv_+ \rVert_{L_p ((0, T), W^1_q (\Omega_+))} \\
& \quad + \lVert \pd_t H_\phi \rVert_{L_p ((0, T), W^1_q (\Omega_+))} \lVert \bv_+ \rVert_{L_\infty ((0, T), W^1_q (\Omega_+))} \\
& \quad + \lVert H_\phi \rVert_{L_\infty ((0, T), W^2_q (\Omega_+))} \Big(\lVert \pd_t \bv_+ \rVert_{L_p ((0, T), L_q (\Omega_+))} + \lVert \bv_+ \rVert_{L_p ((0, T), W^2_q (\Omega_+))}\Big) \\
& \quad + \lVert \varrho_+ \rVert_{L_\infty ((0, T), W^1_q (\Omega_+))} \lVert \varrho_+ \rVert_{L_p ((0, T), W^3_q (\Omega_+))} + \lVert \varrho_+ \rVert_{L_\infty ((0, T), W^2_q (\Omega_+))} \lVert\varrho_+ \rVert_{L_p ((0, T), W^2_q (\Omega_+))} \Big\} \\
&\leq C \varepsilon_T^2,
\end{split}
\end{equation}
where $C$ is a positive constant. Analogously we have
\begin{align}
\label{est-2}
\lVert \ov \bff_- (\bv_-, H_\phi) \rVert_{L_p (\BR, L_q (\Omega_-))} \leq C \varepsilon_T^2.
\end{align}
\par
To handle other nonlinear terms we extend the transformation: $w= x + H_\phi \bn_*$ to $w = x + E_\phi [H_\phi] \bn_*$, where $\bn_*$ is the extension of $\bn$ from $\Gamma$ to $\dot \Omega$ satisfying the estimate $\lVert \bn_* \rVert_{W^2_\infty (\dot \Omega)} \leq C$ with some positive constant $C$. Setting $\widetilde \bm = \nabla (E_\phi [H_\phi] \bn_*)$, by \eqref{2.7*} we see that
\begin{align*}
J_0 (\widetilde \bm) \dv \bv_- + J (\widetilde \bm) \bM_0 (\widetilde \bm) \colon \nabla \bv_- = \dv (J (\widetilde \bm) (\bI + {}^\top\!\bM_0 (\widetilde \bm)) \bv_-).
\end{align*}
Then we set
\begin{align*}
\ov f_d & := \ov f_d (\bv_-, H_\phi) = - \Big(J_0 (\widetilde \bm) \dv \bv_- + (1 + J_0 (\widetilde \bm)) \bM_0 (\widetilde \bm) \colon \nabla \bv_- \Big) = \bM_1 (\widetilde \bm) \nabla E_{v -}[\bv_-],\\
\ov \bF_d & := \ov \bF_d (\bv_-, H_\phi) = - (1 + J_0 (\widetilde \bm)) {}^\top\bM_0 (\widetilde \bm) \bv_- = \bM_1 (\widetilde \bm) E_{v -} [\bv_-],
\end{align*}
and
\begin{align*}
\ov f_d (\bv_-, H_\phi) = f_d (\bv_-, H_\phi), \quad \ov \bF_d (\bv_-, H_\phi) = \bF_d (\bv_-, H_\phi) \quad &\text{for $t \in (0, T)$}, \\
\dv \ov \bF_d (\bv_-, H_\phi) = \ov f_d (\bv_-, H_\phi) \quad &\text{in $\Omega_-$},
\end{align*}
where $\bM_1(\widetilde \bm)$ is a $C^\infty$ function of matrix defined on $\lvert \widetilde \bm \rvert \leq \varepsilon_T$ such that $\bM_1 (0) = 0$. Here, $\bM_1 (\widetilde \bm)$ can be estimated as
\begin{equation}
\label{8.15}
\begin{split}
\lVert \bM_1 (\widetilde \bm) \rVert_{L_\infty (\BR, L_\infty (\dot \Omega))} \leq C \Big(\lVert H_\phi \rVert_{L_\infty ((0, T), W^2_q (\dot \Omega))} + \lVert T_\phi (\cdot) H_{h_0} \rVert_{L_\infty ((0, \infty), W^2_q (\dot \Omega))} \Big) \leq C \varepsilon_T.
\end{split}
\end{equation}
Furthermore, since
\begin{align*}
\pd_t \ov \bF_d (\bv_-, h) = \bM_1 (\widetilde \bm) \pd_t E_{v -} (\bv_-) + \bM_1' (\widetilde \bm) (\pd_t \widetilde \bm) E_{v -} [\bv_-],
\end{align*}
where $\bM'_1 (\widetilde \bm)$ is the derivative of $\bM_1 (\widetilde \bm)$ with respect to $\widetilde \bm$, by \eqref{8.1}, \eqref{8.2}, \eqref{8.5}, \eqref{7.10}, and \eqref{8.15} we have
\begin{equation}
\label{est-3}
\begin{split}
\lVert \pd_t \ov \bF_d (\bv_-, H_\phi) \rVert_{L_p (\BR, L_q (\Omega_-))}
& \leq C \bigg(\lVert H_\phi \rVert_{L_\infty ((0, T), W^2_q (\dot \Omega))} + \lVert T_\phi (\cdot) H_{h_0} \rVert_{L_\infty ((0, \infty), W^2_q (\dot \Omega))} \bigg) \\
& \quad \times \bigg(\lVert \pd_t \bv_- \rVert_{L_p ((0, T), L_q (\Omega_-))} + \lVert \pd_t T_{v -} (\cdot) \bu_{0 -} \rVert_{L_p ((0, \infty), L_q (\Omega_-))} \bigg) \\
& \quad \times \bigg(\lVert \pd_t H_\phi \rVert_{L_\infty ((0, T), W^1_q (\dot \Omega))} + \lVert \pd_t T_\phi (\cdot) H_{h_0} \rVert_{L_p ((0, \infty), W^1_q (\dot \Omega))} \bigg) \\
& \quad \times \bigg(\lVert \bv_- \rVert_{L_\infty ((0, T), W^1_q (\Omega_-))} + \lVert T_{v -} (\cdot) \bu_{0 -} \rVert_{L_\infty ((0, \infty), W^1_q (\Omega_-))} \bigg) \\
& \leq C \varepsilon^2_T.
\end{split}
\end{equation}
Analogously, we obtain the estimate
\begin{align}
\label{9.17}
\lVert \ov f_d (\bv_-, H_\phi) \rVert_{L_p (\BR, W^1_q (\Omega_-))} \leq C \varepsilon_T^2.
\end{align}
From Lemmas~\ref{lem-7.1} and \ref{lem-7.2} and the estimates \eqref{8.2}, \eqref{8.5}, and \eqref{9.8}, we have
\begin{equation}
\label{est-3*}
\begin{split}
\lVert \ov f_d (\bv_-, H_\phi) \rVert_{H^{1/2}_p (\BR, L_q (\Omega_-))} & \leq C \bigg(\lVert H_\phi \rVert_{L_\infty ((0, T), W^2_q (\dot \Omega))} + \lVert T_\phi (\cdot) H_{h_0} \rVert_{L_\infty ((0, \infty), W^2_q (\dot \Omega))} \bigg) \\
& \quad \times \bigg(\lVert \pd_t \bv_- \rVert_{L_p ((0, T), L_q (\Omega_-))} + \lVert \pd_t T_{v -} (\cdot) \bu_{0 -} \rVert_{L_p ((0, \infty), L_q (\Omega_-))} \bigg) \\
& \quad \times \bigg(\lVert \bv_- \rVert_{L_\infty ((0, T), W^1_q (\Omega_-))} + \lVert T_{v -} (\cdot) \bu_{0 -} \rVert_{L_\infty ((0, \infty), W^1_q (\Omega_-))} \bigg) \\
& \leq C \varepsilon^2_T.
\end{split}
\end{equation}
Let us define $\ov f_M (\varrho_+, \bv_+, H_\phi)$ by $\ov f_M (\varrho_+, \bv_+, H_\phi) = f_M (\varrho_+, \bv_+, H_\phi)$ for $t \in (0, T)$ and $f_M (\varrho_+, \bv_+, H_\phi) = 0$ for $t \notin (0, T)$. Then employing the argument above, we have the following estimates:
\begin{align}
\label{est-4}
\lVert \ov f_M \rVert_{L_p (\BR, W^1_q (\Omega_-))} \leq C \varepsilon_T^2.
\end{align}
\par
We next extend the right-hand members of the boundary conditions, $d$, $g$, $f^+_B$, $f^-_B$, $\bh$, and $k_-$, to $t \in \BR$. We define $\ov d (\varrho_+, \bv_+, \bv_-, H_\phi) = d (\varrho_+, \bv_+,\bv_-, H_\phi)$ for $t \in (0, T)$ and $\ov d (\varrho_+, \bv_+, \bv_-, H_\phi) = 0$ for $t \notin (0, T)$. On the other hand, to define $g$, $f^+_B$, $f^+_B$, $\bh$, and $k_-$ we use the extensions $E_{v \pm}$, $E_{\varrho_+}$, and $E_\phi$. Using the argument above and choosing $\varepsilon_T > 0$ suitably small such that $\varepsilon_T \leq \min(\rho_{* +} / 3, 1)$, we have the estimates
\begin{align}
\label{est-5}
\begin{aligned}
\lVert \ov d \rVert_{L_p (\BR, W^{2 - 1/q}_q (\Gamma))} & \leq C \varepsilon_T^2, & \quad \lVert \ov g \rVert_{H^{1/2}_p (\BR, L_q (\dot \Omega))} & \leq C \varepsilon_T^2, & \quad \lVert \ov g \rVert_{L_p (\BR, W^1_q (\dot \Omega))} & \leq C \varepsilon_T^2, \\
\lVert \ov f^+_B \rVert_{L_p (\BR, L_q(\dot \Omega))} & \leq C \varepsilon_T^2, & \quad \lVert \ov f^+_B \rVert_{H^{1/2}_p (\BR, L_q (\dot \Omega))} & \leq C \varepsilon_T^2, & \quad \lVert \ov f^+_B \rVert_{L_p (\BR,W^1_q (\dot \Omega))} & \leq C \varepsilon_T^2, \\
\lVert \ov f^-_B \rVert_{L_p (\BR, L_q (\dot \Omega))} & \leq C \varepsilon_T^2, & \quad \lVert \ov f^-_B \rVert_{H^{1/2}_p (\BR, L_q(\dot \Omega))} & \leq C \varepsilon_T^2, & \quad \lVert \ov f^-_B \rVert_{L_p (\BR, W^1_q (\dot \Omega))} & \leq C \varepsilon_T^2, \\
\lVert \ov \bh \rVert_{W^1_p (\BR, L_q (\dot \Omega))} & \leq C \varepsilon_T^2, & \quad \lVert \ov \bh \rVert_{H^{1/2}_p (\BR, \dot\Omega)} & \leq C \varepsilon_T^2, & \quad \lVert \ov \bh\rVert_{L_p (\BR, W^2_q (\dot \Omega))} & \leq C\varepsilon_T^2, \\
\lVert \ov k_- \rVert_{W^1_p (\BR, L_q (\dot \Omega))} & \leq C \varepsilon_T^2, & \quad \lVert \ov k_- \rVert_{H^{1/2}_p (\BR, \dot \Omega)} & \leq C \varepsilon_T^2, & \quad \lVert \ov k_- \rVert_{L_p (\BR, W^2_q (\dot \Omega))} &\leq C \varepsilon_T^2.
\end{aligned}
\end{align}
Summing up, from Theorem \ref{th-1.5}, \eqref{est-1}, \eqref{est-2}, \eqref{est-3},\eqref{est-3*}, \eqref{est-4}, and \eqref{est-5}, we see that the problem~\eqref{8.7} admits a unique solution
\begin{align*}
\rho_+ \in W^{3, 1}_{q, p} (\Omega_+ \times (0, T)), \enskip \bu_\pm \in W^{2, 1}_{q, p} (\Omega_\pm \times (0, T)), \enskip h \in L_p ((0, T), W^{3 - 1/q}_q (\Gamma)) \cap W^1_p ((0, T), W^{2 - 1/q}_q(\Gamma))	
\end{align*}
with the estimate
\begin{align}
\label{9.20}
\BI_{p, q} (\rho_+, \bu_+, \bu_-, h, 0; (0, T)) \leq C \varepsilon_T^2.
\end{align}
Choosing $\varepsilon_T > 0$ so small that $C \varepsilon_T \leq 1$, we see that $\BI_{p,q} (\rho_+, \bu_+, \bu_-, h, 0; (0, T)) \leq \varepsilon_T$. We define a map $\Phi \colon \BU_{\varepsilon_T} \to\BU_{\varepsilon_T}$ such that $\Phi (\varrho_+, \bv_+, \bv_-, \phi) = (\rho_+, \bu_+, \bu_-, h)$. Then the mapping $\Phi$ is a contraction mapping. In fact, given $(\varrho_{i +}, \bv_{i +}, \bv_{i -},\phi_i) \in \BU_{\varepsilon_T}$ with $i = 1, 2$, we set $(\rho_+, \bu_+, \bu_-, h) = \Phi(\varrho_{i +}, \bv_{i +}, \bv_{i -}, \phi_i)$, and then $\rho_+ = \rho_{2 +} - \rho_{1 +}$, $\bu_\pm =\bu_{2 \pm} - \bu_{1 \pm}$, and $h = h_2 - h_1$ satisfy the following system:
\begin{align}
\left\{\begin{aligned}
\pd_t \rho_+ + \rho_{*+} \dv \bu_+ & = \widetilde f_M (\varrho_+, \bv_+, H_\phi) &\enskip &\text{ in $\Omega_+ \times(0, T)$}, \\
\rho_{*-} \dv \bu_- = \widetilde f_d (\bv_-, H_\phi) & = \rho_{*-} \dv \widetilde \bF_d (\bv_-, H_\phi) &\enskip &\text{ in $\Omega_- \times (0, T)$}, \\
\rho_{*+}\pd_t\bu_+ -\DV \bT_+ (\gamma_1, \gamma_2, \gamma_3, \rho_+, \bu_+) & = \widetilde \bff_+ (\varrho_+, \bv_+, H_\phi) &\enskip &\text{ in $\Omega_+ \times (0, T)$}, \\
\rho_{*-}\pd_t\bu_--\DV \bT_- (\gamma_4, \bu_-, \pi_-) & = \widetilde \bff_- (\bv_-, H_\phi) &\enskip &\text{ in $\Omega_- \times (0, T)$}, \\
\pd_t h - \frac{\langle \rho_{*-} \bu_-, \bn \rangle \rvert_- - \langle \rho_{*+} \bu_+, \bn \rangle \rvert_+}{\rho_{*-} - \rho_{*+}} & = \widetilde d (\varrho_+, \bv_+, \bv_-, H_\phi) &\enskip &\text{ on $\Gamma \times (0, T)$}, \\
\bB (\gamma_1, \gamma_2, \gamma_3, \gamma_4, \rho_+, \bu_+, \bu_-, \pi_-) & = \widetilde \bG (\varrho_+, \bv_+, \bv_-, H_\phi), &\enskip &\text{ on $\Gamma \times (0, T)$}, \\
\bu_+ = 0, \qquad \langle \nabla \rho_+, \bn_+ \rangle & = 0 &\enskip &\text{ on $\Gamma_+ \times (0, T)$}, \\
\bu_- & = 0 &\enskip &\text{ on $\Gamma_- \times (0, T)$}, \\
(\rho_+, \bu_+, \bu_-, h) \rvert_{t = 0} & = (0, 0, 0, 0) &\enskip &\text{ on $\Omega_+ \times \Omega_+ \times \Omega_- \times \Gamma$}.
\end{aligned}\right.
\end{align}
with some $\pi_- \in  L_p((0, T), W^1_q(\Omega_-) + \wh W^1_{q, \Gamma} (\Omega_-))$, where we have set
\begin{alignat*}3
\widetilde F_1 (\varrho_+, \bv_+, \phi) & = F_1(\varrho_{2 +}, \bv_{2 +}, H_{\phi_2}) - F_1 (\varrho_1, \bv_{1 +}, H_{\phi_1}) &\quad &(F_1 \in \{f_M, \bff_+\}), \\
\widetilde F_2 (\bv_-, \phi) & = F_2 (\bv_{2 -}, H_{\phi_2}) - F_2 (\bv_{1 -}, H_{\phi_1}) &\quad &(F_2 \in \{f_d, \bF_d, \bff_-\}), \\
\widetilde d (\varrho_+, \bv_+, \bv_-, \phi) & = d (\varrho_{2 +}, \bv_{2 +}, \bv_{2 -}, \phi_2) - h (\varrho_{1 +}, \bv_{1 +}, \bv_{1 -}, \phi_1), \\
\widetilde F_3 (\varrho_+, \bv_+, \bv_-, \phi) & = F (\varrho_{2 +}, \bv_{2 +}, \bv_{2 -}, H_{\phi_2}) - F (\varrho_1, \bv_{1 +}, \bv_{1 -}, H_{\phi_1}) &\quad &(F_3 \in \{g, f^+_B, f^-_B, \bh\}), \\
\widetilde \bh (\bu_+, \bu_-, \phi) & = \bh (\bu_{2 +}, \bu_{2 -}, H_{\phi_2}) - \bh (\bu_{1 +}, \bu_{1 -}, H_{\phi_1}), \\
\widetilde k_- (\varrho_+, \phi) & = k_- (\varrho_{2 +}, H_{\phi_2}) - k_- (\varrho_{1 +}, H_{\phi_1}).
\end{alignat*}
By the Taylor formula, we denote
{\allowdisplaybreaks \begin{align*}
\widetilde F_1 (\varrho_+, \bv_+, \phi) & = (H_{\phi_2} - H_{\phi_1}) \int_0^1 \frac{\mathrm{d}}{\dtheta} F_1 (\varrho_{2 +}, \bv_{2 +}, \theta H_{\phi_2} + (1 - \theta) H_{\phi_1}) \dtheta \\
& \quad + (\bv_{2 +} - \bv_{1 +}) \int_0^1 \frac{\mathrm{d}}{\dtheta} F_1 (\varrho_{2 +}, \theta \bv_{2 +} + (1 - \theta) \bv_{1 +}, H_{\phi_1}) \dtheta \\
& \quad + (\varrho_{2 +} - \varrho_{1 +}) \int_0^1 \frac{\mathrm{d}}{\dtheta} F_1 (\theta \varrho_{2 +} + (1 - \theta) \varrho_{1 +}, \bv_{1 +}, H_{\phi_1}) \dtheta, \\
\widetilde F_2 (\bv_-, \phi) & = (H_{\phi_2} - H_{\phi_1}) \int_0^1 \frac{\mathrm{d}}{\dtheta} F_2 (\bv_{2 -}, \theta H_{\phi_2} + (1 - \theta) H_{\phi_1}) \dtheta \\
& \quad + (\bv_{2 -} - \bv_{1 -}) \int_0^1 \frac{\mathrm{d}}{\dtheta} F_2 (\theta \bv_{2 -} + (1 - \theta) \bv_{1 -}, H_{\phi_1}) \dtheta, \\
\widetilde d (\varrho_+, \bv_+, \bv_-, \phi) & = (\phi_2 - \phi_1) \int_0^1 \frac{\mathrm{d}}{\dtheta} d (\varrho_{2 +}, \bv_{2 +}, \bv_{2 -}, \theta \phi_2 + (1 - \theta) \phi_1) \dtheta \\
& \quad + (\bv_{2 -} - \bv_{1 -}) \int_0^1 \frac{\mathrm{d}}{\dtheta} d (\varrho_{2 +}, \bv_{2 +}, \theta \bv_{2 -} + (1 - \theta) \bv_{1 -}, \phi_1) \dtheta \\
& \quad + (\bv_{2 +} - \bv_{1 +}) \int_0^1 \frac{\mathrm{d}}{\dtheta} d (\varrho_{2 +}, \theta \bv_{2 +} + (1 - \theta) \bv_{1 +}, \bv_{1 -}, \phi_1) \dtheta \\
& \quad + (\varrho_{2 +} - \varrho_{1 +}) \int_0^1 \frac{\mathrm{d}}{\dtheta} d (\theta \varrho_{2 +} + (1 - \theta) \varrho_{1 +}, \bv_{1 +}, \bv_{1 -}, \phi_1) \dtheta, \\
\widetilde F_3 (\varrho_+, \bv_+, \bv_-, \phi) & = (H_{\phi_2} - H_{\phi_1}) \int_0^1 \frac{\mathrm{d}}{\dtheta} F_3 (\varrho_{2 +}, \bv_{2 +}, \bv_{2 -}, \theta H_{\phi_2} + (1 - \theta) H_{\phi_1})\dtheta \\
& \quad + (\bv_{2 -} - \bv_{1 -}) \int_0^1 \frac{\mathrm{d}}{\dtheta} F_3 (\varrho_{2 +}, \bv_{2 +}, \theta \bv_{2 -} + (1 - \theta) \bv_{1 -}, H_{\phi_1}) \dtheta \\
& \quad + (\bv_{2 +} - \bv_{1 +}) \int_0^1 \frac{\mathrm{d}}{\dtheta} F_3 (\varrho_{2 +}, \theta \bv_{2 +} + (1 - \theta) \bv_{1 +}, \bv_{1 -}, H_{\phi_1}) \dtheta \\
& \quad + (\varrho_{2 +} - \varrho_{1 +}) \int_0^1 \frac{\mathrm{d}}{\dtheta} F_3 (\theta \varrho_{2 +} + (1 - \theta) \varrho_{1 +}, \bv_{1 +}, \bv_{1 -}, H_{\phi_1}) \dtheta, \\
\widetilde \bh (\bv_+, \bv_-, \phi) & = (H_{\phi_2} - H_{\phi_1}) \int_0^1 \frac{\mathrm{d}}{\dtheta} \bh (\bv_{2 +}, \bv_{2 -}, \theta H_{\phi_2} + (1 - \theta) H_{\phi_1})\dtheta \\
& \quad + (\bv_{2 -} - \bv_{1 -}) \int_0^1 \frac{\mathrm{d}}{\dtheta} \bh (\bv_{2 +}, \theta \bv_{2 -} + (1 - \theta) \bv_{1 -}, H_{\phi_1}) \dtheta \\
& \quad + (\bv_{2 +} - \bv_{1 +}) \int_0^1 \frac{\mathrm{d}}{\dtheta} \bh (\theta \bv_{2 +} + (1 - \theta) \bv_{1 +}, \bv_{1 -}, H_{\phi_1}) \dtheta, \\
\widetilde k_- (\varrho_+, \phi) & = (H_{\phi_2} - H_{\phi_1}) \int_0^1 \frac{\mathrm{d}}{\dtheta} k_- (\bv_{2 -}, \theta H_{\phi_2} + (1 - \theta) H_{\phi_1}) \dtheta \\
& \quad + (\bv_{2 -} - \bv_{1 -}) \int_0^1 \frac{\mathrm{d}}{\dtheta} k_- (\theta \bv_{2 -} + (1 - \theta) \bv_{1 -}, H_{\phi_1}) \dtheta.
\end{align*}}\noindent
Since $\varrho_{2 +} - \varrho_{1 +} = \bv_{2 \pm} - \bv_{1 \pm} = \phi_2 - \phi_1 = 0$ at $t = 0$, we extend $\ov f_M$, $\ov \bff_\pm$, $\ov f_d$, $\ov \bF_d$, $\ov d$, $\ov g$, $\ov f^\pm_B$, $\ov \bh$, and $\ov k_-$ to $t \in \BR$ by using the extensions $e_T [\varrho_{2 +} - \varrho_{1 +}]$, $e_T [\bv_{2 \pm} - \bv_{1 \pm}]$, and $e_T [\phi_2 - \phi_1]$. Then, employing the same argument as in proving the estimate \eqref{9.20}, we obtain the estimate
\begin{align*}
& \BI_{p,q} (\rho_{2 +} - \rho_{1 +}, \bu_{2 +} - \bu_{1 +}, \bu_{2 -} - \bu_{1 -}, h_2 - h_1, 0; (0, T)) \\
& \leq C \varepsilon_T \BI_{p,q} (\varrho_{2 +} - \varrho_{1 +}, \bv_{2 +} - \bv_{1 +}, \bv_{2 -} - \bv_{1 -}, \phi_2 - \phi_1, 0; (0, T)).
\end{align*}
Choosing $\varepsilon_T \in (0, 1)$ suitably small such that $C \varepsilon_T \leq 1$, we see that $\Phi$ is a contraction mapping on $\BU_{\varepsilon_T}$. Hence, by the Banach fixed point theorem there exists a \textit{unique} fixed point $(\rho_+, \bu_+, \bu_-, h) \in \BU_{\varepsilon_T}$ of the mapping $\Phi$. Furthermore, these $(\rho_+, \bu_+, \bu_-, h)$ enjoy the system \eqref{eq-1.5} with a suitable pressure term $\pi_- \in L_p ((0, T), W^1_q (\Omega_-) + \wh W^1_{q, \Gamma} (\Omega_-))$. Summing up, we have completed the proof of Theorem \ref{th-1.9}.

\appendix
\section{Explicit formulas of nonlinear terms}
\label{sect-2}
\subsection{Transformation of the mass equations and momentum equations}
\label{sect-2.1}
We first transform problem \eqref{eq-1.1}-\eqref{cond-1.3} to a domain with a fixed interface $\Gamma$, where $\Gamma_t$ is parameterized over $\Gamma$ by means of an unknown height function $h(x,t)$, whose idea is said to be the \textit{Hanzawa transformation}. We emphasize that the Lagrangian transformation is \textit{not} available in the phase transition problem case, because the interface is moved not only by advection but also by the phase flux, see Pr{\"u}ss and Simonett~\cite{PS2016} for further explanations. Furthermore, in general, if the surface tension is present on the moving boundary, the Lagrangian transformation is not well-adapted due to a lack of precise information on the regularity of the free boundary. \par
We first assume that the interface $\Gamma_t$ is given by
\begin{align*}
\Gamma_t = \{w = x + h(x, t) \bn(x) \mid x \in \Gamma \} \qquad (t \in (0, T)),
\end{align*}
where $h(x,t)$ is a small unknown function and $\bn$ is the outer unit normal to $\Gamma$. Although the hypersurface $\Gamma_t$ is unknown, it is possible to assume that the moving interface $\Gamma_t$ is approximated by a real analytic hypersurface $\Gamma$ because the $C^2$-hypersurface $\Gamma_t$ admits a \textit{tubular neighbourhood}. Namely, there exists a positive constant $r_0$ such that the mapping
\begin{equation*}
\Theta \colon \Gamma \times (- r_0, r_0) \to \BR^N, \qquad \Theta (x, r) := x + r \bn(x)
\end{equation*}
is a diffeomorphism, see Pr{\"u}ss and Simonett~\cite[Chapter 2]{PS2016} for further introductions. Let $d_\Gamma$ be the signed distance from $x \in \Omega$ to $\Gamma$, whose magnitude is given by $\lvert d_\Gamma \rvert = \mathrm{dist}\, (x, \Gamma)$. We define that $d_\Gamma$ is strictly negative if and only if $x \in \Omega_+$. To introduce the transformation $\Gamma \mapsto \Gamma_t$, let $H_h (x,t)$ and $\bn_* (x)$ be extensions of $h(x, t)$ and $\bn (x)$ from $\Gamma$ to $\dot \Omega$, respectively. The extension $\bn_*(x)$ is defined to be a sufficiently regular vector field and $H_h (x, t)$ is a small function in the sense of
\begin{align}
\label{2.2}
\sup_{t \in (0, T)} \lVert H_h (\cdot, t) \rVert_{W^1_\infty (\dot \Omega)} \leq \widetilde \varepsilon,
\end{align}
where $\widetilde \varepsilon \in (0, 1)$ is a suitably small given constant such that $\widetilde \varepsilon < r_0 / 6$. Then the \textit{Hanzawa transformation} is defined by
\begin{align}
\label{2.1}
w = x + \chi \left(\frac{3 d_\Gamma}{r_0} \right) H_h (x, t) \bn_* (x)  \colon \enskip \dot \Omega \to \dot \Omega_t,
\end{align}
where $0 \leq \chi (\xi) \leq 1$ is a cut-off function such that $\chi (\xi) = 1$ for $\lvert \xi \rvert \leq 1$ and $\chi (\xi) = 0$ for $\lvert \xi \rvert \geq 2$. Here, the assumption \eqref{2.2} guarantees the injectivity of the transformation~\eqref{2.1} for each $t \in (0, T)$. If the condition \eqref{2.2} holds, we set
\begin{align*}
\dot \Omega_t = \{w = x + \Psi (x, t) \mid x \in \dot \Omega\} \qquad (t \in (0, T)).
\end{align*}
For simplicity of notation, in the following, we may use the symbol $\Psi (x, t) = \chi(3 d_\Gamma / r_0) H_h (x, t) \bn_* (x)$. Let $\pd w/\pd x$ be the Jacobi matrix of the transformation \eqref{2.1}, that is,
\begin{align*}
\frac{\pd w}{\pd x} = \bI + \nabla \Psi(x, t), \quad \nabla \Psi = (\pd_i \Psi_j), \quad \bigg(\pd_i \Psi_j := \frac{\pd \Psi_j}{\pd x_i} \bigg),
\end{align*}
where we have set $\Psi(x, t) = (\Psi_1(x, t), \dots, \Psi_N(x, t))$. If $\widetilde \varepsilon \in (0, 1)$ is suitably small, we see that
\begin{align*}
\bigg(\frac{\pd w}{\pd x} \bigg)^{- 1} = \bI + \sum_{k = 1}^\infty (- (\nabla \Psi(x, t)))^k
\end{align*}
exists. Hence, there exists an $N \times N$ matrix $\bM_0 (\bm)$ of $C^\infty$ functions defined on $\lvert \bm \rvert < \widetilde \varepsilon$ such that $\bM_0 (0) = 0$ and $(\pd w / \pd x)^{- 1} = \bI + \bM_0 (\nabla \Psi(x, t))$. Here and in the following, we write $\bm = (m_{ij})$ and $m_{ij}$ denote the variables corresponding to $\pd_i \Psi_j$ $(1 \leq i, j \leq N)$. From the assumption \eqref{2.2}, we have the estimates
\begin{equation*}
\lVert \bM_0 (\nabla \Psi) \rVert_{L_\infty (\dot \Omega)} \leq C \widetilde \varepsilon, \quad \lVert \nabla^n \bM_0 (\nabla \Psi) \rVert_{L_q (\dot \Omega)} \leq C \lVert \nabla^{1 + n} \Psi \rVert_{L_q (\dot \Omega)} \quad (n = 0, 1, 2, 3),
\end{equation*}
Let $\rho_+ (x, t) = \varrho_+ (x + \Psi(x, t), t) - \rho_{* +}$, $\bu_\pm (x, t) = \bv_\pm (x + \Psi (x, t), t)$, and $\pi_- (x, t) = \fp_- (x + \Psi(x, t), t) - \pi_{* -}$. \par
Let $M_{0ij} (\bm)$ be the $(i, j)$th component of $\bM_0 (\bm)$ and $\nabla_x$ and $\nabla_w$ be the gradient with respect to $x$ and $w$, respectively. In addition, let $\Delta_x = \nabla_x \cdot\nabla_x$ and $\Delta_w = \nabla_w \cdot \nabla_w$ be the Laplace operator with respect to $x$ and $w$, respectively. We see that
\begin{align*}
\nabla_w = (\bI + \bM_0 (\bm)) \nabla_x, \quad \frac{\pd}{\pd w_i} = \sum_{j = 1}^N (\delta_{ij} + M_{0ij} (\bm)) \frac{\pd}{\pd x_j}.
\end{align*}
We then observe that
\begin{gather}
\label{2.4}
\bD_w (\bv_\pm) = \bD_x (\bu_\pm) + \CD_\bD \nabla_x \bu_\pm, \quad (\CD_\bD \nabla_x \bu_\pm)_{ij} = \sum_{k = 1}^N \bigg(M_{0ij} (\bm) \frac{\pd \bu_{i \pm}}{\pd x_k} \bigg), \\
\label{2.5}
\dv_w \bv_\pm = \sum_{j = 1}^N \frac{\pd v_{j \pm}}{\pd w_j} = \sum_{j, k = 1}^N (\delta_{ij} + M_{0ij} (\bm)) \frac{\pd \bu_{j\pm}}{\pd x_k} = \dv_x \bu_\pm + \bM_0 (\bm) \colon \nabla_x \bu_\pm.
\end{gather}
Let $J = J (\bm)$ be the Jacobi matrix of the transformation of \eqref{2.1}. Choosing $\widetilde \varepsilon \in (0, 1)$ small enough, we may assume that $J = 1 + J_0 (\bm)$, where $J_0 (\bm)$ is a $C^\infty$ function defined for $\lvert \bm \rvert < \widetilde \varepsilon$ such that $J_0 (0) = 0$. To obtain the representation formula of $\dv_w \bv_-$, we use the inner product $(\cdot, \cdot)_{\Omega_{t -}}$. Set $\zeta_x (x) = \zeta_w (w)$ for any $\zeta_w \in C^\infty_0 (\Omega_{t -})$. Then we obtain
\begin{align*}
(\dv_w \bv_-, \zeta_w)_{\Omega_{t -}} & = - (\bv_-, \nabla_w \zeta_w)_{\Omega_{t -}} \\
& = - (J (\bm) \bu_-, (\bI + \bM_0) \nabla_x \zeta_x)_{\Omega_-} \\
& = (\dv_x (J (\bm) (\bI + {}^\top\!\bM_0)) \bu_-, \zeta_x)_{\Omega_-} \\
& = (J^{- 1} (\bm) \dv_x (J (\bm) (\bI + {}^\top\!\bM_0) \bu_-), \zeta_w)_{\Omega_{t -}}.
\end{align*}
Summing up, we have
\begin{align*}
\dv_w \bv_- = \dv_x \bu_- + \bM_0 (\bm) \colon \nabla_x \bu_- = J^{- 1} \dv_x (J (\bI + {}^\top\!\bM_0) \bu_-),
\end{align*}
which yields
\begin{align}
\label{2.7*}
J_0 (\bm) \dv_x \bu_- + J (\bm) \bM_0 (\bm) \colon \nabla_x \bu_- = \dv_x (J (\bm) (\bI + {}^\top\!\bM_0) \bu_-).
\end{align}
Setting
\begin{align*}
f_d & = f_d (\bu_-, h) = - (J_0 (\bm) \dv_x \bu_- + (1 + J_0 (\bm)) \bM_0 (\bm) \colon \nabla_x \bu_-), \\
\bF_d & = \bF_d (\bu_-, h) = - (1 + J_0 (\bm)) {}^\top \bM_0 (\bm) \bu_-,
\end{align*}
the divergence-free condition $\dv_w \bu_- = 0$ is equivalent to
\begin{align}
\label{l-1}
\dv_x  \bu_- = f_d = \dv \bF_d \quad \text{in $\Omega$}.
\end{align}
Since
\begin{align*}
\frac{\pd}{\pd t} \Big(v_{i \pm} (x + \Psi(x, t), t) \Big) = \frac{\pd v_{i \pm}}{\pd t} (w, t) + \sum_{j = 1}^N \frac{\pd \Psi_j}{\pd t} \frac{\pd v_{i \pm}}{\pd w_j} (w, t),
\end{align*}
we see that
\begin{align}
\label{2.7}
\frac{\pd v_{i \pm}}{\pd t} (w, t) = \frac{\pd u_{i \pm}}{\pd t} (x, t) - \sum_{j, k = 1}^N \frac{\pd \Psi_j}{\pd t} (\delta_{jk} + M_{0jk} (\bm)) \frac{\pd u_{i \pm}}{\pd x_k} (x, t).
\end{align}
Analogously, we have
\begin{align}
\label{2.8}
\frac{\pd \varrho_+}{\pd t} (w, t) = \frac{\pd \rho_+}{\pd t} (x, t) - \sum_{j, k = 1}^N \frac{\pd \Psi_j}{\pd t} (\delta_{jk} + M_{0jk} (\bm)) \frac{\pd \rho_+}{\pd x_k} (x, t).
\end{align}
Hence, from \eqref{2.5} and \eqref{2.8} the first equation in \eqref{eq-1.1} is transformed into
\begin{align}
\label{l-2}
\pd_t \rho_+ + \rho_{* +} \dv_x \bu_+ = f_M (\rho_+, \bu_+, h),
\end{align}
where we have set
\begin{align*}
f_M (\rho_+, \bu_+, h)  & = \sum_{j, k = 1}^N \frac{\pd \Psi_j}{\pd t} (\delta_{jk} + M_{0jk} (\bm)) \frac{\pd \rho_+}{\pd x_k} + \rho_{* +} \bM_0 (\bm) \colon \nabla_x \bu_+ + \langle \bu_+, (\bI + \bM_0 (\bm)) \nabla_x \rho_+ \rangle.
\end{align*}
Next, by \eqref{2.4} and \eqref{2.5}, we observe that the $i$th component of $\DV_w \BT_+$ can be written in the form of
\begin{equation}
\label{2.10}
\begin{split}
& \sum_{j = 1}^N \frac{\pd}{\pd w_j} \bigg\{\mu_+ \bD_w (\bv_+)_{ij} + \Big((\nu_+ - \mu_+) \dv_w \bv_+ - \fp_+ + \frac{\kappa_+}{2} \lvert \nabla_w \varrho_+ \rvert^2 + \kappa_+ \varrho \Delta_w \varrho_+\Big) \delta_{ij} - \kappa_+ \pd_i \varrho_+ \pd_j \varrho_+ \bigg\} \\
& = \sum_{j,k = 1}^N (\delta_{jk} + M_{0jk} (\bm)) \frac{\pd}{\pd x_k} \left\{\mu_+ \Big(\bD_x (\bu_+)_{ij} + (\CD_\bD (\bm) \nabla_x \bu_+)_{ij} \Big) \right\} \\
& \quad + \sum_{j,k = 1}^N \delta_{ij} (\delta_{jk} + M_{0jk}) \frac{\pd}{\pd x_j} \left\{(\nu_+ - \mu_+) \Big(\dv_x \bu_+ + \bM_0 (\bm) \colon \nabla_x \bu_+ \Big) \right\} \\
& \quad - \sum_{j,k = 1}^N \delta_{ij} (\delta_{jk} + M_{0jk}) \bigg\{\fp'_+ (\rho_{* +}) \frac{\pd \rho_+}{\pd x_j} + \frac{\pd}{\pd x_j} \bigg(\rho_+^2 \int_0^1 (1 - \theta) \fp''_+ (\rho_{* +} + \theta \rho_+) \dtheta \bigg) \bigg\} \\
& \quad + \sum_{j = 1}^N \delta_{ij} (\delta_{jk} + M_{0jk}) \frac{\pd}{\pd x_j} \bigg(\frac{\kappa_+}{2} \lvert (\bI+\bM_0(\bm))\nabla_x \rho \rvert^2 \\
& \qquad + \kappa_+ (\rho_{* +} + \rho_+) \Big(\Delta_x \rho_+ + \dv_x (\bM_0(\bm) \nabla_x \rho_+) + \bM_0 (\bm) \colon \nabla_x ((\bI + \bM_0 (\bm)) \nabla_x \rho_+) \Big) \bigg) \\
& \quad + \sum_{j,k = 1}^N (\delta_{jk} + M_{0jk} (\bm)) \frac{\pd}{\pd x_k} \bigg\{\kappa_+ \bigg(\sum_{l = 1}^N (\delta_{il} + M_{0il}) \frac{\pd \rho_+}{\pd x_l}\bigg) \bigg(\sum_{m = 1}^N (\delta_{im} + M_{0im}) \frac{\pd \rho_+}{\pd x_m}\bigg) \bigg\}.
\end{split}
\end{equation}
Here, we have used the following identity:
\begin{align*}
\fp_+ (\varrho_+) = \fp_+ (\rho_{* +} + \rho_+) = \fp_+ (\rho_{* +}) + \fp'_+ (\rho_{* +}) \rho_+ + \rho_+^2 \int_0^1 (1 - \theta) \fp''_+ (\rho_{* +} + \theta \rho_+) \dtheta,
\end{align*}
where $\fp'_+$ and $\fp''_+$ denote the first and second derivative of $\fp$ with respect to $\varrho_+$, respectively. In the sequel, we set $\pi_{* +} = \fp_+ (\rho_{* +})$. Combining with \eqref{2.7} and \eqref{2.10}, from the second equation of \eqref{eq-1.1}, we have
{\allowdisplaybreaks \begin{align*}
0 & = (\rho_{* +} + \rho_+) \bigg\{\frac{\pd u_{i +}}{\pd t} + \sum_{j,k = 1}^N \bigg(u_{j +} - \frac{\pd \Psi_j}{\pd t}\bigg) (\delta_{jk} + M_{0jk} (\bm)) \frac{\pd u_{i +}}{\pd x_k} \bigg\} \\
& \quad - \sum_{j,k = 1}^N (\delta_{jk} + M_{0jk} (\bm)) \frac{\pd}{\pd x_k} \left\{\mu_+ \Big(\bD_x (\bu_+)_{ij} + (\CD_\bD (\bm) \nabla_x \bu_+)_{ij} \Big) \right\} \\
& \quad - \sum_{j,k = 1}^N \delta_{ij} (\delta_{jk} + M_{0jk}) \frac{\pd}{\pd x_j} \left\{(\nu_+ - \mu_+) \Big(\dv_x \bu_+ + \bM_0 (\bm) \colon \nabla_x \bu_+ \Big) \right\} \\
& \quad + \sum_{j,k = 1}^N \delta_{ij} (\delta_{jk} + M_{0jk}) \bigg\{\fp'_+ (\rho_{* +}) \frac{\pd \rho_+}{\pd x_j} + \frac{\pd}{\pd x_j} \bigg(\rho_+^2 \int_0^1 (1 - \theta) \fp''_+ (\rho_{* +} + \theta \rho_+) \dtheta \bigg) \bigg\} \\
& \quad - \sum_{j = 1}^N \delta_{ij} (\delta_{jk} + M_{0jk}) \frac{\pd}{\pd x_j} \bigg(\frac{\kappa_+}{2} \lvert (\bI+\bM_0(\bm))\nabla_x \rho \rvert^2 \\
& \qquad + \kappa_+ (\rho_{* +} + \rho_+) \Big(\Delta_x \rho_+ + \dv_x (\bM_0(\bm) \nabla_x \rho_+) + \bM_0 (\bm) \colon \nabla_x ((\bI + \bM_0 (\bm)) \nabla_x \rho_+) \Big) \bigg) \\
& \quad - \sum_{j,k = 1}^N (\delta_{jk} + M_{0jk} (\bm)) \frac{\pd}{\pd x_k} \bigg\{\kappa_+ \bigg(\sum_{l = 1}^N (\delta_{il} + M_{0il}) \frac{\pd \rho_+}{\pd x_l}\bigg) \bigg(\sum_{m = 1}^N (\delta_{im} + M_{0im}) \frac{\pd \rho_+}{\pd x_m}\bigg) \bigg\}
\end{align*}}\noindent
for $i = 1, \dots, N$. Hence, we define an $N$-vector of functions $\bff_+ (\rho_+, \bu_+, h)$ by
\begin{align*}
& \bff_+(\rho_+, \bu_+, h) \vert_i \\
& = - \rho_+ \frac{\pd u_{i +}}{\pd t} - (\rho_{* +} + \rho_+) \bigg\{\sum_{j,k = 1}^N \bigg(u_{j +} - \frac{\pd \Psi_j}{\pd t}\bigg) (\delta_{jk} + M_{0jk} (\bm)) \frac{\pd u_{i +}}{\pd x_k} \bigg\} \\
& \quad + \sum_{j,k = 1}^N \left[\delta_{jk} \frac{\pd}{\pd x_k} \Big(\mu_+ (\CD_\bD (\bm) \nabla_x \bu_+)_{ij} \Big) + M_{0jk} (\bm) \frac{\pd}{\pd x_k} \left\{\mu_+ \Big(\bD_x (\bu_+)_{ij} + (\CD_\bD (\bm) \nabla_x \bu_+)_{ij} \Big) \right\} \right] \\
& \quad + \sum_{j,k = 1}^N \delta_{ij} \left[ \delta_{jk} \frac{\pd}{\pd x_j} \Big((\nu_+ - \mu_+) (\bM_0 (\bm) \colon \nabla_x \bu_+) \Big) + M_{0jk} \frac{\pd}{\pd x_j} \left\{(\nu_+ - \mu_+) \Big(\dv_x \bu_+ + \bM_0 (\bm) \colon \nabla_x \bu_+ \Big) \right\} \right] \\
& \quad - \sum_{j,k = 1}^N \delta_{ij} \left[M_{0jk} \fp'_+ (\rho_{* +}) \frac{\pd \rho_+}{\pd x_j} + (\delta_{jk} + M_{0jk}) \frac{\pd}{\pd x_j} \bigg(\rho_+^2 \int_0^1 (1 - \theta) \fp''_+ (\rho_{* +} + \theta \rho_+) \dtheta \bigg)\right] \\
& \quad + \sum_{j = 1}^N \delta_{ij} (\delta_{jk} + M_{0jk}) \frac{\pd}{\pd x_j} \bigg(\frac{\kappa_+}{2} \lvert (\bI+\bM_0(\bm))\nabla_x \rho \rvert^2 \\
& \qquad + \kappa_+ \rho_+ \Big(\Delta_x \rho_+ + \dv_x (\bM_0(\bm) \nabla_x \rho_+) + \bM_0 (\bm) \colon \nabla_x ((\bI + \bM_0 (\bm)) \nabla_x \rho_+) \Big) \bigg) \\
& \quad + \rho_{* +} \sum_{j = 1}^N \delta_{ij} M_{0jk} \frac{\pd}{\pd x_j} \Big(\kappa_+ \Delta_x \rho_+ \Big) \\
& \quad + \sum_{j,k = 1}^N (\delta_{jk} + M_{0jk} (\bm)) \frac{\pd}{\pd x_k} \bigg\{\kappa_+ \bigg(\sum_{l = 1}^N (\delta_{il} + M_{0il}) \frac{\pd \rho_+}{\pd x_l}\bigg) \bigg(\sum_{m = 1}^N (\delta_{im} + M_{0im}) \frac{\pd \rho_+}{\pd x_m}\bigg) \bigg\},
\end{align*}
where $\bff_+ (\rho_+, \bu_+, h) \vert_i$ denotes the $i$th element of $\bff_+ (\rho_+, \bu_+, h)$. We, therefore, see that the second equation in \eqref{eq-1.1} is transformed to
\begin{equation}
\label{l-3}
\rho_{* +} \pd_t \bu_+ - \DV \bT_+ (\gamma_1, \gamma_2, \gamma_3, \bu_+, \rho_+) = \bff_+ (\rho_+, \bu_+, h) \quad \text{in $\Omega_+\times (0,T)$}
\end{equation}
with $\gamma_{* +} = \fp'_+ (\rho_{* +})$. Using the similar argument above, we easily see that the fourth equation in \eqref{eq-1.1} is transformed into
{\allowdisplaybreaks \begin{align*}
0 & = \rho_{* -} \bigg\{\frac{\pd u_{i -}}{\pd t} + \sum_{j,k = 1}^N \bigg(u_{j -} - \frac{\pd \Psi_j}{\pd t}\bigg) (\delta_{jk} + M_{0jk} (\bm)) \frac{\pd u_{i -}}{\pd x_k} \bigg\} \\
& \quad - \sum_{j,k = 1}^N (\delta_{jk} + M_{0jk} (\bm)) \frac{\pd}{\pd x_k} \left\{\mu_+ \Big(\bD_x (\bu_+)_{ij} + (\CD_\bD (\bm) \nabla_x \bu_+)_{ij} \Big) \right\} + \sum_{j = 1}^N (\delta_{ij} + M_{0ij} (\bm)) \frac{\pd}{\pd x_j} \pi_-.
\end{align*}}\noindent
Since $(\bI + \nabla \Psi) (\bI + \bM_0 (\bm)) = (\pd w / \pd x) (\pd x / \pd w) = \bI$, we have
\begin{align*}
\sum_{i = 1}^N (\delta_{ni} + \pd_n \Psi_i) (\delta_{ij} + M_{0ij} (\bm)) = \delta_{nj}
\end{align*}
for each $n = 1, \dots, N$, which yields that
\begin{align*}
0 & = \rho_{* -} \sum_{i = 1}^N (\delta_{ni} + \pd_n \Psi_i) \bigg\{\frac{\pd u_{i -}}{\pd t} + \sum_{j,k = 1}^N \bigg(u_{j -} - \frac{\pd \Psi_j}{\pd t}\bigg) (\delta_{jk} + M_{0jk} (\bm)) \frac{\pd u_{i -}}{\pd x_k} \bigg\} \\
& \quad - \sum_{i,j,k = 1}^N (\delta_{ni} + \pd_n \Psi_i) (\delta_{jk} + M_{0jk} (\bm)) \frac{\pd}{\pd x_k} \left\{\mu_+ \Big(\bD_x (\bu_+)_{ij} + (\CD_\bD (\bm) \nabla_x \bu_+)_{ij} \Big) \right\} + \frac{\pd}{\pd x_n} \pi_-.
\end{align*}
Hence, changing $i$ to $l$ and $n$ to $i$ in the above identity, we arrive at
\begin{align}
\label{l-4}
\rho_{* -} \pd_t \bu_- - \DV \bT_- (\gamma_4, \bu_-, \pi_-) = \bff_- ( \bu_-, h) \quad \text{in $\Omega_- \times (0, T)$},
\end{align}
where we have set
\begin{align*}
& \bff_- (\bu_-, h) \vert_i \\
& = - \rho_{* -} \bigg\{\sum_{j,k = 1}^N \bigg(u_{j -} - \frac{\pd \Psi_j}{\pd t}\bigg) (\delta_{jk} + M_{0jk} (\bm)) \frac{\pd u_{n -}}{\pd x_k} \bigg\} \\
& \quad - \rho_{* -} \sum_{l = 1}^N \pd_i \Psi_l \bigg\{\frac{\pd u_{l -}}{\pd t} + \sum_{j,k = 1}^N \bigg(u_{j -} - \frac{\pd \Psi_j}{\pd t}\bigg) (\delta_{jk} + M_{0jk} (\bm)) \frac{\pd  u_{l -}}{\pd x_k} \bigg\} \\
& \quad + \sum_{j = 1}^N \frac{\pd}{\pd x_j} \Big(\mu_- (\CD_\bD (\bm) \nabla_x \bu_-)_{ij}\Big) + \sum_{j,k = 1}^N M_{0jk} (\bm) \frac{\pd}{\pd x_k} \Big\{\mu_- \Big(\bD_x(\bu_-)_{ij} + (\CD_\bD (\bm)\nabla_x \bu_-)_{ij} \Big)\Big\} \\
&\quad + \sum_{j,k,l = 1}^N \pd_i \Psi_l (\delta_{jk} + M_{0jk} (\bm)) \frac{\pd}{\pd x_k}\Big\{\mu_- \Big(\bD_x (\bu_-)_{lj} + (\CD_\bD (\bm) \nabla_x \bu_-)_{lj} \Big)\Big\}.
\end{align*}
Here, $\bff_- (\bu_-, h) \vert_i$ stands the $i$th component of $\bff_- (\bu_-, h)$. We emphasize that $\bff_-$ is independent of $\pi_-$.

\subsection{Laplace-Beltrami operator}
We start with the following proposition, which was essentially proved by Enomoto and Shibata~\cite[Appendix]{ES2013}, see also Shibata~\cite[Proposition~2.2]{Shi2019}.
\begin{prop}
\label{prop-2.1}\label{prop-6.1}
Let $\Omega_+$ and $\Omega_-$ are uniform $W^{4, 3}_r$ and $W^{4, 2}_r$ domain in $\BR^N$ for $N < r < \infty$, respectively. {\color{blue} Let $\Gamma := \Gamma^0$, $\Gamma_+ := \Gamma^1$, and $\Gamma_- := \Gamma^2$.} Furthermore, let $M_1 \in (0, 1)$ be any given positive number. Then there exist {\color{red} positive constants $M_2 \geq 1$ and $d^i \in (0, 1)$, $i = 0, \dots, 4$, and at most countably many $N$-vector of functions $\Phi_j^0 \in W^4_r (\BR^N)^N$, $\Phi^1_j \in W^3_r (\BR^N)^N$, $\Phi^2_j \in W^2_r (\BR^N)^N$, points $x^i_j \in \Gamma^i$ $(i = 0,1,2)$, $x^3_j \in \Omega_+$, and $x^4_j \in \Omega_-$} such that the following assertions hold:
\begin{enumerate}\renewcommand{\labelenumi}{(\arabic{enumi})}
\item The mappings: $\BR^N\ni x\mapsto {\color{red} \Phi_j^i(x)} \in \BR^N$ {\color{blue} ($i = 0,1,2$, $j \in \BN$)} are bijections of $C^1$-class {\color{blue} such that $\nabla \Phi_j^i = \bA^i_j + \bB^i_j$ and $\nabla (\Phi^i_j)^{- 1} = \bA^i_{j -} + \bB^i_{j -}$, respectively, where $\bA^i_j$ and $\bA^i_{j -}$ are $N \times N$ orthogonal matrices with constant coefficients and $\bB_j^i$ and $\bB_{j -}^i$ are $N \times N$ matrices of $W^{3 - i}_r (\BR^N)$ functions satisfying $\lVert (\bB^i_j, \bB_{j -}^i) \rVert_{L_\infty (\BR^N)} \leq M_1$ and $\lVert (\nabla \bB_j^i, \nabla \bB_{j -}^i) \rVert_{W^{3 - i}_r (\BR^N)} \leq M_2$.}
\item {\color{red}
Set $D^0=\BR^N$, $D^1=\RP$, and $D^2=\RM$.
Then it holds that
\begin{align*}
\Omega=\bigg\{\bigcup_{i=0,1,2}\bigcup_{j=1}^\infty
(\Phi^i_j(D^i)\cap B_{d^i}(x^i_j)) \bigg\}\cup
\bigg\{\bigcup_{i=3,4}\bigcup_{j=1}^\infty B_{d^i}(x^i_j) \bigg\},	
\end{align*}
where $\Phi^i_j(D^i)\cap B_{d^i}(x^i_j)=
\Omega_+\cap B_{d^i}(x^i_j)$ $(i=0,1)$,
$\Phi^i_j(D^i)\cap B_{d^i}(x^i_j)=
\Omega_- \cap B_{d^i}(x^i_j)$ $(i=0,2)$,
$B_{d^3}(x^3_j)\subset \Omega_+$,
$B_{d^4}(x^4_j)\subset \Omega_-$, and
$\Phi^i_j(\RZ)\cap B_{d^i}(x^i_j)
=\Gamma^i\cap B_{d^i}(x^i_j)$ $(i=0,1,2)$.
}
\item {\color{blue} There exists $C^\infty$ functions $\zeta^i_j$ and $\tzeta^i_j$, $(i = 0, \dots, 4, j \in \BN)$ such that
\begin{gather*}
0\leq \zeta^i_j,\tzeta^i_j\leq 1,\quad
\supp\zeta^i_j,\enskip\supp\tzeta^i_j\subset B_{d^i}(x^i_j),\quad
\|(\zeta^i_j,\tzeta^i_j)\|_{W^3_\infty(\BR^N)}\leq c, \\
\tzeta^i_j=1 \enskip \text{ on $\supp\zeta^i_j$},\quad 
\sum_{i=0}^{4}\sum_{j=1}^{\infty}\zeta^i_j=1
\enskip\text{ on $\overline{\Omega}$},\quad
\sum_{j=1}^\infty\zeta^i_j=1\enskip \text{ on $\Gamma^i$ $(i=0,1,2)$}.	
\end{gather*}
Here, $c$ is a positive constant independent of $j\in\BN$.
}
\item There exists an integer $L \geq 2$ such that any $L + 1$ distinct sets of $\{{\color{red} B_{d^i} (x^i_j) \mid i = 0, \dots, 4, \enskip j \in \BN}\}$ have an empty intersection.
\end{enumerate}	
\end{prop}
We next introduce the Laplace-Beltrami operator defined on $\Gamma_t$ and $\Gamma$. In the following, we write $B_m = {\color{red} B_{d^0}} (x_m)$ and $\Gamma_m = {\color{red} \Phi_m^0} (\RZ)$ $(m \in \BN)$ for short. To this end, let $p = \{p_1, \dots, p_N\}$ be a local coordinate system in a neighbourhood of $x_l\in\Gamma$ $(l\in\BN)$ such that
\begin{align}
\Omega_\pm \cap B_m = \{x = {\color{red} \Phi_m^0} (p) \mid p \in \BR^N_\pm \}\cap B_m, \quad \Gamma\cap B_m = \{y = {\color{red} \Phi_m^0} (p',0) \mid (p',0) \in \RZ \} \cap B_m
\end{align}
where $p' = (p_1, \dots, p_{N - 1}) \in \BR^{N - 1}$. By abuse of notation, we let $x_m (p)$ stand for ${\color{red} \Phi_m^0} (p)$. Let $G$ be the first form on $\Gamma_m$ such that
\begin{align*}
G = (g_{ij}), \quad g_{ij} = \frac{\pd x}{\pd p_i} \frac{\pd x}{\pd p_j}, \quad (1 \leq i,j \leq N - 1).
\end{align*}
In addition, let $G^{- 1} = (g^{ij})$ be the inverse matrix of $G$, hence $g_{ij} g^{ij} = \delta_{ij}$. We then define the Laplace-Beltrami operator $\Delta_\Gamma$ on $\Gamma$ by
\begin{align}
\label{2.15}
\Delta_\Gamma f = \frac{1}{\sqrt{\det G}} \sum_{i,j = 1}^{N - 1} \frac{1}{\pd p_i} \bigg(\sqrt{\det G} g^{ij} \frac{\pd f}{\pd p_j} \bigg),
\end{align}
which is defined on $\BR^N$. Notice that by Proposition \ref{prop-2.1} (3), we see that
\begin{align}
\label{2.16}
\bigg\lVert \bigg(\frac{\pd x}{\pd p_i}, g_{ij}, g^{ij}, \sqrt{\det G} \bigg) \bigg\rVert_{W^2_\infty (\BR^N)} \leq C_{M_2}.
\end{align}
Let $H_\Gamma$ be the $(N - 1)$-times mean curvature of $\Gamma$, which is given by $H_\Gamma = \langle \Delta_\Gamma x, \bn \rangle$ on $\Gamma \cap B_0$. Using the symbols defined above, $H_\Gamma$ can be written as
\begin{align*}
H_\Gamma = \sum_{i,j = 1}^{N - 1} g^{ij} \bigg\langle \frac{\pd^2 x}{\pd p_i \pd p_j}, \bn \bigg\rangle \bigg \vert_{p_N = 0} \quad \text{ on $\Gamma \cap B_m$}
\end{align*}
because $\langle \bn, (\pd \bn / \pd p_i)\rangle = 0$ as follows from $\lvert \bn \rvert = 1$. \par
We finally consider the Laplace-Beltrami operator $\Delta_{\Gamma_t}$ on $\Gamma_t$. Recall that $\Gamma_t$ is defined by $\Gamma_t = \{w = x + h(x, t) \bn(x) \mid x \in \Gamma\}$ for $t \in (0, T)$. Let $G_t$ be the first fundamental form on $\Gamma_t$ such that
\begin{align*}
G_t = (g_{t, ij}), \quad g_{t, ij} = \frac{\pd w}{\pd p_i} \cdot \frac{\pd w}{\pd p_j}, \quad (1 \leq i,j \leq N - 1).
\end{align*}
for $w \in \Gamma_t\cap B_m$. Furthermore, let $G_t^{- 1} = (g^{ij}_t)$ be the inverse matrix of $\Gamma_t$. Using these symbols, the Laplace-Beltrami operator $\Delta_{\Gamma_t}$ on $\Gamma_t$ is represented by
\begin{align*}
\Delta_{\Gamma_t} f = \frac{1}{\sqrt{\det G}} \sum_{i,j = 1}^{N - 1} \frac{1}{\pd p_i} \bigg(\sqrt{\det G} g^{ij} \frac{\pd f}{\pd p_j} \bigg) \quad \text{ for $x \in \Gamma \cap B_m$}.
\end{align*}
By the definition of transformation, we have
\begin{align*}
\frac{\pd w}{\pd p_i} = \frac{\pd x}{\pd p_i} + \sum_{k = 1}^N \bigg(\frac{\pd}{\pd x_k}(h(x, t) \bn(x)) \bigg) \frac{\pd x_k}{\pd p_i}.
\end{align*}
Hence, choosing $\widetilde \varepsilon > 0$ sufficiently small in \eqref{2.2} and using \eqref{2.16}, there exist scalar functions $\bG_1 (\bm)$ and $\bG_{2ij} (\bm)$ of $C^2$-class defined on $\BR^N \times B_{\widetilde \varepsilon} (0)$ such that
\begin{gather}
\bG_1 (0) = \bG_{2ij} (0) = 0, \quad \lVert (\bG_1, \bG_{2ij}) \rVert_{W^2_\infty (\BR^N \times B_{\widetilde \varepsilon} (0))} \leq C_{M_2}, \\
\label{2.18}
(\Delta_{\Gamma_t} - \Delta_\Gamma) f = \sum_{i,j = 1}^{N - 1} \left\{\bG_1 (\bm) \frac{\pd}{\pd p_i} \bigg((\sqrt{\det G} g^{ij} + \bG_{2ij} (\bm)) \frac{\pd f}{\pd p_j}\bigg) + \frac{1}{\sqrt{\det G}} \frac{\pd}{\pd p_i} \bigg(\bG_{2ij} (\bm) \frac{\pd f}{\pd p_j} \bigg)\right\}.
\end{gather}
on $\Gamma \cap B_m$ with $\bm = (\nabla (h(x, t) \bn (x))) \circ x (p)$.

\subsection{Transformation of the kinetic equation}
In this subsection, we consider the interface condition~\eqref{eq-1.2}. For this purpose, we first treat the outer unit normal $\bn_t$. Since
\begin{align*}
0 = \langle \bn, dx \rangle = \bigg\langle \bn, \frac{\pd x}{\pd w} dw \bigg\rangle = \langle\bn, (\bI + \bM_0 (\bm)) dx \rangle = \langle(\bI + {}^\top\!\bM_0 (\bm)) \bn, dx \rangle
\end{align*}
on $\Gamma$, the outer unit normals $\bn_t$ and $\bn$ have the following relationship:
\begin{align}
\label{2.19*}
\bn_t = \frac{(\bI + {}^\top\!\bM_0 (\bm)) \bn}{\lvert (\bI + {}^\top\!\bM_0 (\bm)) \bn\rvert}.
\end{align}
Choosing $\widetilde \varepsilon$ in \eqref{2.2} small enough, we observe that there exist an $N$-vector function $\bN (\bm)$ defined on $\BR^N \times B^{N^2}_{\widetilde \varepsilon} (0)$ such that $\bN (0) = 0$, $\lVert \bN \rVert_{W^2_\infty (\BR^N) \times B^{N^2}_{\widetilde \varepsilon} (0)} \leq C_{M_2}$, and 
\begin{align}
\label{2.19}
\bn_t = \bn + \bN (\bm).
\end{align}
Here we have set $B^{N^2}_{\widetilde \varepsilon} (0) = \{\bm \in \BR^{N^2} \mid \lvert \bm \rvert < \widetilde \varepsilon \}$. \par
Since $w = x + h(x, t) \bn (x)$ on $x \in \Gamma$, we have $V_{\Gamma_t} = \langle \pd w / \pd t, \bn_t \rangle = \langle (\pd_t h) \bn, \bn_t \rangle$. Then the kinetic equation can be written in the form of
\begin{align*}
\langle (\pd_t h) \bn, \bn_t \rangle = \frac{\langle \{\rho_{* -} \bu_- \vert_- - (\rho_{* +} + \rho_+) \bu_+ \vert_+\}, \bn_t \rangle}{\rho_{* -} - (\rho_{* +} + \rho_+) \vert_+},
\end{align*}
which, combined with \eqref{2.19}, yields
\begin{align}
\label{l-5}
\pd_t h - - \frac{\langle \rho_{* -} \bu_-, \bn\rangle \vert_- - \langle \rho_{* +} \bu_+, \bn\rangle \vert_+}{\rho_{* -} - \rho_{* +}} = d(\rho_+, \bu_+, \bu_-, h) \quad \text{ on $\Gamma$}
\end{align}
with
\begin{align*}
d (\rho_+, \bu_+, \bu_-, h) & = - \langle \bn, \bN (\bm) \rangle (\pd_t h) + \frac{\rho_+ \Big\vert_+ \Big(\rho_{* -} \langle \bu_-, \bn \rangle \vert_- - \rho_{* +} \langle \bu_+, \bn \rangle \vert_+ \Big)}{(\rho_{* -} - \rho_{* +}) \{\rho_{* -} - (\rho_{* +} + \rho_+) \vert_+ \}} \\
& \quad + \frac{\rho_{* -} \langle \bu_-, \bN (\bm) \rangle \vert_- - \rho_{* +} \langle \bu_+, \bN (\bm) \rangle \vert_+ - \rho_+ \langle \bu_+, \bn + \bN (\bm) \rangle \vert_+}{\rho_{* -} - (\rho_{* +} + \rho_+) \vert_+}.
\end{align*}

\subsection{Transformation of the interface condition}
To transform the third jump condition in \eqref{eq-1.2}, we use the following lemma proven by Shibata and Shimizu~\cite[Lemma 2.1]{SS2011} (cf. Solonnikov~\cite[p.155]{Sol2003}).
\begin{lemm}
\label{lem-2.2}
If $\langle \bn_t, \bn \rangle \neq 0$, then $\bd = 0$ is equivalent to
\begin{align*}
\Pi_\bn \Pi_{\bn_t} \bd = 0 \quad \text{and} \quad \langle \bd, \bn \rangle = 0
\end{align*}	
for any $N$-vector field $\bd$.	
\end{lemm}
By Lemma \ref{lem-2.2}, the third condition in \eqref{eq-1.2} is equivalent to the following two condition:
\begin{align}
\label{2.21}
\Big\lbrac \Pi_\bn \Pi_{\bn_t} (- \j \bv + \mu \bD_w (\bv) \bn_t) \Big\rbrac + \Pi_\bn \Pi_{\bn_t} \Big(\kappa_+ (\nabla_w \varrho_+ \otimes \nabla_w \varrho_+) \bn_t \Big) \Big\vert_+ & = 0 , \\
\label{2.22}
\Big\lbrac - \langle \j \bv, \bn \rangle + \langle \BT \bn_t, \bn \rangle \Big\rbrac - \sigma \langle \Delta_{\Gamma_t} (x + h \bn), \bn \rangle & = 0 
\end{align}
on $\Gamma_t$ for $t \in (0, T)$. From \eqref{2.19}, using $\Pi_{\bn_0} \Pi_{\bn_0} = \Pi_{\bn_0}$, we rewrite \eqref{2.21} as
\begin{equation}
\label{l-6}
\begin{split}
& \Pi_\bn (\mu_- \bD_x (\bu_-) \bn) \vert_- - \Pi_\bn (\mu_+ \bD_x (\bu_+) \bn) \vert_+ \\
& = \Big\lbrac \Pi_\bn (\Pi_\bn - \Pi_{\bn_t}) (\mu_- \bD_x (\bu_-) (\bn + \bN (\bm))) + \Pi_\bn (\mu_- \bD_x (\bu_-) \bN (\bm)) \\
& \qquad - \Pi_\bn \Pi_{\bn_t} (- \j  \bu_- + (\CD_\bD (\bm) \nabla_x \bu_-) (\bn + \bN (\bm))) \Big\rbrac \\
& \quad + \Pi_\bn \Pi_{\bn_t} \Big(\kappa_+ \big((\bI + \bM_0 (\bm)) \nabla_x \rho_+ \big) \otimes \big((\bI + \bM_0 (\bm)) \nabla_x \rho_+ \big) \Big)\Big\vert_+ \\
& =: g(\rho_+, \bu_+, \bu_-, h),
\end{split}
\end{equation}
where we have used the second condition in \eqref{eq-1.2}. Notice that $\j$ is given by
\begin{align}
\label{2.24}
\j = \frac{\langle \bu_- \vert_- - \bu_+ \vert_+, \bn + \bN \rangle}{1 / \rho_{* -} - 1 /(\rho_{* +} + \rho_+) \vert_+}.
\end{align}
By abuse of notation, in the following, let $\j$ be the symbol defined by \eqref{2.24}. Obviously, the quantity $\j$ is determined by $\bn$, $\bN$, $\rho_+$, $\bu_+$, and $\bu_-$. We  next consider the term $\langle \Delta_{\Gamma_t} (x + h \bn), \bn \rangle$. By \eqref{2.18}, we see that
\begin{align*}
\langle \Delta_{\Gamma_t} (x + h \bn), \bn \rangle & = \bigg\langle \bG_1 (\bm) \sum_{i,j = 1}^{N - 1} \frac{\pd}{\pd p_i} \bigg((\sqrt{\det G} g^{ij} + \bG_{2ij} (\bm)) \frac{\pd }{\pd p_j} (x + h \bn) \bigg), \bn \bigg\rangle \\
& \quad + \bigg\langle \frac{1}{\sqrt{\det G}} \sum_{i,j = 1}^{N - 1} \frac{\pd}{\pd p_i} \bigg(\bG_{2ij} (\bm) \frac{\pd }{\pd p_j} (x + h \bn) \bigg), \bn \bigg\rangle + \langle\Delta_\Gamma (x + h \bn), \bn \rangle.
\end{align*}
From $\langle \pd x / \pd p_j, \bn \rangle = 0$ $(j = 1, \dots, N - 1)$, we obtain
\begin{gather*}
\bigg\langle \bG_1 (\bm) \sum_{i,j = 1}^{N - 1} \frac{\pd}{\pd p_i} \bigg((\sqrt{\det G}  g^{ij} + \bG_{2ij} (\bm)) \frac{\pd x}{\pd p_j} \bigg), \bn \bigg\rangle = \bG_1 (\bm) \sum_{i,j = 1}^{N - 1} (\sqrt{\det G} g^{ij} + \bG_{2ij} (\bm)) \bigg\langle \frac{\pd^2 x}{\pd p_i \pd p_j}, \bn \bigg\rangle, \\
\bigg\langle \frac{1}{\sqrt{\det G}} \sum_{i,j = 1}^{N - 1} \frac{\pd}{\pd p_i} \bigg(\bG_{2ij} (\bm) \frac{\pd x}{\pd p_j} \bigg), \bn \bigg\rangle = \frac{1}{\sqrt{\det G}} \sum_{i,j = 1}^{N - 1} \bG_{2ij} (\bm) \bigg\langle \frac{\pd^2 x}{\pd p_i \pd p_j}, \bn \bigg\rangle,
\end{gather*}
while by $\langle \pd \bn / \pd p_j, \bn \rangle = 0$ $(j = 1, \dots, N - 1)$, we observe that
\begin{align*}
& \bigg\langle \bG_1 (\bm) \sum_{i,j = 1}^{N - 1} \frac{\pd}{\pd p_i} \bigg((\sqrt{\det G}g^{ij} + \bG_{2ij} (\bm)) \frac{\pd}{\pd p_j} (h \bn) \bigg), \bn \bigg\rangle \\
& \quad = \bG_1 (\bm) \sum_{i,j = 1}^{N - 1} \bigg\{\frac{\pd}{\pd p_i} \bigg((\sqrt{\det G}g^{ij} + \bG_{2ij} (\bm)) \frac{\pd h}{\pd p_j} \bigg) + (\sqrt{\det G} g^{ij} + \bG_{2ij} (\bm)) \bigg\langle \frac{\pd^2 \bn}{\pd p_i \pd p_i}, \bn \bigg\rangle\bigg\}, \\
& \bigg\langle \frac{1}{\sqrt{\det G}} \sum_{i,j = 1}^{N - 1} \frac{\pd}{\pd p_i} \bigg(\bG_{2ij} (\bm) \frac{\pd}{\pd p_j} (h \bn) \bigg), \bn \bigg\rangle \\
& \quad = \frac{1}{\sqrt{\det G}} \sum_{i,j = 1}^{N - 1} \bigg\{\frac{\pd}{\pd p_i}  \bigg(\bG_{2ij} (\bm) \frac{\pd h}{\pd p_j} \bigg) + \bG_{2ij} (\bm) h \bigg\langle\frac{\pd^2 \bn}{\pd p_i \pd p_j}, \bn \bigg\rangle \bigg\}.
\end{align*}
In addition, using \eqref{2.15} we see that
\begin{align*}
\langle \Delta_\Gamma (h \bn), \bn \rangle = \bigg\langle \frac{1}{\sqrt{\det G}} \sum_{i,j = 1}^{N - 1} \frac{1}{\pd p_i} \bigg\{\sqrt{\det G} g^{ij} \bigg(\frac{\pd h}{\pd p_j} \bn + h\frac{\pd \bn}{\pd p_j} \bigg)\bigg\}, \bn \bigg\rangle = (\langle \Delta_\Gamma \bn, \bn\rangle + \Delta_\Gamma) h,
\end{align*}
where we have used the fact that $\langle\pd \bn/\pd p_i,\bn\rangle=0$
for $i=1,\dots,N-1$. Recalling that $\Delta_\Gamma x = H_\Gamma \bn$, we obtain
\begin{align*}
\langle\Delta_\Gamma (h \bn), \bn \rangle = H_\Gamma + (\langle \Delta_\Gamma \bn, \bn\rangle + \Delta_\Gamma) h.
\end{align*}
Hence, from \eqref{2.19*}, we can rewrite \eqref{2.22} as
\begin{equation}
\label{l-7}
\begin{split}
\langle \bT_- (\gamma_4, \bu_-, \pi_-) \bn, \bn \rangle \vert_- - \langle \bT_+ (\gamma_1, \gamma_2, \gamma_3, \bu_+, \rho_+) \bn, \bn \rangle \vert_+ \\
- \sigma (\langle \Delta_\Gamma \bn, \bn \rangle + \Delta_\Gamma) h & = f^+_B (\rho_+, \bu_+, \bu_-, h)
\end{split}
\end{equation}
on $\Gamma \times (0, T)$, where we have set
{\allowdisplaybreaks \begin{align*}
&f^+_B (\rho_+, \bu_+, \bu_-, h) \\
& = \left(\frac{\langle \bn, \bN(\bm)\rangle}{1 + \langle \bn, \bN(\bm) \rangle} - 1\right) \Big\lbrac \langle - \j \bu, \bn \rangle + \mu \CD_\bD \nabla_x \bu \Big\rbrac \\
&  \quad
\color{red}
+ \frac{\langle \bn, \bN(\bm)\rangle}{1 + \langle \bn, \bN(\bm) \rangle} \bigg[\Big\lbrac 2 \mu \bD_x (\bu) \Big\rbrac - \bigg\{(\nu_+ - \mu_+) \Big(\bM_0 (\bm) \colon \nabla_x \bu_+\Big) - \gamma_{* +} \rho_+ \bigg\} \bigg\vert_+ \bigg] \\
& \quad 
\color{red}
- \left(\frac{\langle \bn, \bN(\bm)\rangle}{1 + \langle \bn, \bN(\bm) \rangle} - 1 \right) \bigg \{- \rho_+^2 \int_0^1 (1 - \theta) \fp''_+ (\rho_{* +} + \theta \rho_+) \dtheta \\
& \quad \quad 
\color{red}
+ \frac{\kappa_+}{2} \lvert (\bI + \bM_0 (\bm)) \nabla_x \rho_+ \rvert^2 + \kappa_+ \rho_+ \Delta \rho_+ \bigg\} \bigg\vert_+ \\
& \quad + \left(\frac{\langle \bn, \bN(\bm)\rangle}{1 + \langle \bn, \bN(\bm) \rangle} - 1\right) \Big(\kappa_+ \big((\bI + \bM_0 (\bm)) \nabla_x \rho_+\big) \otimes \big((\bI + \bM_0 (\bm)) \nabla_x \rho_+\big)\Big) \bigg\vert_+ \\
& \quad - \left(\frac{\langle \bn, \bN(\bm)\rangle}{1 + \langle \bn, \bN(\bm) \rangle} - 1\right) \Big(\bG_{N, 1} (\bm) + \bG_{N, 2} (\bm) \Big) + \frac{\langle \bn, \bN(\bm)\rangle}{1 + \langle \bn, \bN(\bm) \rangle} \Big(H_\Gamma + (\langle \Delta_\Gamma \bn, \bn\rangle + \Delta_\Gamma) h \Big),
\end{align*}}\noindent
where we have set
\begin{align*}
\bG_{N, 1} (\bm) & = \bG_1 (\bm) \sum_{i,j = 1}^{N - 1} \bigg\{\frac{\pd}{\pd p_i} \bigg((\sqrt{\det G}g^{ij} + \bG_{2ij} (\bm)) \frac{\pd h}{\pd p_j} \bigg) + (\sqrt{\det G} g^{ij} + \bG_{2ij} (\bm)) \bigg\langle \frac{\pd^2 (x + \bn)}{\pd p_i \pd p_i}, \bn \bigg\rangle\bigg\}, \\
\bG_{N, 2} (\bm) & = \frac{1}{\sqrt{\det G}} \sum_{i,j = 1}^{N - 1} \bigg\{\frac{\pd}{\pd p_i}  \bigg(\bG_{2ij} (\bm) \frac{\pd h}{\pd p_j} \bigg) + \bG_{2ij} (\bm) \bigg(h \bigg\langle\frac{\pd^2 \bn}{\pd p_i \pd p_j}, \bn \bigg\rangle + \bigg\langle \frac{\pd^2 x}{\pd p_i \pd p_j}, \bn \bigg\rangle \bigg).
\end{align*}
To obtain the representation formula above, we have used the assumption $\pi_{* -} -\pi_{* +} = \sigma H_\Gamma$. Notice that $f^+_B$ is independent of $\pi_-$. \par
The second jump condition in (\ref{eq-1.2}) is equivalent to the conditions
\begin{align}
\label{2.25}
\lbrac \Pi_{\bn_t} \bu\rbrac = 0, \quad \Pi_{\bn_t} (\bu_\pm - \bu_\Gamma) = 0, \quad \lbrac \rho (\bu - \bu_\Gamma) \cdot \bn_t \rbrac = 0.
\end{align}
Since the second and third conditions in \eqref{2.25} have already used for deriving the equations (cf. Watanabe~\cite{Wat2017}), we need the first condition in \eqref{2.25} to derive the linearized problem. Applying the transform~\eqref{2.1}, the rest jump condition in \eqref{eq-1.2} take the following form:
\begin{align}
\label{l-8}
\left\{\begin{aligned}
\frac{1}{\rho_{*-}}\langle \bT_- (\gamma_4, \bu_-, \pi_-)\bn, \bn \rangle \vert_-
- \Big\{\frac{1}{\rho_{*+}} \langle \bT_+ (\gamma_1, \gamma_2, \gamma_3, \bu_+, \rho_+) \bn, \bn \rangle +\gamma_{**}^+ \rho_+ \Big\} \Big\vert_+ & = f^-_B (\rho_+, \bu_+, \bu_-, h), \\
\Pi_\bn \bu_- \vert_- - \Pi_\bn \bu_+ \vert_+ & = \bh (\bu_+, \bu_-, h), \\
\langle \nabla \rho_+, \bn\rangle \vert_+ & = k_- (\rho_+, h)
\end{aligned}\right.
\end{align}
on $\Gamma \times (0, T)$, where we have set
\begin{align*}
& f^-_B (\rho_+, \bu_+, \bu_-, h) \\
& = \frac{\j^2}{2} \bigg(\frac{1}{\rho_{* -}} - \frac{1}{\rho_{* +} + \rho_+} \bigg\vert_- \bigg) + z \int_0^1 \bigg(\frac{\pd^2 \psi_+}{\pd \rho \pd z}(\rho_{* +}, \theta z) + \frac{\pd \psi_+}{\pd z}(\rho_{* +}, \theta z) \bigg) \dtheta \\
& \quad - \frac{1}{\rho_{* -}} \bigg\{\Big\langle \bD_x (\bu_-) \bN (\bm), (\bn + \bN (\bm))\Big\rangle + \Big\langle \bD_x (\bu_-) \bn, \bN (\bm)\Big\rangle \bigg\} \bigg\vert_- \\
& \quad - \frac{1}{\rho_{* -}} \Big\langle \big(\mu \CD_\bD \nabla_x \bu_- \big) (\bn + \bN (\bm)), (\bn + \bN (\bm)) \Big\rangle \Big\vert_- - \frac{\rho_+}{\rho_{* +} (\rho_{* +} + \rho_+)} \langle \bT_+ (\gamma_1, \gamma_2, \gamma_3, \bu_+, \rho_+) \bn, \bn \rangle \Big\vert_+ \\
& \quad + \frac{\mu_+}{\rho_{* +} + \rho_+} \Big(\Big\langle \bD_x (\bu_+) \bN (\bm), (\bn + \bN (\bm)) \Big\rangle + \Big\langle \bD_x (\bu_+) \bn, \bN (\bm) \Big\rangle \Big) \\
& \quad + \frac{1}{\rho_{* +} + \rho_+} \Big\{\Big\langle \big(\mu \CD_\bD \nabla_x \bu\big) (\bn + \bN (\bm)), (\bn + \bN (\bm)) \Big\rangle + (\nu_+ - \mu_+) \big(\bM_0 (\bm) \colon \nabla_x \bu_+ \big) \\
& \qquad + \frac{\kappa_+}{2} \lvert (\bI + \bM_0 (\bm)) \nabla_x \rho_+ \rvert^2 + \rho_{* +} \kappa_+ \Big(\dv_x \big(\bM_0 (\bm) \nabla_x \bu_+ \big) + \bM_0 (\bm) \colon \nabla_x \big((\bI + \bM_0 (\bm)) \nabla_x \bu_+ \big) \Big) \\
& \qquad + \kappa_+ \rho_+ \Big(\dv_x \big((\bI + \bM_0 (\bm)) \nabla_x \bu_+ \big) + \bM_0 (\bm) \colon \nabla_x \big((\bI + \bM_0 (\bm)) \nabla_x \bu_+ \big) \Big) \Big\} \Big\vert_+, \\
& \bh(\bu_+, \bu_-, h) = \Big\lbrac \langle \bu, \bn \rangle \bN (\bm) + \langle \bu, \bN (\bm) \rangle (\bn + \bN (\bm)) \Big\rbrac, \\
& k_- (\rho_+, h) = - \Big\{\Big\langle \bM_0 (\bm) \nabla_x \rho_+, \bn \Big\rangle + \Big\langle (\bI +\bM_0 (\bm)) \nabla_x \rho_+, \bN (\bm) \Big\rangle \Big\}\Big\vert_+.
\end{align*}
with $z = \lvert (\bI + \bM_0 (\bm)) \nabla_x \rho_+ \rvert^2$.	
Here, we have used the Gibbs-Thomson condition \eqref{cond-GT} and the Taylor formula:
\begin{align*}
\psi_+ (\rho_{* +} + \rho_+, z)
& =\psi_+ (\rho_{* +}, 0) + \rho_+ \bigg(\frac{\pd \psi_+}{\pd \varrho_+} (\rho_{* +}, 0) + z\int_0^1 \frac{\pd^2 \psi_+}{\pd \varrho \pd z}(\rho_{* +}, \theta z) \dtheta \bigg) + z \int_0^1 \frac{\pd \psi_+}{\pd z} (\rho_{* +}, \theta z) \dtheta \\
& =: \psi_+ (\rho_{* +}, 0) + \gamma_{**}^+ \rho_+ + z \int_0^1 \bigg(\frac{\pd^2 \psi_+}{\pd \varrho \pd z} (\rho_{* +}, \theta z) + \frac{\pd \psi_+}{\pd z} (\rho_{* +}, \theta z) \bigg)\dtheta.
\end{align*}
Recall that the Helmholtz free energy $\psi_+$ depends on not only the density $\varrho_+$ but also the square of the gradient of density $\lvert \nabla \varrho_+ \rvert^2$ if the compressible fluid is dominated by the Navier-Stokes-Korteweg equations. Summing up, from \eqref{l-1}, \eqref{l-2}, \eqref{l-3}, \eqref{l-4}, \eqref{l-5}, \eqref{l-6}, \eqref{l-7}, and \eqref{l-8}, we have derived \eqref{eq-1.5}.

{\color{blue}
\section*{Acknowledgements}
\noindent
The author is grateful to the anonymous referees for their helpful suggestions and comments which lead to the improvement of this paper.
}

\end{document}

%% file: localized.tex
\subsubsection{Unit outer normal and Laplace-Beltrami operator on a bent interface}
Let $\Phi:\BR^N\ni x\mapsto y\in \BR^N$ be a 
bijective map of $C^1$ class.
Furthermore, let $\Phi^{-1}$ be an inverse map of $\Phi$.
We assume that $\nabla_x\Phi$ and 
$\nabla_y\Phi^{-1}$ can be written in the form $\nabla_x\Phi(x)=\bA+\bB(x)$ and $\nabla_y\Phi^{-1}(\Phi(x))=\bA_{-1}+\bB_{-1}(\Phi(x))$, respectively, where $\bA$ and $\bA_{-1}$ are orthogonal matrices
with constant coefficients satisfying $\det\bA=\det\bA_{-1}=1$
and $\bB(x)$ and $\bB_{-1}(\Phi(x))$ are matrices  
of functions in $W^3_r(\BR^N)$, $N < r < \infty$, such that
\begin{align}\label{cond-5.1}
\|(\bB,\bB_{-1})\|_{L_\infty(\BR^N)}\leq M_1,\quad
\|\nabla(\bB,\bB_{-1})\|_{W^2_r(\BR^N)}\leq M_2.
\end{align}
In the following, we write $\bB=\bB(x)$ and
$\bB_{-1}=\bB_{-1}(\Phi(x))$ for short.
Let 
$a_{ij}$ and $b_{ij}$ be the $(i,j)$
elements of $\bA_{-1}$ and $\bB_{-1}$, respectively.
Since we will chose $M_1$ small enough,
we may assume that $0<M_1\leq 1\leq M_2$ beforehand.
Set $\Omega_\pm^\Phi = \Phi (\BR^N_\pm)$ and 
$\Gamma^\Phi =\Phi(\RZ)$ and $\widetilde{\bn}$ be the unit outer normal to $\Gamma^\Phi$.
We see that $\Gamma^\Phi$ is represented by $\Phi_{-1,N}(y)=0$ with
$\Phi^{-1}=(\Phi_{-1,1},\dots,\Phi_{-1,N})$, which furnishes that
\begin{equation}\label{5.2}
\begin{split}
\widetilde{\bn}(y)=\widetilde{\bn}(\Phi(x))
=&-\frac{(\nabla_y\Phi_{-1,N})\circ\Phi(x)}
{|(\nabla_y\Phi_{-1,N})\circ\Phi(x)|}\\
=&-\frac{(a_{N1}+b_{N1},\dots,a_{NN}+b_{NN})}
{(\sum_{i=1}^{N}(a_{Ni}+b_{Ni})^2)^{1/2}} = \frac{{}^\top \!(\bA_{-1}+\bB_{-1})\bn_0}
{|{}^\top \!(\bA_{-1}+\bB_{-1})\bn_0|}
\end{split}
\end{equation}
with $\bn_0=(0,\dots,0,-1)$.
Since $\sum_{i=1}^{N}(a_{Ni}+b_{Ni})^2=1+\sum_{i=1}^{N}(a_{Ni}b_{Ni}+b^2_{Ni})$,
by (\ref{cond-5.1}) we can write
\begin{align}\label{5.3}
\bigg\{\sum_{i=1}^{N}(a_{Ni}+b_{Ni})\bigg\}^{-1/2}=1+b_0
\end{align}
with $b_0\in W^2_r(\BR^N)$ and small positive constant $M_1$ possessing the following estimates:
\begin{align}
\label{cond-5.4}
\|b_0\|_{L_\infty(\BR^N)}\leq C_N M_1,\quad
\|\nabla b_0\|_{W^2_r(\BR^N)}\leq C_{N,r} M_2.
\end{align}
Let $G=(g_{ij}(x))$ be the first fundamental form of $\Gamma^\Phi$ defined by 
\begin{align*}
g_{ij}(x) = g_{ij}(x',0) & = \bigg\langle\frac{\pd}{\pd x_i}\Phi(x',0),
\frac{\pd}{\pd x_j}\Phi(x',0)\bigg\rangle\nonumber\\
&=\sum_{k=1}^{N}(a_{ki}+b_{ki}(x))(a_{kj}+b_{kj}(x))\nonumber\\
&=\delta_{ij}+\sum_{k=1}^{N}(a_{ki} b_{kj}(x)+a_{kj} b_{ki}(x)
+b_{ki}(x) b_{kj}(x)) =:\delta_{ij}+\tg_{ij}(x',0),\nonumber
\end{align*}
where $\delta_{ij}$ is the Kronecker delta symbol.
Choosing $M_1>0$ suitably small, from (\ref{cond-5.1}) we see that
$G$ is symmetric and positive definite, that is,
the determinant $\BG:=\det G$ is positive.
By (\ref{cond-5.1}), we obtain
\begin{align}\label{est-5.6}
\|\tg_{ij}\|_{L_\infty(\BR^N)}\leq C_N M_1,\quad
\|\nabla\tg_{ij}\|_{W^2_r(\BR^N)}\leq C_{N,r} M_2.
\end{align}
Since $\BG>0$, there exists the inverse matrix $G^{-1}=(g^{ij})$ of $G$.
From (\ref{est-5.6}), we have 
$\sqrt{\BG}=1+\widetilde{\BG}$ and $g^{ij}=\delta_{ij}+\tg^{ij}$ with
\begin{align}\label{est-5.7}
\|(\widetilde{\BG},\tg^{ij})\|_{L_\infty(\BR^N)} \leq C_N M_1,\quad
\|\nabla(\widetilde{\BG},\tg^{ij})\|_{W^2_r(\BR^N)} \leq C_{N,r} M_2.
\end{align}
Then the Laplace-Beltrami operator $\Delta_{\Gamma^\Phi}$ on $\Gamma^\Phi$
is defied for scalar fields by $\Delta_{\Gamma^\Phi} f = \dv_{\Gamma^\Phi} (\nabla_{\Gamma^\Phi} f)$, which can be read as
\begin{align*}
\Delta_{\Gamma^\Phi} f (x')
=\frac{1}{\sqrt{\BG(x',0)}}\sum_{i,j=1}^{N-1}
\frac{\pd}{\pd x_j}
\Big(\sqrt{\BG(x',0)}g^{ij}(x',0)\frac{\pd f(x')}{\pd x_j} \Big)
\end{align*}
in local coordinates. Let $\Delta'_x f(x') =\sum_{j = 1}^{N-1} \pd^2 f(x') / \pd x^2_j$ and
\begin{align*}
\CD'_xf(x')
&=\sum_{i,j=1}^{N-1}\tg^{ij}(x',0)
\frac{\pd^2 f(x')}{\pd x_i\pd x_j}
+\sum_{i,j=1}^{N-1}\left(\frac{\pd \tg^{ij}(x',0)}{\pd x_i}
+\frac{g^{ij}(x',0)}{\sqrt{\BG(x',0)}}\frac{\pd \widetilde{\BG}(x',0)}
{\pd x_i} \right)\frac{\pd f(x')}{\pd x_j}\\
&=:\sum_{i,j=1}^{N-1}\tg^{ij}(x',0)
\frac{\pd^2 f(x')}{\pd x_i\pd x_j}
+\sum_{j=1}^{N-1}\tg^j(x',0)\frac{\pd f(x')}{\pd x_j},
\end{align*}
and then we write 
$\Delta_{\Gamma^\Phi} f(x')=\Delta'_x f(x')+\CD'_x f(x')$.
By (\ref{est-5.7}) we have
\begin{align}\label{est-5.8}
\|\tg^j\|_{L_\infty(\BR^N)}\leq C_{N} M_1\quad
\|\tg^j\|_{W^2_r(\BR^N)}\leq C_{N,r} M_2.
\end{align}

\subsubsection{Construction of solution operators}
Let $\tgamma_k(y)$ be coefficients of equations
that is defined on $\BR^N$ satisfying the conditions:
\begin{align}\label{cond-5.9}
\frac{1}{2} \gamma^-_{k*} \leq \tgamma_k(y)
\leq \frac{3}{2} \gamma^+_{k*}, \quad
|\tgamma_k(y)-\gamma^*_{k}|\leq M_1,\quad
\|\nabla \tgamma_{k}\|_{L_r(\BR^N)}\leq C_r
\end{align}
for $y\in\BR^N$, $k=1,\dots,4$, and the same constants
$\gamma^\pm_{k*}$ as in Assumption \ref{as-star},
where $\gamma^*_{k}$ are some constants with 
$\gamma^-_{k*}\leq \gamma^*_{k}\leq \gamma^+_{k*}$.
Set $\gamma_{k}(x)=\tgamma_{k}(\Phi(x))$ for $k=1,\dots,4$.
Then the following conditions hold:
\begin{align}\label{cond-5.10}
\begin{aligned} 
|\gamma_k(x)-\gamma^*_{k}|&\leq M_1,&\quad
\|\nabla \gamma_{k}\|_{L_r(\RP)}&\leq C_r
&\quad(k=1,\dots,4,\enskip x\in\RP),\\
|\gamma_{4}(x)-\gamma^*_{5}|&\leq M_1,&\quad
\|\nabla \gamma_{4}\|_{L_r(\RM)}&\leq C_r
&\quad(x\in\RM).
\end{aligned}
\end{align}
For $\widetilde{\bu}_{-}\in W^2_q(\Omega_-^\Phi)^N$ let
$\CK^0_1(\widetilde{\bu}_-)$ be a unique solution to the variational problem:
\begin{align*}
(\nabla_y \CK^0_1(\widetilde{\bu}_-),\nabla_y\varphi)_{\Omega_-^\Phi}
=(\DV_y(\tgamma_{4}\bD_y(\widetilde{\bu}_{-}))
-\rho_{*-}\nabla_y\dv\widetilde{\bu}_{-},\nabla_y\varphi)_{\Omega_-^\Phi}
\quad\text{ for any $\varphi\in\hW^1_{q',0}(\Omega_-^\Phi)$}
\end{align*}
such that $\CK^0_1(\widetilde{\bu}_-)
=\tgamma_{4}\langle\bD_y(\widetilde{\bu}_{-})\widetilde{\bn},
\widetilde{\bn}\rangle-\rho_{*-}\dv_y\widetilde{\bu}_-$
on $\Gamma^\Phi$, while for $\tilh\in W^{3-1/q}_q(\Gamma^\Phi)$ let $\CK^0_2(\tilh)$
be a unique solution to the variational problem:
\begin{align*}
(\nabla_y\CK^0_2(\tilh),\nabla_y\varphi)_{\Omega_-^\Phi}=0
\quad\text{ for any $\varphi\in \hW_{q',0}(\Omega_-^\Phi)$}
\end{align*}
such that $\CK^0_2(\tilh)
= -(\rho_{*-}-\rho_{*+})^{-1}\rho_{*-}\sigma\Delta_{\Gamma^\Phi} \tilh$ on $\Gamma^\Phi$.
Under these situations, let us consider
\begin{align}
\label{eq-5.11}
\left\{\begin{aligned}
\lambda \wt \rho_+ + \rho_{*+} \dv \wt \bu_+ & = \wt f_M &\enskip &\text{ in $\Omega_+^\Phi$,} \\
\rho_{*+} \lambda \wt \bu_+ - \DV \bT_{0 +} (\wt \gamma_1, \wt \gamma_2, \wt \gamma_3, \wt \bu_+, \wt \rho_+) & = \wt \bff_+ &\enskip &\text{ in $\Omega_+^\Phi$}, \\
\rho_{*-} \lambda \wt \bu_- - \DV \bT_- (\wt \gamma_4, \wt \bu_-, \CK_1^0 (\wt \bu_-) + \CK_2^0 (\wt h)) & = \wt \bff_- &\enskip &\text{ in $\Omega_-^\Phi$}, \\
\lambda \wt h - \frac{\langle \rho_{* -} \wt \bu_-, \bn \rangle \vert_- - \langle \rho_{* +} \wt \bu_+, \bn \rangle \vert_+} {\rho_{* -} - \rho_{* +}} & = \wt d &\enskip &\text{ on $\Gamma^\Phi$}, \\
\bB_0 (\wt \gamma_1, \wt \gamma_2, \wt \gamma_3, \wt \gamma_4, \wt \rho_+, \wt \bu_+, \wt \bu_-) & = \wt \bG_0 &\enskip &\text{ on $\Gamma^\Phi$}, \\
\rho_{*-} (\dv \wt \bu_-) \vert_- & = \wt g_- &\enskip &\text{ on $\Gamma^\Phi$},
\end{aligned}\right.
\end{align}
where $\wt \bG_0 = {}^\top\! (\wt g, \wt g_+, \wt \bh, \wt k_-)$. The main purpose of this subsection is to prove the following theorem.
\begin{theo}\label{th-5.1}
Let $1 < q < \infty$ and $\rho_{* +} \neq \rho_{* -}$. Assume that $\rho_{* +}$, $\tgamma_1$, $\tgamma_2$, and $\tgamma_3$ satisfy Assumption~\ref{as-star} (c) with $\gamma_k = \tgamma_k$ $(k = 1, 2, 3)$. Then there exists constant $\varepsilon_* \in (0, \pi/2)$ such that for any $\varepsilon\in (\varepsilon_*, \pi/2)$ there exists a constant $\lambda_4 > 0$ with the following assertions valid:
\begin{enumerate}\renewcommand{\labelenumi}{(\arabic{enumi})}
\item For any $\lambda\in \Sigma_{\varepsilon, \lambda_4}$, there exists operators
\begin{align*}
\CA^+_F (\lambda) & \in \Hol (\Sigma_{\varepsilon, \lambda_4}, \CL (\CZ_q (\Omega_+^\Phi, \Omega^\Phi_-, \Gamma^\Phi), W^3_q (\Omega_+^\Phi))), \\
\CB^\pm_F (\lambda) & \in \Hol (\Sigma_{\varepsilon, \lambda_4}, \CL (\CZ_q (\Omega_+^\Phi, \Omega^\Phi_-, \Gamma^\Phi), W^2_q (\Omega_\pm^\Phi)^N)), \\
\CH_F (\lambda) & \in \Hol (\Sigma_{\varepsilon, \lambda_4}, \CL (\CZ_q (\Omega_+^\Phi, \Omega^\Phi_-, \Gamma^\Phi), W^{3 - 1/q}_q (\Gamma^\Phi)))
\end{align*}	
such that for any $\widetilde \bF_Z := (\tf_M, \widetilde \bff_+, \widetilde \bff_-, \td, \tg, \tg_+, \tg_-, \widetilde \bh, \tk_-) \in Z_q (\Omega_+^\Phi, \Omega^\Phi_-, \Gamma^\Phi)$, the problem \eqref{eq-4.5} admits a unique solution $(\trho_+, \widetilde \bu_+, \widetilde \bu_-, \widetilde h)$ defined by $\wt \rho_+ = \CA^+_F (\lambda) G_\lambda (\wt \bF_Z)$, $\wt \bu_\pm = \CB^\pm_F (\lambda) G_\lambda (\wt \bF_Z)$, and $\wt h = \CH_F (\lambda) G_\lambda (\wt \bF_Z)$, where $G_\lambda$ is the operator defined in Theorem~\ref{th-4.4}.
\item For $s = 0, 1$, $i = 0, 1, 2, 3$, $j = 0, 1, 2$, and $k = 0, 1$, there exists a positive constant $c_4$ such that
\begin{equation*}
\begin{split}
\CR_{\CL (\CZ_q (\Omega_+^\Phi, \Omega^\Phi_-, \Gamma^\Phi), W^{3 - i}_q (\Omega_+^\Phi))} (\{(\tau \pd_\tau)^s (\lambda^{i/2} \CA^+_F (\lambda)) \mid \lambda \in \Sigma_{\varepsilon, \lambda_4} \}) & \leq c_4, \\
\CR_{\CL (\CZ_q (\Omega_+^\Phi, \Omega^\Phi_-, \Gamma^\Phi), W^{2 - j}_q (\Omega_\pm^\Phi)^N)} (\{(\tau \pd_\tau)^s (\lambda^{j/2} \CB^\pm_F (\lambda)) \mid \lambda \in \Sigma_{\varepsilon, \lambda_4}\}) & \leq c_4, \\
\CR_{\CL (\CZ_q (\Omega_+^\Phi, \Omega^\Phi_-, \Gamma^\Phi), W^{3 - 1/q - k}_q (\Gamma^\Phi))} (\{(\tau \pd_\tau)^s (\lambda^k \CH_F (\lambda)) \mid \lambda \in \Sigma_{\varepsilon, \lambda_4}\}) & \leq c_4.
\end{split}
\end{equation*}
Here, the constant $c_4$ is independent of $\lambda$.	
\end{enumerate}
\end{theo}
In the following, we prove Theorem \ref{th-5.1}.
By change of variable: $y=\Phi(x)$, we have
\begin{align*}
\frac{\pd}{\pd y_j}=\sum_{k=1}^{N}(a_{kj}+b_{kj}(x))\frac{\pd}{\pd x_j}.
\end{align*}
Accordingly, the variational equation:
\begin{align*}
(\nabla_y \tu,\nabla_y\tvarphi)_{\Omega_-^\Phi}=(\widetilde{\bff},
\nabla_y\tvarphi)_{\Omega_-^\Phi}
\quad\text{ for any $\tvarphi\in\hW^1_{q',0}(\Omega_-^\Phi)$}
\end{align*}
subject to $\tu=\tf$ on $\Gamma^\Phi$ is transformed to the variational equation:
\begin{align}\label{eq-5.13}
(\nabla_x u,\nabla_x\varphi)_{\RM}+(\bG^0\nabla_x u,\nabla_x\varphi)_{\RM}
=(\bg,\nabla_x\varphi)_{\RM} \quad\text{ for any $\varphi\in W^1_{q',0} (\BR^N_-)$}
\end{align}
subject to $u=f$ on $\RZ$, where we have set
$\bg=\BG(\bA_{-1}+\bB_{-1})\bff$, $\bff=\widetilde{\bff}\circ\Phi$,
$f=\tf\circ\Phi$, $u=\tu\circ\Phi$, $\varphi=\tvarphi\circ\Phi$,
and $\bG^0$ is an $N\times N$ matrix 
whose $(k,l)$ component $G^0_{kl}$ is given by
\begin{align}\label{5.14}
G^0_{kl}=\delta_{kl}(\BG-1)
+\sum_{j=1}^{N}(a_{kj}b_{lj}+a_{lj}b_{kj}+b_{kj}b_{lj})\BG.
\end{align}
Since $N<r<\infty$, by the Sobolev embedding theorem, we have
$\|\nabla b_{ij}\|_{L_\infty(\BR^N)}\leq \|\nabla b_{ij}\|_{W^1_r(\BR^N)}$.
Hence, by (\ref{cond-5.1}) we obtain
\begin{equation}\label{est-5.15}
\begin{split} 
\|(\BG-1,G^0_{kl})\|_{L_\infty(\BR^N)}\leq C_N M_1, \quad \|(\nabla\BG,\nabla G^0_{kl})\|_{W^2_r(\BR^N)}\leq C_{N,r} M_2.
\end{split} 
\end{equation}
Choosing $M_1>0$ small enough and using the Banach 
fixed point theorem, we can prove the following lemma.
\begin{lemm}\label{lem-5.2}
Let $1<q<\infty$. Then  there exists an operator $K_2$ with
\begin{align*}
K_2 \in \CL(L_q(\RM)^N\times W^{1-1/q}_q(\RM),W^1_q(\RM)+\hW^1_{q,0}(\RM) )
\end{align*}
such that for any $\bg\in L_q(\RM)^N$ and
$f\in W^{1-1/q}_q(\BR^{N-1})$, $u=K_2(\bg,f)$ is
a unique solution to the variational problem (\ref{eq-5.13})
possessing the estimate:
\begin{align*}
\|\nabla u\|_{L_q(\RM)}\leq
C_{N,q}\Big(\|\bg\|_{L_q(\RM)}+\|f\|_{W^{1-1/q}_q(\BR^{N-1})}\Big).
\end{align*}		
\end{lemm}	

Next, we transform (\ref{eq-5.11}) to a problem 
with the flat interface $\RZ$ by change of variable: $y=\Phi(x)$.
To this end, we set
\begin{alignat*}6
\gamma_k(x)&:=\tgamma_k(\Phi(x))&\quad (k=&1,\dots,4),
&\quad \rho_{+}(x)&:=\trho_{+}(\Phi(x)), &\quad \bu_\pm(x)&:=\widetilde{\bu}_\pm(\Phi(x)),
\\ 
h(x)&:=\tilh(\Phi(x)), &\quad
f_M (x) &:=\tf_M (\Phi(x)),&\quad
\bff_\pm (x) &:=\widetilde{\bff}_\pm(\Phi(x)),&\quad
d(x) &:= \td(\Phi(x)),\\
g(x) &:=\tg(\Phi(x)), &\quad
f^\pm_B (x)&:=\tf^\pm_B (\Phi(x)), &\quad
\bh(x) &:=\widetilde{\bh}(\Phi(x))&\quad
k_-(x)&:=\tk_-(\Phi(x)).
\end{alignat*}
From the change of variable: $y=\Phi(x)$, we see that
\begin{align*}
\frac{\pd}{\pd y_j} = \sum_{k=1}^{N}(a_{kj}+b_{kj})\frac{\pd}{\pd x_k}, \quad \nabla_y={}^\top \!(\bA_{-1}+\bB_{-1})\nabla_x, \quad \frac{\pd^2}{\pd y_j \pd y_k} = \sum_{l,m=1}^{N}a_{lj}a_{mk} \frac{\pd^2}{\pd x_l \pd x_m}+\CD_{jk},
\end{align*}
where $\CD_{jk}$ are second order differential operators defined by
\begin{align*}
\CD_{jk}=\sum_{l,m=1}^{N}
(a_{lj}b_{mk} +a_{mk}b_{lj}+b_{lj}b_{mk})
\frac{\pd^2}{\pd x_l\pd x_m}+\sum_{l,m=1}^N(a_{lj}+b_{lj})
\Big(\frac{\pd}{\pd x_l}b_{mk} \Big) \frac{\pd}{\pd x_m}.
\end{align*}
In addition, the following conditions hold:
\begin{align*}
\Delta_y & = \Delta_x+\sum_{j=1}^{N}\CD_{jj},\quad
\dv_y\widetilde{\bu}_\pm
=\dv_x(\bA_{-1}\bu_\pm)+\bB_{-1}:\nabla_x\bu_\pm,\\
\bD_y(\widetilde{\bu}_\pm) & = (\nabla_x\bu_\pm)(\bA_{-1}+\bB_{-1})+{}^\top \!(\bA_{-1}+\bB_{-1})
{}^\top \!(\nabla_x \bu_\pm),\\
\nabla_y \dv_y\widetilde{\bu}_\pm & = {}^\top \!(\bA_{-1}+\bB_{-1})\nabla_x\dv_x(\bA_{-1}\bu_\pm) +{}^\top \!(\bA_{-1}+\bB_{-1}) \sum_{j,k=1}^N \nabla_x\Big(b_{kj}\frac{\pd}{\pd x_k} u_{j\pm} \Big).
\end{align*}
\par
\textit{\underline{The first line of (\ref{eq-5.11}).}}
We easily see that
\begin{align*}
\lambda \rho_{+}+\rho_{*+}\dv_x(\bA_{-1}\bu_+)+\CR^1(\bu_+)
=f_M \quad\text{in $\RP$,}
\end{align*}
where $\CR^1(\bu_+)$ is a linear operator
such that
$\CR^1(\bu_+)=\rho_{*+} \bB_{-1}:\nabla_x\bu_+$.
\par
\textit{\underline{The second line of (\ref{eq-5.11}).}}
In the following, we write
${}^\top \!\bA_{-1}={}^\top \!(\bA_{-1})$ and
${}^\top \!\bB_{-1}={}^\top \!(\bB_{-1})$ for short.
Since the equation can be written as
\begin{align*}
\rho_{*+}\lambda\widetilde{\bu}_+
-\tgamma_{1}\Delta_y\widetilde{\bu}_+
-\tgamma_{2}\nabla_y\dv_y
\widetilde{\bu}_+
-\rho_{*+}\tgamma_{3}
\nabla_y\Delta_y\trho_+
-\bD_y(\widetilde{{\bu}}_+)
\nabla_y\tgamma_{1}&\\
-(\dv_y\widetilde{\bu}_+)
\nabla_y(\tgamma_{2}-\tgamma_{1})
-(\Delta_y\trho_+)\nabla_y
(\rho_{*+}\tgamma_{3})
&=\widetilde{\bff}_+ \quad\text{ in $\Omega_+$,}
\end{align*}
we have
\begin{align*}
\rho_{*+}\lambda(\bA_{-1}\bu_+)
-\gamma_{1}\Delta_x (\bA_{-1}\bu_+)
-\gamma_{2}\nabla_x\dv_x (\bA^{-1}\bu_+)&\\
-\rho_{*+}\gamma_{3}
\nabla_x\Delta_x\rho_+
+\bA_{-1}\CR^2(\bu_+,\rho_{+})
&=\bA_{-1}\bff_+ \quad 
\text{in $\RP$},\nonumber
\end{align*}
with
\begin{align*}
\CR^2(\bu_+,\rho_{+}) & = -\gamma_{1}\sum_{j=1}^N\CD_{jj}\bu_+-\gamma_{2}{}^\top \!\bA_{-1}
\nabla_x(\bB_{-1}:\nabla_x\bu_+)\\
&\quad-\gamma_{2}{}^\top \!\bB_{-1}\nabla_x
\{\dv_x(\bA_{-1}\bu_+)+\bB_{-1}:\nabla_x\bu_+\}
-\rho_{*+}\gamma_{3}\Delta_x({}^\top \!\bB_{-1}\nabla_x\rho_+)\\
&\quad-\rho_{*+}\gamma_{3}\sum_{j=1}^{N}\CD_{jj}
\{{}^\top \!(\bA_{-1}+\bB_{-1})\nabla_x\rho_{+}\}\\
&\quad-\{(\nabla_x\bu_+)(\bA_{-1}+\bB_{-1})
+{}^\top \!(\bA_{-1}+\bB_{-1})
{}^\top \!(\nabla_x\bu_+) \}
{}^\top \!(\bA_{-1}+\bB_{-1})\nabla_x\gamma_{1}\\
&\quad-\{\dv_x(\bA_{-1}\bu_+)+\bB_{-1}:\nabla_x\bu_+\}
{}^\top \!(\bA_{-1}+\bB_{-1})\nabla_x(\gamma_{2}-\gamma_{1})\\
&\quad-\Big(\Delta_x\rho_{+}+\sum_{j=1}^{N}\CD_{jj}
\rho_{+} \Big)(\bA_{-1}+\bB_{-1})\nabla_x(\rho_{*+}\gamma_{3}).
\end{align*}
\par
\textit{\underline{The third line of (\ref{eq-5.11}).}}
For a $N\times N$ matrix $\bL$,
let ${}^{-\top} \!\bL=({}^\top \!\bL)^{-1}={}^\top \!(\bL^{-1})$.
Since $(\bA_{-1})^{-1}={}^{-\top} \!\bA$, we see that ${}^{-\top} \!(\bA_{-1} + \bB_{-1})
=\bA_{-1}{}^{-\top} \!(\bI+{}^{\top} \!\bA_{-1}\bB_{-1})$. Combining this relation and (\ref{5.2}), we have
\begin{align*}
\bD_y(\widetilde{\bu}_\pm)\widetilde{\bn} & = (1+b_0){}^{\top} \!(\bA_{-1} + \bB_{-1})
\Big\{\bA_{-1}({}^{-\top} \!(\bI+{}^{\top} \!\bA_{-1} \bB_{-1})-\bI)
(\nabla_x\bu_\pm)\bn_0 \\
&\quad+\bA_{-1}{}^{-\top} \!(\bI+{}^{\top} \!\bA_{-1} \bB_{-1})
(\nabla_x\bu_\pm)(\bA_{-1} {}^{\top} \!\bB_{-1}
+\bB_{-1}{}^{\top} \!\bA_{-1}+\bB_{-1}{}^{\top} \!\bB_{-1})\bn_0
+\bD_x(\bA_{-1}\bu_\pm)\bn_0 \Big\}.
\end{align*}
From (\ref{5.2}) and (\ref{5.3}) we see that
\begin{align}\label{5.16}
\bD_y(\widetilde{\bu}_\pm)\widetilde{\bn}
=(1+b_0)({}^{\top} \!\bA_{-1}+{}^{\top} \!\bB_{-1})
\{\bD_x(\bA_{-1}\bu_\pm)\bn_0+\CR^3(\bu_\pm) \}.
\end{align}
with
\begin{align*}
\CR^3(\bu_\pm)
=&\bA_{-1}({}^{-\top} \!(\bI+{}^{\top} \!\bA_{-1}\bB_{-1})-\bI)
(\nabla_x\bu_\pm)\bn_0+{}^{\top} \!(\nabla_x\bu_\pm)
{}^{\top} \!\bB_{-1}\bn_0 \\&
+\bA_{-1}{}^{-\top} \!(\bI+{}^{\top} \!\bA_{-1} \bB_{-1})
(\nabla_x\bu_\pm)(\bA_{-1} {}^{\top} \!\bB_{-1}
+\bB_{-1}{}^{\top} \!\bA_{-1}+\bB_{-1}{}^{\top} \!\bB_{-1})\bn_0.
\end{align*}
Hence, we can write 
$\langle\bD_y(\widetilde{\bu}_\pm)\widetilde{\bn},\widetilde{\bn}\rangle$ as
\begin{align}\label{5.17}
\langle\bD_y(\widetilde{\bu}_\pm)\widetilde{\bn},\widetilde{\bn}\rangle
=\langle\bD_x(\bA_{-1}\bu_\pm)\bn_0,\bn_0\rangle+\CR^4(\bu_\pm),
\end{align}
with
\begin{align*}
\CR^4(\bu_\pm)
=&(2b_0+b_0^2)\langle\bD_x(\bA_{-1}\bu_\pm)\bn_0,\bn_0 \rangle
+(1+b_0)^2\langle{}^\top\!\bA_{-1}\bD_x(\bA_{-1}\bu_\pm) \bn_0,
{}^\top\!\bB_{-1}\bn_0\rangle\\
&+(1+b_0)^2\langle{}^\top\!\bB_{-1}(\bD_x(\bA_{-1}\bu_\pm) + \CR^3(\bu_\pm)),
({}^\top\!\bA_{-1}+{}^\top\!\bB_{-1})\bn_0\rangle + (1+b_0)^2 \langle \CR^3(\bu_\pm), \bn_0\rangle,
\end{align*}
because $\bA_{-1}$ is the orthogonal matrix.
Define 
$\CR^5(\bu_-)=\CR^4(\bu_-)-\bB_{-1}:\nabla_x\bu_-$, then we have
\begin{align*}
\langle\bD_y(\widetilde{\bu}_-)\widetilde{\bn},\widetilde{\bn}\rangle
-\rho_{*-}\dv_y\widetilde{\bu}_-
&=\langle\gamma_{4}\bD_x(\bA_{-1}\bu_-)\bn_0,\bn_0\rangle
-\dv_x(\bA_{-1}\bu)+\CR^4(\bu_-).
\end{align*}
In addition, we see that
\begin{align}\label{5.18}
\DV_y(\tgamma_2 \bD_y(\widetilde{\bu}_-))-\rho_{*-}\nabla_y\dv_y\widetilde{\bu}_- = \gamma_{4}\Delta_x\bu_--\rho_{*-}{}^{\top} \!\bA_{-1}
\nabla_x\dv_x (\bA_{-1}\bu_-) +\CR^6(\bu_-),
\end{align}
where we have set
\begin{align*}
\CR^6(\bu_-)
=&\gamma_{4}\sum_{j=1}^{N}\CD_{jj}\bu_-+\{(\nabla_x\bu_-)(\bA_{-1}+\bB_{-1})
+{}^{\top} \!(\bA_{-1}+\bB_{-1}){}^{\top} \!(\nabla_x\bu_-)\}
{}^{\top} \!(\bA_{-1}+\bB_{-1})\nabla_x\gamma_{4}\\
&-\rho_{*-}{}^{\top} \!\bA_{-1}\nabla_x
(\bB_{-1}:\nabla_x\bu_-) -\rho_{*-}{}^{\top} \!\bB_{-1}\nabla_x
\{\dv_x(\bA_{-1}\bu_-)+\bB_{-1}:\nabla_x\bu_- \}.
\end{align*}
Then by  (\ref{5.18}) we have
\begin{align}\label{5.19}
&\BG(\bA_{-1}+\bB_{-1})
\{\DV_y(\tgamma_2\bD_y(\widetilde{\bu}_-))
-\rho_{*-}\nabla_y\dv_y\widetilde{\bu}_-\}\\
&=\gamma_{4}\Delta_x(\bA_{-1}\bu_-)
-\rho_{*-}\nabla_x\dv_x(\bA_{-1}{\bu}_-)
+\CR^7(\bu_-)\nonumber
\end{align}
with the linear operator $\CR^7(\bu_-)$:
\begin{align*}
\CR^{7}(\bu_-) = \{(\BG-1)\bA_{-1}+\BG\bB_{-1}\}(\gamma_{4}\Delta_x\bu_- -\rho_{*-}{}^{\top}\!\bA_{-1}\nabla_x\dv_x(\bA_{-1}\bu_-)) + \BG(\bA_{-1}+\bB_{-1})\CR^6(\bu_-)
\end{align*}
Let $p_1=\CK^0_1(\widetilde{\bu}_-)\circ\Phi$ and
$p_2=\CK^0_2(\tilh)\circ\Phi$.
From (\ref{eq-5.13}), (\ref{5.18}), and (\ref{5.19}),
$p_1$ and $p_2$ satisfy the following variational equations:
\begin{align*}
(\nabla p_1,\nabla\varphi)_{\RM} + (\bG^0\nabla p_1,\nabla\varphi)_{\RM} = (\gamma_{4}\Delta_x(\bA_{-1}\bu_-) -\rho_{*-}\nabla_x\dv_x(\bA_{-1}{\bu}_-)
+\CR^7(\bu_-),\varphi)_{\RM}\nonumber
\end{align*}
for any $\varphi\in W^1_{q',0}(\RM)$ subject to
$p_1
=\langle\gamma_{4}\bD_x(\bA_{-1}\bu_-)\bn_0,\bn_0\rangle
-\rho_{*-}\dv_x(\bA_{-1}\bu)+\CR^4(\bu_-)$ 
on $\RZ$, and 
\begin{align*}
&(\nabla p_2,\nabla\varphi)_{\RM}
+(\bG^0\nabla p_2,\nabla\varphi)_{\RM}=0
\quad\text{for any $\varphi\in W^1_{q',0}(\RM)$}
\end{align*}
subject to 
$p_2
=-(\rho_{*-}-\rho_{*+})^{-1}\rho_{*-}\sigma\Delta'_x h
-(\rho_{*-}-\rho_{*+})^{-1}\rho_{*-}\sigma \CD'_x h$
on $\RZ$, respectively.
Hence, the third line of (\ref{eq-5.11}) is rewritten as follows:
\begin{align*}
\rho_{*-}\lambda\bu_- -\gamma_{4}\Delta_x\bu_- 
-({}^{\top} \!\bA_{-1}+{}^{\top} \!\bB_{-1})\nabla_x(p_1+p_2)+\CR^8(\bu_-)
=\bff_-
\end{align*}
with
\begin{align*} 
\CR^8(\bu_-)
=&-\gamma_{4}\sum_{j=1}^{N}\CD_{jj}\bu_-
-\{(\nabla_x\bu_-)(\bA_{-1}+\bB_{-1})
+{}^{\top} \!(\bA_{-1}+\bB_{-1}){}^{\top} \!(\nabla_x\bu_-) \}
{}^{\top} \!(\bA_{-1}+\bB_{-1})\nabla_x\gamma_{4}.
\end{align*}
From (\ref{eq-4.2}), (\ref{eq-4.3}), and Lemma \ref{lem-5.2}, we obtain
\begin{equation}\label{5.20}
\begin{split} 
p_1=&\CK_{F1}(\bA_-\bu_-)
+K_2(\CR^7(\bu_-)-\bG^0\nabla_x\CK_{F1}(\bu_-), \CR^5(\bu_-)),\\
p_2=&\CK_{F2}(h)-K_2(\bG^0\nabla_x\CK_{F2}(h),
(\rho_{*-}-\rho_{*+})^{-1}\rho_{*-}\sigma\CD'_xh).
\end{split} 
\end{equation}
Accordingly, we arrive at
\begin{align*}
\rho_{*-} \lambda(\bA_{-1}\bu_-) -\gamma_{4}\Delta_x(\bA_{-1}\bu_-)
-\nabla_x(\CK_{F1}(\bA_-\bu_-)+\CK_{F2}(h) )+\CR^9(\bu_-,h)
&=\bA_{-1}\bff_-, 
\end{align*}
where we have set
\begin{align*}
\CR^9(\bu_-,h) & =-\nabla_x\Big\{K_2(\CR^7(\bu_-)-\bG^0 \nabla_x\CK_{F1}(\bu_-), \CR^5(\bu_-))\\
&\quad-K_2(\bG^0\nabla_x\CK_{F2}(h),
(\rho_{*-}-\rho_{*+})^{-1}\rho_{*-}\sigma\CD'_xh)\Big\}\\
&\quad-\bA_{-1}{}^{\top} \!\bB_{-1}\nabla_x
\Big\{\CK_{F1}(\bA_-\bu_-)+K_2(\CR^7(\bu_-)
-\bG^0\nabla_x\CK_{F1}(\bu_-),\CR^5(\bu_-))\\
&\quad+\CK_{F2}(h)-K_2(\bG^0\nabla_x\CK_{F2}(h),
(\rho_{*-}-\rho_{*+})^{-1}\rho_{*-}\sigma\CD'_xh)\Big\}
+\bA_{-1}\CR^8(\bu_-).
\end{align*}
\par
\textit{\underline{The fourth line of (\ref{eq-5.11}).}}
By (\ref{5.2}) and (\ref{5.3}), we have
\begin{align*}
\lambda h
-\frac{\langle\rho_{*-}\bA_{-1}\bu_-,\bn_0\rangle|_-
-\langle\rho_{*+}\bA_{-1}\bu_+,\bn_0\rangle|_+}
{\rho_{*-}-\rho_{*+}}+\CR^{10}(\bu_+,\bu_-)=d
\end{align*}
with
\begin{align*}
\CR^{10}(\bu_+,\bu_-)
=&\frac{\rho_{*-}}{\rho_{*-}-\rho_{*+}}
\{\langle\bB_{-1}\bu_-,\bn_0\rangle
+b_0\langle(\bA_{-1}+\bB_{-1})\bu_-,\bn_0\rangle \}\Big|_-\\
&-\frac{\rho_{*+}}{\rho_{*-}-\rho_{*+}}
\{\langle\bB_{-1}\bu_+,\bn_0\rangle
+b_0\langle(\bA_{-1}+\bB_{-1})\bu_+,\bn_0\rangle\}\Big|_+.
\end{align*}
\par
\textit{\underline{The fifth line of (\ref{eq-5.11}).}}
By (\ref{5.17}) and (\ref{5.20}), we have
\begin{align*}
{}^\top\!\bA_{-1}\tgamma_{4}\bD_y(\widetilde{\bu}_-)\widetilde{\bn}
&=\gamma_{4}\bD_x(\bA_{-1}\bu_-)\bn_0+\CR^{11}(\bu_{-}),\\
{}^\top\!\bA_{-1}\rho_{*+}\bD_y(\widetilde{\bu}_+)\widetilde{\bn}
&=\gamma_1\bD_x(\bA_{-1}\bu_+)\bn_0+\CR^{12}(\bu_+),
\end{align*}
where we have set
\begin{align*}
\CR^{11}(\bu_-)
=&\gamma_{4}\Big((b_0\bI+(1+b_0){}^\top\!\bA_{-1}{}^{\top}\!\bB_{-1})
\bD_x(\bA_{-1}\bu_-)\bn_0+(1+b_0)(\bI+{}^\top\!\bA_{-1}
{}^\top\!\bB_{-1})\CR^3(\bu_{-})\Big)\Big|_-,\\
\CR^{12}(\bu_+)
=&\gamma_1\Big((b_0\bI+(1+b_0){}^\top\!\bA_{-1}{}^{\top}\!\bB_{-1})
\bD_x(\bA_{-1}\bu_+)\bn_0
+(1+b_0)(\bI+{}^\top\!\bA_{-1}{}^\top\!\bB_{-1})\CR^3(\bu_{-})\Big)\Big|_+.
\end{align*}
Hence, from (\ref{5.2}), (\ref{5.3}), and (\ref{5.17}) we obtain
\begin{align*}
\Pi_{\bn_0}(\gamma_{4}\bD_x(\bA_{-1}\bu_{-})\bn_0)|_-
-\Pi_{\bn_0}(\gamma_{1}\bD_x(\bA_{-1}\bu_{+})\bn_0)|_-
+\CR^{13}(\bu_{+},\bu_{-})=\bA_{-1}g,
\end{align*}
where we have set
\begin{align*}
\CR^{13}(\bu_{+},\bu_{-})
=&\CR^{11}(\bu_{-})|_- -\langle\gamma_{4}\bD_x(\bA_{-1}\bu_{-})\bn_0,\bn_0\rangle
(b_0\bI+(1+b_0)\bA_{-1}{}^\top\!\bB_{-1})|_-\\
&-\CR^3(\bu_{-})(1+b_0)(\bI+\bA_{-1}{}^\top\!\bB_{-1})|_-
-\CR^{12}(\bu_{+})|_+\\
&-\langle\gamma_{1}\bD_x(\bA_{-1}\bu_{+})\bn_0,\bn_0\rangle
(b_0\bI+(1+b_0)\bA_{-1}{}^\top\!\bB_{-1})|_+\\
&-\CR^3(\bu_{+})(1+b_0)(\bI+\bA_{-1}{}^\top\!\bB_{-1})|_+.
\end{align*}
\par
\textit{\underline{The sixth line of (\ref{eq-5.11}).}}
We easily see that
\begin{align*}
\rho_{*-}\dv_x(\bA_{-1}\bu_{-})|_-+\CR^{14}(\bu_{-})=g_-\quad\text{on $\RZ$},
\end{align*}
where $\CR^{14}(\bu_{-})$ is a linear function such that
$\CR^{14}(\bu_{-})=\rho_{*-}(\bB_{-1}:\nabla_x\bu_{-})|_-$.
\par
\textit{\underline{The seventh line of (\ref{eq-5.11}).}}
From (\ref{5.17}) we rewrite the sixth line of (\ref{eq-5.11}) as
\begin{equation*}
\begin{split} 
&\langle(\tgamma_{1}\bD_y(\widetilde{\bu}_+)
+(\tgamma_{2}-\tgamma_{1})\dv_y\widetilde{\bu}_+\bI 
+\rho_{*+}\tgamma_{2}\Delta_y\trho_+\bI)\widetilde{\bn},\widetilde{\bn}\rangle|_+\\
&=\{\langle\gamma_{1}\bD_x(\bA_{-1}\bu_+)\bn_0,\bn_0\rangle
+(\gamma_{2}-\gamma_{1})\dv_x(\bA_{-1}\bu_+)+\rho_{*+}\gamma_{2}\Delta_x\rho_+\}|_+ 
+\CR^{15}(\bu_+,\rho_{+})
\end{split}
\end{equation*}
with
\begin{align*}
\CR^{15}(\bu_+,\rho_+)
& = \Big\{(1+b_0)^2\bigg(\bB_{-1}:\nabla\bu_+
+\rho_{*+}\gamma_{3}\sum_{j=1}^N\CD_{jj}\rho_+\bigg)\\
& \quad + (2b_0+b_0^2)(\dv_x(\bA_{-1}\bu_+)+
\rho_{*+}\gamma_{3}\Delta_x\rho_+)\Big\}\Big|_+.
\end{align*}
Then the seventh line of (\ref{eq-5.11}) can be rewritten as
\begin{equation*}
\begin{split}
\{\langle\gamma_{1}\bD_x(\bA_{-1}\bu_+)\bn_0,\bn_0\rangle
+(\gamma_{2}-\gamma_{1})\dv_x(\bA_{-1}\bu_+)
+\rho_{*+}\gamma_{3}\Delta_x\rho_+\}|_+&\\
-\frac{\rho_{*+}\sigma}{\rho_{*-}-\rho_{*+}}\Delta'_x h
+\CR^{16}(\bu_+,\rho_{+},h)&=g_+
\end{split}
\end{equation*}
with
\begin{align*}
\CR^{16}(\bu_+,\rho_{+},h)
=\CR^{15}(\bu_{+},\rho_{+})-\frac{\rho_{*+}\sigma}{\rho_{*-}-\rho_{*+}}\CD'_x h.
\end{align*}
\par
\textit{\underline{The eighth line of (\ref{eq-5.11}).}}
From (\ref{5.2}) and (\ref{5.3}) we see that
\begin{align*}
\bA_{-1}\Pi_{\wt \bn} \widetilde{\bu}_\pm
=\Pi_{\bn_0} (\bA_{-1} \bu_\pm)+\CR^{17}(\bu_\pm)
\end{align*}
with
\begin{align*}
\CR^{17}(\bu_\pm) = (2b_0+b_0^2)\Pi_{\bn_0}(\bA_{-1}\bu_{-})
+(1+b_0)^2\{\langle(\bA_{-1}+\bB_{-1})\bu_\pm,\bn_0\rangle
\bA_{-1}{}^\top\!\bB_{-1}\bn_0 + \langle\bB_{-1}\bu_\pm,\bn_0\rangle\bn_0 \}.
\end{align*}
Here, we use the fact that
$\langle\bu_\pm,{}^\top\!(\bA_{-1}+\bB_{-1})\bn_0\rangle
=\langle(\bA_{-1}+\bB_{-1})\bu_\pm,\bn_0\rangle$.
If we define an operator
$\CR^{18}(\bu_+,\bu_-)$ as
$\CR^{18}(\bu_+,\bu_-)=\CR^{17}(\bu_-)|_- -\CR^{17}(\bu_+)|_+$, 
we have
\begin{align*}
\Pi_{\bn_0}(\bA_{-1}\bu_-)|_- -\Pi_{\bn_0}(\bA_{-1}\bu_+)|_+
+\CR^{18}(\bu_+,\bu_-) =\bA_{-1}\bh.
\end{align*}
\par
\textit{\underline{The ninth line of (\ref{eq-5.11}).}}
By (\ref{5.2})  there exists a linear operator
$\CR^{19}(\rho_{+})$ such that
\begin{gather*}
\langle\nabla_x\rho_{+},\bn_0\rangle+\CR^{19}(\rho_{+}) = k_-, \\
\CR^{19}(\rho_+) = (1+b_0)(\langle{}^\top\!\bB_{-1}\nabla_x\rho_{+},
({}^\top\!\bA_{-1}+{}^\top\!\bB_{-1})\bn_0\rangle
+\langle{}^\top\!\bB_{-1}\nabla_x\rho_{+},{}^\top\!\bA_{-1}\bn_0\rangle) + b_0 \langle\nabla_x\rho_{+},\bn_0\rangle
\end{gather*}

We now set
\begin{align*}
\bu^\bA_\pm=&\bA_{-1}\bu_\pm, \quad
\bff_+^\bA=\bA_{-1}\bff_+,\quad
\bff_-^\bA=\bA_{-1}\bff_-,\quad
\bh^\bA=\bA_{-1}\bh,\quad g^\bA=\bA_{-1} g.
\end{align*}
Summing up, noting
$\bu_\pm = (\bA_{-1})^{-1}\bu^\bA_\pm = {}^{\top} \!\bA_{-1}\bu^\bA_\pm$,
from the argumentation mentioned above, we observe
\begin{align}\label{eq-5.22}
\left\{\begin{aligned}
\lambda \rho_+ + \rho_{*+} \dv \bu_+^\bA + R_1 (\bu_+^\bA) & = f_M &\enskip &\text{ in $\RP$,} \\
\rho_{*+} \lambda \bu_+^\bA - \DV \bT_{0 +} (\gamma_{10}, \gamma_{20}, \gamma_{30}, \bu_+^\bA, \rho_+) + R_2 (\rho_+, \bu_+^\bA) & = \bff_+^\bA &\enskip &\text{ in $\RP$}, \\
\rho_{*-} \lambda \bu_-^\bA - \DV \bT_- (\gamma_{40}, \bu_-, \CK_1 (\bu_-^\bA) + \CK_2 (h)) + R_3 (\bu_-^\bA, h) & = \bff_-^\bA &\enskip &\text{ in $\RM$}, \\
\lambda h - \frac{\langle \rho_{* -} \bu_-^\bA, \bn \rangle \vert_- - \langle \rho_{* +} \bu_+^\bA, \bn \rangle \vert_+} {\rho_{* -} - \rho_{* +}} + R_4 (\bu_+^\bA, \bu_-^\bA) & = d &\enskip &\text{ on $\RZ$}, \\
\bB_0 (\gamma_{10}, \gamma_{20}, \gamma_{30}, \gamma_{40}, \rho_+, \bu_+^\bA, \bu_-^\bA) + R (\rho_+, \bu_+^\bA, \bu_-^\bA, h) & = \bG_0^\bA &\enskip &\text{ on $\RZ$}, \\
\rho_{*-} (\dv \bu_-^\bA) \vert_- + R_7 (\bu_-^\bA) & = g_- &\enskip &\text{ on $\RZ$},
\end{aligned}\right.
\end{align}
where
\begin{align*}
\bG_0^\bA = & {}^\top\! (g^\bA, g_+, \bh, k_-), \\
R (\rho_+, \bu_+^\bA, \bu_-^\bA, h) = & {}^\top\! (R_5 (\bu^\bA_+,\bu^\bA_-), R_6 (\bu^\bA_+,\rho_{+},h), R_8 (\bu^\bA_+,\bu^\bA_-), R_9 (\rho_+)), \\
R_1(\bu^\bA_+)=&\CR^1({}^{\top} \!\bA_{-1}\bu^\bA_+),\\
R_2(\bu^\bA_+,\rho_+)=
&\bA_{-1}\CR^2({}^{\top} \!\bA_{-1}\bu^\bA_+,\rho_{+}) 
-(\gamma_{1}-\gamma_{10})\Delta\bu^\bA_+
- (\gamma_{2}-\gamma_{20})\nabla\dv\bu^\bA_+
-\rho_{*+}(\gamma_{3}-\gamma_{30})\Delta\nabla\rho_+,\\
R_3(\bu^\bA_-,h) =&\CR^9({}^{\top} \!\bA_{-1}\bu^\bA_-,h)
-(\gamma_{4}-\gamma_{40})\Delta\bu^\bA_- ,\\
R_4(\bu^\bA_+,\bu^\bA_-)=
&\CR^{10}({}^{\top} \!\bA_{-1}\bu^\bA_+,{}^{\top} \!\bA_{-1}\bu^\bA_-),\\
R_5(\bu^\bA_+,\bu^\bA_-)=
&\CR^{13}({}^{\top} \!\bA_{-1}\bu^\bA_+,{}^{\top} \!\bA_{-1}\bu^\bA_-)
+\Pi_{\bn_0}\{(\gamma_{4}-\gamma_{40})\bD(\bu_-^\bA)|_- 
-(\gamma_{1}-\gamma_{10})\bD(\bu_-^\bA)|_+ \},\\
R_6(\bu^\bA_+,\rho_{+},h)=
&\CR^{16}({}^\top\!\bA_{-1}\bu_{+}^\bA,\rho_{+},h)
+\{\langle(\gamma_{1}-\gamma_{10})\bD(\bu_+^\bA)\bn_0,\bn_0\rangle\\
&+((\gamma_{2}-\gamma_{1})-(\gamma_{20}-\gamma_{10}))\dv(\bu_+^\bA)
+\rho_{*+}(\gamma_{3}-\gamma_{30})\Delta\rho_+\}|_+,\\
R_7(\bu^\bA_-)=
&\CR^{14}({}^\top\!\bA_{-1}\bu_{-}^\bA) ,\quad
R_8(\bu^\bA_+,\bu^\bA_-)=
\CR^{18}({}^{\top} \!\bA_{-1}\bu^\bA_+,{}^{\top} \!\bA_{-1}\bu^\bA_-),
\quad R_9(\rho_+)=\CR^{20}(\rho_+).
\end{align*}
To estimate $R_i$, $i=1,\dots,9$, we use the following lemma proved by Shibata~\cite[Lemma 2.4]{Shi2013}).
\begin{lemm}\label{lem-5.3}
Let $1< q< r<\infty$, $N<r<\infty$, and $D\in\{\RP,\RM,\dot{\BR}^N \}$.
In addition, let $f\in L_r(D)$, and $g\in W^1_q(D)$.
Then there exist a positive constant $C_{N,q,r}$ such that
\begin{align*}
\|fg\|_{L_q(D)}\leq C_{N,q,r}
\|f\|_{L_r(D)}\|g\|^{1-N/r}_{L_q(D)}
\|\nabla g\|^{N/r}_{L_q(D)}.
\end{align*}
Furthermore, for any $\delta>0$ we have
\begin{align*}
\|fg\|_{L_q(D)}\leq 
\delta \|\nabla g\|_{L_q(D)}
+C_{N,q,r}\delta^{-N/(r-N)}
\|f\|^{r/(r-N)}_{L_r(D)}
\|g\|_{L_q(D)}.
\end{align*}
\end{lemm}
From the representation formulas of $\CR^i$ $(i=1,\dots,19)$,
we see that the following estimates hold:
{\allowdisplaybreaks \begin{align*}
\|R_1(\bu^\bA_+)\|_{W^1_q(\RP)}\leq
& C_\gamma (M_1+\delta) \|\nabla^2\bu^\bA_+\|_{L_q(\RP)} 
+ C_{M_2,\delta}\|\nabla\bu^\bA_+\|_{L_q(\RP)}, \\
\|R_1(\bu^\bA_+)\|_{L_q(\RP)}\leq
& C_\gamma M_1 \|\nabla \bu^\bA_+\|_{L_q(\RP)},   \\
\|R_2(\bu^\bA_+,\rho_+)\|_{L_q(\RP)}\leq
&C_\gamma M_1 \Big(|\lambda|\|\bu^\bA_+\|_{L_q(\RP)}
+\|\nabla^3\rho_{+} \|_{L_q(\RM)}\Big)  \\
&+C_\gamma (M_1+\delta)\|\nabla^2\bu^\bA_+\|_{L_q(\RP)} 
+ C_{M_2,\delta}\Big(\|\bu^\bA_+\|_{W^1_q(\RP)}
+\|\rho_+ \|_{W^2_q(\RM)}\Big),   \\
\|R_3(\bu^\bA_-,h) \|_{L_q(\RM)}\leq
&C_\gamma (M_1+\delta) \Big(\|\nabla^2\bu^\bA_-\|_{L_q(\RM)}
+\|h\|_{W^{3-1/q}_q(\BR^{N-1})}\Big)  \\
&+C_{M_2,\delta}\Big(\|\bu^\bA_-\|_{W^1_q(\RM)}
+\|h\|_{W^{2-1/q}_q(\BR^{N-1})}\Big),  \\
\|R_4(\bu^\bA_+,\bu^\bA_-)\|_{W^2_q(\dot{\BR}^N)}\leq
&C_\gamma M_1\Big(\|\nabla^2\bu^\bA_+\|_{L_q(\RP)} 
+\|\nabla^2\bu^\bA_-\|_{L_q(\RM)}\Big)  \\ 
&+ C_{M_2,\delta}\Big(\|\bu^\bA_+\|_{W^1_q(\RP)}
+ \|\bu^\bA_-\|_{W^1_q(\RM)}\Big),   \\
\|R_5(\bu^\bA_+,\bu^\bA_-)\|_{W^1_q(\dot{\BR}^N)}\leq&
C_\gamma (M_1+\delta)\Big( \|\nabla^2\bu^\bA_+\|_{L_q(\RP)}
+ \|\nabla^2\bu^\bA_-\|_{L_q(\RM)} \Big)  \\
&+ C_{M_2,\delta}\Big(\|\nabla\bu^\bA_+\|_{L_q(\RP)}
+\|\nabla\bu^\bA_-\|_{L_q(\RM)}\Big),  \\
\|R_6(\bu^\bA_+,\rho_{+},h)\|_{W^1_q(\dot{\BR}^N)}\leq&
C_\gamma (M_1+\delta) \Big(\|\nabla^2\bu^\bA_+\|_{L_q(\RP)}
+\|\nabla^3\rho_{+}\|_{L_q(\RP)} +\|h\|_{W^{3-1/q}_q(\BR^{N-1})}\Big)  \\
&+C_{M_2,\delta}\Big(\|\bu^\bA_+\|_{W^1_q(\RP)}
+\|\rho_+\|_{W^2_q(\RP)} +\|h\|_{W^{2-1/q}_q(\BR^{N-1})}\Big),  \\
\|R_7(\bu^\bA_-)\|_{W^1_q(\dot{\BR}^N)}\leq&
C_\gamma (M_1+\delta) \|\nabla^2\bu^\bA_-\|_{L_q(\RM)} 
+ C_{M_2,\delta}\|\nabla\bu^\bA_-\|_{L_q(\RM)},  \\
\|R_8(\bu^\bA_+,\bu^\bA_-)\|_{W^2_q(\dot{\BR}^N)}\leq&
C_\gamma(M_1+\delta)\Big(\|\nabla^2\bu^\bA_+\|_{L_q(\RP)}
\|\nabla^2\bu^\bA_-\|_{L_q(\RM)}\Big)  \\
&+C_{M_2,\delta}\Big(\|\bu^\bA_+\|_{W^1_q(\RP)}+\|\bu^\bA_-\|_{W^1_q(\RM)}\Big),  \\
\|R_9(\rho_+)\|_{W^2_q(\RP)}\leq&
C_\gamma(M_1+\delta)\|\nabla^3\rho_+\|_{L_q(\RP)}
+C_{M_2,\delta}\|\rho_+\|_{W^2_q(\RP)}.  
\end{align*}}\noindent
Here, $C_\gamma$ are positive constant depending only on $N$, $q$, $r$,
and $\gamma^\pm_{k*}$ $(k=1,\dots,4)$
while $C_{M_2,\delta}$ are positive constant depending only on
$M_2$, $\delta$, $N$, $q$, $r$, and $\gamma^\pm_{k*}$ $(k=1,\dots,4)$.
In fact, from (\ref{cond-5.1}), (\ref{cond-5.4}), 
(\ref{cond-5.9}), and (\ref{cond-5.10}) 
we see that $R_2(\bu^\bA_+,\rho_+)$ has the following form:
\begin{align*}
R_2(\bu^\bA_+,\rho_+)=\lambda\CM_1\bu^\bA_++\CM_2\nabla^2\bu_+^\bA
+\CM_3\nabla^3\rho_+ +\CN_1\nabla\bu^\bA_+ +\CN_2\nabla\rho_+
\end{align*}
with some matrices of functions $\CM_1$, $\CM_2$, $\CM_3$, 
$\CN_1$, and $\CN_2$ possessing the estimates:
\begin{align*}
\|(\CM_1,\CM_2,\CM_3)\|_{L_\infty(\RP)}&\leq C_\gamma M_1,\quad
\|(\CN_1,\CN_2)\|_{L_r(\RP)} \leq C_{M_2,\delta}.
\end{align*}
Hence, by Lemma \ref{lem-5.3}, we have the estimate for $R_2(\bu^\bA_+,\rho_+)$.
Analogously, we have the estimate for $R_1(\bu^\bA_+)$, 
$R_4(\bu^\bA_+,\bu^\bA_-)$, $R_5(\bu^\bA_+,\bu^\bA_-)$, $R_7(\bu^\bA_-)$,
$R_8(\bu_+^\bA,\bu_-^\bA)$, and $R_9(\rho_+)$.
Furthermore, using Lemmas~\ref{lem-5.2} and \ref{lem-5.3} and the estimates (\ref{est-4.4}), (\ref{cond-5.1}),~(\ref{est-5.7}),~(\ref{est-5.8}), 
we obtain the estimate for $R_3(\bu^\bA_-,h)$ and $R_6(\bu^\bA_+,\bu^\bA_-,\rho_{+},h)$.
\par
Let $\CA^+_{F0}(\lambda)$, $\CB^\pm_{F0}(\lambda)$, and $\CH_{F0}(\lambda)$
be the solution operators of problem (\ref{eq-5.22}) given in Theorem~\ref{th-4.4}, and set
$\rho_+ = \CA^+ _{F0}(\lambda)G_\lambda(\bF_Z)$, $\bu^\bA_+ = \CB^\pm_{F0} (\lambda) G_\lambda(\bF_Z)$,
$h = \CH_{F0} (\lambda) G_\lambda(\bF_Z)$.
Then the problem (\ref{eq-5.22}) can be rewritten as follows:
\begin{align*}
\left\{\begin{aligned}
\lambda \rho_+ + \rho_{*+} \dv \bu_+^\bA & = f_M - V_1(\lambda)G_\lambda(\bF_Z) &\enskip &\text{ in $\RP$,} \\
\rho_{*+} \lambda \bu_+^\bA - \DV \bT_{0 +} (\gamma_{10}, \gamma_{20}, \gamma_{30}, \bu_+^\bA, \rho_+) & = \bff_+^\bA - V_2 (\lambda)G_\lambda(\bF_Z) &\enskip &\text{ in $\RP$}, \\
\rho_{*-} \lambda \bu_-^\bA - \DV \bT_- (\gamma_{40}, \bu_-, \CK_1 (\bu_-^\bA) + \CK_2 (h)) & = \bff_-^\bA - V_3 (\lambda) G_\lambda(\bF_Z) &\enskip &\text{ in $\RM$}, \\
\lambda h - \frac{\langle \rho_{* -} \bu_-^\bA, \bn \rangle \vert_- - \langle \rho_{* +} \bu_+^\bA, \bn \rangle \vert_+} {\rho_{* -} - \rho_{* +}} & = d - V_4 (\lambda)G_\lambda(\bF_Z) &\enskip &\text{ on $\RZ$}, \\
\bB_0 (\gamma_{10}, \gamma_{20}, \gamma_{30}, \gamma_{40}, \rho_+, \bu_+^\bA, \bu_-^\bA)& = \bG_0^\bA - V (\lambda) G_\lambda (\bF_Z) &\enskip &\text{ on $\RZ$}, \\
\rho_{*-} (\dv \bu_-^\bA) \vert_- & = g_- - V_7 (\lambda) G_\lambda (\bF_Z) &\enskip &\text{ on $\RZ$},
\end{aligned}\right.
\end{align*}
where we have set
{\allowdisplaybreaks
\begin{align*}
V (\lambda)G_\lambda(\bF_Z) = & {}^\top\!(V_5(\lambda)G_\lambda(\bF_Z), V_6(\lambda)G_\lambda(\bF_Z), V_8(\lambda)G_\lambda(\bF_Z), V_9(\lambda)G_\lambda(\bF_Z)), \\
V_1(\lambda)G_\lambda(\bF_Z)=
&R_1(\CB^+_{F0}(\lambda)G_\lambda(\bF_Z)), \\
V_2(\lambda)G_\lambda(\bF_Z)=
&R_2(\CB^+_{F0}(\lambda)G_\lambda(\bF_Z),\CA^+_{F0}(\lambda)G_\lambda(\bF_Z)),\\
V_3(\lambda)G_\lambda(\bF_Z)=
&R_3(\CB^+_{F0}(\lambda)G_\lambda(\bF_Z),\CH_{F0}(\lambda)G_\lambda(\bF_Z)), \\
V_4(\lambda)G_\lambda(\bF_Z)=
&R_4(\CB^+_{F0}(\lambda)G_\lambda(\bF_Z),\CB^-_{F0}(\lambda)G_\lambda(\bF_Z)),\\
V_5(\lambda)G_\lambda(\bF_Z)=
&R_5(\CB^+_{F0}(\lambda)G_\lambda(\bF_Z),\CB^-_{F0}(\lambda)G_\lambda(\bF_Z)), \\
V_6(\lambda)G_\lambda(\bF_Z)=
&R_6(\CB^+_{F0}(\lambda)G_\lambda(\bF_Z),\CA^+_{F0}(\lambda)G_\lambda(\bF_Z),\CH_{F0}(\lambda)G_\lambda(\bF_Z)),  \\
V_7(\lambda)G_\lambda(\bF_Z)=
&R_7(\CB^-_{F0}(\lambda)G_\lambda(\bF_Z)),\\
V_8(\lambda)G_\lambda(\bF_Z)=
&R_8(\CB^+_{F0}(\lambda)G_\lambda(\bF_Z),\CB^-_{F0}(\lambda)G_\lambda(\bF_Z)),\\
V_9(\lambda)G_\lambda(\bF_Z)=
&R_9 (\CA^+_{F0}(\lambda)G_\lambda(\bF_Z)).
\end{align*}}\noindent
Let $\bQ(\lambda)$ and $\CQ(\lambda)$ be the operator such that
{\allowdisplaybreaks
\begin{align*}
\bQ(\lambda)\bF=
&(V_1(\lambda)G_\lambda(\bF_Z),\dots,V_9(\lambda)G_\lambda(\bF_Z))\\
\CQ(\lambda)G=&( \lambda^{1/2}V_1(\lambda)G, \nabla V_1(\lambda)G,
V_2(\lambda)G, V_3(\lambda)G, V_4(\lambda)G, \lambda^{1/2}V_5(\lambda)G,
\nabla V_5(\lambda)G, \lambda^{1/2} V_6(\lambda)G, \\
&\nabla V_6(\lambda)G, \lambda^{1/2}V_7(\lambda)G,\nabla V_7(\lambda)G,
V_7(\lambda)G,\lambda V_8(\lambda)G, \lambda^{1/2}\nabla V_8(\lambda)G,
\nabla^2 V_8(\lambda)G,\\
&\lambda V_9(\lambda)G, \lambda^{1/2}\nabla V_9(\lambda)G,\nabla^2 V_9(\lambda)G),
\end{align*}}\noindent
respectively. We easily see that
\begin{align}\label{5.27}
G_\lambda(\bQ(\lambda)(\bF_Z))=\CQ(\lambda)G_\lambda(\bF_Z).
\end{align}
From the definitions of $\CR$-boundedness, Lemma \ref{lem-4.2},
Theorem \ref{th-4.4}, we have
\begin{align*}
\CR_{\CL(\CZ_q(\RP,\RM,\RZ))}(\{(\tau\pd_\tau)^s\CQ(\lambda)\mid
\lambda\in\Sigma_{\varepsilon,\lambda_4} \})
\leq C_\gamma(M_1+\delta)+ C_{M_2,\delta}\lambda_4^{-1/2}
\qquad (s=0,1)
\end{align*}
for any $\lambda_4\geq \max(1,\lambda_1)$ with some positive constant $C$
independent of $\lambda_4$.
We now choose $M_1>0$ and $\delta>0$ so small that
$C_\gamma M_1\leq 1/4$ and $C_{M_2,\delta}\leq 1/4$,
and then we choose $\lambda_4>\max(1,\lambda_1)$ so large that
$C_{M_2,\delta}\lambda_4^{-1/2}\leq 1/4$.
Then we arrive at
\begin{align}\label{est-5.28}
\CR_{\CL(\CZ_q(\RP,\RM,\RZ))}(\{(\tau\pd_\tau)^s\CQ(\lambda)\mid
\lambda\in\Sigma_{\varepsilon,\lambda_4} \})
\leq \frac34 \qquad (s=0,1),
\end{align}
which, combined with (\ref{5.27}), furnishes that
\begin{align}\label{est-5.29}
\|G_\lambda(\bQ(\lambda)(\bF_Z))\|_{\CZ_q(\RP,\RM,\RZ)}
\leq \frac34 \|G_\lambda(\bF_Z)\|_{\CZ_q(\RP,\RM,\RZ)}.
\end{align}
Since $\|G_\lambda(\bF_Z)\|_{\CZ_q(\RP,\RM,\RZ)}
=\|\bF_Z\|_{Z_q(\RP,\RM,\RZ)}$
gives us equivalent norms for $\lambda\neq 0$, by (\ref{est-5.29}) the operator $\bQ(\lambda)$ is a contraction map from $Z_q(\RP,\RM,\RZ)$ into itself, so that there exists an
inverse operator $(\bI-\bQ(\lambda))^{-1}$ of $\bI-\bQ(\lambda)$ in
$\CL(\CZ_q(\RP,\RM,\RZ))$ for any
$\lambda\in\Sigma_{\varepsilon,\lambda_4}$ .
Hence, 
\begin{align*}
\rho_+ = \CA^+_{F0}(\lambda)G_\lambda((\bI-\bQ(\lambda))^{-1}\bF_Z), \enskip                         \bu^\bA_\pm = \CB^\pm_{F0}(\lambda)G_\lambda((\bI-\bQ(\lambda))^{-1}\bF_Z), \enskip                                                                                                                                                             h = \CH_{F0}(\lambda)G_\lambda((\bI-\bQ(\lambda))^{-1}\bF_Z)
\end{align*}
are unique solutions of problem (\ref{eq-5.22}).
On the other hand, from (\ref{est-5.28}) we see that
$(\bI-\CQ(\lambda))^{-1}=\sum_{j=1}^{\infty}(\CQ(\lambda))^j$
exists in $\CL(\CZ_q(\RP,\RM,\RZ))$ enjoying the estimate
\begin{align*}
\CR_{\CL(\CZ_q(\RP,\RM,\RZ))}(\{(\tau\pd_\tau\tau)^s
(\bI-\CQ(\lambda))^{-1}\mid \lambda\in\Sigma_{\varepsilon,\lambda_4} \})
\leq c\qquad(s=0,1)
\end{align*}
with some positive constant $c$.
Since
\begin{align}\label{5.31}
G_\lambda(\bI-\bQ(\lambda))^{-1}
=\sum_{j=1}^{\infty}G_\lambda (\bQ(\lambda))^j
=\bigg(\sum_{j=1}^{\infty}(\CQ(\lambda))^j\bigg)G_\lambda
=(\bI-\CQ(\lambda))^{-1}G_\lambda
\end{align}
as follows from (\ref{5.27}), if we define operator families
$\widetilde{\CA}^+(\lambda)$, $\widetilde{\CB}^\pm(\lambda)$
$\widetilde{\CH}(\lambda)$ by
\begin{align*}
\widetilde{\CA}^+ (\lambda) = \CA^+_{F0}(\lambda)(\bI-\CQ(\lambda))^{-1},\quad
\widetilde{\CB}^\pm (\lambda) = \CB^\pm_{F0}(\lambda)(\bI-\CQ(\lambda))^{-1},\quad
\widetilde{\CH}(\lambda) = \CH_{F0}(\lambda)(\bI-\CQ(\lambda))^{-1},
\end{align*}
then by (\ref{5.31}) we have
\begin{align*}
\rho_+=\widetilde{\CA}^+(\lambda)G_\lambda(\bF_Z),\quad
\bu^\bA_\pm=\widetilde{\CB}^\pm(\lambda)G_\lambda(\bF_Z),\quad
h=\widetilde{\CH}(\lambda)G_\lambda(\bF_Z).
\end{align*}
In addition, by Lemma \ref{lem-4.2}, Theorem \ref{th-4.4}, 
and (\ref{est-5.29}), we obtain
\begin{equation}\label{est-5.32}
\begin{split} 
\CR_{\CL(\CZ_q(\RP,\RM,\RZ),W^{3-i}_q(\BR^N_\pm))}(\{(\tau\pd_\tau)^s
(\lambda^{i/2}\widetilde{\CA}^+(\lambda))
\mid \lambda\in\Sigma_{\varepsilon,\lambda_4} \})&\leq c_4,\\
\CR_{\CL(\CZ_q(\RP,\RM,\RZ),W^{2-j}_q(\RP)^N)}(\{(\tau\pd_\tau)^s
(\lambda^{j/2}\widetilde{\CB}^\pm(\lambda))
\mid \lambda\in\Sigma_{\varepsilon,\lambda_4} \})&\leq c_4,\\
\CR_{\CL(\CZ_q(\RP,\RM,\RZ),W^{3-1/q-k}_q(\RZ))}(\{(\tau\pd_\tau)^s
(\lambda^{k}\widetilde{\CH}(\lambda))
\mid \lambda\in\Sigma_{\varepsilon,\lambda_4} \})&\leq c_4
\end{split}
\end{equation}
for $s=0,1$, $i=0,1,2,3$, $j=0,1,2$, and $k=0,1$
with some positive constant $c_4$.

Recalling that
\begin{alignat*}6
\rho_+ & = \trho_+\circ\Phi, &\quad \bu^\bA_\pm & = \bA_{-1}\widetilde{\bu}\circ\Phi, &\quad h & = \tilh\circ\Phi, &\quad
f_M & = \tf_M \circ\Phi, &\quad \bff_\pm & = \bA_{-1}\widetilde{\bff}_\pm\circ\Phi, \\
d & = \td \circ\Phi &\quad g^\bA & = \bA_{-1}\tg\circ\Phi,&\quad g_\pm & = \tg_\pm\circ\Phi, 
&\quad \bh & = \bA_{-1}\widetilde{\bh}\circ\Phi, 
&\quad k_- & = \widetilde{k}_- \circ\Phi,
\end{alignat*}
we define operators $\CA^+_F(\lambda)$, $\CB^\pm_{F}(\lambda)$, $\CH_{F}(\lambda)$ such that
\begin{align*}
\CA^+_{F} (G) & = \Big[\widetilde{\CA}^+(\lambda) G^\Phi \Big]\circ\Phi^{-1}, \enskip \CB^\pm_{F} (G) = \Big[{}^\top \!\bA_{-1} \widetilde{\CB}^\pm(\lambda) G^\Phi \Big]\circ\Phi^{-1}, \enskip \CH_{F} (G) = \Big[\widetilde{\CH}(\lambda) G^\Phi \Big]\circ\Phi^{-1}, \\
G & := (G_1,\dots,G_{18}), \\
G^\Phi & := (G_1\circ\Phi,G_2\circ\Phi,\bA_{-1}G_3\circ\Phi,\bA_{-1}G_4\circ\Phi,
G_5\circ\Phi,\bA_{-1}G_6\circ\Phi,\bA_{-1}G_7\circ\Phi,\bA_{-1}G_8\circ\Phi,G_9\circ\Phi,\\
& \quad G_{10}\circ\Phi,G_{11}\circ\Phi,G_{12}\circ\Phi,\bA_{-1}G_{13}\circ\Phi, \bA_{-1}G_{14}\circ\Phi,\bA_{-1}G_{15}\circ\Phi,G_{16}\circ\Phi, G_{17}\circ\Phi,G_{18}\circ\Phi).
\end{align*}
Then by (\ref{est-5.32}) we see that $\CA^+_{F}(\lambda)$,
$\CB^\pm_{F}(\lambda)$, and $\CH_{F}(\lambda)$ satisfy the
required properties in Theorem \ref{th-5.1}.
Hence, we complete the proof of Theorem \ref{th-5.1}.

\subsection{Proof of Theorem \ref{TH-3.1}}\label{sect-6}
\subsubsection{Some preparation for the proof of Theorem \ref{TH-3.1}}
Since $\gamma_{k}(x)$ $(k=1,\dots,4)$ are
uniformly continuous functions defined on $\BR^N$
satisfying the Assumption \ref{as-star} (2), choosing $d^i>0$
small if necessary, we may assume that
$|\gamma_{k}(x)-\gamma_{k}(x^i_j)|\leq M_1$ for any
$x\in B_{d^i}(x^i_j)$ with $i=0,\dots,4$ and $j\in\BN$.
Furthermore, we choose $M_2$ so large that $\|\nabla\gamma_k\|_{L_r(B_{d^i}(x^i_j))}\leq M_2$ $(k=1,\dots,4)$. Thanks to the choices of $M_1$, $M_2$, and $d^i$, we may assume that
\begin{align}\label{cond-6.1}
\begin{gathered}
\frac12 \gamma^-_{k*}\leq \gamma_k(x^i_j)\leq \frac32 \gamma^+_{k*},
\quad (x\in B_{d^i}(x^i_j))\\
|\gamma_k(x)-\gamma_k(x^i_j)|\leq M_1,\quad
\|\nabla \gamma_k\|_{L_r(B_{d^i}(x^i_j))}\leq M_2
\end{gathered}
\end{align}
for $k=1,2,3,4$. \par
In the following, to simplify notation, we write $B^i_j$ instead of $B_{d^i}(x^i_j)$.
From the finite intersection property stated in 
Proposition \ref{prop-6.1} (4), we see that for any 
$1\leq s<\infty$ there exists a positive constant $C_{s,L}$
such that for any $f\in L_s(D)$ with an open set $D$ of
$\BR^N$ and for $i=0,\dots,4$,
\begin{align}\label{cond-6.2}
\bigg(\sum_{j=1}^\infty\|f\|^s_{L_s(D\cap B^i_j)} \bigg)^{1/s}
\leq C_{s,L}\|f\|_{L_s(D)}.
\end{align}
In fact, when $1\leq s<\infty$, we see that
\begin{align*}
\sum_{j=1}^\infty\|f\|^s_{L_s(D\cap B^i_j)} = \int_{D}\bigg(\sum_{j=1}^{\infty}\chi_{B^i_j}(x) \bigg)|f(x)|^s \dx \leq \bigg\|\sum_{j=1}^{\infty}\chi_{B^i_j}(x) \bigg\|_{L_\infty(\BR^N)}
\|f\|^s_{L_s(D)}\leq L\|f\|^s_{L_s(D)}.
\end{align*}
For simplicity of notation, we write
$D^0_{j\pm}=\Phi^0_j(\BR^N_\pm)$, $D^0_j=D^0_{j+}\cup D^0_{j-}$,
$D^1_j=\Phi^1_{j}(\RP)$,
$D^2_j=\Phi^2_{j}(\RM)$ $D^3_j=D^4_j=\BR^N$,
and $\Gamma^i_j=\Phi^i_j(\RZ)$ $(i=0,1,2)$ for short.
Let $n\in \BN_0$, $f\in W^n_q(\Omega)$, and let $\eta^i_j$ be functions in
$C^\infty_0(B^i_j)$ with $\|\eta^i_j\|_{W^n_\infty(\BR^n)}\leq c_1$ for some
positive constant $c_1$ independent of $j\in\BN$.
Since $\Omega\cap B^i_j=D^i_j\cap B^i_j$, from (\ref{cond-6.2}) we have
\begin{align}\label{est-6.4}
\sum_{j=1}^{\infty}\|\eta^i_j f\|^q_{W^n_q(D^i_j)}\leq C_q\|f\|_{W^n_q(\Omega)}
\end{align}
According to Shibata~\cite[Proposition 9.5.4, 9.5.5]{Shi2016a},
we know the following two propositions,
which we will use for defining the infinite sum of $\CR$-bounded operator families
defined on $D^i_j$.
\begin{prop}\label{prop-6.4}
Let $1<q<\infty$, and $n=2,3$. Then  the following assertions hold:\\
(1) There exist extension maps 
$\bT^n_j: W^{n-1/q}_q(\Gamma_j^0)\to W^n_q(D^0_{j-})$
such that for $h\in W^{n-1/q}_q(\Gamma^0_j)$, $\bT^n_j h=h$ on $\Gamma^n_j$ and
$\|\bT^n_j h\|_{W^n_q(D^0_{j-})}\leq C\|h\|_{W^{n-1/q}_q(\Gamma^0_j)}$ with some
positive constant $C$ independent of $j\in\BN$.\\
(2) There exist extension maps 
$\bT^n_\Gamma: W^{n-1/q}_q(\Gamma)\to W^n_q(\Omega_{-})$
such that for $h\in W^{n-1/q}_q(\Gamma)$, $\bT^n_\Gamma h=h$ on $\Gamma$ and
$\|\bT^n_j h\|_{W^n_q(\Omega_{-})}\leq C\|h\|_{W^{n-1/q}_q(\Gamma)}$ with some
positive constant $C$ independent of $j\in\BN$.	
\end{prop}
\begin{prop}\label{prop-6.5}
Let $1<q<\infty$ and $n=2,3$, and let $\eta_j\in C^\infty_0(B^0_j)$ $(j\in \BN)$
with $\|\eta_j\|_{W^n_\infty(\BR^N)}\leq c_7$ for some positive constant $c_7$
independent of $j\in\BN$.
Then  we have the following assertions:\\
(1) Let $f_j$ $(j\in\BN)$ be functions in $W^{n-1/q}_q(\Gamma_j^0)$ such that
$\sum_{j=1}^{\infty}\|f_j\|^q_{W^{n-1/q}_q(\Gamma^0_j)} < \infty$,
and then the infinite sum $\sum_{j=1}^{\infty}\eta_j f_j$ convergence to some
$f\in W^{n-1/q}_q(\Gamma)$ strongly in $W^{n-1/q}_q(\Gamma)$ and
\begin{align*}
\|f\|_{W^{n-1/q}_q(\Gamma)}
\leq C_q \bigg(\sum_{j=1}^{\infty}\|f_j\|^q_{W^{n-1/q}_q(\Gamma^0_j)} \bigg)^{1/q}.
\end{align*}
(2) For any $h\in W^{n-1/q}_q(\Gamma)$, we have $\sum_{j=1}^{\infty}\|\eta_j h\|^q_{W^{n-1/q}_q(\Gamma^0_j)} \leq C \|h\|^q_{W^{n-1/q}_q(\Gamma)}$.
\end{prop}

\subsubsection{Local solutions}
In view of (\ref{cond-6.1}), we define the functions $\gamma^i_{jk}$ by
\begin{align*}
\gamma^i_{jk}(x)=(\gamma_k(x)-\gamma_k(x^i_j))
\tzeta^i_j(x)+\gamma_k(x^i_j)\quad (k=1,\dots,4).
\end{align*}
Note that $0\leq \tzeta^i_j\leq 1$ and
$\|\nabla\tzeta^i_j\|_{L_\infty(\BR^N)}\leq c_0$.
From (\ref{cond-6.1}) we obtain
\begin{equation*}
\begin{split}
\frac12 \gamma^-_{k*} &\leq \gamma^0_{jk}(x) \leq 
\frac32 \gamma^+_{k*}
\quad(x\in D^0_{j+},\enskip k=1,2,3),\\
\frac12 \gamma^-_{4*} &\leq \gamma^0_{j4}(x) \leq 
\frac32 \gamma^+_{4*}
\quad(x\in D^0_{j-}),\\
\frac12 \gamma^-_{k*} &\leq \gamma^i_{jk}(x) \leq 
\frac32 \gamma^+_{k*}
\quad(x\in D^i_{j},\enskip i=1,2,3,4,\enskip k=1,2,3,4),\\
\|\gamma^0_{jk}-\gamma^0_{jk}(x^0_j)\|_{L_\infty(D^0_{j+})}&\leq M_1,\quad
\|\nabla \gamma^0_{jk}\|_{L_r(D^0_{j+})}\leq C_{M_2}
\quad(k=1,2,3),\\
\|\gamma^0_{j5}-\gamma^0_{j4}(x^0_j)\|_{L_\infty(D^0_{j-})}&\leq M_1,\quad
\|\nabla \gamma^0_{j4}\|_{L_r(D^0_{j-})}\leq C_{M_2},\\
\|\gamma^i_{jk}-\gamma^i_{jk}(x^i_j)\|_{L_\infty(D^i_j)}&\leq M_1,\quad
\|\nabla \gamma^i_{jk}\|_{L_r(D^i_j)}\leq C_{M_2}
\quad(i=1,2,3,4,\enskip k=1,2,3,4),
\end{split}
\end{equation*}
In addition, since $\tzeta^i_j(x)=1$ on $\supp\zeta^i_j$, we have $\gamma^i_{jk}(x)=\gamma_k(x)$ for $x \in \supp \zeta^i_j$, $k=1,2,3,4$. For $\bF_Y \in Y_q(\Omega_+,\Omega_-,\Gamma)$, we consider the following equations:
\allowdisplaybreaks{
\begin{align}
\label{eq-6.7}
\left\{\begin{aligned}
\lambda\rho_{j+}^0+\rho_{*+}\dv\bu_{j+}^0
&=\tzeta^0_j f_M &\quad &\text{ in $D^0_{j+}$,}\\
\rho_{*+}\lambda\bu_{j+}^0
-\DV \bT_{0 +} (\gamma_{j1}^0, \gamma_{j2}^0, \gamma_{j3}^0, \bu_{j+}^0, \rho_{j+}^0) & =\tzeta^0_j \bff_+ &\quad&\text{ in $D^0_{j+}$}, \\
\rho_{*-}\lambda\bu_{j-}^0
-\DV \bT_- (\gamma_{j4}^0, \bu_{j-}^0, \CK^0_{j1}(\bu_{j-}^0)+\CK_{j2}^0(h_j)) &=\tzeta^0_j \bff_- 
&\quad&\text{ in $D^0_{j-}$},\\
\lambda h_j-\frac{\langle\rho_{*-}\bu_{j-}^0,\bn^0_j\rangle|_- 
-\langle\rho_{*+}\bu_{j+}^0,\bn^0_j\rangle|_+}{\rho_{*-}-\rho_{*+}}
&=\tzeta^0_j d &\quad&\text{ on $\Gamma^0_j$},\\
\bB_0 (\gamma_{j1}^0, \gamma_{j2}^0, \gamma_{j3}^0, \gamma_{j4}^0, \rho_{j +}^0, \bu_{j +}^0, \bu_{j -}^0) &=\tzeta^0_j \bG_0 &\quad&\text{ on $\Gamma^0_j$},\\
\langle \bT_- (\gamma^0_{j4}, \bu_{j -}^0, h_j) \bn^0_j, \bn^0_j\rangle\Big|_- - \frac{\rho_{*-}\sigma}{\rho_{*-}-\rho_{*+}}\Delta_{\Gamma^0_j} h_j &=\tzeta^0_j g_-&\quad&\text{ on $\Gamma^0_j$},\\
\end{aligned}\right. \\
\label{eq-6.8}
\left\{\begin{aligned}
\lambda\rho_{j}^1+\rho_{*+}\dv\bu_{j}^1
&=\tzeta^1_j f_M &\quad &\text{ in $D^1_{j}$,}\\
\rho_{*+}\lambda\bu_{j}^1
-\DV \bT_{0 +} (\gamma_{j1}^1, \gamma_{j2}^1, \gamma_{j3}^1, \bu_{j}^1, \rho_{j}^1) &=\tzeta^1_j \bff_+ &\quad&\text{ in $D^1_{j}$}, \\
\bu_{j}^1|_+=0,\quad \langle\nabla\rho_{j}^1,\bn_j^1\rangle|_+&=\tzeta^1_j k_+
&\quad&\text{ in $\Gamma^1_{j}$},
\end{aligned}\right. \\
\label{eq-6.9}
\left\{\begin{aligned}
\rho_{*-}\lambda\bu_{j}^2
-\DV \bT_- (\gamma_{j1}^2, \bu_{j}^2, \CK^2_{j1}(\bu_{j}^2)) &=\tzeta^2_j \bff_-
&\quad&\text{ in $D^2_{j}$},\\
\bu_{j}^2|_-&=0&\quad&\text{ in $\Gamma^2_{j}$},
\end{aligned}\right. \\
\label{eq-6.10}
\left\{\begin{aligned}
\lambda\rho_{j}^3+\rho_{*+}\dv\bu_{j}^3
&=\tzeta^3_j f_M &\quad &\text{ in $D^3_{j}$,}\\
\rho_{*+}\lambda\bu_{j}^3
-\DV \bT_{0+} (\gamma_{j1}^3, \gamma_{j2}^3, \bu_{j}^3, \rho_{j}^3) &=\tzeta^3_j \bff_+ &\quad&\text{ in $D^3_{j}$}, \\
\end{aligned}\right. \\
\label{eq-6.11}
\quad\rho_{*-}\lambda\bu_{j}^4 - \DV \bT_- (\gamma_{j1}^4, \bu_{j}^4, \CK^4_{j1}(\bu_{j}^4)) = \tzeta^4_j \bff_-
\quad \text{ in $D^4_{j}$},
\end{align}}Here, $\bn_j^i$ $(i=0,1)$ denote
the unit outer normal to $\Gamma_i^i$ and
$\Delta_{\Gamma^0_j}$ the Laplace-Beltrami operator on $\Gamma^0_j$.
We may assume that $\bn_j^i$ are defined in
$\BR^N$ and satisfy the estimates:
\begin{align*}
\|\bn_j^i\|_{L_\infty(\BR^N)}\leq C_N,\quad
\|\nabla \bn_j^i \|_{W^{2-i}_r(\BR^N)}\leq C_{N,r}M_2 
\quad\text{for $i=0,1$, $j=1,2,\dots$},
\end{align*}
while we can write $\Delta_{\Gamma^0_j}$ as $\Delta_{\Gamma^0_j} f=\Delta' f+\CD_{\Gamma^0_j}$ with $\CD_{\Gamma^0_j} f=\sum_{i,j=1}^{N-1}g^0_{ij}\pd_i\pd_j f
+\sum_{i=1}^{N-1}g^0_i\pd_i f$ in local coordinates, where $g^0_{ij}$ and $g^0_j$ can be bounded by $
\|g^0_{ij}\|_{L_\infty(\BR^N)}\leq C_N M_1$ and $\|(\nabla g^0_{ij},g^0_j)\|_{W^2_r(\BR^N)}\leq C_{N,r}M_2$, respectively. In addition, the functions $\CK^i_{j1}$, $\CK^0_{j2}$ $(i=0,2,4\enskip j\in \BN)$ 
are given as follows:
For $\bu_{j-}^0\in W^2_q(D^0_{j-})$, 
let $\CK^0_{j1}(\bu_{j-}^0)\in W^1_q(D^0_{j-})+\hW^1_q(D^0_{j-})$ be a unique
solution of the variational equation:
\begin{align}\label{eq-6.13}
(\nabla \CK_{j1}^0(\bu_{j-}^0),\nabla\varphi)_{D^0_{j-}}
=(\DV(\rho_{*+} \bD(\bu_{j-}^0))-\rho_{*-}\nabla\dv \bu_{j-}^0,
\nabla\varphi)_{D^0_{j-}}
\end{align}
for any $\varphi\in\hW^1_{q',0}(D^0_{j-})$
subject to $\CK_{j1}^0(\bu_{j-}^0)
=\langle\rho_{*+}\bD(\bu_{j-}^0)\bn_j^0,\bn_j^0\rangle
-\rho_{*-}\dv\bu_{j-}^0$ on $\Gamma_j^0$, while for $h_j\in W^{3-1/q}_q(D^0_j)$,
let $\CK_{j2}^0\in W^1_q(D^0_{j-})+\hW^1_{q,0}(D^0_{j-})$ be a unique solution
of the variational equation:
\begin{align}\label{eq-6.14}
(\nabla\CK_{j2}^0(h_j),\nabla\varphi)_{D^0_j}=0
\qquad\text{for $\varphi\in\hW^1_{q',0}(D^0_{j-})$}
\end{align}
subject to $\CK_{j2}^0(h_j)
=-(\rho_{*-}-\rho_{*+})^{-1}\rho_{*-}\sigma\Delta_{\Gamma^0_j}h_j$
on $\Gamma_j^0$.
For $\bu_{j}^i\in W^2_q(D^i_{j})$ $(i=2,4)$,
let $\CK^i_{j1}(\bu_{j}^i)\in W^1_q(D^i_{j})+\hW^1_q(D^i_{j})$ be a unique
solution of the variational equation:
\begin{align}\label{eq-6.15}
(\nabla \CK_{j1}^i(\bu_{j}^i),\nabla\varphi)_{D^i_{j}}
=(\DV(\rho_{*+} \bD(\bu_{j}^i))-\rho_{*-}\nabla\dv \bu_{j}^i,
\nabla\varphi)_{D^i_{j}}
\end{align}
for any $\varphi\in\hW^1_{q',0}(D^i_{j})$.

Choosing $0 < M_1 < 1$ suitably small, we have a unique existence of solutions
to (\ref{eq-6.13}), (\ref{eq-6.14}), and~(\ref{eq-6.15}) possessing the estimates:
\begin{equation*}
\begin{split} 
\|\nabla\CK_{j1}^0(\bu_{j-}^0)\|_{L_q(D^0_{j-})}
&\leq c_8\|\bu_{j-}^0\|_{W^1_q(D^0_{j-})},\\
\|\nabla\CK_{j2}^0(h_{j})\|_{L_q(D^0_{j-})}
&\leq c_8\|h_{j}\|_{W^{3-1/q}_q(\Gamma^0_{j-})},\\
\|\nabla\CK_{j1}^i(\bu_{j}^i)\|_{L_q(D^i_{j})}
&\leq c_8\|\bu_{j}^i\|_{W^1_q(D^0_{j})}
\end{split}
\end{equation*}
for $i=2,4$ and $j\in\BN$ with some constant $c_8$ independent of $j\in\BN$.
For $i=1,2,3,4$ and $j\in\BN$ define function spaces $\widetilde{Y}_q(i,j)$
and $\widetilde{\CY}_q(i,j)$ by
\begin{align*}
\widetilde{Y}_q(i,j)=
\begin{cases}
Z_q(D^0_{j+},D^0_{j-},\Gamma^0_j),\\
W^1_q(D^1_j)\times L_q(D^1_j)^N\times W^2_q(D^1_j),\\
L_q(D^2_j)^N,\\
W^1_q(D^3_j)\times L_q(D^3_j)^N,\\
L_q(D^4_j)^N,
\end{cases}
\quad
\widetilde{\CY}_q(i,j)=
\begin{cases}
\CZ_q(D^0_{j+},D^0_{j-},\Gamma^0_j),\\
W^1_q(D^1_j)\times L_q(D^1_j)^N\times W^2_q(D^1_j),\\
L_q(D^2_j)^N,\\
W^1_q(D^3_j)\times L_q(D^3_j)^N,\\
L_q(D^4_j)^N,
\end{cases}
\end{align*}
respectively.
Then, according to the results of analysis in the whole space, the whole/half spaces with a bent interface, see Saito~\cite{Sai2019}, Shibata~\cite{Shi2014}, and Theorem \ref{th-5.1}, there exists a constant
$\lambda_6\geq \max(\lambda_4,\lambda_5)$ and operator families
\begin{alignat*}3
\CA^i_{j +}(\lambda)
&\in\Hol(\Sigma_{\varepsilon_0,\lambda_0},
\CL(\widetilde{\CY}_q(i,j), W^3_q(D^i_{j\pm})))
&\quad &(i=0,1,3),\\
\CB^i_{j \pm}(\lambda)
&\in\Hol(\Sigma_{\varepsilon_0,\lambda_0},\CL(\widetilde{\CY}_q(i,j),W^2_q(D^i_j) ))
&\quad &(i=0,1,2,3,4), \\
\CH_{j}(\lambda)
&\in\Hol(\Sigma_{\varepsilon_0,\lambda_0},\CL(\CZ_q(D^0_{j+},D^0_{j-},
\Gamma^0_{j}), W^{3-1/q}_q(\Gamma^0_{j})))
\end{alignat*}
such that problem (\ref{eq-6.7}), (\ref{eq-6.8}), (\ref{eq-6.9}),
(\ref{eq-6.10}), (\ref{eq-6.11}) admit unique solutions
\begin{align*}
\rho_{j+}^0 & = \CA^0_{j+}(\lambda) G_\lambda (\bF^0_j), \quad \bu_{j\pm}^0 = \CB^0_{j\pm}(\lambda) G_\lambda(\bF^0_j),
\quad h_j = \CH_{j}(\lambda) G_\lambda(\bF^0_j),\\
\rho_{j}^1 & = \CA^1_{j +} (\lambda)(\lambda^{1/2}\tzeta^1_j f_M, \nabla\tzeta^1_j f_M, \tzeta^1_j \bff_+,
\tzeta^1_j k_+), \quad \bu_{j}^1 = \CB^1_{j +} (\lambda)(
\lambda^{1/2}\tzeta^1_j f_M,\nabla\tzeta^1_j f_M, \tzeta^1_j \bff_+,
\tzeta^1_j k_+), \\
\rho_{j}^3 & = \CA^3_{j +} (\lambda)(
\lambda^{1/2}\tzeta^3_j f_M,\nabla\tzeta^3_j f_M,\tzeta^3_j \bff_+), \quad
\bu_{j}^2 =\CB^2_{j -}(\lambda)\tzeta^2_j \bff_-,\quad
\bu_{j}^3 = \CB^3_{j +}(\lambda)(
\lambda^{1/2}\tzeta^3_j f_M,\nabla\tzeta^3_j f_M,\tzeta^3_j \bff_+), \\
\bu_{j}^4 & = \CB^4_{j -} (\lambda)\tzeta^4_j\bff_-
\end{align*}
for $j\in\BN$. Here, we have set $\bF^0_j = \tzeta^0_j \bF_Z$.
Furthermore, we have the estimates
\begin{equation}\label{6.17}
\begin{split}
\CR_{\CL(\widetilde{\CY}_q(i,j),W^{3-l}_q(D^i_{j+}))}
(\{(\tau\pd_\tau)^s(\lambda^{l/2} \CA^i_{j+}(\lambda))
\mid \lambda \in \Sigma_{\varepsilon,\lambda_6}\})\leq& \omega_{0}
\quad(i=0,1,3),\\
\CR_{\CL(\widetilde{\CY}_q(i,j),W^{2-m}_q(D^i_{j\pm})^N)}	
(\{(\tau\pd_\tau)^s(\lambda^{m/2} \CB^i_{j\pm}(\lambda))
\mid \lambda \in \Sigma_{\varepsilon,\lambda_6}\})\leq& \omega_{0}
\quad(i=0,1,2,3,4),\\
\CR_{\CL(\CZ_{q}(D^0_{j+},D^0_{j-},\Gamma^0_j),W^{3-1/q-n}_q(\Gamma^0_j))}
(\{(\tau\pd_\tau)^s(\lambda^n \CH_j (\lambda))
\mid \lambda  \in \Sigma_{\varepsilon,\lambda_6}\})\leq& \omega_{0}
\end{split}
\end{equation}
for $s=0,1$, $j\in\BN$, $l=0,1,2,3$, $m=0,1,2$, and $n=0,1$
with some positive constant $\omega_0$ independent of $j\in\BN$.
Since the $\CR$-boundedness implies the uniform boundedness,
by \eqref{6.17} for any $j\in\BN$ and $\lambda\in\Sigma_{\varepsilon_0,\lambda_6}$ there exists a constant $C$ independent of $j\in\BN$ such that
\begin{equation}\label{6.18}
\begin{split}
&\|(\lambda \bu_{j+}^0,\lambda^{1/2}\nabla\bu_{j+}^0,
\nabla^2\bu_{j+}^0,\lambda^{3/2}\rho_{j+}^0,\lambda\nabla\rho_{j+}^0,
\lambda^{1/2}\nabla^2\rho_{j+}^0,\nabla^3\rho_{j+})\|_{L_q(D^0_{j+})}\\
& +\|(\lambda \bu_{j-}^0,\lambda^{1/2}\nabla\bu_{j-}^0,
\nabla^2\bu_{j-}^0)\|_{L_q(D^0_{j-})}
+\|(\lambda h_j,\nabla h_j)\|_{W^{2-1/q}_q(\Gamma^0_j)} \leq C \lVert G_\lambda (\bF^0_j) \rVert_{\CZ_q (D^0_{j +}, D^0_{j -}, \Gamma_j^0)}, \\
&\|(\lambda \bu_{j}^1,\lambda^{1/2}\nabla\bu_{j}^1,\nabla^2\bu_{j}^1,
\lambda^{3/2}\rho_{j}^1,\lambda\nabla\rho_{j}^1,\lambda^{1/2}\nabla^2\rho_{j}^1,
\nabla^3\rho_{j}^1)\|_{L_q(D^i_j)}\\
&\quad \leq C\|(\lambda^{1/2}\tzeta^1_j f_M,\nabla\tzeta^1_j f_M,
\tzeta^1_j\bff_+,\lambda\tzeta^1_j k_+,\lambda^{1/2}\nabla\tzeta^1_j k_+,
\nabla^2\tzeta^1_j k_+)\|_{L_q(D^1_j)},\\
&\|(\lambda \bu_{j}^i,\lambda^{1/2}\nabla\bu_{j}^i,\nabla^2\bu_{j}^i)\|_{L_q(D^i_j)}
\leq C\|\tzeta^i_j\bff_-\|_{L_q(D^i_j)} \qquad(i=2,4),\\
&\|(\lambda \bu_{j}^3,\lambda^{1/2}\nabla\bu_{j}^3,\nabla^2\bu_{j}^3,
\lambda^{3/2}\rho_{j}^3,\lambda\nabla\rho_{j}^3,\lambda^{1/2}\nabla^2\rho_{j}^3,
\nabla^3\rho_{j}^3)\|_{L_q(D^3_j)} \leq C\|\lambda^{1/2}\tzeta^3_j f_M,\nabla\tzeta^3_j f_M,
\tzeta^3_j\bff_+\|_{L_q(D^3_j)}.
\end{split}
\end{equation}

\subsubsection{Construction of a parametrix}
For 
$\bF_Y \in Y_q(\Omega_+,\Omega_-,\Gamma)$ we consider the reduced equations (\ref{eq-3.8}).
We use the notations defined in Proposition \ref{prop-6.1}.
Let us define $\wh \rho_+$, $\wh \bu_\pm$, and $\wh h$ by
\begin{align}\label{6.19}
\begin{aligned}
\wh \rho_+&=\sum_{j=1}^{\infty}\bigg(\zeta^0_j\rho^0_{j+}
+\sum_{i=1,3}\zeta^i_j\rho_{j}^i\bigg),&\quad 
\wh \bu_+&=\sum_{j=1}^\infty\bigg(\zeta^0_j\bu^0_{j+}
+\sum_{i=1,3}\zeta^i_j\bu_{j}^i\bigg),\\
\wh \bu_-&=\sum_{j=1}^\infty\bigg(\zeta^0_j\bu^0_{j-}
+\sum_{i=2,4}\zeta^i_j\bu_{j}^i\bigg),
&\quad \wh h&=\sum_{j=1}^{\infty}\zeta^0_j h_j.
\end{aligned}
\end{align}
The right-hand sides of (\ref{6.19}) convergence strongly to $\wh \rho_+\in W^3_q(\Omega_+)$,
$\wh \bu_\pm \in W^2_q(\Omega_\pm)$, and $\wh h\in W^{3-1/q}_q(\Gamma)$, respectively.
In fact, from (\ref{est-6.4}) and \eqref{6.18} we see that
\begin{align*}
\sum_{j=1}^{\infty}\|h_j\|_{W^{3-1/q}(\Gamma_j^0)} \leq C\sum_{j=1}^{\infty} \lVert G_\lambda (\bF^0_j) \rVert_{\CZ_q (D^0_{j+}, D^0_{j-}, \Gamma^0_j)} \leq C\lVert G_\lambda (\bF_Z) \rVert_{\CZ_q (\Omega_+, \Omega_-, \Gamma)} <\infty.
\end{align*}
Then, by Proposition \ref{prop-6.5}, the sum $\sum_{j=1}^{\infty}\zeta^0_j h_j$
convergences to $\wh h$ in $W^{3-1/q}_q(\Gamma)$ strongly and
\begin{align*}
\|\wh h\|_{W^{3-1/q}_q(\Gamma)}\leq C_q \lVert F_\lambda (\bF_Y) \rVert_{\CY_q (\Omega_+, \Omega_-, \Gamma)}.
\end{align*}
Analogously, the infinite sum
\begin{align*}
\sum_{j=1}^{\infty}\bigg(\zeta^0_j\rho^0_{j+}+\sum_{i=1,3}\zeta^i_j\rho_{j}^i\bigg), \quad \sum_{j=1}^\infty\bigg(\zeta^0_j\bu^0_{j+}+\sum_{i=1,3}\zeta^i_j\bu_{j}^i\bigg),\quad
\sum_{j=1}^\infty\bigg(\zeta^0_j\bu^0_{j-}+\sum_{i=2,4}\zeta^i_j\bu_{j}^i\bigg)
\end{align*}
strongly convergence to $\wh \rho_+\in W^3_q(\Omega_+)$ 
and $\wh \bu_\pm \in W^2_q(\Omega_\pm)$, respectively,
and the following estimates hold:
\begin{align*}
\|\wh \rho_+\|_{W^3_q(\Omega_+)} + \|\wh \bu_\pm\|_{W^2_q(\Omega_\pm)}\leq&
C_q \lVert F_\lambda (\bF_Y) \rVert_{\CY_q (\Omega_+, \Omega_-, \Gamma)}.
\end{align*}
Since $\bn=\bn^0_j$ on $\supp \zeta^0_j\cap \Gamma$ and
$\bn_+=\bn^1_j$ on $\supp \zeta^1_j\cap \Gamma_+$,
we see that $(\wh \rho_+, \wh \bu_+, \wh \bu_-, \wh h)$ satisfies
\begin{align}\label{6.20}
\left\{\begin{aligned}
\lambda \wh \rho_+ + \rho_{*+} \dv \wh \bu_+ & = \wh f_M - \bV^1 (\lambda) \bF_Y &\enskip &\text{ in $\Omega_+$,} \\
\rho_{*+} \lambda \wh \bu_+ - \DV \bT_{0 +} (\gamma_1, \gamma_2, \gamma_3, \wh \bu_+, \wh \rho_+) & = \wh \bff_+ - \bV^2 (\lambda) \bF_Y &\enskip &\text{ in $\Omega_+$}, \\
\rho_{*-} \lambda \wh \bu_- - \DV \bT_- (\gamma_4, \wh \bu_-, \CK_1 (\wh \bu_-) + \CK_2 (\wh h)) & = \wh \bff_- - \bV^3 (\lambda) \bF_Y &\enskip &\text{ in $\Omega_-$}, \\
\lambda \wh h - \frac{\langle \rho_{*-} \wh \bu_-, \bn \rangle \vert_- - \langle \rho_{* +} \wh \bu_+, \bn \rangle\vert_+} {\rho_{* -} - \rho_{* +}} & = \wh d &\enskip &\text{ on $\Gamma$}, \\
\bB_0 (\gamma_1, \gamma_2, \gamma_3, \gamma_4, \wh \rho_+, \wh \bu_+, \wh \bu_-) & = \wh \bG_0 - \ov \bV (\lambda) \bF_Y &\enskip &\text{ on $\Gamma$}, \\
\langle \bT_- (\gamma_4, \wh \bu_-, \CK_1 (\wh \bu_-) + \CK_2 (\wh h)) \bn, \bn \rangle \Big\vert_- - \frac{\rho_{*-} \sigma}{\rho_{*-} - \rho_{*+}} \Delta_\Gamma \wh h & = \wh g_- - \bV^5 (\lambda) \bF_Y &\enskip &\text{ on $\Gamma$}, \\
\wh \bu_+ = 0, \qquad \langle \nabla \rho_+, \bn_+ \rangle & = \wh k_+ - \bV^8(\lambda)\bF_Y &\enskip & \text{ on $\Gamma_+$}, \\
\wh \bu_- & = 0 &\enskip & \text{ on $\Gamma_-$}.
\end{aligned}\right.
\end{align}
where we have set
\begin{align}\label{6.21}
\ov \bV(\lambda)\bF_Y= & {}^\top\! (\bV^4(\lambda)\bF_Y, \bV^6(\lambda)\bF_Y, 0, \bV^7(\lambda)\bF_Y), \\
\bV^1(\lambda)\bF_Y=&-\sum_{j=1}^{\infty}\Big(
\rho_{*+}(\dv(\zeta^0_j\bu_{j+}^0)-\zeta^0_j\dv(\bu^0_{j+}))
+ \sum_{i=1,3} \rho_{*+}(\dv(\zeta^i_j\bu_{j}^i)-\zeta^i_j\dv(\bu_{j}^i))\Big), \nonumber\\
\bV^2(\lambda)\bF_Y=&\sum_{j=1}^{\infty}\Big\{
\DV(\gamma_{j1}^0\bD(\zeta^0_j\bu_{j+}^0))-\zeta^0_j\DV(\gamma^0_{j1}\bD(\bu^0_j))
+\DV((\gamma^0_{j2}-\gamma^0_{j1})\dv(\zeta^0_j\bu_{j+}^0)\bI)
\nonumber\\
&-\zeta^0_j\DV((\gamma^0_{j2}-\gamma^0_{j1})\dv\bu_{j+}^0\bI)
+\DV(\rho_{*+}\gamma^0_{j3}\Delta(\zeta^0_j\rho_{j+}^0)\bI)
-\zeta^0_j\DV(\rho_{*+}\gamma^0_{j3}\Delta\rho_{j+}^0\bI)
\nonumber\\
&+\sum_{i=1,3}\Big(
\DV(\gamma_{j1}^i\bD(\zeta^i_j\bu_{j}^i))-\zeta^i_j\DV(\gamma^i_{j1}\bD(\bu^i_j))
+\DV((\gamma^i_{j2}-\gamma^i_{j1})\dv(\zeta^i_j\bu_{j+}^i)\bI)
\nonumber\\
&-\zeta^i_j\DV((\gamma^i_{j2}-\gamma^i_{j1})\dv\bu_{j+}^i\bI)
+\DV(\rho_{*+}\gamma^i_{j3}\Delta(\zeta^i_j\rho_{j+}^i)\bI)
-\zeta^i_j\DV(\rho_{*+}\gamma^i_{j3}\Delta\rho_{j+}^i\bI)\Big)\Big\},
\nonumber\\
\bV^3(\lambda)\bF_Y=&
\sum_{j=1}^{\infty}\Big\{\DV(\gamma_{j1}^0\bD(\zeta^0_j\bu_{j+}^0))
-\zeta^0_j\DV(\gamma^0_{j1}\bD(\bu^0_j))
\nonumber\\
&+\sum_{i=2,4}\Big( \DV(\gamma_{j1}^i\bD(\zeta^i_j\bu_{j}^i))
-\zeta^i_j\DV(\gamma^i_{j1}\bD(\bu^i_j))\Big)\Big\}
\nonumber\\
&+\sum_{j=1}^{\infty}\Big\{\nabla \CK_1(\zeta^0_j\bu_{j-}^0) 
- \zeta^0_j\nabla \CK_{j1}^0(\bu_{j-}^0)
+\sum_{i=2,4}\Big(\nabla \CK_1(\zeta^i_j\bu_{j-}^i)
-\zeta^i_j\nabla\CK_{j1}^i(\bu_{j-}^i) \Big)\Big\}
\nonumber\\
&+\sum_{j=1}^{\infty}\Big(\nabla \CK_2(\zeta^0_j h_{j}) 
- \zeta^0_j\nabla \CK_{j2}^0 (h_{j})\Big),
\nonumber\\
\bV^4(\lambda)\bF_Y=&
\rho_{*+}(\bD(\zeta^0_j\bu_{j-}^0)- \zeta^0_j\bD(\bu_{j-}^0))
-\rho_{*+}(\langle\bD(\zeta^0_j\bu_{j-}^0)\bn_j^0,\bn_j^0\rangle
- \langle\zeta^0_j\bD(\bu_{j-}^0)\bn_j^0,\bn_j^0\rangle)\bn_j^0
\nonumber\\
&-\gamma_{j1}^0(\bD(\zeta^0_j\bu_{j+}^0)- \zeta^0_j\bD(\bu_{j+}^0))
-\gamma_{j1}^0(\langle\bD(\zeta^0_j\bu_{j+}^0)\bn_j^0,\bn_j^0\rangle 
- \langle\zeta^0_j\bD(\bu_{j+}^0)\bn_j^0,\bn_j^0\rangle)\bn_j^0,
\nonumber\\
\bV^5(\lambda)\bF_Y=&-\sum_{j=1}^{\infty}\Big\{
\rho_{*+}(\langle\bD(\zeta^0_j\bu_{j-}^0)\bn_j^0,\bn_j^0\rangle-
\langle\zeta^0_j\bD(\bu_{j-}^0)\bn_j^0,\bn_j^0\rangle)
\nonumber\\
&-\rho_{*+}(\langle\bD(\zeta^0_j\bu^0_{j-})\bn_j^0,\bn_j^0\rangle
-\langle\zeta^0_j\bD(\bu_{j-}^0)\bn_j^0,\bn_j^0\rangle) - \rho_{*-}(\nabla\dv(\zeta^0_j\bu_{j-}^0)-\zeta^0_j\nabla\dv\bu_{j-}^0)
\Big\},\nonumber\\
\bV^6(\lambda)\bF_Y=&-\sum_{j=1}^{\infty}\Big\{
\gamma_{j1}^0(\langle\bD(\zeta^0_j\bu_{j+}^0)\bn_j^0,\bn_j^0\rangle
-\langle\zeta^0_j\bD(\bu_{j+}^0)\bn_j^0,\bn_j^0\rangle)
\nonumber\\
&-(\gamma_{j2}^0-\gamma_{j1}^0)(\dv(\zeta^0_j\bu_{j+}^0)
-\zeta^0_j\dv(\bu_{j+}^0))-\rho_{*+}\gamma_{j3}^0(\Delta(\zeta^0_j\rho_{j+}^0)
-\zeta^0_j\Delta\rho_{j+}^0)
\nonumber\\
&-\frac{\rho_{*+}}{\rho_{*-}-\rho_{*+}}(\Delta_{\Gamma^0_j}
(\zeta^0_j h_j) -\zeta^0_j\Delta_{\Gamma_j^0} h_j )\Big\},
\nonumber\\
\bV^7(\lambda)\bF_Y=&-\sum_{j=1}^{\infty}\langle\nabla(\zeta^0_j\rho_{j+}^0)
-\zeta^0_j\nabla\rho_{j+}^0, \bn_j^0\rangle,
\nonumber\\
\bV^8(\lambda)\bF_Y=&-\sum_{j=1}^{\infty}\langle\nabla(\zeta^1_j\rho_{j}^0)
-\zeta^1_j\nabla\rho_{j}^0, \bn_j^1\rangle.
\nonumber
\end{align} 
To obtain the representation formulas above we have used the fact that
\begin{align}\label{5.22}
\nabla \CK_1(\wh \bu_-)
=\sum_{j=1}^\infty\bigg(\nabla\CK_1(\zeta^0_j\bu^0_{j-})
+\sum_{i=2,4}\nabla\CK_1(\zeta^i_j\bu_{j}^i)\bigg),\quad
\nabla \CK_2(\wh h)=\sum_{j=1}^\infty\nabla\CK_2(\zeta^0_j h_j),
\end{align}
which follow from the strong convergence in (\ref{6.19}) and the continuity
of $K_1$ in Definition \ref{def-1.2}.
In view of (\ref{6.19}), we define the operators $\CA^+ (\lambda)$, $\CB^\pm (\lambda)$, $\CH(\lambda)$ acting on $F:=(F_1,\dots,F_{21})\in \CY_q(\Omega_{+},\Omega_{-},\Gamma)$ by
\begin{align*}
\CA^+(\lambda)F & = \sum_{j=1}^\infty
\Big(\zeta^0_j \CA^0_{j+}(\lambda)(F)
+\zeta^1_j \CA^1_j(\lambda) \tzeta^1_j (F_1,F_2,F_3,F_{19},F_{20},F_{21}) 
+\zeta^3_j \CA^3_j(\lambda) \tzeta^3_j (F_1,F_2,F_3) \Big),\\
\CB^+(\lambda)F & =
\sum_{j=1}^\infty\Big(\zeta^0_j \CB^0_{j+}(\lambda) (F)
+\zeta^1_j \CB^1_j(\lambda) \tzeta^1_j (F_1,F_2,F_3,F_{19},F_{20},F_{21}) 
+\zeta^3_j \CB^3_j(\lambda) \tzeta^3_j (F_1,F_2,F_3) \Big),\\
\CB^-(\lambda)F & = \sum_{j=1}^\infty
\bigg(\zeta^0_j \CB^0_{j-}(\lambda) (F)
+\sum_{i=2,4}\zeta^i_j \CB^i_j(\lambda) \tzeta^i_j F_4 \bigg),\\
\CH(\lambda)F & = \sum_{j=1}^{\infty}
\zeta^0_j \CH_j(\lambda) \tzeta^0_j (F).
\end{align*}
Then  we see that
\begin{equation}\label{6.23}
\begin{split}
&\CA^+ (\lambda)
\in\Hol(\Sigma_{\varepsilon_0,\lambda_0},\CL(\CY_q(\Omega_{+},\Omega_{-},
\Gamma), W^3_q(\Omega_{\pm}))),\\
&\CB^\pm(\lambda)
\in\Hol(\Sigma_{\varepsilon_0,\lambda_0},\CL(\CY_q(\Omega_{+},\Omega_{-},
\Gamma), W^2_q(\Omega_{+}))),\\
&\CH(\lambda)
\in\Hol(\Sigma_{\varepsilon_0,\lambda_0},\CL(\CY_q(\Omega_{+},\Omega_{-},
\Gamma), W^{3-1/q}_q(\Gamma))),\\
&\CR_{\CL(\CY_{q}(\Omega_{+},\Omega_{-},\Gamma),W^{3-l}_q(\Omega_+))}
(\{(\tau\pd_\tau)^s(\lambda^{l/2} \CA^+(\lambda))
\mid \lambda \in \Sigma_{\varepsilon,\lambda_0}\})\leq C_q\omega_{0},\\
&\CR_{\CL(\CY_{q}(\Omega_{+},\Omega_{-},\Gamma),W^{2-m}_q(\Omega_\pm)^N)}	
(\{(\tau\pd_\tau)^s(\lambda^{m/2} \CB^\pm(\lambda))
\mid \lambda \in \Sigma_{\varepsilon,\lambda_0}\})\leq C_q\omega_{0},\\
&\CR_{\CL(\CY_{q}(\Omega_{+},\Omega_{-},\Gamma),W^{3-1/q-n}_q(\Gamma))}
(\{(\tau\pd_\tau)^s(\lambda^n \CH (\lambda))
\mid \lambda  \in \Sigma_{\varepsilon,\lambda_0}\})\leq C_q\omega_{0}
\end{split}
\end{equation}
for $s=0,1$, $l=0,1,2,3$, $m=0,1,2$, and $n=0,1$.
In fact, for any $m\in\BN$, ${\lambda}_{k=1}^m\subset
\Sigma_{\varepsilon,\lambda_0}$, $\{E_k\}_{k=1}^m\subset
\CY_q(\Omega_+,\Omega_-,\Gamma)$ and $\{a_k\}_{k=1}^m\subset\BC$,
from (\ref{6.17}) we have
\begin{align*}
\bigg\|\sum_{k=1}^{m}a_k\lambda^n_k \CH_j(\lambda_k)\tzeta^0_j E_k \bigg\|_
{W^{3-1/q-n}_q(\Gamma^0_j)}
&\leq \sum_{k=1}^{m}|a_k| \|\lambda^n_k \CH_j(\lambda_k) 
\tzeta^0_j E_k\|_{W^{3-1/q-n}_q(\Gamma^0_j)}\\
&\leq  \omega_0\sum_{k=1}^{m}|a_k| 
\|\tzeta^0_j E_k\|_{\CY(D^0_{j+},D^0_{j-},\Gamma_j^0)}
\end{align*}
for $n=0,1$.
Hence, applying the same argument as in the proof of Proposition 5.3 
in~\cite{Shi2014}, by (\ref{est-6.4}) and Proposition \ref{prop-6.5}
we obtain
\begin{gather*}
\CH(\lambda)\in \Hol(\Sigma_{\varepsilon,\lambda_0},
\CL(\CY_q(\Omega_+,\Omega_-,\Gamma),W^{3-1/q-n}_q(\Gamma))),\\
\bigg\|\sum_{k=1}^{m}a_k \lambda^n_k\CH(\lambda_k)E_k\bigg\|^q_
{W^{3-1/q-n}_q(\Gamma)}
\leq C_q\sum_{j=1}^{\infty}\bigg\|\sum_{k=1}^{m}a_k\lambda^n_k\zeta^0_j
\CH_j(\lambda_k)\tzeta^0_j E_k\bigg\|^q_{W^{3-1/q-n}_q(\Gamma_j^0)}
\end{gather*}
for $n=0,1$.
By (\ref{6.17}), the monotone convergence theorem for Lebesgue integral,
(\ref{est-6.4}), and Proposition~\ref{prop-6.5} (2), we have
\begin{align*}
\int_{0}^{1}\bigg\|\sum_{k=1}^{m}r_k(t)\lambda^n_k\CH(\lambda_k) E_k \bigg\|^q_{W^{3-1/q-n}_q(\Gamma)}\dt & \leq C_q\sum_{j=1}^{\infty}\int_0^1 \bigg\|\sum_{k=1}^{m}r_k(t)\lambda^n_k
\CH_j(\lambda_k)\tzeta^0_j E_k\bigg\|^q_{W^{3-1/q-n}_q(\Gamma_j^1)} \dt\\
&\leq C_q\omega_0\int_{0}^{1}\sum_{j=1}^{\infty}
\bigg\|\tzeta^0_j\sum_{k=1}^{m}r_k(t)E_k\bigg\|^q
_{\CY_q(D^0_{j+},D_{j-}^0,\Gamma_j^0)}\dt\\
&\leq C_q\omega_0\int^1_0\bigg\|\sum_{k=1}^{m}r_k(t) E_k\bigg\|^q_
{\CY_q(\Omega_+,\Omega_-,\Gamma)}\dt\qquad(n=0,1),
\end{align*}
which completes the proof of the assertion for $\CH(\lambda)$.
Analogously, we obtain the assertions for $\CA^+ (\lambda)$ and $\CB^\pm (\lambda)$. Finally, combining Propositions \ref{prop-6.4} and \ref{prop-6.5}, we obtain the following lemma.
\begin{lemm}\label{lem-6.6}
Let $\lambda_6$ and $\omega_0$ be the same constants as in (\ref{6.17}).
In addition, let $\bV^i(\lambda)$ $(i=1,\dots,8)$ be the operators defined
in (\ref{6.21}), and set
\begin{align*}
\bV(\lambda)\bF_Y=(\bV^1(\lambda),\bV^2(\lambda),\bV^3(\lambda),0,\bV^4(\lambda),
\bV^5(\lambda),\bV^6(\lambda),0,\bV^7(\lambda),\bV^8(\lambda),0).
\end{align*}
Then  there exists an operator family
$\CV(\lambda)\in \Hol(\Sigma_{\varepsilon,\lambda_*},
\CL(\CY_q(\Omega_+,\Omega_-,\Gamma)))$ such that
\begin{equation}\label{6.25}
\begin{split}
&F_\lambda(\bV(\lambda)\bF_Y)=\CV(\lambda)F_\lambda(\bF_Y)
\qquad(\bF_Y\in Y_q(\Omega_+,\Omega_-,\Gamma)),\\
&\CR_{\CL(\CY_{q}(\Omega_{+},\Omega_{-},\Gamma))}(\{(\tau\pd_\tau )^s
\CV(\lambda)\mid \lambda\in \Sigma_{\varepsilon,\lambda_*}\})
\leq (\sigma_*+C_{\sigma_*}\lambda^{-1/2}_*)\omega_0\qquad(s=0,1)
\end{split}
\end{equation}
for any $\sigma_*>0$ and $\lambda_*\geq \max(\lambda_6,1)$,
where $C_{\sigma_*}$ is some positive constant independent of
$\lambda_*$.
\end{lemm}

\subsubsection{Completion of the proof of Theorem~\ref{TH-3.1}}
From Lemma \ref{lem-6.6} we have
\begin{align*}
\|F_\lambda(\bV(\lambda)\bF_Y) \|_{\CY_q(\Omega_{+},\Omega_{-},\Gamma)} = \|\CV(\lambda)F_\lambda(\bF_Y)\|_{\CY_q(\Omega_+,\Omega_-,\Gamma)} \leq (\sigma_*+C_{\sigma_*}\lambda_*^{-1/2})\omega_0
\|F_\lambda(\bF_Y)\|_{\CY_q(\Omega_+,\Omega_{-},\Gamma)}.
\end{align*}
We choose $\sigma_*$ very small and $\lambda_*$ sufficiently large such that
$0<\sigma_*\omega_0\leq 1/4$ and $C_{\sigma_*}\lambda^{-1/2}_* \omega_0\leq 1/4$.
Then, we have $\|F_\lambda(\bV(\lambda)\bF_Y) \|_{\CY_q(\Omega_{+},\Omega_{-},\Gamma)}
\leq 2^{- 1} \|F_\lambda(\bF_Y)\|_{\CY_q(\Omega_+,\Omega_{-},\Gamma)}$ for any $\lambda\in\Sigma_{\varepsilon,\lambda_*}$,
so that $(I-\bV(\lambda))^{-1}=\sum_{j=1}^{\infty}(\bV(\lambda))^j$
exists in $\CL(\CY_q(\Omega_+,\Omega_-,\Gamma))$.
In addition, from (\ref{6.25}),
$(I-\CV(\lambda))^{-1}=\sum_{j=1}^{\infty}(F_\lambda\CV(\lambda))^j$
exists and
\begin{align}\label{6.*}
\CR_{\CL(\CY_q(\Omega_+,\Omega_-,\Gamma))}
(\{(\tau\pd_\tau\tau)^s(I-\CV(\lambda))^{-1}
\mid \lambda\in\Sigma_{\varepsilon,\lambda_*} \})\leq 2
\end{align}
for $s=0,1$.
By (\ref{6.20}) and (\ref{6.23}),
\begin{align*}
\wh \rho_+ = \CA^+ (\lambda)F_\lambda(I-\bV(\lambda))^{-1}\bF_Y, \enskip \wh \bu_\pm = \CB^\pm (\lambda)F_\lambda(I-\bV(\lambda))^{-1}\bF_Y, \enskip \wh h = \CH(\lambda)F_\lambda(I-\bV(\lambda))^{-1}\bF_Y
\end{align*}
are solutions of problem (\ref{eq-3.8}).
Since $F_\lambda\bV(\lambda)=\CV(\lambda)F_\lambda$, we have
\begin{align*}
F_\lambda(I-\bV(\lambda))^{-1}=F_\lambda\sum_{j=1}^{\infty}(\bV(\lambda))^j
=\sum_{j=1}^{\infty}(\CV(\lambda))^jF_\lambda
=(I-\CV(\lambda))^{-1}F_\lambda.
\end{align*}
Then, we rewrite the representation of $\wh \rho_{+}$, $\wh \bu_{\pm}$, $\wh h$ as follows:
\begin{align*}
\wh \rho_+ = \CA^+ (\lambda)(I-\CV(\lambda))^{-1}F_\lambda(\bF_Y), \enskip \wh \bu_\pm = \CB^\pm (\lambda)(I-\CV(\lambda))^{-1}F_\lambda(\bF_Y), \enskip \wh h = \CH(\lambda)(I-\CV(\lambda))^{-1}F_\lambda(\bF_Y).
\end{align*}
Hence, setting $P^+(\lambda)=\CA^+(\lambda)(I-\CV(\lambda))^{-1}$, $U^\pm(\lambda)=\CB^\pm(\lambda)(I-\CV(\lambda))^{-1}$, and $H(\lambda)=\CH(\lambda)(I-\CV(\lambda))^{-1}$,
by (\ref{6.23}), (\ref{6.*}), and Lemma \ref{lem-4.2}, we see that $P^+(\lambda)$,  $U^\pm(\lambda)$, and $H(\lambda)$ satisfy the properties stated in Theorem~\ref{TH-3.1}.